%% file: Arnold_Inv_Sum.tex
\documentclass[a4paper]{amsart}
\usepackage{amsthm, amsmath, amssymb, url, color, enumerate, bm}
\usepackage{graphicx}
\definecolor{Black}{cmyk}{0,0,0,1}
\usepackage[varg]{txfonts}

\newcommand{\I}{\mathcal{I}}
\newcommand{\R}{\mathbb{R}}
\newcommand{\Z}{\mathbb{Z}}

\newcommand{\abs}[1]{\left\lvert #1\right\rvert}
\newcommand{\ang}[1]{\left\langle#1\right\rangle}
\newcommand{\GI}{\mathcal{GI}}

\newcommand{\ind}{{\rm ind}}
\newcommand{\Int}{{\rm Int}\,}

\newcommand{\St}{\mathit{St}}

\theoremstyle{definition}
\newtheorem{lem}{Lemma}[section]
\newtheorem{defn}[lem]{Definition}
\newtheorem{prop}[lem]{Proposition}
\newtheorem{rem}[lem]{Remark}
\newtheorem{thm}[lem]{Theorem}
\newtheorem{ex}[lem]{Example}
\newtheorem{cor}[lem]{Corollary}
\newtheorem{notation}[lem]{Notation}

\numberwithin{equation}{section}
\numberwithin{figure}{section}

\title[Generalized connected sum formula for the Arnold invariants]{Generalized connected sum formula for the Arnold invariants of generic plane curves}
\author{Keiichi Sakai \and Ryutaro Sugiyama}
\thanks{KS is partially supported by JSPS KAKENHI Grant Number 16K05144}
\address{Faculty of Mathematics, Shinshu University, 3-1-1 Asahi, Matsumoto, Nagano 390-8621, Japan}
\email{ksakai@math.shinshu-u.ac.jp}
\email{15sm108b@shinshu-u.ac.jp}
\date{\today}

\begin{document}
\maketitle

\begin{abstract}
We define the \emph{generalized connected sum} for generic closed plane curves, generalizing the \emph{strange sum} defined by Arnold, and completely describe how the Arnold invariants $J^{\pm}$ and $\St$ behave under the generalized connected sums.
\end{abstract}

\section{Introduction}
A \emph{generic plane curve} is a $C^{\infty}$-immersion $c\colon M\looparrowright\R^2$, where $M$ is either a closed connected interval or $S^1$, without self-tangency points, nor multiple points of multiplicity higher or equal to three.
We call the image $C=c(M)\subset\R^2$ an \emph{unoriented} generic plane curve.
The classification problem of generic \emph{closed} plane curves $c\colon S^1\looparrowright\R^2$ was posed in \cite{Arnold} and is known to be quite complicated.
Arnold defined in \cite{Arnold} three invariants $J^+,J^-$ and $\St$ (the last one is called the \emph{strangeness}) of unoriented generic closed plane curves (Theorems~\ref{thm:J^+}--\ref{thm:St}).
These invariants form a basis of the space of degree one Vassiliev type invariants of generic closed plane curves.
Although they are not complete invariants, geometric understanding of them would be important,
with high-dimensional / high-order analogues in mind.

Let $C_0,C_1$ be any unoriented generic closed plane curves in general position and let $\Gamma$ be a \emph{bridge} between them, namely a generic plane curve that emanates from a generic point (i.e.\ non-double point) on $C_0$ and terminates at a generic point on $C_1$ creating no triple point nor tangency point with $C_0$ and $C_1$.
The \emph{generalized connected sum} of $C_0,C_1$ along $\Gamma$ is  obtained by performing a ``surgery'' on $C_0\cup C_1$ along $\Gamma$ (Definition~\ref{def:generalized_connected_sum}).
This generalizes the \emph{connected-sum} and the \emph{strange sum} \cite[\S3]{Arnold} of \emph{separated} closed curves.

In \cite{Arnold} it has been shown that $\St$ is additive under the strange sums and $J^{\pm}$ is additive under the connected sums (Theorem~\ref{thm:Arnold_additivity}).
Afterward Mendes de Jesus and Romero Fuster \cite{MendesRomero} partially generalized Arnold's result to the case that $C_0\cap C_1$ may be non-empty but $\Int\Gamma\cap(C_0\cup C_1)=\emptyset$ and $\Gamma$ has no double point (Theorem~\ref{thm:MendesRomero}).
They use the numbers $T^{\pm}$ and $T^{\St}_{\Gamma}$ which count with signs the tangencies / triple points that occur in a separating homotopy of $C_0,C_1$ (Definitions~\ref{def:T^pm}, \ref{def:T^St}).
Though these numbers occupy important places in connected sum formulas, any explicit way to calculate them has not been known.

Our aim is to establish complete generalized connected sum formulas for $X=J^{\pm},\St$ in terms of $X(C_i)$ ($i=0,1$) and some geometric data of $C_i$ and $\Gamma$ (Theorems~\ref{thm:main_J}, \ref{thm:main_St}).
Our formulas also include $T^{\pm}$ and $T^{\St}_{\Gamma}$, and we give algorithms to compute them (Propositions~\ref{prop:T^pm}, \ref{prop:T^St}).

We remark that the numbers of double points of generalized connected sums are large in general, and our formulas would be essentially applied to compute the Arnold invariants of generic curves that are not in the tables in \cite[\S10 and appendix by Aicardi]{Arnold}.
It would also be an interesting problem to find a family of generic closed curves whose Arnold invariants can be computed in a systematic way using our formulas.
It would also be curious to understand our formulas using geometric formulas due to Viro and Shmakovich (see \cite{Shumakovich95}).



\input{2_Terminologies}
\input{3_Strange_Sum}
\input{4_General_formula}
\input{5_T}

\section*{Acknowledgment}
The authors express their sincere appreciation to Mendes de Jesus and Romero Fuster for answering questions and giving fruitful comments on the draft of this manuscript, and to the referees for careful reviewing.

\end{document}

%% file: 2_Terminologies.tex
\section{Terminologies}
All the curves that we will consider are plane curves.
The Arnold invariants $J^{\pm}$ and $\St$ assign integers to (unoriented) \emph{generic closed curves}.
In this section we recall the axiomatic definition of these invariants following \cite{Arnold}.

\subsection{Generic curves}
A \emph{regular closed curve} is a smooth immersion $S^1:=\R/\Z\looparrowright\R^2$.
A regular closed curve is endowed with a natural orientation.
In the following, unless otherwise stated, regular closed curves (maps) $c\colon S^1\looparrowright\R^2$ are denoted by lowercase letters and their images $C:=c(S^1)\subset\R^2$ are denoted by capital letters.
We call $C$ an \emph{unoriented regular closed curve}.

Let $\I$ be the space of regular closed curves (maps) equipped with the Whitney $C^{\infty}$-topology.
A \emph{regular homotopy} of regular closed curves is a continuous path $[0,1]\to\I$, or equivalently a continuous family $[0,1]\times S^1\to\R^2$ of regular closed curves.
Regular homotopy for unoriented closed curves is defined in an obvious way.

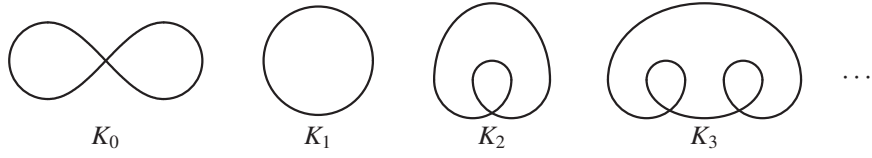
\begin{figure}
\begin{center}
\input{standard_curves}
\end{center}
\caption{The standard curves}
\label{fig:standard_curves}
\end{figure}

\begin{defn}\label{def:rotation_number}
The \emph{rotation number} of $c\in\I$ is defined by
\[
 r(c):=\deg\Bigl(\frac{c'}{\abs{c'}}\colon S^1\to S^1\Bigr)\in\Z.
\]
\end{defn}
It is known \cite{Whitney} that $c_0,c_1\in\I$ can be connected to each other by a regular homotopy in $\I$ if and only if $r(c_0)=r(c_1)$.
Thus the set of regular closed curves up to regular homotopy, namely $\pi_0(\I)$, is identified with $\Z$ via the rotation number.
It can be seen that the rotation number of (a parametrization of) the \emph{standard curve} $K_n$ shown in Figure~\ref{fig:standard_curves} is $\pm n$,
and hence any $c\in\I$ can be transformed to a parametrization of $K_{\abs{r(c)}}$ by a regular homotopy.

\begin{defn}
A \emph{double point} of $c\in\I$ is a point $c(s)=c(t)$ for some distinct $s,t\in S^1$.
A \emph{self-tangency point} of $c$ is a double point $c(s)=c(t)$,
such that $c'(s)=kc'(t)$ for some $k\ne 0$.
If $k>0$ then the self-tangency point is said to be \emph{direct},
and \emph{inverse} otherwise.
A \emph{triple point} of $c$ is a point $c(s)=c(t)=c(u)$,
where $s,t,u\in S^1$ are all distinct.

Define $\Sigma^d,\Sigma^i,\Sigma^t\subset\I$ to be the subspaces consisting of $c\in\I$ with respectively at least one direct self-tangency point, one inverse self-tangency point, and one triple point.
The union $\Sigma:=\Sigma^d\cup\Sigma^i\cup\Sigma^t$ is called the \emph{discriminant set}.
Any element of the complement $\GI:=\I\setminus\Sigma$ is called a \emph{generic closed curve}.
We say two generic closed curves are \emph{equivalent as generic closed curves} if they are joined by a regular homotopy which does not intersect $\Sigma$.
\end{defn}

\begin{rem}
Double / triple points of \emph{unoriented} generic curves are defined in obvious ways.
Moreover ``direct / inverse self-tangency points of \emph{unoriented} curves'' make sense because, if $p\in\R^2$ is a direct (resp.\ an inverse) self-tangency point of $c$, then $p$ is also a direct (resp.\ an inverse) self-tangency point of $c\circ\psi$ for any diffeomorphism $\psi\colon S^1\to S^1$.
\end{rem}

All the multiple points of a generic curve are transverse double points, and there is no triple point nor self-tangency point.
The space $\GI$ of generic closed curves is open dense in $\I$.

The subspace of $\Sigma^d$ (resp.~$\Sigma^i$) consisting of curves with exactly one direct (resp.~inverse) \emph{generic self-tangency point}, namely at the point the second derivatives of two branches of the curve are not equal, is an open dense subspace.
Similarly the subspace of $\Sigma^t$ consisting of curves with exactly one \emph{generic triple point}, namely at the point any two of three branches of the curve are transverse is open dense.
These ``most generic parts'' of $\Sigma$ are of codimension one in $\I$.

\begin{defn}
We say a regular homotopy $[0,1]\to\I$ is \emph{generic} if it transversely intersects $\Sigma$ in a finite set contained in the most generic part of $\Sigma$.
\end{defn}

If $h\colon [0,1]\to\I$ is a generic regular homotopy then $h(s)\in\GI$ for all but finitely many $s_1,\dots,s_k\in[0,1]$,
and for any $i=1,\dots,k$ and any small $\epsilon>0$, the curves $h(s_i+\epsilon)$ and $h(s_i-\epsilon)\in\GI$ are distinct as generic closed curves.
For $i=1,\dots,k$ the regular closed curve $h(s_i)$ has exactly one generic self-tangency point or generic triple point.

\begin{figure}
\begin{center}
\input{tangency_-}
\caption{Getting across an inverse self-tangency}
\label{fig:tangency_-}
\end{center}
\end{figure}
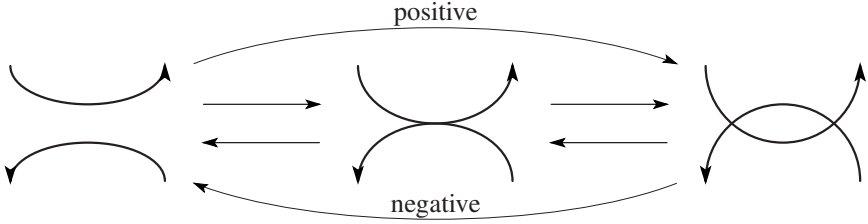

\begin{defn}\label{def:signs}
Define the \emph{co-orientation} of the most generic part of $\Sigma$ as follows;
we say a generic regular homotopy $h\colon [-\epsilon,+\epsilon]\to\I$, $h(0)\in\Sigma^d$ (resp.~$h(0)\in\Sigma^i$), \emph{gets across} $\Sigma^d$ (resp.~$\Sigma^i$) \emph{positively} if the number of double points of $h(\epsilon)$ is greater than that of $h(-\epsilon)$ by two (see Figures~\ref{fig:tangency_+}, \ref{fig:tangency_-}).
We say a generic regular homotopy $h\colon [-\epsilon,\epsilon]\to\I$, $h(0)\in\Sigma^t$, \emph{gets across} $\Sigma^t$ \emph{positively} if $(-1)^q=+1$, where $q$ is the number of edges of the \emph{newborn triangle} of $h(\epsilon)$ near the triple point (see Figures~\ref{fig:newborn_triangle_+}, \ref{fig:newborn_triangle_-}) whose orientations are compatible with the cyclic order of the edges determined by the parameter of $h(\epsilon)$.
\end{defn}
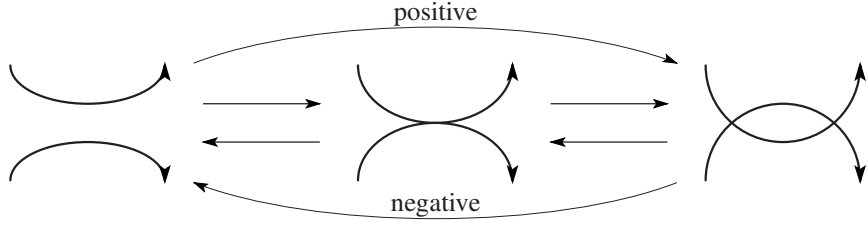
\begin{figure}
\begin{center}
\input{tangency_+}
\caption{Getting across a direct self-tangency}
\label{fig:tangency_+}
\end{center}
\end{figure}

\begin{ex}\label{ex:signs_triplepoint}
In both Figures~\ref{fig:newborn_triangle_+}, \ref{fig:newborn_triangle_-} the edges of the newborn triangles are given the counterclockwise order by the parameter of $h(\epsilon)$.
In the case of Figure~\ref{fig:newborn_triangle_+} the orientations of edges $1,2$ of newborn triangle are compatible with the counterclockwise order of the edges, and hence $q=2$.
Thus in Figure~\ref{fig:newborn_triangle_+} the homotopy gets positively across $\Sigma^t$.
In the case of Figure~\ref{fig:newborn_triangle_-} the orientations of all the edges $1,2,3$ of the newborn triangle are compatible with the order of the edges, and hence $q=3$.
Thus the homotopy gets negatively across $\Sigma^t$.
\end{ex}\begin{figure}
\begin{center}
\input{newborn_triangle_+}
\caption{Getting across a triple point positively ($q=2$)}
\label{fig:newborn_triangle_+}
\end{center}
\end{figure}
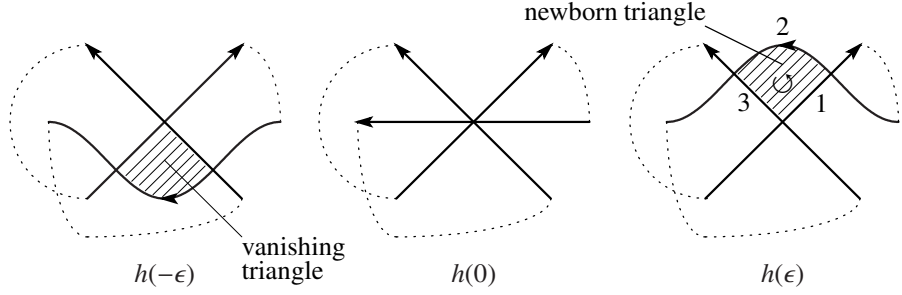
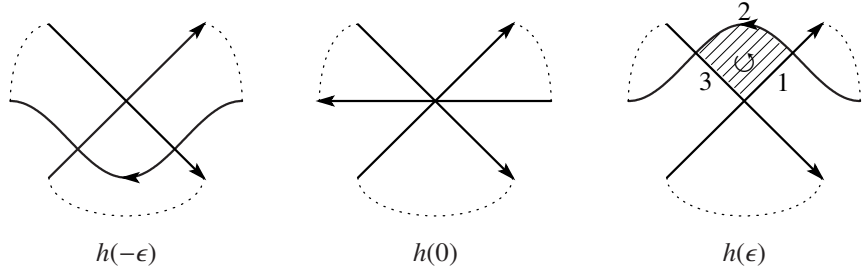
\begin{figure}
\begin{center}
\input{newborn_triangle_-}
\caption{Getting across a triple point negatively ($q=3$)}
\label{fig:newborn_triangle_-}
\end{center}
\end{figure}

\begin{rem}
The co-orientation of $\Sigma^t$ is well-defined, namely cyclic permutations of the order of edges of newborn triangle do not change $q$.
Moreover the definition of the co-orientation of $\Sigma$ is independent of the orientation of the curve;
for example reversing the orientations of curves in Figures~\ref{fig:newborn_triangle_+}, \ref{fig:newborn_triangle_-} does not change the parity of $q$.
Therefore it makes sense saying that ``a generic regular homotopy of \emph{unoriented} closed curves gets across a self-tangency / a triple point \emph{positively} (or \emph{negatively})''.
Here we say a subset $C\subset\R^2$ is an \emph{unoriented generic closed curve} if it is the image of some $c\in\GI$.
\end{rem}

\begin{rem}
If the sign of a newborn triangle is $\pm 1$, then the sign of the corresponding \emph{vanishing triangle} (see Figure~\ref{fig:newborn_triangle_+}) is respectively $\mp 1$.
\end{rem}

\subsection{The Arnold invariants}
The Arnold invariants are the functions $\GI/\mathrm{Diff}(S^1)\to\Z$ characterized by the following Theorems.

\begin{thm}[\cite{Arnold}]\label{thm:J^+}
There exists a unique invariant $J^+$ of unoriented generic closed curves under generic homotopies getting across no direct self-tangency,
satisfying
\begin{itemize}
\item
	$J^+(K_0)=0$ and $J^+(K_n)=-2(n-1)$ ($n\ge 1$) for the standard curves $K_n$ (see Figure~\ref{fig:standard_curves}), and
\item
	if $C_0$ is transformed to $C_1$ through a positive direct self-tangency point, then
	\[
	 J^+(C_1)-J^+(C_0)=2.
	\]
\end{itemize}
\end{thm}

Notice that, after getting through a direct self-tangency negatively, the value of $J^+$ decreases by $2$.
The above properties enable us to compute $J^+(C)$ for any unoriented generic closed curves since any $C$ can be transformed to one of $K_n$'s by a generic regular homotopy \cite{Whitney}, and in fact $J^+(C)$ is independent of the choice of generic homotopies since the most generic part of $\Sigma$ is ``consistently co-oriented''.
See \cite[\S2]{Arnold} for details.
Similarly $J^-(C)$ and $\St(C)$ below are well defined and can be computed in principle.

\begin{thm}[\cite{Arnold}]\label{thm:J^-}
There exists a unique invariant $J^-$ of unoriented generic closed curves under generic homotopies getting across no inverse self-tangency, satisfying
\begin{itemize}
\item
	$J^-(K_0)=-1$ and $J^-(K_n)=-3(n-1)$ ($n\ge 1$), and
\item
	if $C_0$ is transformed to $C_1$ through a positive inverse self-tangency point, then
	\[
	 J^-(C_1)-J^-(C_0)=-2.
	\]
\end{itemize}
\end{thm}

\begin{thm}[\cite{Arnold}]\label{thm:St}
There exists a unique invariant $\St$ of unoriented generic closed curves under generic homotopies getting across no triple point, satisfying
\begin{itemize}
\item
	$\St(K_0)=0$ and $\St(K_n)=n-1$ ($n\ge 1$), and
\item
	if $C_0$ is transformed to $C_1$ through a positive triple point, then
	\[
	 \St(C_1)-\St(C_0)=1.
	\]
\end{itemize}
\end{thm}

\subsection{Generalized connected sum}
Let $C_0,C_1$ be unoriented generic closed curves.

\begin{defn}
The image $\Gamma:=\gamma([0,1])$ of a generic curve $\gamma\colon [0,1]\to\R^2$ is a \emph{bridge between $C_0$ and $C_1$} if $b_i:=\gamma(i)\in C_i$ ($i=0,1$) and $\Gamma$ is generic with respect to $C_0$ and $C_1$ (see Figure~\ref{fig:generalized_sum}), that means
\begin{enumerate}[(1)]
\item
	$\Gamma\cap D(C_i)=\emptyset$ ($i=0,1$), where $D(C_i)$ is the set of double points of $C_i$,
\item
	$D(\Gamma)\cap C_i=\emptyset$ ($i=0,1$), and
\item
	$\Gamma$ intersects $C_0,C_1$ transversely (transversality at $b_0$ and $b_1$ is also required).
\end{enumerate}
\end{defn}

\begin{defn}
We say that two generic closed curves are \emph{in general position} if there exists no tangency nor triple point involving both curves.
\end{defn}

\begin{defn}\label{def:generalized_connected_sum}
Let $C_0,C_1$ be unoriented generic closed curves in general position and $\Gamma$ a bridge between them with endpoints $b_i\in C_i$ ($i=0,1$).
Define the \emph{generalized connected sum} $C_0+_{\Gamma}C_1$ along $\Gamma$ as in Figure~\ref{fig:generalized_sum};
\begin{figure}
\begin{center}
\input{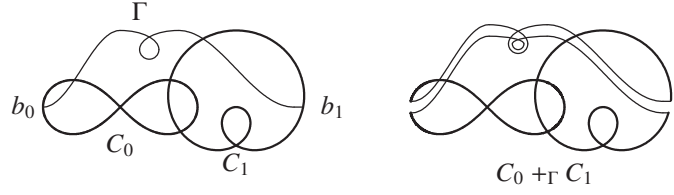}
\caption{An example of generalized connected sum}
\label{fig:generalized_sum}
\end{center}
\end{figure}
first remove from $C_i$ ($i=0,1$) a small interval $I_i$ satisfying $b_i\in I_i$, $I_i\cap\Gamma=\{b_i\}$, $I_i\cap D(C_i)=\emptyset$ and $I_i\cap C_{1-i}=\emptyset$.
Next replace $\Gamma$ with a pair of generic curves $\Gamma_{\pm}$ that are parallel to the original $\Gamma$ and connect the boundary points of $I_0$ and $I_1$, and that are close to $\Gamma$ so that $\Gamma_{\pm}\cap(C_0\cup C_1)$ can be naturally identified with $\Gamma\cap(C_0\cup C_1)$.
Define $C_0+_{\Gamma}C_1$ as the curve $(\Gamma\setminus(I_0\sqcup I_1))\cup(\Gamma_+\cup\Gamma_-)$ (with its corners smoothed).
\end{defn}

$C_0+_{\Gamma}C_1$ is a generic closed curve if $C_0,C_1$ are generic and in general position, and if the parallel curves $\Gamma_{\pm}$ are taken to be sufficiently close to $\Gamma$.

\begin{defn}\label{def:compatible_ori}
Let $C_0,C_1$ be unoriented generic closed curves and $\Gamma$ a bridge between them.
We say parameters $c_i$ ($i=0,1$) of $C_i$ are \emph{compatible with respect to $\Gamma$} if there exists a parameter $c$ of $C_0+_{\Gamma}C_1$ that determines the same orientation of $C_i$ as that by $c_i$ (see Figure~\ref{fig:compatible_ori}) .
\end{defn}
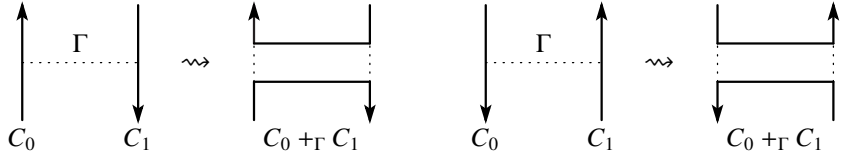
\begin{figure}
\begin{center}
\input{compatible_ori}
\end{center}
\caption{Compatible orientations}
\label{fig:compatible_ori}
\end{figure}
If $c_0$ and $c_1$ are compatible then so are $-c_0$ and $-c_1$, where in general $(-c)(t):=c(1-t)$ is the curve with opposite orientation to $c$.
Any quantities in the following are independent of the choices.
Unless otherwise stated, we always choose compatible parameters $c_i$ of $C_i$ ($i=0,1$) with respect to bridge $\Gamma$.

Note that the complement $\R^2\setminus C$ of a generic closed curve $C$ is decomposed into finitely many connected components, exactly one of which is unbounded.

\begin{defn}
For a generic closed curve $C$, we denote the unbounded component of $\R^2\setminus C$ by $R^{\infty}$.
We say (unoriented) generic closed curves $C_0,C_1$ are \emph{separated} if $C_0\subset R^{\infty}_1$ and $C_1\subset R^{\infty}_0$,
where $R^{\infty}_i\subset\R^2\setminus C_i$ ($i=0,1$) is the unbounded component.
\end{defn}

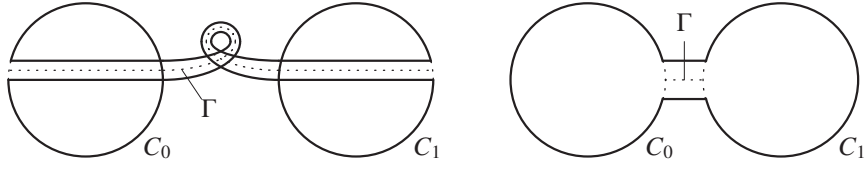
\begin{figure}
\begin{center}
\input{strange_connected_sum}
\caption{Strange sum and connected sum}
\label{fig:strange_connected_sum}
\end{center}
\end{figure}

When $C_0$ and $C_1$ are separated, $C_0+_{\Gamma}C_1$ is \emph{strange sum} introduced in \cite[\S3]{Arnold}.
If moreover $\Int\Gamma\cap(C_0\cup C_1)=\emptyset$ (where $\Int\Gamma:=\Gamma\setminus\partial\Gamma$) and $D(\Gamma)=\emptyset$, then $C_0+_{\Gamma}C_1$ is \emph{connected sum} introduced in \cite[\S3]{Arnold}.
See Figure~\ref{fig:strange_connected_sum}.

Our aim is to generalize the following additivity formulas for the Arnold invariants (and their generalizations given in \cite{MendesRomero}, see Theorem~\ref{thm:MendesRomero}) to generalized connected sums.

\begin{thm}[\cite{Arnold}]\label{thm:Arnold_additivity}
\begin{enumerate}[(1)]
\item
	If $C_0+_{\Gamma}C_1$ is a connected sum, then
	$J^{\pm}(C_0+_{\Gamma}C_1)=J^{\pm}(C_0)+J^{\pm}(C_1)$.
\item
	If $C_0+_{\Gamma}C_1$ is a strange sum, then $\St(C_0+_{\Gamma}C_1)=\St(C_0)+\St(C_1)$.
\end{enumerate}
\end{thm}

%% file: standard_curves.tex
{\unitlength 0.1in%
\begin{picture}(41.6500,6.2700)(7.0000,-26.2700)%
\put(12.0000,-27.0000){\makebox(0,0){$K_0$}}%
\put(32.0000,-27.0000){\makebox(0,0){$K_2$}}%
%
\special{pn 13}%
\special{ar 2300 2300 283 283 0.0000000 6.2831853}%
\put(23.0000,-27.0000){\makebox(0,0){$K_1$}}%
\put(43.0000,-27.0000){\makebox(0,0){$K_3$}}%
\put(51.0000,-24.0000){\makebox(0,0){$\cdots$}}%
\special{pn 13}%
\special{pa 900 2500}%
\special{pa 910 2500}%
\special{pa 915 2499}%
\special{pa 920 2499}%
\special{pa 925 2498}%
\special{pa 930 2498}%
\special{pa 940 2496}%
\special{pa 945 2494}%
\special{pa 955 2492}%
\special{pa 960 2490}%
\special{pa 965 2489}%
\special{pa 985 2481}%
\special{pa 990 2478}%
\special{pa 995 2476}%
\special{pa 1000 2473}%
\special{pa 1005 2471}%
\special{pa 1025 2459}%
\special{pa 1030 2455}%
\special{pa 1040 2449}%
\special{pa 1050 2441}%
\special{pa 1055 2438}%
\special{pa 1080 2418}%
\special{pa 1085 2413}%
\special{pa 1090 2409}%
\special{pa 1095 2404}%
\special{pa 1100 2400}%
\special{pa 1105 2395}%
\special{pa 1110 2391}%
\special{pa 1120 2381}%
\special{pa 1125 2377}%
\special{pa 1160 2342}%
\special{pa 1165 2336}%
\special{pa 1185 2316}%
\special{pa 1190 2310}%
\special{pa 1210 2290}%
\special{pa 1215 2284}%
\special{pa 1235 2264}%
\special{pa 1240 2258}%
\special{pa 1275 2223}%
\special{pa 1280 2219}%
\special{pa 1290 2209}%
\special{pa 1295 2205}%
\special{pa 1300 2200}%
\special{pa 1305 2196}%
\special{pa 1310 2191}%
\special{pa 1315 2187}%
\special{pa 1320 2182}%
\special{pa 1345 2162}%
\special{pa 1350 2159}%
\special{pa 1360 2151}%
\special{pa 1370 2145}%
\special{pa 1375 2141}%
\special{pa 1395 2129}%
\special{pa 1400 2127}%
\special{pa 1405 2124}%
\special{pa 1410 2122}%
\special{pa 1415 2119}%
\special{pa 1435 2111}%
\special{pa 1440 2110}%
\special{pa 1445 2108}%
\special{pa 1455 2106}%
\special{pa 1460 2104}%
\special{pa 1470 2102}%
\special{pa 1475 2102}%
\special{pa 1480 2101}%
\special{pa 1485 2101}%
\special{pa 1490 2100}%
\special{pa 1500 2100}%
\special{fp}%
\special{pn 13}%
\special{pa 900 2100}%
\special{pa 910 2100}%
\special{pa 915 2101}%
\special{pa 920 2101}%
\special{pa 925 2102}%
\special{pa 930 2102}%
\special{pa 940 2104}%
\special{pa 945 2106}%
\special{pa 955 2108}%
\special{pa 960 2110}%
\special{pa 965 2111}%
\special{pa 985 2119}%
\special{pa 990 2122}%
\special{pa 995 2124}%
\special{pa 1000 2127}%
\special{pa 1005 2129}%
\special{pa 1025 2141}%
\special{pa 1030 2145}%
\special{pa 1040 2151}%
\special{pa 1050 2159}%
\special{pa 1055 2162}%
\special{pa 1080 2182}%
\special{pa 1085 2187}%
\special{pa 1090 2191}%
\special{pa 1095 2196}%
\special{pa 1100 2200}%
\special{pa 1105 2205}%
\special{pa 1110 2209}%
\special{pa 1120 2219}%
\special{pa 1125 2223}%
\special{pa 1160 2258}%
\special{pa 1165 2264}%
\special{pa 1185 2284}%
\special{pa 1190 2290}%
\special{pa 1210 2310}%
\special{pa 1215 2316}%
\special{pa 1235 2336}%
\special{pa 1240 2342}%
\special{pa 1275 2377}%
\special{pa 1280 2381}%
\special{pa 1290 2391}%
\special{pa 1295 2395}%
\special{pa 1300 2400}%
\special{pa 1305 2404}%
\special{pa 1310 2409}%
\special{pa 1315 2413}%
\special{pa 1320 2418}%
\special{pa 1345 2438}%
\special{pa 1350 2441}%
\special{pa 1360 2449}%
\special{pa 1370 2455}%
\special{pa 1375 2459}%
\special{pa 1395 2471}%
\special{pa 1400 2473}%
\special{pa 1405 2476}%
\special{pa 1410 2478}%
\special{pa 1415 2481}%
\special{pa 1435 2489}%
\special{pa 1440 2490}%
\special{pa 1445 2492}%
\special{pa 1455 2494}%
\special{pa 1460 2496}%
\special{pa 1470 2498}%
\special{pa 1475 2498}%
\special{pa 1480 2499}%
\special{pa 1485 2499}%
\special{pa 1490 2500}%
\special{pa 1500 2500}%
\special{fp}%
%
\special{pn 13}%
\special{ar 900 2300 200 200 1.5707963 4.7123890}%
%
\special{pn 13}%
\special{ar 1500 2300 200 200 4.7123890 1.5707963}%
%
\special{pn 13}%
\special{ar 3100 2400 200 200 6.2831853 3.1415927}%
%
\special{pn 13}%
\special{ar 3300 2400 200 200 6.2831853 3.1415927}%
%
\special{pn 13}%
\special{ar 3200 2400 100 100 3.1415927 6.2831853}%
%
\special{pn 13}%
\special{ar 3200 2400 300 400 3.1415927 6.2831853}%
%
\special{pn 13}%
\special{ar 4000 2400 200 200 6.2831853 3.1415927}%
%
\special{pn 13}%
\special{ar 4600 2400 200 200 6.2831853 3.1415927}%
%
\special{pn 13}%
\special{ar 4300 2400 300 200 6.2831853 3.1415927}%
%
\special{pn 13}%
\special{ar 4100 2400 100 100 3.1415927 6.2831853}%
%
\special{pn 13}%
\special{ar 4500 2400 100 100 3.1415927 6.2831853}%
%
\special{pn 13}%
\special{ar 4300 2400 500 400 3.1415927 6.2831853}%
\end{picture}}%

%% file: tangency_+.tex
\unitlength 0.1in
\begin{picture}( 44.0000, 11.3500)(  8.0000,-12.0000)
%
{\color[named]{Black}{%
\special{pn 13}%
\special{ar 1200 1000 400 200  3.1415927  6.2831853}%
}}%
%
{\color[named]{Black}{%
\special{pn 13}%
\special{pa 1600 990}%
\special{pa 1600 1000}%
\special{fp}%
\special{sh 1}%
\special{pa 1600 1000}%
\special{pa 1620 934}%
\special{pa 1600 948}%
\special{pa 1580 934}%
\special{pa 1600 1000}%
\special{fp}%
}}%
%
{\color[named]{Black}{%
\special{pn 13}%
\special{ar 1200 400 400 200  6.2831853  3.1415927}%
}}%
%
{\color[named]{Black}{%
\special{pn 13}%
\special{pa 1600 410}%
\special{pa 1600 400}%
\special{fp}%
\special{sh 1}%
\special{pa 1600 400}%
\special{pa 1580 468}%
\special{pa 1600 454}%
\special{pa 1620 468}%
\special{pa 1600 400}%
\special{fp}%
}}%
%
{\color[named]{Black}{%
\special{pn 8}%
\special{pa 1800 600}%
\special{pa 2400 600}%
\special{fp}%
\special{sh 1}%
\special{pa 2400 600}%
\special{pa 2334 580}%
\special{pa 2348 600}%
\special{pa 2334 620}%
\special{pa 2400 600}%
\special{fp}%
}}%
%
{\color[named]{Black}{%
\special{pn 8}%
\special{pa 2400 800}%
\special{pa 1800 800}%
\special{fp}%
\special{sh 1}%
\special{pa 1800 800}%
\special{pa 1868 820}%
\special{pa 1854 800}%
\special{pa 1868 780}%
\special{pa 1800 800}%
\special{fp}%
}}%
\put(30.0000,-1.3000){\makebox(0,0){positive}}%
%
{\color[named]{Black}{%
\special{pn 13}%
\special{ar 3000 400 400 300  6.2831853  3.1415927}%
}}%
%
{\color[named]{Black}{%
\special{pn 13}%
\special{pa 3400 414}%
\special{pa 3400 400}%
\special{fp}%
\special{sh 1}%
\special{pa 3400 400}%
\special{pa 3380 468}%
\special{pa 3400 454}%
\special{pa 3420 468}%
\special{pa 3400 400}%
\special{fp}%
}}%
%
{\color[named]{Black}{%
\special{pn 13}%
\special{ar 3000 1000 400 300  3.1415927  6.2831853}%
}}%
%
{\color[named]{Black}{%
\special{pn 13}%
\special{pa 3400 988}%
\special{pa 3400 1000}%
\special{fp}%
\special{sh 1}%
\special{pa 3400 1000}%
\special{pa 3420 934}%
\special{pa 3400 948}%
\special{pa 3380 934}%
\special{pa 3400 1000}%
\special{fp}%
}}%
%
{\color[named]{Black}{%
\special{pn 8}%
\special{pa 4200 800}%
\special{pa 3600 800}%
\special{fp}%
\special{sh 1}%
\special{pa 3600 800}%
\special{pa 3668 820}%
\special{pa 3654 800}%
\special{pa 3668 780}%
\special{pa 3600 800}%
\special{fp}%
}}%
%
{\color[named]{Black}{%
\special{pn 8}%
\special{pa 3600 600}%
\special{pa 4200 600}%
\special{fp}%
\special{sh 1}%
\special{pa 4200 600}%
\special{pa 4134 580}%
\special{pa 4148 600}%
\special{pa 4134 620}%
\special{pa 4200 600}%
\special{fp}%
}}%
%
{\color[named]{Black}{%
\special{pn 13}%
\special{ar 4800 400 400 400  6.2831853  3.1415927}%
}}%
%
{\color[named]{Black}{%
\special{pn 13}%
\special{pa 5200 416}%
\special{pa 5200 400}%
\special{fp}%
\special{sh 1}%
\special{pa 5200 400}%
\special{pa 5180 468}%
\special{pa 5200 454}%
\special{pa 5220 468}%
\special{pa 5200 400}%
\special{fp}%
}}%
%
{\color[named]{Black}{%
\special{pn 13}%
\special{ar 4800 1000 400 400  3.1415927  6.2831853}%
}}%
%
{\color[named]{Black}{%
\special{pn 13}%
\special{pa 5200 986}%
\special{pa 5200 1000}%
\special{fp}%
\special{sh 1}%
\special{pa 5200 1000}%
\special{pa 5220 934}%
\special{pa 5200 948}%
\special{pa 5180 934}%
\special{pa 5200 1000}%
\special{fp}%
}}%
%
{\color[named]{Black}{%
\special{pn 4}%
\special{ar 3000 700 1600 500  3.8163336  5.6084444}%
}}%
%
{\color[named]{Black}{%
\special{pn 4}%
\special{pa 4236 382}%
\special{pa 4250 388}%
\special{fp}%
\special{sh 1}%
\special{pa 4250 388}%
\special{pa 4196 344}%
\special{pa 4200 368}%
\special{pa 4180 380}%
\special{pa 4250 388}%
\special{fp}%
}}%
%
{\color[named]{Black}{%
\special{pn 4}%
\special{ar 3000 700 1600 500  0.6747409  2.4668517}%
}}%
%
{\color[named]{Black}{%
\special{pn 4}%
\special{pa 1766 1018}%
\special{pa 1752 1012}%
\special{fp}%
\special{sh 1}%
\special{pa 1752 1012}%
\special{pa 1804 1058}%
\special{pa 1800 1034}%
\special{pa 1820 1020}%
\special{pa 1752 1012}%
\special{fp}%
}}%
\put(30.0000,-11.3000){\makebox(0,0){negative}}%
\end{picture}%

%% file: tangency_-.tex
\unitlength 0.1in
\begin{picture}( 44.0000, 11.3500)(  8.0000,-12.0000)
%
{\color[named]{Black}{%
\special{pn 13}%
\special{ar 1200 1000 400 200  3.1415927  6.2831853}%
}}%
%
{\color[named]{Black}{%
\special{pn 13}%
\special{pa 800 990}%
\special{pa 800 1000}%
\special{fp}%
\special{sh 1}%
\special{pa 800 1000}%
\special{pa 820 934}%
\special{pa 800 948}%
\special{pa 780 934}%
\special{pa 800 1000}%
\special{fp}%
}}%
%
{\color[named]{Black}{%
\special{pn 13}%
\special{ar 1200 400 400 200  6.2831853  3.1415927}%
}}%
%
{\color[named]{Black}{%
\special{pn 13}%
\special{pa 1600 410}%
\special{pa 1600 400}%
\special{fp}%
\special{sh 1}%
\special{pa 1600 400}%
\special{pa 1580 468}%
\special{pa 1600 454}%
\special{pa 1620 468}%
\special{pa 1600 400}%
\special{fp}%
}}%
%
{\color[named]{Black}{%
\special{pn 8}%
\special{pa 1800 600}%
\special{pa 2400 600}%
\special{fp}%
\special{sh 1}%
\special{pa 2400 600}%
\special{pa 2334 580}%
\special{pa 2348 600}%
\special{pa 2334 620}%
\special{pa 2400 600}%
\special{fp}%
}}%
%
{\color[named]{Black}{%
\special{pn 8}%
\special{pa 2400 800}%
\special{pa 1800 800}%
\special{fp}%
\special{sh 1}%
\special{pa 1800 800}%
\special{pa 1868 820}%
\special{pa 1854 800}%
\special{pa 1868 780}%
\special{pa 1800 800}%
\special{fp}%
}}%
\put(30.0000,-1.3000){\makebox(0,0){positive}}%
%
{\color[named]{Black}{%
\special{pn 13}%
\special{ar 3000 400 400 300  6.2831853  3.1415927}%
}}%
%
{\color[named]{Black}{%
\special{pn 13}%
\special{pa 3400 414}%
\special{pa 3400 400}%
\special{fp}%
\special{sh 1}%
\special{pa 3400 400}%
\special{pa 3380 468}%
\special{pa 3400 454}%
\special{pa 3420 468}%
\special{pa 3400 400}%
\special{fp}%
}}%
%
{\color[named]{Black}{%
\special{pn 13}%
\special{ar 3000 1000 400 300  3.1415927  6.2831853}%
}}%
%
{\color[named]{Black}{%
\special{pn 13}%
\special{pa 2600 988}%
\special{pa 2600 1000}%
\special{fp}%
\special{sh 1}%
\special{pa 2600 1000}%
\special{pa 2620 934}%
\special{pa 2600 948}%
\special{pa 2580 934}%
\special{pa 2600 1000}%
\special{fp}%
}}%
%
{\color[named]{Black}{%
\special{pn 8}%
\special{pa 4200 800}%
\special{pa 3600 800}%
\special{fp}%
\special{sh 1}%
\special{pa 3600 800}%
\special{pa 3668 820}%
\special{pa 3654 800}%
\special{pa 3668 780}%
\special{pa 3600 800}%
\special{fp}%
}}%
%
{\color[named]{Black}{%
\special{pn 8}%
\special{pa 3600 600}%
\special{pa 4200 600}%
\special{fp}%
\special{sh 1}%
\special{pa 4200 600}%
\special{pa 4134 580}%
\special{pa 4148 600}%
\special{pa 4134 620}%
\special{pa 4200 600}%
\special{fp}%
}}%
%
{\color[named]{Black}{%
\special{pn 13}%
\special{ar 4800 400 400 400  6.2831853  3.1415927}%
}}%
%
{\color[named]{Black}{%
\special{pn 13}%
\special{pa 5200 416}%
\special{pa 5200 400}%
\special{fp}%
\special{sh 1}%
\special{pa 5200 400}%
\special{pa 5180 468}%
\special{pa 5200 454}%
\special{pa 5220 468}%
\special{pa 5200 400}%
\special{fp}%
}}%
%
{\color[named]{Black}{%
\special{pn 13}%
\special{ar 4800 1000 400 400  3.1415927  6.2831853}%
}}%
%
{\color[named]{Black}{%
\special{pn 13}%
\special{pa 4400 986}%
\special{pa 4400 1000}%
\special{fp}%
\special{sh 1}%
\special{pa 4400 1000}%
\special{pa 4420 934}%
\special{pa 4400 948}%
\special{pa 4380 934}%
\special{pa 4400 1000}%
\special{fp}%
}}%
%
{\color[named]{Black}{%
\special{pn 4}%
\special{ar 3000 700 1600 500  3.8163336  5.6084444}%
}}%
%
{\color[named]{Black}{%
\special{pn 4}%
\special{pa 4236 382}%
\special{pa 4250 388}%
\special{fp}%
\special{sh 1}%
\special{pa 4250 388}%
\special{pa 4196 344}%
\special{pa 4200 368}%
\special{pa 4180 380}%
\special{pa 4250 388}%
\special{fp}%
}}%
%
{\color[named]{Black}{%
\special{pn 4}%
\special{ar 3000 700 1600 500  0.6747409  2.4668517}%
}}%
%
{\color[named]{Black}{%
\special{pn 4}%
\special{pa 1766 1018}%
\special{pa 1752 1012}%
\special{fp}%
\special{sh 1}%
\special{pa 1752 1012}%
\special{pa 1804 1058}%
\special{pa 1800 1034}%
\special{pa 1820 1020}%
\special{pa 1752 1012}%
\special{fp}%
}}%
\put(30.0000,-11.3000){\makebox(0,0){negative}}%
\end{picture}%

%% file: newborn_triangle_+.tex
{\unitlength 0.1in%
\begin{picture}(46.0000,13.6500)(4.0000,-17.3500)%
\special{pn 13}%
\special{pa 600 1000}%
\special{pa 610 1000}%
\special{pa 615 1001}%
\special{pa 620 1001}%
\special{pa 625 1002}%
\special{pa 630 1002}%
\special{pa 640 1004}%
\special{pa 645 1006}%
\special{pa 655 1008}%
\special{pa 660 1010}%
\special{pa 665 1011}%
\special{pa 685 1019}%
\special{pa 690 1022}%
\special{pa 695 1024}%
\special{pa 700 1027}%
\special{pa 705 1029}%
\special{pa 725 1041}%
\special{pa 730 1045}%
\special{pa 740 1051}%
\special{pa 750 1059}%
\special{pa 755 1062}%
\special{pa 780 1082}%
\special{pa 785 1087}%
\special{pa 790 1091}%
\special{pa 795 1096}%
\special{pa 800 1100}%
\special{pa 805 1105}%
\special{pa 810 1109}%
\special{pa 820 1119}%
\special{pa 825 1123}%
\special{pa 860 1158}%
\special{pa 865 1164}%
\special{pa 885 1184}%
\special{pa 890 1190}%
\special{pa 910 1210}%
\special{pa 915 1216}%
\special{pa 935 1236}%
\special{pa 940 1242}%
\special{pa 975 1277}%
\special{pa 980 1281}%
\special{pa 990 1291}%
\special{pa 995 1295}%
\special{pa 1000 1300}%
\special{pa 1005 1304}%
\special{pa 1010 1309}%
\special{pa 1015 1313}%
\special{pa 1020 1318}%
\special{pa 1045 1338}%
\special{pa 1050 1341}%
\special{pa 1060 1349}%
\special{pa 1070 1355}%
\special{pa 1075 1359}%
\special{pa 1095 1371}%
\special{pa 1100 1373}%
\special{pa 1105 1376}%
\special{pa 1110 1378}%
\special{pa 1115 1381}%
\special{pa 1135 1389}%
\special{pa 1140 1390}%
\special{pa 1145 1392}%
\special{pa 1155 1394}%
\special{pa 1160 1396}%
\special{pa 1170 1398}%
\special{pa 1175 1398}%
\special{pa 1180 1399}%
\special{pa 1185 1399}%
\special{pa 1190 1400}%
\special{pa 1200 1400}%
\special{fp}%
\special{pn 13}%
\special{pa 1200 1400}%
\special{pa 1210 1400}%
\special{pa 1215 1399}%
\special{pa 1220 1399}%
\special{pa 1225 1398}%
\special{pa 1230 1398}%
\special{pa 1240 1396}%
\special{pa 1245 1394}%
\special{pa 1255 1392}%
\special{pa 1260 1390}%
\special{pa 1265 1389}%
\special{pa 1285 1381}%
\special{pa 1290 1378}%
\special{pa 1295 1376}%
\special{pa 1300 1373}%
\special{pa 1305 1371}%
\special{pa 1325 1359}%
\special{pa 1330 1355}%
\special{pa 1340 1349}%
\special{pa 1350 1341}%
\special{pa 1355 1338}%
\special{pa 1380 1318}%
\special{pa 1385 1313}%
\special{pa 1390 1309}%
\special{pa 1395 1304}%
\special{pa 1400 1300}%
\special{pa 1405 1295}%
\special{pa 1410 1291}%
\special{pa 1420 1281}%
\special{pa 1425 1277}%
\special{pa 1460 1242}%
\special{pa 1465 1236}%
\special{pa 1485 1216}%
\special{pa 1490 1210}%
\special{pa 1510 1190}%
\special{pa 1515 1184}%
\special{pa 1535 1164}%
\special{pa 1540 1158}%
\special{pa 1575 1123}%
\special{pa 1580 1119}%
\special{pa 1590 1109}%
\special{pa 1595 1105}%
\special{pa 1600 1100}%
\special{pa 1605 1096}%
\special{pa 1610 1091}%
\special{pa 1615 1087}%
\special{pa 1620 1082}%
\special{pa 1645 1062}%
\special{pa 1650 1059}%
\special{pa 1660 1051}%
\special{pa 1670 1045}%
\special{pa 1675 1041}%
\special{pa 1695 1029}%
\special{pa 1700 1027}%
\special{pa 1705 1024}%
\special{pa 1710 1022}%
\special{pa 1715 1019}%
\special{pa 1735 1011}%
\special{pa 1740 1010}%
\special{pa 1745 1008}%
\special{pa 1755 1006}%
\special{pa 1760 1004}%
\special{pa 1770 1002}%
\special{pa 1775 1002}%
\special{pa 1780 1001}%
\special{pa 1785 1001}%
\special{pa 1790 1000}%
\special{pa 1800 1000}%
\special{fp}%
\put(12.0000,-18.0000){\makebox(0,0){$h(-\epsilon)$}}%
\put(28.0000,-18.0000){\makebox(0,0){$h(0)$}}%
\put(44.0000,-18.0000){\makebox(0,0){$h(\epsilon)$}}%
\special{pn 13}%
\special{pa 4400 600}%
\special{pa 4410 600}%
\special{pa 4415 601}%
\special{pa 4420 601}%
\special{pa 4425 602}%
\special{pa 4430 602}%
\special{pa 4440 604}%
\special{pa 4445 606}%
\special{pa 4455 608}%
\special{pa 4460 610}%
\special{pa 4465 611}%
\special{pa 4485 619}%
\special{pa 4490 622}%
\special{pa 4495 624}%
\special{pa 4500 627}%
\special{pa 4505 629}%
\special{pa 4525 641}%
\special{pa 4530 645}%
\special{pa 4540 651}%
\special{pa 4550 659}%
\special{pa 4555 662}%
\special{pa 4580 682}%
\special{pa 4585 687}%
\special{pa 4590 691}%
\special{pa 4595 696}%
\special{pa 4600 700}%
\special{pa 4605 705}%
\special{pa 4610 709}%
\special{pa 4620 719}%
\special{pa 4625 723}%
\special{pa 4660 758}%
\special{pa 4665 764}%
\special{pa 4685 784}%
\special{pa 4690 790}%
\special{pa 4710 810}%
\special{pa 4715 816}%
\special{pa 4735 836}%
\special{pa 4740 842}%
\special{pa 4775 877}%
\special{pa 4780 881}%
\special{pa 4790 891}%
\special{pa 4795 895}%
\special{pa 4800 900}%
\special{pa 4805 904}%
\special{pa 4810 909}%
\special{pa 4815 913}%
\special{pa 4820 918}%
\special{pa 4845 938}%
\special{pa 4850 941}%
\special{pa 4860 949}%
\special{pa 4870 955}%
\special{pa 4875 959}%
\special{pa 4895 971}%
\special{pa 4900 973}%
\special{pa 4905 976}%
\special{pa 4910 978}%
\special{pa 4915 981}%
\special{pa 4935 989}%
\special{pa 4940 990}%
\special{pa 4945 992}%
\special{pa 4955 994}%
\special{pa 4960 996}%
\special{pa 4970 998}%
\special{pa 4975 998}%
\special{pa 4980 999}%
\special{pa 4985 999}%
\special{pa 4990 1000}%
\special{pa 5000 1000}%
\special{fp}%
\special{pn 13}%
\special{pa 3800 1000}%
\special{pa 3810 1000}%
\special{pa 3815 999}%
\special{pa 3820 999}%
\special{pa 3825 998}%
\special{pa 3830 998}%
\special{pa 3840 996}%
\special{pa 3845 994}%
\special{pa 3855 992}%
\special{pa 3860 990}%
\special{pa 3865 989}%
\special{pa 3885 981}%
\special{pa 3890 978}%
\special{pa 3895 976}%
\special{pa 3900 973}%
\special{pa 3905 971}%
\special{pa 3925 959}%
\special{pa 3930 955}%
\special{pa 3940 949}%
\special{pa 3950 941}%
\special{pa 3955 938}%
\special{pa 3980 918}%
\special{pa 3985 913}%
\special{pa 3990 909}%
\special{pa 3995 904}%
\special{pa 4000 900}%
\special{pa 4005 895}%
\special{pa 4010 891}%
\special{pa 4020 881}%
\special{pa 4025 877}%
\special{pa 4060 842}%
\special{pa 4065 836}%
\special{pa 4085 816}%
\special{pa 4090 810}%
\special{pa 4110 790}%
\special{pa 4115 784}%
\special{pa 4135 764}%
\special{pa 4140 758}%
\special{pa 4175 723}%
\special{pa 4180 719}%
\special{pa 4190 709}%
\special{pa 4195 705}%
\special{pa 4200 700}%
\special{pa 4205 696}%
\special{pa 4210 691}%
\special{pa 4215 687}%
\special{pa 4220 682}%
\special{pa 4245 662}%
\special{pa 4250 659}%
\special{pa 4260 651}%
\special{pa 4270 645}%
\special{pa 4275 641}%
\special{pa 4295 629}%
\special{pa 4300 627}%
\special{pa 4305 624}%
\special{pa 4310 622}%
\special{pa 4315 619}%
\special{pa 4335 611}%
\special{pa 4340 610}%
\special{pa 4345 608}%
\special{pa 4355 606}%
\special{pa 4360 604}%
\special{pa 4370 602}%
\special{pa 4375 602}%
\special{pa 4380 601}%
\special{pa 4385 601}%
\special{pa 4390 600}%
\special{pa 4400 600}%
\special{fp}%
%
\special{pn 8}%
\special{pn 8}%
\special{pa 1600 600}%
\special{pa 1609 600}%
\special{fp}%
\special{pa 1645 610}%
\special{pa 1652 614}%
\special{fp}%
\special{pa 1684 637}%
\special{pa 1690 643}%
\special{fp}%
\special{pa 1715 673}%
\special{pa 1720 680}%
\special{fp}%
\special{pa 1740 714}%
\special{pa 1743 722}%
\special{fp}%
\special{pa 1759 757}%
\special{pa 1762 765}%
\special{fp}%
\special{pa 1774 803}%
\special{pa 1776 811}%
\special{fp}%
\special{pa 1785 849}%
\special{pa 1787 857}%
\special{fp}%
\special{pa 1793 896}%
\special{pa 1794 904}%
\special{fp}%
\special{pa 1798 943}%
\special{pa 1799 952}%
\special{fp}%
\special{pa 1800 991}%
\special{pa 1800 1000}%
\special{fp}%
%
\special{pn 13}%
\special{pa 800 1400}%
\special{pa 1600 600}%
\special{fp}%
\special{sh 1}%
\special{pa 1600 600}%
\special{pa 1539 633}%
\special{pa 1562 638}%
\special{pa 1567 661}%
\special{pa 1600 600}%
\special{fp}%
%
\special{pn 13}%
\special{pa 1200 1400}%
\special{pa 1190 1400}%
\special{fp}%
\special{sh 1}%
\special{pa 1190 1400}%
\special{pa 1257 1420}%
\special{pa 1243 1400}%
\special{pa 1257 1380}%
\special{pa 1190 1400}%
\special{fp}%
%
\special{pn 13}%
\special{pa 1600 1400}%
\special{pa 800 600}%
\special{fp}%
\special{sh 1}%
\special{pa 800 600}%
\special{pa 833 661}%
\special{pa 838 638}%
\special{pa 861 633}%
\special{pa 800 600}%
\special{fp}%
%
\special{pn 13}%
\special{pa 4000 1400}%
\special{pa 4800 600}%
\special{fp}%
\special{sh 1}%
\special{pa 4800 600}%
\special{pa 4739 633}%
\special{pa 4762 638}%
\special{pa 4767 661}%
\special{pa 4800 600}%
\special{fp}%
%
\special{pn 13}%
\special{pa 4800 1400}%
\special{pa 4000 600}%
\special{fp}%
\special{sh 1}%
\special{pa 4000 600}%
\special{pa 4033 661}%
\special{pa 4038 638}%
\special{pa 4061 633}%
\special{pa 4000 600}%
\special{fp}%
%
\special{pn 8}%
\special{pn 8}%
\special{pa 4800 600}%
\special{pa 4809 600}%
\special{fp}%
\special{pa 4841 608}%
\special{pa 4847 611}%
\special{fp}%
\special{pa 4875 629}%
\special{pa 4881 634}%
\special{fp}%
\special{pa 4904 658}%
\special{pa 4908 663}%
\special{fp}%
\special{pa 4927 691}%
\special{pa 4931 698}%
\special{fp}%
\special{pa 4947 729}%
\special{pa 4950 736}%
\special{fp}%
\special{pa 4963 768}%
\special{pa 4966 776}%
\special{fp}%
\special{pa 4976 811}%
\special{pa 4978 818}%
\special{fp}%
\special{pa 4986 854}%
\special{pa 4988 862}%
\special{fp}%
\special{pa 4993 899}%
\special{pa 4994 906}%
\special{fp}%
\special{pa 4998 944}%
\special{pa 4999 952}%
\special{fp}%
\special{pa 5000 991}%
\special{pa 5000 1000}%
\special{fp}%
%
\special{pn 13}%
\special{pa 4400 600}%
\special{pa 4390 600}%
\special{fp}%
\special{sh 1}%
\special{pa 4390 600}%
\special{pa 4457 620}%
\special{pa 4443 600}%
\special{pa 4457 580}%
\special{pa 4390 600}%
\special{fp}%
\put(46.0000,-9.0000){\makebox(0,0){$1$}}%
\put(42.0000,-9.0000){\makebox(0,0){$3$}}%
\put(44.0000,-5.0000){\makebox(0,0){$2$}}%
%
\special{pn 8}%
\special{pn 8}%
\special{pa 800 1400}%
\special{pa 792 1400}%
\special{fp}%
\special{pa 755 1398}%
\special{pa 748 1397}%
\special{fp}%
\special{pa 712 1390}%
\special{pa 704 1389}%
\special{fp}%
\special{pa 669 1378}%
\special{pa 661 1375}%
\special{fp}%
\special{pa 628 1361}%
\special{pa 621 1358}%
\special{fp}%
\special{pa 589 1339}%
\special{pa 582 1335}%
\special{fp}%
\special{pa 552 1314}%
\special{pa 546 1309}%
\special{fp}%
\special{pa 519 1284}%
\special{pa 513 1279}%
\special{fp}%
\special{pa 489 1251}%
\special{pa 484 1245}%
\special{fp}%
\special{pa 462 1215}%
\special{pa 459 1208}%
\special{fp}%
\special{pa 441 1176}%
\special{pa 438 1169}%
\special{fp}%
\special{pa 423 1135}%
\special{pa 421 1128}%
\special{fp}%
\special{pa 411 1092}%
\special{pa 409 1085}%
\special{fp}%
\special{pa 403 1048}%
\special{pa 402 1041}%
\special{fp}%
\special{pa 400 1004}%
\special{pa 400 996}%
\special{fp}%
\special{pa 402 959}%
\special{pa 403 952}%
\special{fp}%
\special{pa 409 916}%
\special{pa 410 908}%
\special{fp}%
\special{pa 421 872}%
\special{pa 423 865}%
\special{fp}%
\special{pa 437 831}%
\special{pa 441 824}%
\special{fp}%
\special{pa 459 792}%
\special{pa 463 785}%
\special{fp}%
\special{pa 484 755}%
\special{pa 489 749}%
\special{fp}%
\special{pa 513 722}%
\special{pa 519 715}%
\special{fp}%
\special{pa 546 691}%
\special{pa 552 686}%
\special{fp}%
\special{pa 582 665}%
\special{pa 588 660}%
\special{fp}%
\special{pa 621 643}%
\special{pa 627 639}%
\special{fp}%
\special{pa 661 624}%
\special{pa 669 622}%
\special{fp}%
\special{pa 704 612}%
\special{pa 712 610}%
\special{fp}%
\special{pa 748 603}%
\special{pa 755 602}%
\special{fp}%
\special{pa 792 600}%
\special{pa 800 600}%
\special{fp}%
%
\special{pn 8}%
\special{pn 8}%
\special{pa 4010 1400}%
\special{pa 4002 1400}%
\special{fp}%
\special{pa 3965 1397}%
\special{pa 3957 1397}%
\special{fp}%
\special{pa 3921 1390}%
\special{pa 3913 1388}%
\special{fp}%
\special{pa 3879 1378}%
\special{pa 3872 1376}%
\special{fp}%
\special{pa 3839 1362}%
\special{pa 3832 1358}%
\special{fp}%
\special{pa 3801 1341}%
\special{pa 3795 1337}%
\special{fp}%
\special{pa 3766 1317}%
\special{pa 3761 1313}%
\special{fp}%
\special{pa 3735 1290}%
\special{pa 3730 1286}%
\special{fp}%
\special{pa 3707 1262}%
\special{pa 3703 1256}%
\special{fp}%
\special{pa 3681 1228}%
\special{pa 3677 1222}%
\special{fp}%
\special{pa 3659 1192}%
\special{pa 3655 1185}%
\special{fp}%
\special{pa 3641 1153}%
\special{pa 3638 1146}%
\special{fp}%
\special{pa 3626 1112}%
\special{pa 3624 1104}%
\special{fp}%
\special{pa 3616 1068}%
\special{pa 3615 1061}%
\special{fp}%
\special{pa 3611 1024}%
\special{pa 3610 1016}%
\special{fp}%
\special{pa 3611 979}%
\special{pa 3611 971}%
\special{fp}%
\special{pa 3615 934}%
\special{pa 3617 926}%
\special{fp}%
\special{pa 3625 891}%
\special{pa 3627 884}%
\special{fp}%
\special{pa 3639 850}%
\special{pa 3642 843}%
\special{fp}%
\special{pa 3657 812}%
\special{pa 3660 806}%
\special{fp}%
\special{pa 3679 775}%
\special{pa 3683 770}%
\special{fp}%
\special{pa 3703 744}%
\special{pa 3708 737}%
\special{fp}%
\special{pa 3731 714}%
\special{pa 3736 709}%
\special{fp}%
\special{pa 3762 686}%
\special{pa 3768 682}%
\special{fp}%
\special{pa 3796 662}%
\special{pa 3802 658}%
\special{fp}%
\special{pa 3833 641}%
\special{pa 3839 638}%
\special{fp}%
\special{pa 3872 624}%
\special{pa 3879 622}%
\special{fp}%
\special{pa 3913 612}%
\special{pa 3921 610}%
\special{fp}%
\special{pa 3957 604}%
\special{pa 3965 602}%
\special{fp}%
\special{pa 4002 600}%
\special{pa 4010 600}%
\special{fp}%
%
\special{pn 8}%
\special{pn 8}%
\special{pa 800 1600}%
\special{pa 792 1600}%
\special{fp}%
\special{pa 756 1586}%
\special{pa 750 1580}%
\special{fp}%
\special{pa 722 1553}%
\special{pa 718 1547}%
\special{fp}%
\special{pa 697 1514}%
\special{pa 693 1507}%
\special{fp}%
\special{pa 677 1472}%
\special{pa 673 1464}%
\special{fp}%
\special{pa 660 1428}%
\special{pa 658 1420}%
\special{fp}%
\special{pa 646 1383}%
\special{pa 644 1375}%
\special{fp}%
\special{pa 635 1337}%
\special{pa 633 1329}%
\special{fp}%
\special{pa 626 1291}%
\special{pa 624 1282}%
\special{fp}%
\special{pa 617 1244}%
\special{pa 616 1236}%
\special{fp}%
\special{pa 611 1197}%
\special{pa 610 1189}%
\special{fp}%
\special{pa 606 1150}%
\special{pa 606 1142}%
\special{fp}%
\special{pa 603 1103}%
\special{pa 603 1095}%
\special{fp}%
\special{pa 601 1056}%
\special{pa 601 1047}%
\special{fp}%
\special{pa 600 1008}%
\special{pa 600 1000}%
\special{fp}%
\put(44.0000,-8.0000){\makebox(0,0){$\circlearrowleft$}}%
%
\special{pn 13}%
\special{pa 2400 1400}%
\special{pa 3200 600}%
\special{fp}%
\special{sh 1}%
\special{pa 3200 600}%
\special{pa 3139 633}%
\special{pa 3162 638}%
\special{pa 3167 661}%
\special{pa 3200 600}%
\special{fp}%
%
\special{pn 13}%
\special{pa 3200 1400}%
\special{pa 2400 600}%
\special{fp}%
\special{sh 1}%
\special{pa 2400 600}%
\special{pa 2433 661}%
\special{pa 2438 638}%
\special{pa 2461 633}%
\special{pa 2400 600}%
\special{fp}%
%
\special{pn 13}%
\special{pa 3400 1000}%
\special{pa 2200 1000}%
\special{fp}%
\special{sh 1}%
\special{pa 2200 1000}%
\special{pa 2267 1020}%
\special{pa 2253 1000}%
\special{pa 2267 980}%
\special{pa 2200 1000}%
\special{fp}%
%
\special{pn 8}%
\special{pn 8}%
\special{pa 3200 600}%
\special{pa 3208 600}%
\special{fp}%
\special{pa 3241 608}%
\special{pa 3248 612}%
\special{fp}%
\special{pa 3275 629}%
\special{pa 3280 634}%
\special{fp}%
\special{pa 3303 658}%
\special{pa 3308 663}%
\special{fp}%
\special{pa 3327 691}%
\special{pa 3331 697}%
\special{fp}%
\special{pa 3348 730}%
\special{pa 3351 737}%
\special{fp}%
\special{pa 3364 770}%
\special{pa 3366 777}%
\special{fp}%
\special{pa 3377 812}%
\special{pa 3379 820}%
\special{fp}%
\special{pa 3387 856}%
\special{pa 3388 864}%
\special{fp}%
\special{pa 3394 901}%
\special{pa 3395 909}%
\special{fp}%
\special{pa 3398 946}%
\special{pa 3399 954}%
\special{fp}%
\special{pa 3400 992}%
\special{pa 3400 1000}%
\special{fp}%
%
\special{pn 8}%
\special{pn 8}%
\special{pa 2400 1400}%
\special{pa 2392 1400}%
\special{fp}%
\special{pa 2355 1397}%
\special{pa 2347 1396}%
\special{fp}%
\special{pa 2311 1390}%
\special{pa 2303 1388}%
\special{fp}%
\special{pa 2268 1377}%
\special{pa 2261 1375}%
\special{fp}%
\special{pa 2228 1361}%
\special{pa 2221 1357}%
\special{fp}%
\special{pa 2189 1340}%
\special{pa 2182 1336}%
\special{fp}%
\special{pa 2152 1314}%
\special{pa 2146 1309}%
\special{fp}%
\special{pa 2119 1284}%
\special{pa 2113 1279}%
\special{fp}%
\special{pa 2089 1251}%
\special{pa 2084 1245}%
\special{fp}%
\special{pa 2063 1214}%
\special{pa 2058 1208}%
\special{fp}%
\special{pa 2041 1176}%
\special{pa 2038 1169}%
\special{fp}%
\special{pa 2024 1135}%
\special{pa 2021 1128}%
\special{fp}%
\special{pa 2011 1093}%
\special{pa 2009 1085}%
\special{fp}%
\special{pa 2003 1049}%
\special{pa 2002 1041}%
\special{fp}%
\special{pa 2000 1004}%
\special{pa 2000 996}%
\special{fp}%
\special{pa 2002 959}%
\special{pa 2003 951}%
\special{fp}%
\special{pa 2009 915}%
\special{pa 2011 907}%
\special{fp}%
\special{pa 2021 872}%
\special{pa 2024 865}%
\special{fp}%
\special{pa 2038 831}%
\special{pa 2041 824}%
\special{fp}%
\special{pa 2058 792}%
\special{pa 2063 786}%
\special{fp}%
\special{pa 2084 755}%
\special{pa 2089 749}%
\special{fp}%
\special{pa 2113 721}%
\special{pa 2119 716}%
\special{fp}%
\special{pa 2146 691}%
\special{pa 2152 686}%
\special{fp}%
\special{pa 2182 664}%
\special{pa 2189 660}%
\special{fp}%
\special{pa 2221 643}%
\special{pa 2228 639}%
\special{fp}%
\special{pa 2261 625}%
\special{pa 2269 622}%
\special{fp}%
\special{pa 2304 612}%
\special{pa 2311 610}%
\special{fp}%
\special{pa 2347 604}%
\special{pa 2355 603}%
\special{fp}%
\special{pa 2392 600}%
\special{pa 2400 600}%
\special{fp}%
%
\special{pn 4}%
\special{pa 4520 640}%
\special{pa 4290 870}%
\special{fp}%
\special{pa 4550 670}%
\special{pa 4320 900}%
\special{fp}%
\special{pa 4580 700}%
\special{pa 4350 930}%
\special{fp}%
\special{pa 4620 720}%
\special{pa 4380 960}%
\special{fp}%
\special{pa 4480 620}%
\special{pa 4260 840}%
\special{fp}%
\special{pa 4420 620}%
\special{pa 4230 810}%
\special{fp}%
\special{pa 4380 600}%
\special{pa 4200 780}%
\special{fp}%
\special{pa 4270 650}%
\special{pa 4170 750}%
\special{fp}%
%
\special{pn 8}%
\special{pa 4000 500}%
\special{pa 4400 700}%
\special{fp}%
\put(40.0000,-5.0000){\makebox(0,0)[rb]{newborn triangle}}%
%
\special{pn 8}%
\special{pn 8}%
\special{pa 2400 1600}%
\special{pa 2392 1599}%
\special{fp}%
\special{pa 2359 1588}%
\special{pa 2354 1584}%
\special{fp}%
\special{pa 2329 1561}%
\special{pa 2324 1556}%
\special{fp}%
\special{pa 2304 1526}%
\special{pa 2300 1519}%
\special{fp}%
\special{pa 2284 1488}%
\special{pa 2280 1481}%
\special{fp}%
\special{pa 2267 1448}%
\special{pa 2264 1441}%
\special{fp}%
\special{pa 2253 1406}%
\special{pa 2250 1399}%
\special{fp}%
\special{pa 2241 1364}%
\special{pa 2239 1356}%
\special{fp}%
\special{pa 2231 1320}%
\special{pa 2229 1313}%
\special{fp}%
\special{pa 2223 1277}%
\special{pa 2221 1269}%
\special{fp}%
\special{pa 2216 1233}%
\special{pa 2215 1225}%
\special{fp}%
\special{pa 2210 1188}%
\special{pa 2209 1180}%
\special{fp}%
\special{pa 2206 1144}%
\special{pa 2205 1136}%
\special{fp}%
\special{pa 2203 1099}%
\special{pa 2202 1091}%
\special{fp}%
\special{pa 2201 1053}%
\special{pa 2201 1045}%
\special{fp}%
\special{pa 2200 1008}%
\special{pa 2200 1000}%
\special{fp}%
%
\special{pn 8}%
\special{pn 8}%
\special{pa 4000 1600}%
\special{pa 3992 1600}%
\special{fp}%
\special{pa 3961 1588}%
\special{pa 3955 1585}%
\special{fp}%
\special{pa 3932 1565}%
\special{pa 3928 1560}%
\special{fp}%
\special{pa 3908 1533}%
\special{pa 3905 1527}%
\special{fp}%
\special{pa 3888 1498}%
\special{pa 3885 1492}%
\special{fp}%
\special{pa 3871 1459}%
\special{pa 3868 1452}%
\special{fp}%
\special{pa 3856 1417}%
\special{pa 3854 1409}%
\special{fp}%
\special{pa 3843 1373}%
\special{pa 3841 1366}%
\special{fp}%
\special{pa 3833 1329}%
\special{pa 3831 1322}%
\special{fp}%
\special{pa 3824 1284}%
\special{pa 3822 1276}%
\special{fp}%
\special{pa 3817 1239}%
\special{pa 3815 1231}%
\special{fp}%
\special{pa 3811 1194}%
\special{pa 3810 1186}%
\special{fp}%
\special{pa 3806 1148}%
\special{pa 3806 1140}%
\special{fp}%
\special{pa 3803 1102}%
\special{pa 3802 1094}%
\special{fp}%
\special{pa 3801 1055}%
\special{pa 3801 1047}%
\special{fp}%
\special{pa 3800 1008}%
\special{pa 3800 1000}%
\special{fp}%
%
\special{pn 8}%
\special{pn 8}%
\special{pa 1600 1400}%
\special{pa 1599 1408}%
\special{fp}%
\special{pa 1584 1440}%
\special{pa 1579 1445}%
\special{fp}%
\special{pa 1551 1469}%
\special{pa 1545 1473}%
\special{fp}%
\special{pa 1513 1491}%
\special{pa 1506 1494}%
\special{fp}%
\special{pa 1472 1508}%
\special{pa 1465 1512}%
\special{fp}%
\special{pa 1430 1523}%
\special{pa 1422 1526}%
\special{fp}%
\special{pa 1387 1536}%
\special{pa 1379 1538}%
\special{fp}%
\special{pa 1343 1547}%
\special{pa 1335 1548}%
\special{fp}%
\special{pa 1299 1557}%
\special{pa 1291 1558}%
\special{fp}%
\special{pa 1255 1564}%
\special{pa 1247 1566}%
\special{fp}%
\special{pa 1211 1571}%
\special{pa 1203 1573}%
\special{fp}%
\special{pa 1166 1578}%
\special{pa 1158 1579}%
\special{fp}%
\special{pa 1122 1583}%
\special{pa 1114 1584}%
\special{fp}%
\special{pa 1077 1587}%
\special{pa 1069 1588}%
\special{fp}%
\special{pa 1032 1592}%
\special{pa 1024 1592}%
\special{fp}%
\special{pa 988 1594}%
\special{pa 980 1595}%
\special{fp}%
\special{pa 943 1597}%
\special{pa 935 1597}%
\special{fp}%
\special{pa 898 1599}%
\special{pa 890 1599}%
\special{fp}%
\special{pa 853 1600}%
\special{pa 845 1600}%
\special{fp}%
\special{pa 808 1600}%
\special{pa 800 1600}%
\special{fp}%
%
\special{pn 8}%
\special{pn 8}%
\special{pa 3200 1400}%
\special{pa 3199 1407}%
\special{fp}%
\special{pa 3185 1438}%
\special{pa 3181 1444}%
\special{fp}%
\special{pa 3154 1467}%
\special{pa 3147 1472}%
\special{fp}%
\special{pa 3116 1489}%
\special{pa 3109 1492}%
\special{fp}%
\special{pa 3076 1507}%
\special{pa 3068 1510}%
\special{fp}%
\special{pa 3034 1522}%
\special{pa 3026 1524}%
\special{fp}%
\special{pa 2991 1535}%
\special{pa 2983 1537}%
\special{fp}%
\special{pa 2948 1546}%
\special{pa 2940 1547}%
\special{fp}%
\special{pa 2904 1555}%
\special{pa 2896 1557}%
\special{fp}%
\special{pa 2860 1564}%
\special{pa 2852 1565}%
\special{fp}%
\special{pa 2815 1571}%
\special{pa 2807 1572}%
\special{fp}%
\special{pa 2770 1577}%
\special{pa 2763 1578}%
\special{fp}%
\special{pa 2726 1583}%
\special{pa 2718 1584}%
\special{fp}%
\special{pa 2681 1587}%
\special{pa 2673 1588}%
\special{fp}%
\special{pa 2635 1591}%
\special{pa 2627 1592}%
\special{fp}%
\special{pa 2590 1594}%
\special{pa 2582 1595}%
\special{fp}%
\special{pa 2545 1597}%
\special{pa 2537 1597}%
\special{fp}%
\special{pa 2499 1598}%
\special{pa 2491 1599}%
\special{fp}%
\special{pa 2454 1600}%
\special{pa 2446 1600}%
\special{fp}%
\special{pa 2408 1600}%
\special{pa 2400 1600}%
\special{fp}%
%
\special{pn 8}%
\special{pn 8}%
\special{pa 4800 1400}%
\special{pa 4799 1407}%
\special{fp}%
\special{pa 4786 1437}%
\special{pa 4782 1443}%
\special{fp}%
\special{pa 4758 1464}%
\special{pa 4753 1468}%
\special{fp}%
\special{pa 4723 1486}%
\special{pa 4716 1489}%
\special{fp}%
\special{pa 4684 1504}%
\special{pa 4676 1507}%
\special{fp}%
\special{pa 4642 1519}%
\special{pa 4635 1522}%
\special{fp}%
\special{pa 4600 1532}%
\special{pa 4592 1534}%
\special{fp}%
\special{pa 4556 1544}%
\special{pa 4548 1546}%
\special{fp}%
\special{pa 4512 1554}%
\special{pa 4504 1555}%
\special{fp}%
\special{pa 4467 1562}%
\special{pa 4460 1564}%
\special{fp}%
\special{pa 4423 1570}%
\special{pa 4415 1571}%
\special{fp}%
\special{pa 4378 1576}%
\special{pa 4370 1577}%
\special{fp}%
\special{pa 4332 1582}%
\special{pa 4324 1583}%
\special{fp}%
\special{pa 4286 1587}%
\special{pa 4278 1587}%
\special{fp}%
\special{pa 4240 1591}%
\special{pa 4232 1591}%
\special{fp}%
\special{pa 4194 1594}%
\special{pa 4186 1595}%
\special{fp}%
\special{pa 4148 1597}%
\special{pa 4140 1597}%
\special{fp}%
\special{pa 4101 1598}%
\special{pa 4093 1599}%
\special{fp}%
\special{pa 4055 1600}%
\special{pa 4047 1600}%
\special{fp}%
\special{pa 4008 1600}%
\special{pa 4000 1600}%
\special{fp}%
%
\special{pn 4}%
\special{pa 1290 1110}%
\special{pa 1070 1330}%
\special{fp}%
\special{pa 1320 1140}%
\special{pa 1100 1360}%
\special{fp}%
\special{pa 1350 1170}%
\special{pa 1140 1380}%
\special{fp}%
\special{pa 1380 1200}%
\special{pa 1190 1390}%
\special{fp}%
\special{pa 1410 1230}%
\special{pa 1270 1370}%
\special{fp}%
\special{pa 1260 1080}%
\special{pa 1030 1310}%
\special{fp}%
\special{pa 1230 1050}%
\special{pa 1000 1280}%
\special{fp}%
%
\special{pn 4}%
\special{pa 1200 1200}%
\special{pa 1600 1600}%
\special{fp}%
\put(16.0000,-16.0000){\makebox(0,0)[lt]{vanishing}}%
\put(16.0000,-18.5000){\makebox(0,0)[lb]{triangle}}%
\end{picture}}%

%% file: newborn_triangle_-.tex
{\unitlength 0.1in%
\begin{picture}(44.0000,12.7000)(6.0000,-17.3500)%
\special{pn 13}%
\special{pa 600 1000}%
\special{pa 610 1000}%
\special{pa 615 1001}%
\special{pa 620 1001}%
\special{pa 625 1002}%
\special{pa 630 1002}%
\special{pa 640 1004}%
\special{pa 645 1006}%
\special{pa 655 1008}%
\special{pa 660 1010}%
\special{pa 665 1011}%
\special{pa 685 1019}%
\special{pa 690 1022}%
\special{pa 695 1024}%
\special{pa 700 1027}%
\special{pa 705 1029}%
\special{pa 725 1041}%
\special{pa 730 1045}%
\special{pa 740 1051}%
\special{pa 750 1059}%
\special{pa 755 1062}%
\special{pa 780 1082}%
\special{pa 785 1087}%
\special{pa 790 1091}%
\special{pa 795 1096}%
\special{pa 800 1100}%
\special{pa 805 1105}%
\special{pa 810 1109}%
\special{pa 820 1119}%
\special{pa 825 1123}%
\special{pa 860 1158}%
\special{pa 865 1164}%
\special{pa 885 1184}%
\special{pa 890 1190}%
\special{pa 910 1210}%
\special{pa 915 1216}%
\special{pa 935 1236}%
\special{pa 940 1242}%
\special{pa 975 1277}%
\special{pa 980 1281}%
\special{pa 990 1291}%
\special{pa 995 1295}%
\special{pa 1000 1300}%
\special{pa 1005 1304}%
\special{pa 1010 1309}%
\special{pa 1015 1313}%
\special{pa 1020 1318}%
\special{pa 1045 1338}%
\special{pa 1050 1341}%
\special{pa 1060 1349}%
\special{pa 1070 1355}%
\special{pa 1075 1359}%
\special{pa 1095 1371}%
\special{pa 1100 1373}%
\special{pa 1105 1376}%
\special{pa 1110 1378}%
\special{pa 1115 1381}%
\special{pa 1135 1389}%
\special{pa 1140 1390}%
\special{pa 1145 1392}%
\special{pa 1155 1394}%
\special{pa 1160 1396}%
\special{pa 1170 1398}%
\special{pa 1175 1398}%
\special{pa 1180 1399}%
\special{pa 1185 1399}%
\special{pa 1190 1400}%
\special{pa 1200 1400}%
\special{fp}%
\special{pn 13}%
\special{pa 1200 1400}%
\special{pa 1210 1400}%
\special{pa 1215 1399}%
\special{pa 1220 1399}%
\special{pa 1225 1398}%
\special{pa 1230 1398}%
\special{pa 1240 1396}%
\special{pa 1245 1394}%
\special{pa 1255 1392}%
\special{pa 1260 1390}%
\special{pa 1265 1389}%
\special{pa 1285 1381}%
\special{pa 1290 1378}%
\special{pa 1295 1376}%
\special{pa 1300 1373}%
\special{pa 1305 1371}%
\special{pa 1325 1359}%
\special{pa 1330 1355}%
\special{pa 1340 1349}%
\special{pa 1350 1341}%
\special{pa 1355 1338}%
\special{pa 1380 1318}%
\special{pa 1385 1313}%
\special{pa 1390 1309}%
\special{pa 1395 1304}%
\special{pa 1400 1300}%
\special{pa 1405 1295}%
\special{pa 1410 1291}%
\special{pa 1420 1281}%
\special{pa 1425 1277}%
\special{pa 1460 1242}%
\special{pa 1465 1236}%
\special{pa 1485 1216}%
\special{pa 1490 1210}%
\special{pa 1510 1190}%
\special{pa 1515 1184}%
\special{pa 1535 1164}%
\special{pa 1540 1158}%
\special{pa 1575 1123}%
\special{pa 1580 1119}%
\special{pa 1590 1109}%
\special{pa 1595 1105}%
\special{pa 1600 1100}%
\special{pa 1605 1096}%
\special{pa 1610 1091}%
\special{pa 1615 1087}%
\special{pa 1620 1082}%
\special{pa 1645 1062}%
\special{pa 1650 1059}%
\special{pa 1660 1051}%
\special{pa 1670 1045}%
\special{pa 1675 1041}%
\special{pa 1695 1029}%
\special{pa 1700 1027}%
\special{pa 1705 1024}%
\special{pa 1710 1022}%
\special{pa 1715 1019}%
\special{pa 1735 1011}%
\special{pa 1740 1010}%
\special{pa 1745 1008}%
\special{pa 1755 1006}%
\special{pa 1760 1004}%
\special{pa 1770 1002}%
\special{pa 1775 1002}%
\special{pa 1780 1001}%
\special{pa 1785 1001}%
\special{pa 1790 1000}%
\special{pa 1800 1000}%
\special{fp}%
\put(12.0000,-18.0000){\makebox(0,0){$h(-\epsilon)$}}%
\put(28.0000,-18.0000){\makebox(0,0){$h(0)$}}%
\put(44.0000,-18.0000){\makebox(0,0){$h(\epsilon)$}}%
\special{pn 13}%
\special{pa 4400 600}%
\special{pa 4410 600}%
\special{pa 4415 601}%
\special{pa 4420 601}%
\special{pa 4425 602}%
\special{pa 4430 602}%
\special{pa 4440 604}%
\special{pa 4445 606}%
\special{pa 4455 608}%
\special{pa 4460 610}%
\special{pa 4465 611}%
\special{pa 4485 619}%
\special{pa 4490 622}%
\special{pa 4495 624}%
\special{pa 4500 627}%
\special{pa 4505 629}%
\special{pa 4525 641}%
\special{pa 4530 645}%
\special{pa 4540 651}%
\special{pa 4550 659}%
\special{pa 4555 662}%
\special{pa 4580 682}%
\special{pa 4585 687}%
\special{pa 4590 691}%
\special{pa 4595 696}%
\special{pa 4600 700}%
\special{pa 4605 705}%
\special{pa 4610 709}%
\special{pa 4620 719}%
\special{pa 4625 723}%
\special{pa 4660 758}%
\special{pa 4665 764}%
\special{pa 4685 784}%
\special{pa 4690 790}%
\special{pa 4710 810}%
\special{pa 4715 816}%
\special{pa 4735 836}%
\special{pa 4740 842}%
\special{pa 4775 877}%
\special{pa 4780 881}%
\special{pa 4790 891}%
\special{pa 4795 895}%
\special{pa 4800 900}%
\special{pa 4805 904}%
\special{pa 4810 909}%
\special{pa 4815 913}%
\special{pa 4820 918}%
\special{pa 4845 938}%
\special{pa 4850 941}%
\special{pa 4860 949}%
\special{pa 4870 955}%
\special{pa 4875 959}%
\special{pa 4895 971}%
\special{pa 4900 973}%
\special{pa 4905 976}%
\special{pa 4910 978}%
\special{pa 4915 981}%
\special{pa 4935 989}%
\special{pa 4940 990}%
\special{pa 4945 992}%
\special{pa 4955 994}%
\special{pa 4960 996}%
\special{pa 4970 998}%
\special{pa 4975 998}%
\special{pa 4980 999}%
\special{pa 4985 999}%
\special{pa 4990 1000}%
\special{pa 5000 1000}%
\special{fp}%
\special{pn 13}%
\special{pa 3800 1000}%
\special{pa 3810 1000}%
\special{pa 3815 999}%
\special{pa 3820 999}%
\special{pa 3825 998}%
\special{pa 3830 998}%
\special{pa 3840 996}%
\special{pa 3845 994}%
\special{pa 3855 992}%
\special{pa 3860 990}%
\special{pa 3865 989}%
\special{pa 3885 981}%
\special{pa 3890 978}%
\special{pa 3895 976}%
\special{pa 3900 973}%
\special{pa 3905 971}%
\special{pa 3925 959}%
\special{pa 3930 955}%
\special{pa 3940 949}%
\special{pa 3950 941}%
\special{pa 3955 938}%
\special{pa 3980 918}%
\special{pa 3985 913}%
\special{pa 3990 909}%
\special{pa 3995 904}%
\special{pa 4000 900}%
\special{pa 4005 895}%
\special{pa 4010 891}%
\special{pa 4020 881}%
\special{pa 4025 877}%
\special{pa 4060 842}%
\special{pa 4065 836}%
\special{pa 4085 816}%
\special{pa 4090 810}%
\special{pa 4110 790}%
\special{pa 4115 784}%
\special{pa 4135 764}%
\special{pa 4140 758}%
\special{pa 4175 723}%
\special{pa 4180 719}%
\special{pa 4190 709}%
\special{pa 4195 705}%
\special{pa 4200 700}%
\special{pa 4205 696}%
\special{pa 4210 691}%
\special{pa 4215 687}%
\special{pa 4220 682}%
\special{pa 4245 662}%
\special{pa 4250 659}%
\special{pa 4260 651}%
\special{pa 4270 645}%
\special{pa 4275 641}%
\special{pa 4295 629}%
\special{pa 4300 627}%
\special{pa 4305 624}%
\special{pa 4310 622}%
\special{pa 4315 619}%
\special{pa 4335 611}%
\special{pa 4340 610}%
\special{pa 4345 608}%
\special{pa 4355 606}%
\special{pa 4360 604}%
\special{pa 4370 602}%
\special{pa 4375 602}%
\special{pa 4380 601}%
\special{pa 4385 601}%
\special{pa 4390 600}%
\special{pa 4400 600}%
\special{fp}%
%
\special{pn 13}%
\special{pa 800 1400}%
\special{pa 1600 600}%
\special{fp}%
\special{sh 1}%
\special{pa 1600 600}%
\special{pa 1539 633}%
\special{pa 1562 638}%
\special{pa 1567 661}%
\special{pa 1600 600}%
\special{fp}%
%
\special{pn 13}%
\special{pa 800 600}%
\special{pa 1600 1400}%
\special{fp}%
\special{sh 1}%
\special{pa 1600 1400}%
\special{pa 1567 1339}%
\special{pa 1562 1362}%
\special{pa 1539 1367}%
\special{pa 1600 1400}%
\special{fp}%
%
\special{pn 13}%
\special{pa 1200 1400}%
\special{pa 1190 1400}%
\special{fp}%
\special{sh 1}%
\special{pa 1190 1400}%
\special{pa 1257 1420}%
\special{pa 1243 1400}%
\special{pa 1257 1380}%
\special{pa 1190 1400}%
\special{fp}%
%
\special{pn 13}%
\special{pa 4000 1400}%
\special{pa 4800 600}%
\special{fp}%
\special{sh 1}%
\special{pa 4800 600}%
\special{pa 4739 633}%
\special{pa 4762 638}%
\special{pa 4767 661}%
\special{pa 4800 600}%
\special{fp}%
%
\special{pn 13}%
\special{pa 4000 600}%
\special{pa 4800 1400}%
\special{fp}%
\special{sh 1}%
\special{pa 4800 1400}%
\special{pa 4767 1339}%
\special{pa 4762 1362}%
\special{pa 4739 1367}%
\special{pa 4800 1400}%
\special{fp}%
%
\special{pn 13}%
\special{pa 4400 600}%
\special{pa 4390 600}%
\special{fp}%
\special{sh 1}%
\special{pa 4390 600}%
\special{pa 4457 620}%
\special{pa 4443 600}%
\special{pa 4457 580}%
\special{pa 4390 600}%
\special{fp}%
\put(44.0000,-8.0000){\makebox(0,0){$\circlearrowleft$}}%
\put(46.0000,-9.0000){\makebox(0,0){$1$}}%
\put(42.0000,-9.0000){\makebox(0,0){$3$}}%
\put(44.0000,-5.3000){\makebox(0,0){$2$}}%
%
\special{pn 8}%
\special{pn 8}%
\special{pa 1600 600}%
\special{pa 1609 600}%
\special{fp}%
\special{pa 1645 610}%
\special{pa 1652 614}%
\special{fp}%
\special{pa 1684 637}%
\special{pa 1690 643}%
\special{fp}%
\special{pa 1715 673}%
\special{pa 1720 680}%
\special{fp}%
\special{pa 1740 714}%
\special{pa 1743 722}%
\special{fp}%
\special{pa 1759 757}%
\special{pa 1762 765}%
\special{fp}%
\special{pa 1774 803}%
\special{pa 1776 811}%
\special{fp}%
\special{pa 1785 849}%
\special{pa 1787 857}%
\special{fp}%
\special{pa 1793 896}%
\special{pa 1794 904}%
\special{fp}%
\special{pa 1798 943}%
\special{pa 1799 952}%
\special{fp}%
\special{pa 1800 991}%
\special{pa 1800 1000}%
\special{fp}%
%
\special{pn 8}%
\special{pn 8}%
\special{pa 600 1000}%
\special{pa 600 991}%
\special{fp}%
\special{pa 601 952}%
\special{pa 602 944}%
\special{fp}%
\special{pa 606 905}%
\special{pa 607 896}%
\special{fp}%
\special{pa 613 858}%
\special{pa 615 849}%
\special{fp}%
\special{pa 624 811}%
\special{pa 626 803}%
\special{fp}%
\special{pa 638 765}%
\special{pa 641 758}%
\special{fp}%
\special{pa 656 721}%
\special{pa 660 714}%
\special{fp}%
\special{pa 680 680}%
\special{pa 684 674}%
\special{fp}%
\special{pa 709 644}%
\special{pa 715 638}%
\special{fp}%
\special{pa 747 615}%
\special{pa 754 611}%
\special{fp}%
\special{pa 792 600}%
\special{pa 800 600}%
\special{fp}%
%
\special{pn 8}%
\special{pn 8}%
\special{pa 1600 1400}%
\special{pa 1600 1408}%
\special{fp}%
\special{pa 1590 1444}%
\special{pa 1587 1451}%
\special{fp}%
\special{pa 1566 1481}%
\special{pa 1560 1487}%
\special{fp}%
\special{pa 1532 1512}%
\special{pa 1525 1516}%
\special{fp}%
\special{pa 1493 1536}%
\special{pa 1486 1540}%
\special{fp}%
\special{pa 1451 1556}%
\special{pa 1443 1558}%
\special{fp}%
\special{pa 1408 1571}%
\special{pa 1400 1573}%
\special{fp}%
\special{pa 1364 1583}%
\special{pa 1356 1584}%
\special{fp}%
\special{pa 1319 1591}%
\special{pa 1311 1592}%
\special{fp}%
\special{pa 1273 1597}%
\special{pa 1265 1597}%
\special{fp}%
\special{pa 1227 1599}%
\special{pa 1219 1600}%
\special{fp}%
\special{pa 1181 1600}%
\special{pa 1173 1600}%
\special{fp}%
\special{pa 1135 1597}%
\special{pa 1127 1597}%
\special{fp}%
\special{pa 1089 1592}%
\special{pa 1081 1591}%
\special{fp}%
\special{pa 1044 1584}%
\special{pa 1036 1582}%
\special{fp}%
\special{pa 1000 1573}%
\special{pa 992 1571}%
\special{fp}%
\special{pa 956 1558}%
\special{pa 949 1556}%
\special{fp}%
\special{pa 914 1540}%
\special{pa 907 1536}%
\special{fp}%
\special{pa 874 1516}%
\special{pa 868 1512}%
\special{fp}%
\special{pa 840 1487}%
\special{pa 834 1481}%
\special{fp}%
\special{pa 813 1451}%
\special{pa 810 1443}%
\special{fp}%
\special{pa 800 1408}%
\special{pa 800 1400}%
\special{fp}%
%
\special{pn 8}%
\special{pn 8}%
\special{pa 4800 600}%
\special{pa 4809 600}%
\special{fp}%
\special{pa 4841 608}%
\special{pa 4847 611}%
\special{fp}%
\special{pa 4875 629}%
\special{pa 4881 634}%
\special{fp}%
\special{pa 4904 658}%
\special{pa 4908 663}%
\special{fp}%
\special{pa 4927 691}%
\special{pa 4931 698}%
\special{fp}%
\special{pa 4947 729}%
\special{pa 4950 736}%
\special{fp}%
\special{pa 4963 768}%
\special{pa 4966 776}%
\special{fp}%
\special{pa 4976 811}%
\special{pa 4978 818}%
\special{fp}%
\special{pa 4986 854}%
\special{pa 4988 862}%
\special{fp}%
\special{pa 4993 899}%
\special{pa 4994 906}%
\special{fp}%
\special{pa 4998 944}%
\special{pa 4999 952}%
\special{fp}%
\special{pa 5000 991}%
\special{pa 5000 1000}%
\special{fp}%
%
\special{pn 8}%
\special{pn 8}%
\special{pa 3800 1000}%
\special{pa 3800 992}%
\special{fp}%
\special{pa 3801 953}%
\special{pa 3802 945}%
\special{fp}%
\special{pa 3805 907}%
\special{pa 3806 899}%
\special{fp}%
\special{pa 3812 861}%
\special{pa 3814 853}%
\special{fp}%
\special{pa 3822 817}%
\special{pa 3824 809}%
\special{fp}%
\special{pa 3835 774}%
\special{pa 3838 767}%
\special{fp}%
\special{pa 3850 735}%
\special{pa 3853 728}%
\special{fp}%
\special{pa 3870 696}%
\special{pa 3874 690}%
\special{fp}%
\special{pa 3890 666}%
\special{pa 3894 660}%
\special{fp}%
\special{pa 3917 636}%
\special{pa 3922 631}%
\special{fp}%
\special{pa 3950 613}%
\special{pa 3957 609}%
\special{fp}%
\special{pa 3992 600}%
\special{pa 4000 600}%
\special{fp}%
%
\special{pn 8}%
\special{pn 8}%
\special{pa 4800 1400}%
\special{pa 4800 1408}%
\special{fp}%
\special{pa 4792 1441}%
\special{pa 4789 1447}%
\special{fp}%
\special{pa 4772 1474}%
\special{pa 4767 1480}%
\special{fp}%
\special{pa 4744 1502}%
\special{pa 4740 1506}%
\special{fp}%
\special{pa 4714 1524}%
\special{pa 4708 1528}%
\special{fp}%
\special{pa 4678 1544}%
\special{pa 4671 1547}%
\special{fp}%
\special{pa 4641 1560}%
\special{pa 4635 1562}%
\special{fp}%
\special{pa 4602 1573}%
\special{pa 4595 1575}%
\special{fp}%
\special{pa 4560 1583}%
\special{pa 4553 1585}%
\special{fp}%
\special{pa 4517 1591}%
\special{pa 4509 1592}%
\special{fp}%
\special{pa 4473 1597}%
\special{pa 4465 1597}%
\special{fp}%
\special{pa 4428 1599}%
\special{pa 4420 1600}%
\special{fp}%
\special{pa 4382 1600}%
\special{pa 4374 1600}%
\special{fp}%
\special{pa 4337 1597}%
\special{pa 4329 1597}%
\special{fp}%
\special{pa 4294 1593}%
\special{pa 4286 1592}%
\special{fp}%
\special{pa 4252 1586}%
\special{pa 4245 1584}%
\special{fp}%
\special{pa 4211 1576}%
\special{pa 4203 1574}%
\special{fp}%
\special{pa 4169 1563}%
\special{pa 4161 1560}%
\special{fp}%
\special{pa 4129 1547}%
\special{pa 4122 1544}%
\special{fp}%
\special{pa 4092 1528}%
\special{pa 4086 1524}%
\special{fp}%
\special{pa 4059 1505}%
\special{pa 4054 1501}%
\special{fp}%
\special{pa 4031 1478}%
\special{pa 4027 1473}%
\special{fp}%
\special{pa 4010 1444}%
\special{pa 4007 1439}%
\special{fp}%
\special{pa 4000 1408}%
\special{pa 4000 1400}%
\special{fp}%
%
\special{pn 13}%
\special{pa 2400 1400}%
\special{pa 3200 600}%
\special{fp}%
\special{sh 1}%
\special{pa 3200 600}%
\special{pa 3139 633}%
\special{pa 3162 638}%
\special{pa 3167 661}%
\special{pa 3200 600}%
\special{fp}%
%
\special{pn 13}%
\special{pa 2400 600}%
\special{pa 3200 1400}%
\special{fp}%
\special{sh 1}%
\special{pa 3200 1400}%
\special{pa 3167 1339}%
\special{pa 3162 1362}%
\special{pa 3139 1367}%
\special{pa 3200 1400}%
\special{fp}%
%
\special{pn 13}%
\special{pa 3400 1000}%
\special{pa 2200 1000}%
\special{fp}%
\special{sh 1}%
\special{pa 2200 1000}%
\special{pa 2267 1020}%
\special{pa 2253 1000}%
\special{pa 2267 980}%
\special{pa 2200 1000}%
\special{fp}%
%
\special{pn 8}%
\special{pn 8}%
\special{pa 3200 600}%
\special{pa 3208 600}%
\special{fp}%
\special{pa 3241 608}%
\special{pa 3248 612}%
\special{fp}%
\special{pa 3275 629}%
\special{pa 3280 634}%
\special{fp}%
\special{pa 3303 658}%
\special{pa 3308 663}%
\special{fp}%
\special{pa 3327 691}%
\special{pa 3331 697}%
\special{fp}%
\special{pa 3348 730}%
\special{pa 3351 737}%
\special{fp}%
\special{pa 3364 770}%
\special{pa 3366 777}%
\special{fp}%
\special{pa 3377 812}%
\special{pa 3379 820}%
\special{fp}%
\special{pa 3387 856}%
\special{pa 3388 864}%
\special{fp}%
\special{pa 3394 901}%
\special{pa 3395 909}%
\special{fp}%
\special{pa 3398 946}%
\special{pa 3399 954}%
\special{fp}%
\special{pa 3400 992}%
\special{pa 3400 1000}%
\special{fp}%
%
\special{pn 8}%
\special{pn 8}%
\special{pa 2200 1000}%
\special{pa 2200 992}%
\special{fp}%
\special{pa 2201 955}%
\special{pa 2202 947}%
\special{fp}%
\special{pa 2205 910}%
\special{pa 2206 903}%
\special{fp}%
\special{pa 2212 867}%
\special{pa 2213 859}%
\special{fp}%
\special{pa 2221 824}%
\special{pa 2223 816}%
\special{fp}%
\special{pa 2232 782}%
\special{pa 2235 775}%
\special{fp}%
\special{pa 2247 742}%
\special{pa 2250 735}%
\special{fp}%
\special{pa 2266 703}%
\special{pa 2270 696}%
\special{fp}%
\special{pa 2290 666}%
\special{pa 2294 661}%
\special{fp}%
\special{pa 2319 634}%
\special{pa 2324 630}%
\special{fp}%
\special{pa 2353 611}%
\special{pa 2359 608}%
\special{fp}%
\special{pa 2392 600}%
\special{pa 2400 600}%
\special{fp}%
%
\special{pn 8}%
\special{pn 8}%
\special{pa 3200 1400}%
\special{pa 3200 1408}%
\special{fp}%
\special{pa 3191 1442}%
\special{pa 3188 1449}%
\special{fp}%
\special{pa 3169 1477}%
\special{pa 3164 1482}%
\special{fp}%
\special{pa 3140 1506}%
\special{pa 3134 1510}%
\special{fp}%
\special{pa 3105 1529}%
\special{pa 3098 1533}%
\special{fp}%
\special{pa 3065 1550}%
\special{pa 3058 1553}%
\special{fp}%
\special{pa 3023 1566}%
\special{pa 3016 1568}%
\special{fp}%
\special{pa 2982 1578}%
\special{pa 2974 1580}%
\special{fp}%
\special{pa 2938 1588}%
\special{pa 2930 1589}%
\special{fp}%
\special{pa 2893 1594}%
\special{pa 2886 1595}%
\special{fp}%
\special{pa 2849 1599}%
\special{pa 2841 1599}%
\special{fp}%
\special{pa 2803 1600}%
\special{pa 2795 1600}%
\special{fp}%
\special{pa 2757 1599}%
\special{pa 2750 1598}%
\special{fp}%
\special{pa 2713 1595}%
\special{pa 2705 1594}%
\special{fp}%
\special{pa 2668 1589}%
\special{pa 2660 1587}%
\special{fp}%
\special{pa 2625 1580}%
\special{pa 2617 1578}%
\special{fp}%
\special{pa 2583 1568}%
\special{pa 2575 1565}%
\special{fp}%
\special{pa 2541 1552}%
\special{pa 2534 1550}%
\special{fp}%
\special{pa 2501 1533}%
\special{pa 2495 1529}%
\special{fp}%
\special{pa 2465 1509}%
\special{pa 2459 1505}%
\special{fp}%
\special{pa 2434 1481}%
\special{pa 2430 1476}%
\special{fp}%
\special{pa 2412 1449}%
\special{pa 2409 1442}%
\special{fp}%
\special{pa 2400 1408}%
\special{pa 2400 1400}%
\special{fp}%
%
\special{pn 4}%
\special{pa 4480 620}%
\special{pa 4260 840}%
\special{fp}%
\special{pa 4520 640}%
\special{pa 4290 870}%
\special{fp}%
\special{pa 4550 670}%
\special{pa 4320 900}%
\special{fp}%
\special{pa 4580 700}%
\special{pa 4350 930}%
\special{fp}%
\special{pa 4620 720}%
\special{pa 4380 960}%
\special{fp}%
\special{pa 4420 620}%
\special{pa 4230 810}%
\special{fp}%
\special{pa 4380 600}%
\special{pa 4200 780}%
\special{fp}%
\special{pa 4270 650}%
\special{pa 4170 750}%
\special{fp}%
\end{picture}}%

%% file: compatible_ori.tex
{\unitlength 0.1in%
\begin{picture}(43.8200,6.3500)(6.4500,-12.3500)%
%
\special{pn 13}%
\special{pa 2000 800}%
\special{pa 2000 600}%
\special{fp}%
\special{sh 1}%
\special{pa 2000 600}%
\special{pa 1980 667}%
\special{pa 2000 653}%
\special{pa 2020 667}%
\special{pa 2000 600}%
\special{fp}%
%
\special{pn 13}%
\special{pa 2000 1200}%
\special{pa 2000 1000}%
\special{fp}%
%
\special{pn 8}%
\special{pa 2000 1000}%
\special{pa 2000 800}%
\special{dt 0.045}%
%
\special{pn 8}%
\special{pa 800 900}%
\special{pa 1400 900}%
\special{dt 0.045}%
%
\special{pn 13}%
\special{pa 2000 800}%
\special{pa 2600 800}%
\special{fp}%
%
\special{pn 13}%
\special{pa 2600 600}%
\special{pa 2600 800}%
\special{fp}%
%
\special{pn 13}%
\special{pa 2600 1000}%
\special{pa 2600 1200}%
\special{fp}%
\special{sh 1}%
\special{pa 2600 1200}%
\special{pa 2620 1133}%
\special{pa 2600 1147}%
\special{pa 2580 1133}%
\special{pa 2600 1200}%
\special{fp}%
%
\special{pn 13}%
\special{pa 2000 1000}%
\special{pa 2600 1000}%
\special{fp}%
\put(8.0000,-13.0000){\makebox(0,0){$C_0$}}%
\put(14.0000,-13.0000){\makebox(0,0){$C_1$}}%
%
\special{pn 8}%
\special{pa 2600 1000}%
\special{pa 2600 800}%
\special{dt 0.045}%
\put(11.0000,-8.0000){\makebox(0,0){$\Gamma$}}%
%
\special{pn 13}%
\special{pa 4400 1000}%
\special{pa 4400 1200}%
\special{fp}%
\special{sh 1}%
\special{pa 4400 1200}%
\special{pa 4420 1133}%
\special{pa 4400 1147}%
\special{pa 4380 1133}%
\special{pa 4400 1200}%
\special{fp}%
%
\special{pn 13}%
\special{pa 4400 600}%
\special{pa 4400 800}%
\special{fp}%
%
\special{pn 8}%
\special{pa 4400 1000}%
\special{pa 4400 800}%
\special{dt 0.045}%
%
\special{pn 8}%
\special{pa 3200 900}%
\special{pa 3800 900}%
\special{dt 0.045}%
%
\special{pn 13}%
\special{pa 5000 1200}%
\special{pa 5000 1000}%
\special{fp}%
%
\special{pn 8}%
\special{pa 5000 1000}%
\special{pa 5000 800}%
\special{dt 0.045}%
%
\special{pn 13}%
\special{pa 5000 800}%
\special{pa 5000 600}%
\special{fp}%
\special{sh 1}%
\special{pa 5000 600}%
\special{pa 4980 667}%
\special{pa 5000 653}%
\special{pa 5020 667}%
\special{pa 5000 600}%
\special{fp}%
\put(35.0000,-8.0000){\makebox(0,0){$\Gamma$}}%
\put(32.0000,-13.0000){\makebox(0,0){$C_0$}}%
\put(38.0000,-13.0000){\makebox(0,0){$C_1$}}%
%
\special{pn 13}%
\special{pa 4400 1000}%
\special{pa 5000 1000}%
\special{fp}%
%
\special{pn 13}%
\special{pa 4400 800}%
\special{pa 5000 800}%
\special{fp}%
%
\special{pn 13}%
\special{pa 800 1200}%
\special{pa 800 600}%
\special{fp}%
\special{sh 1}%
\special{pa 800 600}%
\special{pa 780 667}%
\special{pa 800 653}%
\special{pa 820 667}%
\special{pa 800 600}%
\special{fp}%
%
\special{pn 13}%
\special{pa 1400 600}%
\special{pa 1400 1200}%
\special{fp}%
\special{sh 1}%
\special{pa 1400 1200}%
\special{pa 1420 1133}%
\special{pa 1400 1147}%
\special{pa 1380 1133}%
\special{pa 1400 1200}%
\special{fp}%
\put(47.0000,-13.0000){\makebox(0,0){$C_0+_{\Gamma}C_1$}}%
\put(17.0000,-9.0000){\makebox(0,0){$\rightsquigarrow$}}%
%
\special{pn 13}%
\special{pa 3800 1200}%
\special{pa 3800 600}%
\special{fp}%
\special{sh 1}%
\special{pa 3800 600}%
\special{pa 3780 667}%
\special{pa 3800 653}%
\special{pa 3820 667}%
\special{pa 3800 600}%
\special{fp}%
%
\special{pn 13}%
\special{pa 3200 600}%
\special{pa 3200 1200}%
\special{fp}%
\special{sh 1}%
\special{pa 3200 1200}%
\special{pa 3220 1133}%
\special{pa 3200 1147}%
\special{pa 3180 1133}%
\special{pa 3200 1200}%
\special{fp}%
\put(23.0000,-13.0000){\makebox(0,0){$C_0+_{\Gamma}C_1$}}%
\put(41.0000,-9.0000){\makebox(0,0){$\rightsquigarrow$}}%
\end{picture}}%

%% file: strange_connected_sum.tex
\unitlength 0.1in
\begin{picture}( 44.0000,  8.0000)(  8.0000,-12.0000)
%
{\color[named]{Black}{%
\special{pn 13}%
\special{ar 1200 800 400 400  3.3865713  3.1415927}%
}}%
%
{\color[named]{Black}{%
\special{pn 13}%
\special{ar 2600 800 400 400  6.2831853  6.0382066}%
}}%
%
{\color[named]{Black}{%
\special{pn 13}%
\special{ar 1600 600 400 200  6.2831853  1.5707963}%
}}%
%
{\color[named]{Black}{%
\special{pn 13}%
\special{ar 1900 600 100 100  3.1415927  6.2831853}%
}}%
%
{\color[named]{Black}{%
\special{pn 13}%
\special{ar 2200 600 400 200  1.5707963  3.1415927}%
}}%
%
{\color[named]{Black}{%
\special{pn 13}%
\special{pa 1600 800}%
\special{pa 800 800}%
\special{fp}%
}}%
%
{\color[named]{Black}{%
\special{pn 13}%
\special{pa 2200 800}%
\special{pa 3000 800}%
\special{fp}%
}}%
%
{\color[named]{Black}{%
\special{pn 13}%
\special{ar 1600 600 350 100  6.2831853  1.5707963}%
}}%
%
{\color[named]{Black}{%
\special{pn 13}%
\special{ar 2200 600 350 100  1.5707963  3.1415927}%
}}%
%
{\color[named]{Black}{%
\special{pn 13}%
\special{ar 1900 600 50 50  3.1415927  6.2831853}%
}}%
%
{\color[named]{Black}{%
\special{pn 13}%
\special{pa 816 700}%
\special{pa 1600 700}%
\special{fp}%
}}%
%
{\color[named]{Black}{%
\special{pn 13}%
\special{pa 2200 700}%
\special{pa 2986 700}%
\special{fp}%
}}%
%
{\color[named]{Black}{%
\special{pn 8}%
\special{pa 1975 600}%
\special{pa 1975 607}%
\special{pa 1974 608}%
\special{fp}%
\special{pa 1960 642}%
\special{pa 1960 643}%
\special{pa 1959 643}%
\special{pa 1958 644}%
\special{pa 1958 645}%
\special{pa 1955 648}%
\special{fp}%
\special{pa 1927 674}%
\special{pa 1927 674}%
\special{pa 1926 674}%
\special{pa 1924 675}%
\special{pa 1921 678}%
\special{pa 1920 678}%
\special{fp}%
\special{pa 1887 697}%
\special{pa 1887 697}%
\special{pa 1885 697}%
\special{pa 1884 698}%
\special{pa 1882 699}%
\special{pa 1881 699}%
\special{pa 1880 700}%
\special{fp}%
\special{pa 1843 714}%
\special{pa 1842 715}%
\special{pa 1840 715}%
\special{pa 1839 716}%
\special{pa 1837 716}%
\special{pa 1836 717}%
\special{fp}%
\special{pa 1798 727}%
\special{pa 1797 728}%
\special{pa 1795 728}%
\special{pa 1793 729}%
\special{pa 1791 729}%
\special{pa 1790 729}%
\special{fp}%
\special{pa 1752 737}%
\special{pa 1751 737}%
\special{pa 1749 738}%
\special{pa 1745 738}%
\special{pa 1744 738}%
\special{fp}%
\special{pa 1705 744}%
\special{pa 1703 744}%
\special{pa 1700 745}%
\special{pa 1697 745}%
\special{fp}%
\special{pa 1657 748}%
\special{pa 1654 748}%
\special{pa 1652 749}%
\special{pa 1649 749}%
\special{fp}%
\special{pa 1609 750}%
\special{pa 1601 750}%
\special{fp}%
}}%
%
{\color[named]{Black}{%
\special{pn 8}%
\special{pa 2200 750}%
\special{pa 2192 750}%
\special{fp}%
\special{pa 2155 749}%
\special{pa 2147 749}%
\special{fp}%
\special{pa 2111 746}%
\special{pa 2109 746}%
\special{pa 2107 745}%
\special{pa 2103 745}%
\special{fp}%
\special{pa 2067 740}%
\special{pa 2063 740}%
\special{pa 2062 739}%
\special{pa 2060 739}%
\special{fp}%
\special{pa 2024 733}%
\special{pa 2023 732}%
\special{pa 2020 732}%
\special{pa 2019 731}%
\special{pa 2017 731}%
\special{fp}%
\special{pa 1983 722}%
\special{pa 1981 722}%
\special{pa 1980 721}%
\special{pa 1977 721}%
\special{pa 1976 720}%
\special{pa 1976 720}%
\special{fp}%
\special{pa 1943 709}%
\special{pa 1942 709}%
\special{pa 1941 708}%
\special{pa 1939 708}%
\special{pa 1938 707}%
\special{pa 1937 707}%
\special{pa 1936 706}%
\special{fp}%
\special{pa 1904 692}%
\special{pa 1904 692}%
\special{pa 1903 692}%
\special{pa 1902 691}%
\special{pa 1901 691}%
\special{pa 1899 689}%
\special{pa 1898 689}%
\special{pa 1898 689}%
\special{fp}%
\special{pa 1869 670}%
\special{pa 1869 670}%
\special{pa 1868 670}%
\special{pa 1865 667}%
\special{pa 1864 667}%
\special{pa 1864 666}%
\special{pa 1864 666}%
\special{fp}%
\special{pa 1840 642}%
\special{pa 1840 642}%
\special{pa 1840 641}%
\special{pa 1839 641}%
\special{pa 1839 640}%
\special{pa 1837 638}%
\special{pa 1837 637}%
\special{pa 1836 637}%
\special{fp}%
\special{pa 1825 607}%
\special{pa 1825 607}%
\special{pa 1825 600}%
\special{fp}%
}}%
%
{\color[named]{Black}{%
\special{pn 8}%
\special{pa 1825 600}%
\special{pa 1825 592}%
\special{fp}%
\special{pa 1835 563}%
\special{pa 1835 562}%
\special{pa 1836 562}%
\special{pa 1836 560}%
\special{pa 1837 560}%
\special{pa 1837 559}%
\special{pa 1838 558}%
\special{pa 1838 557}%
\special{fp}%
\special{pa 1860 537}%
\special{pa 1860 537}%
\special{pa 1860 536}%
\special{pa 1861 536}%
\special{pa 1862 535}%
\special{pa 1863 535}%
\special{pa 1864 534}%
\special{pa 1865 534}%
\special{pa 1865 534}%
\special{fp}%
\special{pa 1897 525}%
\special{pa 1904 525}%
\special{fp}%
\special{pa 1935 534}%
\special{pa 1935 534}%
\special{pa 1936 534}%
\special{pa 1937 535}%
\special{pa 1938 535}%
\special{pa 1939 536}%
\special{pa 1940 536}%
\special{pa 1940 537}%
\special{pa 1941 537}%
\special{fp}%
\special{pa 1962 558}%
\special{pa 1962 558}%
\special{pa 1963 559}%
\special{pa 1963 560}%
\special{pa 1964 560}%
\special{pa 1964 562}%
\special{pa 1965 562}%
\special{pa 1965 563}%
\special{pa 1965 563}%
\special{fp}%
\special{pa 1975 592}%
\special{pa 1975 600}%
\special{fp}%
}}%
%
{\color[named]{Black}{%
\special{pn 8}%
\special{pa 810 750}%
\special{pa 1600 750}%
\special{dt 0.045}%
}}%
%
{\color[named]{Black}{%
\special{pn 8}%
\special{pa 2200 750}%
\special{pa 2990 750}%
\special{dt 0.045}%
}}%
%
{\color[named]{Black}{%
\special{pn 4}%
\special{pa 1700 750}%
\special{pa 1800 900}%
\special{fp}%
}}%
\put(18.0000,-9.0000){\makebox(0,0)[lt]{$\Gamma$}}%
%
{\color[named]{Black}{%
\special{pn 13}%
\special{ar 3800 800 400 400  0.2449787  6.0382066}%
}}%
%
{\color[named]{Black}{%
\special{pn 13}%
\special{ar 4800 800 400 400  3.3865713  2.8966140}%
}}%
%
{\color[named]{Black}{%
\special{pn 13}%
\special{pa 4190 900}%
\special{pa 4410 900}%
\special{fp}%
}}%
%
{\color[named]{Black}{%
\special{pn 13}%
\special{pa 4190 700}%
\special{pa 4410 700}%
\special{fp}%
}}%
%
{\color[named]{Black}{%
\special{pn 8}%
\special{pa 4200 800}%
\special{pa 4400 800}%
\special{dt 0.045}%
}}%
\put(43.0000,-5.0000){\makebox(0,0){$\Gamma$}}%
%
{\color[named]{Black}{%
\special{pn 4}%
\special{pa 4300 800}%
\special{pa 4300 600}%
\special{fp}%
}}%
\put(15.0000,-11.0000){\makebox(0,0)[lt]{$C_0$}}%
\put(29.0000,-11.0000){\makebox(0,0)[lt]{$C_1$}}%
\put(51.0000,-11.0000){\makebox(0,0)[lt]{$C_1$}}%
\put(41.0000,-11.0000){\makebox(0,0)[lt]{$C_0$}}%
%
{\color[named]{Black}{%
\special{pn 8}%
\special{pa 800 800}%
\special{pa 800 792}%
\special{fp}%
\special{pa 802 755}%
\special{pa 803 753}%
\special{pa 803 750}%
\special{pa 804 748}%
\special{fp}%
\special{pa 810 711}%
\special{pa 811 708}%
\special{pa 812 704}%
\special{fp}%
}}%
%
{\color[named]{Black}{%
\special{pn 8}%
\special{pa 2988 703}%
\special{pa 2988 704}%
\special{pa 2989 706}%
\special{pa 2989 707}%
\special{pa 2990 709}%
\special{pa 2990 710}%
\special{fp}%
\special{pa 2996 747}%
\special{pa 2996 747}%
\special{pa 2997 748}%
\special{pa 2997 754}%
\special{pa 2997 754}%
\special{fp}%
\special{pa 3000 791}%
\special{pa 3000 799}%
\special{fp}%
}}%
%
{\color[named]{Black}{%
\special{pn 8}%
\special{pa 4190 700}%
\special{pa 4191 701}%
\special{pa 4191 704}%
\special{pa 4192 705}%
\special{pa 4192 708}%
\special{fp}%
\special{pa 4199 747}%
\special{pa 4199 747}%
\special{pa 4200 748}%
\special{pa 4200 755}%
\special{pa 4200 755}%
\special{fp}%
\special{pa 4203 796}%
\special{pa 4203 804}%
\special{fp}%
\special{pa 4200 845}%
\special{pa 4200 845}%
\special{pa 4200 852}%
\special{pa 4199 853}%
\special{fp}%
\special{pa 4192 892}%
\special{pa 4192 895}%
\special{pa 4191 896}%
\special{pa 4191 899}%
\special{pa 4190 900}%
\special{fp}%
}}%
%
{\color[named]{Black}{%
\special{pn 8}%
\special{pa 4410 900}%
\special{pa 4409 899}%
\special{pa 4409 896}%
\special{pa 4408 895}%
\special{pa 4408 892}%
\special{fp}%
\special{pa 4401 853}%
\special{pa 4401 853}%
\special{pa 4400 852}%
\special{pa 4400 846}%
\special{pa 4399 845}%
\special{fp}%
\special{pa 4397 804}%
\special{pa 4397 796}%
\special{fp}%
\special{pa 4400 755}%
\special{pa 4400 748}%
\special{pa 4401 747}%
\special{pa 4401 747}%
\special{fp}%
\special{pa 4408 708}%
\special{pa 4408 705}%
\special{pa 4409 704}%
\special{pa 4409 701}%
\special{pa 4410 700}%
\special{fp}%
}}%
\end{picture}%

%% file: 3_Strange_Sum.tex
\section{A strange sum formula for $J^{\pm}$}
Strange sum formula for $J^{\pm}$ is lacking in Theorem~\ref{thm:Arnold_additivity}.
In this section we partially fill the lack in some special cases.
This is the first step to the generalized connected sum formulas.

Let $c\in\GI$ and put $C:=c(S^1)\subset\R^2$.

\begin{defn}\label{def:region_index}
For any connected component $R\subset\R^2\setminus C$,
define the \emph{index of $c$ with respect to $R$} by $\ind^R(c):=\deg\varphi\in\Z$, where $\varphi\colon S^1\to S^1$ is defined by
\[
 \varphi(t):=\frac{c(t)-p}{\abs{c(t)-p}}
\]
for a fixed point $p\in R$.
\end{defn}

Notice that $\ind^R(c)$ is independent of the choice of the point $p\in R$, and is not changed by any regular homotopy in which $c$ never pass through $p$.
It is not difficult to see $\ind^{R^{\infty}}(c)=0$.
The index $\ind^R(c)$ is also called the ``Alexander numbering'' of $R$.

The index $\ind^R(c)$ \emph{does} depend on the orientation of $c$, indeed $\ind^R(-c)=-\ind^R(c)$.
For an unoriented curve $C$, a specified normal vector to $C$ can play a role as an orientation.

\begin{defn}\label{def:normal_index}
Let $v$ be a non-zero vector that emanates from $x=c(t_0)\in C\setminus D(C)$ and that is transverse to $C$.
We say a parameter $c$ of $C$ is \emph{compatible with} $v$ if the ordered pair $\ang{v,c'(t_0)}$ is a positive basis of $\R^2$ (see Figure~\ref{fig:normal_ind_sign}).
\begin{figure}
\begin{center}
\input{normal_ind_sign}
\caption{The orientation induced by a normal vector $v$}
\label{fig:normal_ind_sign}
\end{center}
\end{figure}
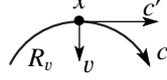

For a connected component $R\subset\R^2\setminus C$, define the \emph{index of $C$ with respect to $(v,R)$} as follows;
choose a compatible parameter $c\colon S^1\looparrowright\R^2$ of $C$ with $v$, and define
\[
 \ind^R_v(C):=\ind^R(c)\in\Z.
\]
Let $R_v$ be the unique component of $\R^2\setminus C$ that contains $x$ in its boundary and that is pointed by the vector $v$ (see Figure~\ref{fig:normal_ind_sign}).
We denote $\ind_v^{R_v}(C)$ by simply $\ind_v(C)$.
\end{defn}
This is a well-defined invariant of \emph{unoriented} generic curves $C$ with specified transverse vector $v$ and connected component $R$ of $\R^2\setminus C$.

In general $\ind^R_v(C)$ can be defined even if $v$ and $R$ do not relate to one another.
In the following however we exclusively consider the cases that a bridge $\Gamma=\gamma([0,1])$ ($\gamma(i)\in C_i$) between generic closed curves $C_0,C_1$ is given, and consider $\ind_v^R(C)$ for $v=(-1)^i\gamma'(i)$ ($i=0,1$) and
\begin{itemize}
\item
	$R=R_{(-1)^i\gamma'(i)}$, or
\item
	$R=R_i$ is the component of $\R^2\setminus C_i$ containing $\gamma(1-i)$.
\end{itemize}

\begin{notation}\label{notation}
In the above setting we put $b_i:=\gamma(i)\in C_i$ and $v_i:=(-1)^i\gamma'(i)$ for $i=0,1$.
See Figure~\ref{fig:bv}
\end{notation}
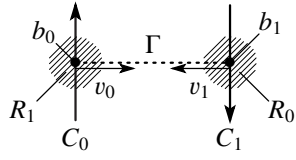
\begin{figure}
\begin{center}
\input{bv}
\end{center}
\caption{Endpoints $b_i\in C_i$ of $\Gamma$ and normal vectors $v_i$ determined by $\Gamma$}
\label{fig:bv}
\end{figure}

Now we are ready to state a strange sum formula for $J^{\pm}$ under some conditions.

\begin{prop}\label{prop:strange_J}
Let $C_0,C_1$ be \emph{separated} unoriented generic closed curves and $\Gamma$ a bridge between them (note that $C_0+_{\Gamma}C_1$ is a strange sum).
Suppose that $\abs{D(\Gamma)}=n_{\Gamma}$ and we can choose a parameter $\gamma\colon [0,1]\to\R^2$ of $\Gamma$ such that, if we let $\Gamma_0:=\gamma([0,1/3])$ and $\Gamma_1:=\gamma([2/3,1])$, then
\begin{equation}\label{eq:strange_assumption}
 \Gamma\cap C_i=\Gamma_i\cap C_i\quad(i=0,1),\quad
 D(\Gamma)=D(\Gamma_0)\sqcup D(\Gamma_1)
\end{equation}
(see Figure~\ref{fig:ex_3_1}).
\begin{figure}
\begin{center}
\input{ex_3_1}
\end{center}
\caption{Strange sums that do / do not satisfy \eqref{eq:strange_assumption}}
\label{fig:ex_3_1}
\end{figure}
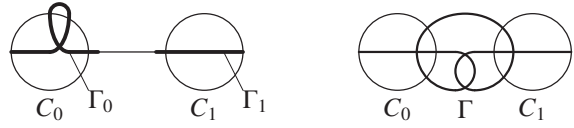
Then we have
\[
 J^{\pm}(C_0+_{\Gamma}C_1)
 =J^{\pm}(C_0)+J^{\pm}(C_1)+\ind_{v_0}(C_0)+\ind_{v_1}(C_1)\pm 2n_{\Gamma}\pm\abs{\Int\Gamma\cap(C_0\cup C_1)}.
\]
\end{prop}

\begin{rem}
In fact Proposition~\ref{prop:strange_J} holds without the assumption \eqref{eq:strange_assumption}.
See Remark~\ref{rem:how_to_recover_formulas}.
\end{rem}

For the proof we need some observations.

\begin{defn}\label{def:same_opposite}
Let $c_0,c_1\in\GI$ be in general position, and let $A\subset\R^2$ be a generic curve such that $A$ is transverse to $C_i$ and $A\cap D(C_i)=\emptyset$ ($i=0,1$).
For any distinct $x\in A\cap C_i$ and $y\in A\cap C_j$ (possibly $i=j$), we say \emph{$c_i$ and $c_j$ point the same side of $A$ at $x$ and $y$} if $c_i$ and $c_j$ are not compatible with respect to the part of $A$ between $x$ and $y$ (see Definition~\ref{def:compatible_ori}).
If otherwise we say \emph{$c_i$ and $c_j$ point the opposite sides of $A$ at $x$ and $y$}.
See Figure~\ref{fig:same_opposite}.
\begin{figure}
\begin{center}
\input{same_opposite}
\end{center}
\caption{$c_i$ and $c_j$ point the same side of $A$ at $x$ and $y$, while $c_i$ and $c_k$ point the opposite sides of $A$ at $x$ and $z$}
\label{fig:same_opposite}
\end{figure}
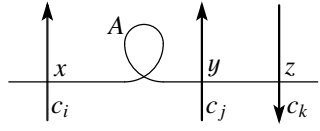
\end{defn}

In the case $i=j$, if $c_i$ points the same side of $A$ at $x$ and $y$, then so does $-c_i$.
Thus it makes sense saying ``an unoriented curve $C$ points the same side of $A$ at $x,y\in A\cap C$''.

\begin{defn}\label{def:sign_Gamma_cap_C}
Let $C_0,C_1$ be generic closed curves and $\Gamma$ a bridge between them.
For $x\in\Int\Gamma\cap C_i$ ($i=0,1$), define $s(x):=+1$ (resp.\ $-1$) if $C_i$ points the same side (resp.\ the opposite sides) of $\Gamma$ at $x$ and $b_i$.
\end{defn}


\begin{lem}\label{lem:adjacent_indices}
Let $C$ be a generic closed curve and $A=\alpha([0,1])$ a generic curve such that
\begin{itemize}
\item
	$A$ is transverse to $C$ and $A\cap D(C)=D(A)\cap C=\emptyset$, and
\item
	$\alpha(0)\in C$, and $\alpha(1)$ is in a connected component $R$ of $\R^2\setminus C$.
\end{itemize}
Let $l^{\pm}$ be the numbers of the points $x\in\Int A\cap C$ such that $s(x)=\pm 1$.
Then $\ind^R_{\alpha'(0)}(C)-\ind_{\alpha'(0)}(C)=l^--l^+$.
\end{lem}

\begin{rem}
We can also state Lemma~\ref{lem:adjacent_indices} as $\ind^R_{\alpha'(0)}(C)-\ind_{\alpha'(0)}(C)=-\sum_{x\in\Int A\cap C}s(x)$.
We often use Lemma~\ref{lem:adjacent_indices} in this form in the proof of Theorem~\ref{thm:main_St}.
\end{rem}

\begin{proof}[Proof of Lemma~\ref{lem:adjacent_indices}]
To compute $\ind^{(R)}_{\alpha'(0)}(C)$ we choose a compatible parameter $c$ of $C$ with $\alpha'(0)$ (see Definition~\ref{def:compatible_ori}).
Let $c(t)\in A\cap C$ and let $R^l,R^r\subset\R^2\setminus C$ be adjacent components to $c(t)$ such that $R^l$ is on the left side of $c$ (see Figure~\ref{fig:adjacent}).
\begin{figure}
\begin{center}
\input{adjacent}
\caption{Comparing $\ind^{R^l}(c)$ with $\ind^{R^r}(c)$ (corners in the right figure should be smoothed)}
\label{fig:adjacent}
\end{center}
\end{figure}
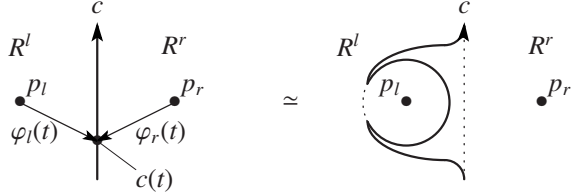
For any points $p_*\in R^*$ ($*=l,r$), we have $\ind^{R^*}(c)=\deg\varphi_*$, where $\varphi_*\colon S^1\to S^1$ is defined by $\varphi_*(t):=(c(t)-p_*)/\abs{c(t)-p_*}$ (see Definition~\ref{def:region_index}).
Transforming $c$ as in Figure~\ref{fig:adjacent}, we see that $\varphi_l(t)$ ``rotates one more than $\varphi_r(t)$ does'' in the counterclockwise direction, and hence
\begin{equation}\label{eq:left_index_bigger}
 \ind^{R^l}(c)=\ind^{R^r}(c)+1.
\end{equation}

Let $\xi\in A\cap C$ be a point such that $s(\xi)=+1$.
Then $\alpha$ gets across $c$ at $\xi$ from the left (consider the case $x=\alpha(0)$ and $y=\xi$ in Figure~\ref{fig:same_opposite}).
Thus by \eqref{eq:left_index_bigger} the index of $c$ with respect to components of $\R^2\setminus C$ decreases by one when we get across $c$ at $\xi$ along $\alpha$.
Similarly we see that the index increases by one when we get across $c$ at $\eta\in A\cap C$ with $s(\eta)=-1$.
Therefore $\ind^R(c)-\ind^{R_{\alpha'(0)}}(c)=l^--l^+$.
We have $\ind^R(c)=\ind^R_{\alpha'(0)}(C)$ and $\ind^{R_{\alpha'(0)}}(c)=\ind_{\alpha'(0)}(C)$ because $c$ is compatible with $\alpha'(0)$.
\end{proof}

The following is proved by inspection (see Figure~\ref{fig:compatible_ori}).

\begin{lem}\label{lem:compatible_ori}
Let $C_0,C_1$ be generic closed curves in general position and $\Gamma$ a bridge between them.
Let $c_0,c_1$ be compatible orientations with respect to $\Gamma$.
Then $c_0$ is compatible with $v_0$ if and only if $c_1$ is compatible with $v_1$.
\end{lem}

\begin{proof}[Proof of Proposition~\ref{prop:strange_J}]
By Lemma~\ref{lem:compatible_ori} we may choose compatible parameters $c_i$ of $C_i$ ($i=0,1$) with respect to both $\Gamma$ and $v_i$.

Transform $C_i$ along $\Gamma_i$ as shown in Figure~\ref{fig:move_along_Gamma} to obtain $\hat{C}_i$ so that $C_0+_{\Gamma}C_1=\hat{C}_0+_{\hat{\Gamma}}\hat{C}_1$, where the right hand side is a connected sum.
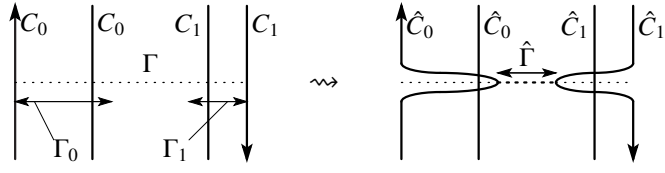
\begin{figure}[htb]
\centering
\input{move_along_Gamma}
\caption{The pushing-appendix along $\Gamma$}
\label{fig:move_along_Gamma}
\end{figure}
In \cite{Arnold} such a transformation is called \emph{pushing-appendix} along $\Gamma$.
It is enough to compute $J^{\pm}(\hat{C}_i)$ since by Theorem~\ref{thm:Arnold_additivity}
\begin{equation}\label{eq:connected_sum_formula_applied}
 J^{\pm}(C_0+_{\Gamma}C_1)
 =J^{\pm}(\hat{C}_0+_{\hat{\Gamma}}\hat{C}_1)
 =J^{\pm}(\hat{C}_0)+J^{\pm}(\hat{C}_1).
\end{equation}

At some finitely many instances the above pushing-appendix gets across self-tangencies,
that derive from
\begin{enumerate}[(i)]
\item
	the double points of $\Gamma$, and
\item
	the intersections of $C_i$ with $\Gamma$.
\end{enumerate}

For (i), define $n_i:=\abs{D(\Gamma_i)}$ ($i=0,1$).
Then $n_0+n_1=n_{\Gamma}$ by \eqref{eq:strange_assumption}.
Corresponding to each double point of $\Gamma$, the above pushing-appendix gets across a pair of a positive direct self-tangency and a positive inverse self-tangency (see Figure~\ref{fig:tangency_from_double_point}).
Thus the self-tangencies derived from (i) add $\pm 2n_i$ to $J^{\pm}(C_i)$.
\begin{figure}
\centering
\input{move_along_Gamma1}
\caption{Self-tangencies derived from a double point of $\Gamma$}
\label{fig:tangency_from_double_point}
\end{figure}
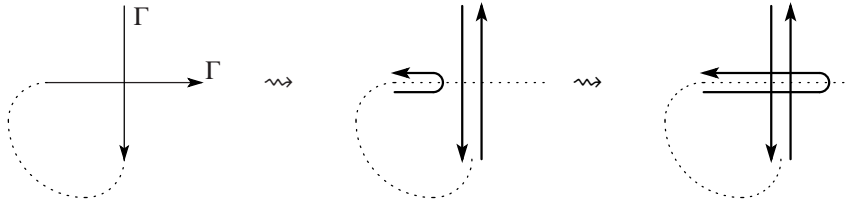

For (ii), we see in Figure~\ref{label:pushing-appendix_+} that, if $x\in\Int\Gamma\cap C_0$ is such that $s(x)=+1$ (resp.\ $s(x)=-1$), then the pushing-appendix gets across positively a direct (resp.\ an inverse) self-tangency at $x$.
\begin{figure}
\begin{center}
\input{move_along_Gamma2}
\end{center}
\caption{Self-tangencies derived from a point in $\Gamma\cap C_0$}
\label{label:pushing-appendix_+}
\end{figure}
Thus self-tangencies derived from (ii) add $\pm 2l^{\pm}_0$ to $J^{\pm}(C_0)$, where $l^{\pm}_0$ are the numbers of points $x\in\Int\Gamma_0\cap C_0$ such that $s(x)=\pm 1$.

By \eqref{eq:strange_assumption} we may suppose that $\gamma(1/3)\in\partial\Gamma_0$ is in the unbounded component of $\R^2\setminus C_0$, whose index is $0$.
Thus $\ind_{v_0}(C_0)=l^+_0-l^-_0$ by Lemma~\ref{lem:adjacent_indices} and by our choice of $c_0$.
Therefore
\begin{align}\label{eq:difference_J_strange_0}
 J^{\pm}(\hat{C}_0)
 &=J^{\pm}(C_0)\pm 2n_0\pm 2l^{\pm}_0 \\
 &=J^{\pm}(C_0)\pm 2n_0+(l^+_0-l^-_0)\pm (l^+_0+l^-_0)\notag \\
 &=J^{\pm}(C_0)\pm 2n_0+\ind_{v_0}(C_0)\pm\abs{\Int\Gamma_0\cap C_0}. \notag
\end{align}
Since $c_1$ is the induced orientation by $v_1$, a similar calculation to \eqref{eq:difference_J_strange_0} shows
\begin{equation}\label{eq:difference_J_strange_1}
 J^{\pm}(\hat{C}_1)=J^{\pm}(C_1)\pm 2n_1+\ind_{v_1}(C_1)\pm\abs{\Int\Gamma_1\cap C_1}.
\end{equation}
Substituting \eqref{eq:difference_J_strange_0}, \eqref{eq:difference_J_strange_1} to \eqref{eq:connected_sum_formula_applied} and using $n_{\Gamma}=n_0+n_1$ and $\abs{\Int\Gamma_0\cap C_0}+\abs{\Int\Gamma_1\cap C_1}=\abs{\Int\Gamma\cap(C_0\cup C_1)}$ (this is because of \eqref{eq:strange_assumption}), we complete the proof.
\end{proof}

%% file: normal_ind_sign.tex
{\unitlength 0.1in%
\begin{picture}(7.5800,3.8600)(6.4200,-8.2100)%
%
\special{pn 13}%
\special{ar 1000 1000 400 400 3.6052403 5.8195377}%
%
\special{pn 13}%
\special{pa 1351 808}%
\special{pa 1358 821}%
\special{fp}%
\special{sh 1}%
\special{pa 1358 821}%
\special{pa 1344 753}%
\special{pa 1333 774}%
\special{pa 1309 772}%
\special{pa 1358 821}%
\special{fp}%
\put(8.0000,-8.0000){\makebox(0,0){$R_v$}}%
%
\special{pn 8}%
\special{pa 1000 600}%
\special{pa 1000 800}%
\special{fp}%
\special{sh 1}%
\special{pa 1000 800}%
\special{pa 1020 733}%
\special{pa 1000 747}%
\special{pa 980 733}%
\special{pa 1000 800}%
\special{fp}%
\put(10.2000,-8.0000){\makebox(0,0)[lt]{$v$}}%
\put(14.0000,-8.0000){\makebox(0,0)[lb]{$c$}}%
\put(10.0000,-6.0000){\makebox(0,0){$\bullet$}}%
\put(10.0000,-5.0000){\makebox(0,0){$x$}}%
%
\special{pn 8}%
\special{pa 1000 600}%
\special{pa 1400 600}%
\special{fp}%
\special{sh 1}%
\special{pa 1400 600}%
\special{pa 1333 580}%
\special{pa 1347 600}%
\special{pa 1333 620}%
\special{pa 1400 600}%
\special{fp}%
\put(13.3000,-5.7000){\makebox(0,0)[lb]{$c'$}}%
\end{picture}}%

%% file: bv.tex
{\unitlength 0.1in%
\begin{picture}(15.3500,6.2700)(8.6500,-17.2700)%
%
\special{pn 13}%
\special{pa 1400 1700}%
\special{pa 1400 1100}%
\special{fp}%
\special{sh 1}%
\special{pa 1400 1100}%
\special{pa 1380 1167}%
\special{pa 1400 1153}%
\special{pa 1420 1167}%
\special{pa 1400 1100}%
\special{fp}%
\put(14.0000,-18.0000){\makebox(0,0){$C_0$}}%
\put(14.0000,-14.0000){\makebox(0,0){$\bullet$}}%
%
\special{pn 13}%
\special{pa 2200 1100}%
\special{pa 2200 1700}%
\special{fp}%
\special{sh 1}%
\special{pa 2200 1700}%
\special{pa 2220 1633}%
\special{pa 2200 1647}%
\special{pa 2180 1633}%
\special{pa 2200 1700}%
\special{fp}%
\put(22.0000,-14.0000){\makebox(0,0){$\bullet$}}%
\put(13.0000,-13.0000){\makebox(0,0)[rb]{$b_0$}}%
\put(23.5000,-12.5000){\makebox(0,0)[lb]{$b_1$}}%
%
\special{pn 13}%
\special{pa 1400 1400}%
\special{pa 2200 1400}%
\special{dt 0.045}%
\put(18.0000,-13.0000){\makebox(0,0){$\Gamma$}}%
\put(22.0000,-18.0000){\makebox(0,0){$C_1$}}%
%
\special{pn 8}%
\special{pa 1400 1430}%
\special{pa 1700 1430}%
\special{fp}%
\special{sh 1}%
\special{pa 1700 1430}%
\special{pa 1633 1410}%
\special{pa 1647 1430}%
\special{pa 1633 1450}%
\special{pa 1700 1430}%
\special{fp}%
\put(15.0000,-15.0000){\makebox(0,0)[lt]{$v_0$}}%
%
\special{pn 8}%
\special{pa 2200 1430}%
\special{pa 1900 1430}%
\special{fp}%
\special{sh 1}%
\special{pa 1900 1430}%
\special{pa 1967 1450}%
\special{pa 1953 1430}%
\special{pa 1967 1410}%
\special{pa 1900 1430}%
\special{fp}%
\put(21.0000,-15.0000){\makebox(0,0)[rt]{$v_1$}}%
%
\special{pn 4}%
\special{pa 2305 1315}%
\special{pa 2115 1505}%
\special{fp}%
\special{pa 2320 1330}%
\special{pa 2135 1515}%
\special{fp}%
\special{pa 2330 1350}%
\special{pa 2155 1525}%
\special{fp}%
\special{pa 2335 1375}%
\special{pa 2175 1535}%
\special{fp}%
\special{pa 2335 1405}%
\special{pa 2205 1535}%
\special{fp}%
\special{pa 2330 1440}%
\special{pa 2240 1530}%
\special{fp}%
\special{pa 2290 1300}%
\special{pa 2100 1490}%
\special{fp}%
\special{pa 2270 1290}%
\special{pa 2090 1470}%
\special{fp}%
\special{pa 2255 1275}%
\special{pa 2075 1455}%
\special{fp}%
\special{pa 2230 1270}%
\special{pa 2065 1435}%
\special{fp}%
\special{pa 2205 1265}%
\special{pa 2065 1405}%
\special{fp}%
\special{pa 2175 1265}%
\special{pa 2065 1375}%
\special{fp}%
\put(24.0000,-16.0000){\makebox(0,0)[lt]{$R_0$}}%
%
\special{pn 8}%
\special{pa 2250 1450}%
\special{pa 2400 1600}%
\special{fp}%
%
\special{pn 8}%
\special{pa 2350 1250}%
\special{pa 2200 1400}%
\special{fp}%
%
\special{pn 8}%
\special{pa 1300 1300}%
\special{pa 1400 1400}%
\special{fp}%
%
\special{pn 4}%
\special{pa 1505 1315}%
\special{pa 1315 1505}%
\special{fp}%
\special{pa 1520 1330}%
\special{pa 1335 1515}%
\special{fp}%
\special{pa 1530 1350}%
\special{pa 1355 1525}%
\special{fp}%
\special{pa 1535 1375}%
\special{pa 1375 1535}%
\special{fp}%
\special{pa 1535 1405}%
\special{pa 1405 1535}%
\special{fp}%
\special{pa 1530 1440}%
\special{pa 1440 1530}%
\special{fp}%
\special{pa 1490 1300}%
\special{pa 1300 1490}%
\special{fp}%
\special{pa 1470 1290}%
\special{pa 1290 1470}%
\special{fp}%
\special{pa 1455 1275}%
\special{pa 1275 1455}%
\special{fp}%
\special{pa 1430 1270}%
\special{pa 1265 1435}%
\special{fp}%
\special{pa 1405 1265}%
\special{pa 1265 1405}%
\special{fp}%
\special{pa 1375 1265}%
\special{pa 1265 1375}%
\special{fp}%
\put(12.0000,-16.0000){\makebox(0,0)[rt]{$R_1$}}%
%
\special{pn 8}%
\special{pa 1350 1450}%
\special{pa 1200 1600}%
\special{fp}%
\end{picture}}%

%% file: ex_3_1.tex
{\unitlength 0.1in%
\begin{picture}(29.0000,4.8500)(7.0000,-10.3500)%
%
\special{pn 8}%
\special{ar 900 800 200 200 0.0000000 6.2831853}%
%
\special{pn 20}%
\special{ar 900 600 100 200 6.2831853 1.5707963}%
%
\special{pn 20}%
\special{ar 1000 600 100 200 1.5707963 3.1415927}%
%
\special{pn 20}%
\special{ar 950 600 50 50 3.1415927 6.2831853}%
%
\special{pn 20}%
\special{pa 1450 800}%
\special{pa 1900 800}%
\special{fp}%
%
\special{pn 20}%
\special{pa 900 800}%
\special{pa 700 800}%
\special{fp}%
%
\special{pn 8}%
\special{ar 1700 800 200 200 0.0000000 6.2831853}%
%
\special{pn 8}%
\special{ar 2700 800 200 200 0.0000000 6.2831853}%
%
\special{pn 13}%
\special{ar 3000 900 100 100 4.7123890 1.5707963}%
%
\special{pn 13}%
\special{ar 3000 800 200 200 1.5707963 3.1415927}%
%
\special{pn 13}%
\special{ar 3100 900 100 100 1.5707963 4.7123890}%
%
\special{pn 13}%
\special{ar 3100 800 200 200 6.2831853 1.5707963}%
%
\special{pn 13}%
\special{ar 3050 800 250 200 3.1415927 6.2831853}%
%
\special{pn 13}%
\special{pa 3000 800}%
\special{pa 2500 800}%
\special{fp}%
%
\special{pn 8}%
\special{ar 3400 800 200 200 0.0000000 6.2831853}%
%
\special{pn 13}%
\special{pa 3100 800}%
\special{pa 3600 800}%
\special{fp}%
%
\special{pn 20}%
\special{pa 1000 800}%
\special{pa 1150 800}%
\special{fp}%
%
\special{pn 8}%
\special{pa 1150 800}%
\special{pa 1450 800}%
\special{fp}%
\put(9.0000,-11.0000){\makebox(0,0){$C_0$}}%
\put(17.0000,-11.0000){\makebox(0,0){$C_1$}}%
\put(34.0000,-11.0000){\makebox(0,0){$C_1$}}%
\put(27.0000,-11.0000){\makebox(0,0){$C_0$}}%
\put(11.0000,-10.0000){\makebox(0,0)[lt]{$\Gamma_0$}}%
\put(19.0000,-10.0000){\makebox(0,0)[lt]{$\Gamma_1$}}%
\put(30.5000,-11.0000){\makebox(0,0){$\Gamma$}}%
%
\special{pn 4}%
\special{pa 1100 1000}%
\special{pa 1000 800}%
\special{fp}%
%
\special{pn 4}%
\special{pa 1900 1000}%
\special{pa 1800 800}%
\special{fp}%
\end{picture}}%

%% file: same_opposite.tex
{\unitlength 0.1in%
\begin{picture}(16.0000,6.0000)(8.0000,-10.0000)%
%
\special{pn 13}%
\special{pa 1000 1000}%
\special{pa 1000 400}%
\special{fp}%
\special{sh 1}%
\special{pa 1000 400}%
\special{pa 980 467}%
\special{pa 1000 453}%
\special{pa 1020 467}%
\special{pa 1000 400}%
\special{fp}%
%
\special{pn 13}%
\special{pa 1800 1000}%
\special{pa 1800 400}%
\special{fp}%
\special{sh 1}%
\special{pa 1800 400}%
\special{pa 1780 467}%
\special{pa 1800 453}%
\special{pa 1820 467}%
\special{pa 1800 400}%
\special{fp}%
%
\special{pn 13}%
\special{pa 2200 400}%
\special{pa 2200 1000}%
\special{fp}%
\special{sh 1}%
\special{pa 2200 1000}%
\special{pa 2220 933}%
\special{pa 2200 947}%
\special{pa 2180 933}%
\special{pa 2200 1000}%
\special{fp}%
%
\special{pn 8}%
\special{ar 1400 600 200 200 6.2831853 1.5707963}%
%
\special{pn 8}%
\special{ar 1600 600 200 200 1.5707963 3.1415927}%
%
\special{pn 8}%
\special{ar 1500 600 100 100 3.1415927 6.2831853}%
%
\special{pn 8}%
\special{pa 800 800}%
\special{pa 1400 800}%
\special{fp}%
%
\special{pn 8}%
\special{pa 1600 800}%
\special{pa 2400 800}%
\special{fp}%
\put(10.2000,-9.0000){\makebox(0,0)[lt]{$c_i$}}%
\put(18.2000,-9.0000){\makebox(0,0)[lt]{$c_j$}}%
\put(22.4000,-9.0000){\makebox(0,0)[lt]{$c_k$}}%
\put(14.0000,-5.5000){\makebox(0,0)[rb]{$A$}}%
\put(10.3000,-7.7000){\makebox(0,0)[lb]{$x$}}%
\put(18.3000,-7.7000){\makebox(0,0)[lb]{$y$}}%
\put(22.3000,-7.7000){\makebox(0,0)[lb]{$z$}}%
\end{picture}}%

%% file: adjacent.tex
{\unitlength 0.1in%
\begin{picture}(32.1000,9.6500)(3.2000,-14.0000)%
%
{\color{Black}{%
\special{pn 13}%
\special{pa 1200 1400}%
\special{pa 1200 600}%
\special{fp}%
\special{sh 1}%
\special{pa 1200 600}%
\special{pa 1180 667}%
\special{pa 1200 653}%
\special{pa 1220 667}%
\special{pa 1200 600}%
\special{fp}%
}}%
\put(12.0000,-5.0000){\makebox(0,0){{\color{Black}{$c$}}}}%
\put(8.0000,-10.0000){\makebox(0,0){{\color{Black}{$\bullet$}}}}%
\put(8.3000,-9.7000){\makebox(0,0)[lb]{{\color{Black}{$p_l$}}}}%
\put(8.0000,-7.0000){\makebox(0,0){{\color{Black}{$R^l$}}}}%
\put(16.0000,-10.0000){\makebox(0,0){{\color{Black}{$\bullet$}}}}%
\put(16.3000,-9.7000){\makebox(0,0)[lb]{{\color{Black}{$p_r$}}}}%
\put(16.0000,-7.0000){\makebox(0,0){{\color{Black}{$R^r$}}}}%
\put(22.0000,-10.0000){\makebox(0,0){{\color{Black}{$\simeq$}}}}%
\put(12.0000,-12.0000){\makebox(0,0){{\color{Black}{$\bullet$}}}}%
\put(14.0000,-13.5000){\makebox(0,0)[lt]{{\color{Black}{$c(t)$}}}}%
%
{\color{Black}{%
\special{pn 8}%
\special{pa 800 1000}%
\special{pa 1200 1200}%
\special{fp}%
\special{sh 1}%
\special{pa 1200 1200}%
\special{pa 1149 1152}%
\special{pa 1152 1176}%
\special{pa 1131 1188}%
\special{pa 1200 1200}%
\special{fp}%
}}%
\put(10.0000,-11.0000){\makebox(0,0)[rt]{{\color{Black}{$\varphi_l(t)$}}}}%
%
{\color{Black}{%
\special{pn 8}%
\special{pa 1600 1000}%
\special{pa 1200 1200}%
\special{fp}%
\special{sh 1}%
\special{pa 1200 1200}%
\special{pa 1269 1188}%
\special{pa 1248 1176}%
\special{pa 1251 1152}%
\special{pa 1200 1200}%
\special{fp}%
}}%
\put(14.0000,-11.0000){\makebox(0,0)[lt]{{\color{Black}{$\varphi_r(t)$}}}}%
%
{\color{Black}{%
\special{pn 4}%
\special{pa 1400 1350}%
\special{pa 1200 1200}%
\special{fp}%
}}%
\put(28.0000,-10.0000){\makebox(0,0){{\color{Black}{$\bullet$}}}}%
%
{\color{Black}{%
\special{pn 13}%
\special{ar 2800 1000 224 224 3.6052403 2.6779450}%
}}%
\put(31.0000,-5.0000){\makebox(0,0){{\color{Black}{$c$}}}}%
\put(25.0000,-7.0000){\makebox(0,0){{\color{Black}{$R^l$}}}}%
\put(35.0000,-7.0000){\makebox(0,0){{\color{Black}{$R^r$}}}}%
\put(27.7000,-9.7000){\makebox(0,0)[rb]{{\color{Black}{$p_l$}}}}%
\put(35.0000,-10.0000){\makebox(0,0){{\color{Black}{$\bullet$}}}}%
\put(35.3000,-9.7000){\makebox(0,0)[lb]{{\color{Black}{$p_r$}}}}%
%
{\color{Black}{%
\special{pn 8}%
\special{pa 3100 1400}%
\special{pa 3100 600}%
\special{dt 0.045}%
}}%
%
{\color{Black}{%
\special{pn 13}%
\special{ar 3000 900 400 200 3.1415927 4.7123890}%
}}%
%
{\color{Black}{%
\special{pn 13}%
\special{ar 3000 600 100 100 6.2831853 1.5707963}%
}}%
%
{\color{Black}{%
\special{pn 13}%
\special{pa 3099 615}%
\special{pa 3100 600}%
\special{fp}%
\special{sh 1}%
\special{pa 3100 600}%
\special{pa 3076 665}%
\special{pa 3096 653}%
\special{pa 3116 668}%
\special{pa 3100 600}%
\special{fp}%
}}%
%
{\color{Black}{%
\special{pn 13}%
\special{ar 3000 1100 400 200 1.5707963 3.1415927}%
}}%
%
{\color{Black}{%
\special{pn 13}%
\special{ar 3000 1400 100 100 4.7123890 6.2831853}%
}}%
%
{\color{Black}{%
\special{pn 8}%
\special{pn 8}%
\special{pa 2600 1100}%
\special{pa 2596 1093}%
\special{fp}%
\special{pa 2586 1065}%
\special{pa 2584 1059}%
\special{fp}%
\special{pa 2577 1025}%
\special{pa 2577 1018}%
\special{fp}%
\special{pa 2577 983}%
\special{pa 2577 975}%
\special{fp}%
\special{pa 2583 943}%
\special{pa 2585 936}%
\special{fp}%
\special{pa 2597 905}%
\special{pa 2600 900}%
\special{fp}%
}}%
\end{picture}}%

%% file: move_along_Gamma.tex
{\unitlength 0.1in%
\begin{picture}(32.8200,8.0000)(5.4500,-12.0000)%
%
\special{pn 13}%
\special{pa 600 1200}%
\special{pa 600 400}%
\special{fp}%
\special{sh 1}%
\special{pa 600 400}%
\special{pa 580 467}%
\special{pa 600 453}%
\special{pa 620 467}%
\special{pa 600 400}%
\special{fp}%
%
\special{pn 8}%
\special{pa 600 800}%
\special{pa 1800 800}%
\special{dt 0.045}%
\put(13.0000,-7.0000){\makebox(0,0){$\Gamma$}}%
\put(7.0000,-5.0000){\makebox(0,0){$C_0$}}%
\put(22.0000,-8.0000){\makebox(0,0){$\rightsquigarrow$}}%
%
\special{pn 13}%
\special{pa 1000 1200}%
\special{pa 1000 400}%
\special{fp}%
\put(11.0000,-5.0000){\makebox(0,0){$C_0$}}%
%
\special{pn 13}%
\special{pa 2600 700}%
\special{pa 2600 400}%
\special{fp}%
\special{sh 1}%
\special{pa 2600 400}%
\special{pa 2580 467}%
\special{pa 2600 453}%
\special{pa 2620 467}%
\special{pa 2600 400}%
\special{fp}%
%
\special{pn 13}%
\special{pa 2600 1200}%
\special{pa 2600 900}%
\special{fp}%
%
\special{pn 13}%
\special{pa 3000 1200}%
\special{pa 3000 400}%
\special{fp}%
\put(27.0000,-5.0000){\makebox(0,0){$\hat{C}_0$}}%
\put(31.0000,-5.0000){\makebox(0,0){$\hat{C}_0$}}%
%
\special{pn 13}%
\special{pa 3800 900}%
\special{pa 3800 1200}%
\special{fp}%
\special{sh 1}%
\special{pa 3800 1200}%
\special{pa 3820 1133}%
\special{pa 3800 1147}%
\special{pa 3780 1133}%
\special{pa 3800 1200}%
\special{fp}%
%
\special{pn 13}%
\special{pa 3600 1200}%
\special{pa 3600 400}%
\special{fp}%
\put(39.0000,-5.0000){\makebox(0,0){$\hat{C}_1$}}%
%
\special{pn 13}%
\special{pa 1800 400}%
\special{pa 1800 1200}%
\special{fp}%
\special{sh 1}%
\special{pa 1800 1200}%
\special{pa 1820 1133}%
\special{pa 1800 1147}%
\special{pa 1780 1133}%
\special{pa 1800 1200}%
\special{fp}%
%
\special{pn 13}%
\special{pa 1600 1200}%
\special{pa 1600 400}%
\special{fp}%
\put(15.0000,-5.0000){\makebox(0,0){$C_1$}}%
\put(19.0000,-5.0000){\makebox(0,0){$C_1$}}%
%
\special{pn 13}%
\special{pa 3800 700}%
\special{pa 3800 400}%
\special{fp}%
%
\special{pn 13}%
\special{ar 2800 700 200 50 1.5707963 3.1415927}%
\put(35.0000,-5.0000){\makebox(0,0){$\hat{C}_1$}}%
%
\special{pn 13}%
\special{ar 2800 900 200 50 3.1415927 4.7123890}%
%
\special{pn 13}%
\special{ar 2800 800 300 50 4.7123890 1.5707963}%
%
\special{pn 13}%
\special{ar 3600 700 200 50 6.2831853 1.5707963}%
%
\special{pn 13}%
\special{ar 3600 900 200 50 4.7123890 6.2831853}%
%
\special{pn 13}%
\special{ar 3600 800 200 50 1.5707963 4.7123890}%
%
\special{pn 13}%
\special{pa 3100 800}%
\special{pa 3400 800}%
\special{dt 0.045}%
\put(32.5000,-6.5000){\makebox(0,0){$\hat{\Gamma}$}}%
%
\special{pn 8}%
\special{pa 3400 800}%
\special{pa 3800 800}%
\special{dt 0.045}%
%
\special{pn 4}%
\special{pa 3100 750}%
\special{pa 3400 750}%
\special{fp}%
\special{sh 1}%
\special{pa 3400 750}%
\special{pa 3333 730}%
\special{pa 3347 750}%
\special{pa 3333 770}%
\special{pa 3400 750}%
\special{fp}%
%
\special{pn 8}%
\special{pa 3400 750}%
\special{pa 3100 750}%
\special{fp}%
\special{sh 1}%
\special{pa 3100 750}%
\special{pa 3167 770}%
\special{pa 3153 750}%
\special{pa 3167 730}%
\special{pa 3100 750}%
\special{fp}%
%
\special{pn 4}%
\special{pa 600 900}%
\special{pa 1100 900}%
\special{fp}%
\special{sh 1}%
\special{pa 1100 900}%
\special{pa 1033 880}%
\special{pa 1047 900}%
\special{pa 1033 920}%
\special{pa 1100 900}%
\special{fp}%
%
\special{pn 4}%
\special{pa 1100 900}%
\special{pa 600 900}%
\special{fp}%
\special{sh 1}%
\special{pa 600 900}%
\special{pa 667 920}%
\special{pa 653 900}%
\special{pa 667 880}%
\special{pa 600 900}%
\special{fp}%
%
\special{pn 4}%
\special{pa 1500 900}%
\special{pa 1800 900}%
\special{fp}%
\special{sh 1}%
\special{pa 1800 900}%
\special{pa 1733 880}%
\special{pa 1747 900}%
\special{pa 1733 920}%
\special{pa 1800 900}%
\special{fp}%
%
\special{pn 4}%
\special{pa 1800 900}%
\special{pa 1500 900}%
\special{fp}%
\special{sh 1}%
\special{pa 1500 900}%
\special{pa 1567 920}%
\special{pa 1553 900}%
\special{pa 1567 880}%
\special{pa 1500 900}%
\special{fp}%
\put(8.0000,-11.0000){\makebox(0,0)[lt]{$\Gamma_0$}}%
\put(15.0000,-11.0000){\makebox(0,0)[rt]{$\Gamma_1$}}%
%
\special{pn 4}%
\special{pa 1500 1100}%
\special{pa 1700 900}%
\special{fp}%
%
\special{pn 4}%
\special{pa 800 1100}%
\special{pa 700 900}%
\special{fp}%
%
\special{pn 8}%
\special{pa 2600 800}%
\special{pa 3100 800}%
\special{dt 0.045}%
\end{picture}}%

%% file: move_along_Gamma1.tex
{\unitlength 0.1in%
\begin{picture}(44.0000,10.3000)(4.0000,-14.0000)%
\put(10.5000,-5.0000){\makebox(0,0)[lb]{$\Gamma$}}%
\put(14.2000,-7.8000){\makebox(0,0)[lb]{$\Gamma$}}%
\put(18.0000,-8.0000){\makebox(0,0){$\rightsquigarrow$}}%
%
\special{pn 8}%
\special{pa 600 800}%
\special{pa 1400 800}%
\special{fp}%
\special{sh 1}%
\special{pa 1400 800}%
\special{pa 1333 780}%
\special{pa 1347 800}%
\special{pa 1333 820}%
\special{pa 1400 800}%
\special{fp}%
%
\special{pn 8}%
\special{pa 1000 400}%
\special{pa 1000 1200}%
\special{fp}%
\special{sh 1}%
\special{pa 1000 1200}%
\special{pa 1020 1133}%
\special{pa 1000 1147}%
\special{pa 980 1133}%
\special{pa 1000 1200}%
\special{fp}%
%
\special{pn 8}%
\special{pn 8}%
\special{pa 1000 1200}%
\special{pa 1000 1208}%
\special{fp}%
\special{pa 995 1244}%
\special{pa 993 1251}%
\special{fp}%
\special{pa 981 1285}%
\special{pa 978 1292}%
\special{fp}%
\special{pa 958 1322}%
\special{pa 954 1328}%
\special{fp}%
\special{pa 928 1354}%
\special{pa 922 1359}%
\special{fp}%
\special{pa 891 1378}%
\special{pa 885 1381}%
\special{fp}%
\special{pa 851 1393}%
\special{pa 844 1395}%
\special{fp}%
\special{pa 808 1400}%
\special{pa 800 1400}%
\special{fp}%
%
\special{pn 8}%
\special{pn 8}%
\special{pa 800 1400}%
\special{pa 792 1400}%
\special{fp}%
\special{pa 756 1398}%
\special{pa 748 1397}%
\special{fp}%
\special{pa 712 1390}%
\special{pa 704 1389}%
\special{fp}%
\special{pa 670 1378}%
\special{pa 662 1376}%
\special{fp}%
\special{pa 629 1361}%
\special{pa 622 1358}%
\special{fp}%
\special{pa 590 1340}%
\special{pa 583 1336}%
\special{fp}%
\special{pa 553 1315}%
\special{pa 547 1310}%
\special{fp}%
\special{pa 520 1285}%
\special{pa 514 1280}%
\special{fp}%
\special{pa 490 1253}%
\special{pa 485 1247}%
\special{fp}%
\special{pa 464 1217}%
\special{pa 460 1210}%
\special{fp}%
\special{pa 442 1179}%
\special{pa 438 1172}%
\special{fp}%
\special{pa 425 1138}%
\special{pa 422 1131}%
\special{fp}%
\special{pa 412 1095}%
\special{pa 410 1088}%
\special{fp}%
\special{pa 403 1052}%
\special{pa 402 1044}%
\special{fp}%
\special{pa 400 1008}%
\special{pa 400 1000}%
\special{fp}%
%
\special{pn 8}%
\special{pn 8}%
\special{pa 400 1000}%
\special{pa 400 992}%
\special{fp}%
\special{pa 405 956}%
\special{pa 407 949}%
\special{fp}%
\special{pa 419 915}%
\special{pa 423 908}%
\special{fp}%
\special{pa 441 878}%
\special{pa 446 872}%
\special{fp}%
\special{pa 471 847}%
\special{pa 477 842}%
\special{fp}%
\special{pa 508 823}%
\special{pa 515 819}%
\special{fp}%
\special{pa 549 807}%
\special{pa 556 805}%
\special{fp}%
\special{pa 592 800}%
\special{pa 600 800}%
\special{fp}%
%
\special{pn 13}%
\special{pa 2400 850}%
\special{pa 2600 850}%
\special{fp}%
%
\special{pn 13}%
\special{ar 2600 800 50 50 4.7123890 1.5707963}%
%
\special{pn 8}%
\special{pn 8}%
\special{pa 2200 1000}%
\special{pa 2200 991}%
\special{fp}%
\special{pa 2206 953}%
\special{pa 2208 946}%
\special{fp}%
\special{pa 2220 913}%
\special{pa 2223 906}%
\special{fp}%
\special{pa 2242 877}%
\special{pa 2247 871}%
\special{fp}%
\special{pa 2271 847}%
\special{pa 2277 842}%
\special{fp}%
\special{pa 2306 823}%
\special{pa 2312 820}%
\special{fp}%
\special{pa 2346 808}%
\special{pa 2353 806}%
\special{fp}%
\special{pa 2391 800}%
\special{pa 2400 800}%
\special{fp}%
%
\special{pn 8}%
\special{pn 8}%
\special{pa 2600 1400}%
\special{pa 2592 1400}%
\special{fp}%
\special{pa 2555 1398}%
\special{pa 2547 1397}%
\special{fp}%
\special{pa 2512 1390}%
\special{pa 2504 1388}%
\special{fp}%
\special{pa 2470 1378}%
\special{pa 2462 1375}%
\special{fp}%
\special{pa 2429 1362}%
\special{pa 2422 1358}%
\special{fp}%
\special{pa 2390 1340}%
\special{pa 2383 1336}%
\special{fp}%
\special{pa 2354 1316}%
\special{pa 2348 1311}%
\special{fp}%
\special{pa 2320 1285}%
\special{pa 2314 1280}%
\special{fp}%
\special{pa 2290 1252}%
\special{pa 2285 1246}%
\special{fp}%
\special{pa 2264 1217}%
\special{pa 2259 1210}%
\special{fp}%
\special{pa 2242 1179}%
\special{pa 2239 1172}%
\special{fp}%
\special{pa 2225 1139}%
\special{pa 2222 1131}%
\special{fp}%
\special{pa 2212 1097}%
\special{pa 2210 1089}%
\special{fp}%
\special{pa 2203 1053}%
\special{pa 2203 1045}%
\special{fp}%
\special{pa 2200 1008}%
\special{pa 2200 1000}%
\special{fp}%
%
\special{pn 8}%
\special{pn 8}%
\special{pa 2800 1200}%
\special{pa 2800 1209}%
\special{fp}%
\special{pa 2795 1246}%
\special{pa 2792 1255}%
\special{fp}%
\special{pa 2779 1288}%
\special{pa 2777 1294}%
\special{fp}%
\special{pa 2757 1323}%
\special{pa 2753 1329}%
\special{fp}%
\special{pa 2728 1353}%
\special{pa 2723 1358}%
\special{fp}%
\special{pa 2693 1377}%
\special{pa 2686 1381}%
\special{fp}%
\special{pa 2654 1393}%
\special{pa 2646 1395}%
\special{fp}%
\special{pa 2609 1400}%
\special{pa 2600 1400}%
\special{fp}%
%
\special{pn 13}%
\special{pa 2750 400}%
\special{pa 2750 1200}%
\special{fp}%
\special{sh 1}%
\special{pa 2750 1200}%
\special{pa 2770 1133}%
\special{pa 2750 1147}%
\special{pa 2730 1133}%
\special{pa 2750 1200}%
\special{fp}%
%
\special{pn 13}%
\special{pa 2850 1200}%
\special{pa 2850 400}%
\special{fp}%
\special{sh 1}%
\special{pa 2850 400}%
\special{pa 2830 467}%
\special{pa 2850 453}%
\special{pa 2870 467}%
\special{pa 2850 400}%
\special{fp}%
\put(34.0000,-8.0000){\makebox(0,0){$\rightsquigarrow$}}%
%
\special{pn 8}%
\special{pa 2400 800}%
\special{pa 3200 800}%
\special{dt 0.045}%
%
\special{pn 13}%
\special{pa 2600 750}%
\special{pa 2400 750}%
\special{fp}%
\special{sh 1}%
\special{pa 2400 750}%
\special{pa 2467 770}%
\special{pa 2453 750}%
\special{pa 2467 730}%
\special{pa 2400 750}%
\special{fp}%
%
\special{pn 8}%
\special{pn 8}%
\special{pa 3800 1000}%
\special{pa 3800 992}%
\special{fp}%
\special{pa 3804 958}%
\special{pa 3806 951}%
\special{fp}%
\special{pa 3816 923}%
\special{pa 3818 916}%
\special{fp}%
\special{pa 3833 890}%
\special{pa 3837 885}%
\special{fp}%
\special{pa 3855 863}%
\special{pa 3860 858}%
\special{fp}%
\special{pa 3883 838}%
\special{pa 3889 834}%
\special{fp}%
\special{pa 3914 820}%
\special{pa 3919 817}%
\special{fp}%
\special{pa 3951 806}%
\special{pa 3958 804}%
\special{fp}%
\special{pa 3992 800}%
\special{pa 4000 800}%
\special{fp}%
%
\special{pn 8}%
\special{pn 8}%
\special{pa 4200 1400}%
\special{pa 4192 1400}%
\special{fp}%
\special{pa 4156 1398}%
\special{pa 4148 1397}%
\special{fp}%
\special{pa 4113 1390}%
\special{pa 4106 1389}%
\special{fp}%
\special{pa 4072 1379}%
\special{pa 4065 1376}%
\special{fp}%
\special{pa 4034 1364}%
\special{pa 4027 1361}%
\special{fp}%
\special{pa 3998 1345}%
\special{pa 3992 1342}%
\special{fp}%
\special{pa 3964 1323}%
\special{pa 3959 1319}%
\special{fp}%
\special{pa 3933 1298}%
\special{pa 3928 1293}%
\special{fp}%
\special{pa 3906 1271}%
\special{pa 3901 1265}%
\special{fp}%
\special{pa 3880 1240}%
\special{pa 3876 1234}%
\special{fp}%
\special{pa 3857 1206}%
\special{pa 3854 1200}%
\special{fp}%
\special{pa 3838 1171}%
\special{pa 3835 1164}%
\special{fp}%
\special{pa 3823 1133}%
\special{pa 3820 1126}%
\special{fp}%
\special{pa 3811 1093}%
\special{pa 3809 1086}%
\special{fp}%
\special{pa 3803 1052}%
\special{pa 3802 1044}%
\special{fp}%
\special{pa 3800 1008}%
\special{pa 3800 1000}%
\special{fp}%
%
\special{pn 8}%
\special{pn 8}%
\special{pa 4400 1200}%
\special{pa 4400 1208}%
\special{fp}%
\special{pa 4396 1242}%
\special{pa 4394 1250}%
\special{fp}%
\special{pa 4383 1282}%
\special{pa 4380 1288}%
\special{fp}%
\special{pa 4365 1314}%
\special{pa 4361 1319}%
\special{fp}%
\special{pa 4342 1341}%
\special{pa 4338 1345}%
\special{fp}%
\special{pa 4316 1363}%
\special{pa 4311 1366}%
\special{fp}%
\special{pa 4284 1381}%
\special{pa 4279 1384}%
\special{fp}%
\special{pa 4249 1394}%
\special{pa 4242 1395}%
\special{fp}%
\special{pa 4208 1400}%
\special{pa 4200 1400}%
\special{fp}%
%
\special{pn 13}%
\special{pa 4350 400}%
\special{pa 4350 1200}%
\special{fp}%
\special{sh 1}%
\special{pa 4350 1200}%
\special{pa 4370 1133}%
\special{pa 4350 1147}%
\special{pa 4330 1133}%
\special{pa 4350 1200}%
\special{fp}%
%
\special{pn 13}%
\special{pa 4450 1200}%
\special{pa 4450 400}%
\special{fp}%
\special{sh 1}%
\special{pa 4450 400}%
\special{pa 4430 467}%
\special{pa 4450 453}%
\special{pa 4470 467}%
\special{pa 4450 400}%
\special{fp}%
%
\special{pn 13}%
\special{pa 4000 850}%
\special{pa 4600 850}%
\special{fp}%
%
\special{pn 13}%
\special{pa 4600 750}%
\special{pa 4000 750}%
\special{fp}%
\special{sh 1}%
\special{pa 4000 750}%
\special{pa 4067 770}%
\special{pa 4053 750}%
\special{pa 4067 730}%
\special{pa 4000 750}%
\special{fp}%
%
\special{pn 13}%
\special{ar 4600 800 50 50 4.7123890 1.5707963}%
%
\special{pn 8}%
\special{pa 4000 800}%
\special{pa 4800 800}%
\special{dt 0.045}%
\end{picture}}%

%% file: move_along_Gamma2.tex
{\unitlength 0.1in%
\begin{picture}(37.3500,9.2700)(-1.1500,-13.2700)%
%
\special{pn 13}%
\special{pa 600 900}%
\special{pa 600 1200}%
\special{fp}%
%
\special{pn 13}%
\special{pa 600 700}%
\special{pa 600 400}%
\special{fp}%
\special{sh 1}%
\special{pa 600 400}%
\special{pa 580 467}%
\special{pa 600 453}%
\special{pa 620 467}%
\special{pa 600 400}%
\special{fp}%
%
\special{pn 8}%
\special{pa 600 800}%
\special{pa 1400 800}%
\special{dt 0.045}%
%
\special{pn 13}%
\special{pa 1200 1200}%
\special{pa 1200 400}%
\special{fp}%
\special{sh 1}%
\special{pa 1200 400}%
\special{pa 1180 467}%
\special{pa 1200 453}%
\special{pa 1220 467}%
\special{pa 1200 400}%
\special{fp}%
\put(14.2000,-7.8000){\makebox(0,0)[lb]{$\Gamma$}}%
\put(7.0000,-5.0000){\makebox(0,0){$c_0$}}%
\put(13.0000,-5.0000){\makebox(0,0){$c_0$}}%
\put(12.0000,-8.0000){\makebox(0,0){$\bullet$}}%
\put(12.2000,-8.2000){\makebox(0,0)[lt]{$x$}}%
%
\special{pn 8}%
\special{pa 2800 800}%
\special{pa 3600 800}%
\special{dt 0.045}%
%
\special{pn 13}%
\special{pa 3400 400}%
\special{pa 3400 1200}%
\special{fp}%
\special{sh 1}%
\special{pa 3400 1200}%
\special{pa 3420 1133}%
\special{pa 3400 1147}%
\special{pa 3380 1133}%
\special{pa 3400 1200}%
\special{fp}%
\put(31.0000,-14.0000){\makebox(0,0){the opposite sides of $\Gamma\iff$inverse}}%
\put(34.2000,-8.2000){\makebox(0,0)[lt]{$x$}}%
\put(34.0000,-8.0000){\makebox(0,0){$\bullet$}}%
\put(29.0000,-5.0000){\makebox(0,0){$c_0$}}%
\put(35.0000,-5.0000){\makebox(0,0){$c_0$}}%
\put(36.2000,-7.8000){\makebox(0,0)[lb]{$\Gamma$}}%
%
\special{pn 13}%
\special{ar 700 800 300 50 4.7123890 1.5707963}%
%
\special{pn 13}%
\special{ar 700 900 100 50 3.1415927 4.7123890}%
%
\special{pn 13}%
\special{ar 700 700 100 50 1.5707963 3.1415927}%
%
\special{pn 13}%
\special{pa 2800 900}%
\special{pa 2800 1200}%
\special{fp}%
%
\special{pn 13}%
\special{ar 2900 900 100 50 3.1415927 4.7123890}%
%
\special{pn 13}%
\special{ar 2900 800 300 50 4.7123890 1.5707963}%
%
\special{pn 13}%
\special{pa 2800 700}%
\special{pa 2800 400}%
\special{fp}%
\special{sh 1}%
\special{pa 2800 400}%
\special{pa 2780 467}%
\special{pa 2800 453}%
\special{pa 2820 467}%
\special{pa 2800 400}%
\special{fp}%
%
\special{pn 13}%
\special{ar 2900 700 100 50 1.5707963 3.1415927}%
%
\special{pn 8}%
\special{pa 1000 700}%
\special{pa 1300 700}%
\special{fp}%
\special{sh 1}%
\special{pa 1300 700}%
\special{pa 1233 680}%
\special{pa 1247 700}%
\special{pa 1233 720}%
\special{pa 1300 700}%
\special{fp}%
%
\special{pn 8}%
\special{pa 3200 700}%
\special{pa 3500 700}%
\special{fp}%
\special{sh 1}%
\special{pa 3500 700}%
\special{pa 3433 680}%
\special{pa 3447 700}%
\special{pa 3433 720}%
\special{pa 3500 700}%
\special{fp}%
\put(9.0000,-14.0000){\makebox(0,0){the same side of $\Gamma\iff$ direct}}%
\put(8.0000,-10.0000){\makebox(0,0)[lt]{$\gamma(0)$}}%
%
\special{pn 4}%
\special{pa 800 1000}%
\special{pa 600 800}%
\special{fp}%
%
\special{pn 4}%
\special{pa 3000 1000}%
\special{pa 2800 800}%
\special{fp}%
\put(30.0000,-10.0000){\makebox(0,0)[lt]{$\gamma(0)$}}%
\end{picture}}%

%% file: 4_General_formula.tex
\section{The generalized connected sum formula}
\subsection{The statement}
For any $c_0,c_1\in\GI$ in general position, there exists a generic homotopy $h\colon[0,1]\to\GI$ starting from $h(0)=c_0$ such that
\begin{itemize}
\item
	$h(s)\not\in\Sigma$ for any $s\in[0,1]$ (hence $h$ is essentially a parallel translation),
\item
	$h(1)$ and $c_1$ are separated,
\item
	except for finitely many instances $h(s)$ and $c_1$ are always in general position, and at the instances there exists exactly one generic tangency or generic triple point involving both $h(s)$ and $c_1$.
\end{itemize}
We call such a homotopy a \emph{generic separating homotopy} of $c_0$ from $c_1$.

\begin{defn}[\cite{MendesRomero}]\label{def:T^pm}
For $c_0,c_1\in\GI$ in general position, let $T^+(c_0,c_1)$ be the number of direct tangencies between $c_0$ and $c_1$ that occur in a generic separating homotopy of $c_0$ from $c_1$.
Each direct tangency is counted as $+1$ (resp.\ $-1$) if it is positive (resp.\ negative).
Similarly define $T^-(c_0,c_1)$ as the number of inverse tangencies, counted with signs, between $c_0$ and $c_1$ that occur in a generic separating homotopy of $c_0$ from $c_1$.

For unoriented generic closed curves $C_0$ and $C_1$ in general position and a bridge $\Gamma$ between them, we choose compatible parameters $c_i$ of $C_i$ with respect to $\Gamma$ and define $T^{\pm}_{\Gamma}(C_0,C_1):=T^{\pm}(c_0,c_1)$.
\end{defn}

\begin{defn}[\cite{MendesRomero}]\label{def:T^St}
Let $C_0$ and $C_1$ be unoriented generic closed curves in general position and $\Gamma$ a bridge between them.
Define $T^{\St}_{\Gamma}(C_0,C_1)$ as the number of triple points occurring in a generic separating homotopy of $C_0$ from $C_1$, counted with signs.
Here the sign $(-1)^q$ of each newborn triangle is determined by extending the homotopy to that of $C_0+_{\Gamma}C_1$ in any way.
\end{defn}

\begin{rem}
Well-definedness of $T^{\pm}$ and $T^{\St}_{\Gamma}$ is proved in \cite{MendesRomero}.
Topological descriptions of them given in \S\ref{s:T} can be seen as alternative proofs.

$T^{\pm}_{\Gamma}(C_0,C_1)$ and $T^{\St}_{\Gamma}(C_0,C_1)$ are independent of the choice of compatible parameters.
To define $T^{\St}_{\Gamma}$ parameters $c_0,c_1$ are not enough;
to determine a cyclic order of the edges of a newborn triangle we only need (the endpoints of) a bridge $\Gamma$ since not all the three edges of a newborn triangle are contained in a single curve $C_i$, $i=0$ or $1$.
\end{rem}

Now we can state the generalized connected sum formula in the most general form, up to computations of $T^{\pm}_{\Gamma}$ and $T^{\St}_{\Gamma}$.

\begin{thm}\label{thm:main_J}
Let $C_0,C_1$ be generic closed curves in general position and $\Gamma$ a bridge between them, with $n_{\Gamma}$ double points.
Let $R_i\subset\R^2\setminus C_i$ ($i=0,1$) be the connected component that contains $b_{1-i}$ (see Notation~\ref{notation}).
Then
\begin{multline*}
 J^{\pm}(C_0+_{\Gamma}C_1)=J^{\pm}(C_0)+J^{\pm}(C_1)\mp 2T^{\pm}_{\Gamma}(C_0,C_1)+\ind_{v_0}(C_0)+\ind_{v_1}(C_1)\\
 +\ind^{R_0}_{v_0}(C_0)+\ind^{R_1}_{v_1}(C_1)\pm 2n_{\Gamma}\pm\abs{\Int\Gamma\cap(C_0\cup C_1)}.
\end{multline*}
\end{thm}
\begin{thm}\label{thm:main_St}
Let $C_0,C_1$ and $\Gamma$ be as in Theorem~\ref{thm:main_J}.
Then
\[
 \St(C_0+_{\Gamma}C_1)=\St(C_0)+\St(C_1)-T^{\St}_{\Gamma}(C_0,C_1)+2\,\ind_{v_0}(C_0)\ind_{v_1}(C_1)
 -2\sum s(x)s(y),
\]
where the sum in the right hand side is taken over all the pairs $(x,y)$, $x\in\Int\Gamma\cap C_0$ and $y\in\Int\Gamma\cap C_1$,
such that $b_0,x,y$ and $b_1$ are placed on $\Gamma$ in this order.
\end{thm}

\begin{rem}\label{rem:how_to_recover_formulas}
Theorem~\ref{thm:main_J} recovers and extends Proposition~\ref{prop:strange_J} to any strange sums, as the special case that $T^{\pm}_{\Gamma}(C_0,C_1)=0$ and $\ind^{R_i}_{v_i}(C_i)=0$;
in Proposition~\ref{prop:strange_J}, $C_0,C_1$ are separated and $R_i$ is the unbounded component $R^{\infty}_i$ of $\R^2\setminus C_i$.
If moreover $C_0+_{\Gamma}C_1$ is a connected sum, then Arnold's formula for $J^{\pm}$ (Theorem~\ref{thm:Arnold_additivity} (1)) is recovered;
in this case $n_{\Gamma}=\abs{\Int\Gamma\cap(C_0\cup C_1)}=0$ and $\ind_{v_i}(C_i)=0$ because $R_{v_i}=R^{\infty}_i$.

Similarly Theorem~\ref{thm:main_St} recovers Arnold's formula for $\St$ (Theorem~\ref{thm:Arnold_additivity} (2));
in this case $T^{\St}_{\Gamma}(C_0,C_1)=0$, and $b_0,x,y,b_1$ are always placed on $\Gamma$ in this order for any $x\in\Int\Gamma\cap C_0$ and $y\in\Int\Gamma\cap C_1$.
Thus
\begin{align*}
 \sum s(x)s(y)&=\sum_{x\in\Int\Gamma\cap C_0}s(x)\sum_{y\in\Int\Gamma\cap C_1}s(y)\\
 &=\bigl(\ind_{v_0}(C_0)-\ind^{R^{\infty}_0}_{v_0}(C_0)\bigr)\bigl(\ind_{v_1}(C_1)-\ind^{R^{\infty}_1}_{v_1}(C_1)\bigr)\\
 &=\ind_{v_0}(C_0)\ind_{v_1}(C_1).
\end{align*}
The second equality follows from Lemma~\ref{lem:adjacent_indices}, and the last one holds since $\ind^{R^{\infty}_i}(C_i)=0$.
\end{rem}

\subsection{Mendes de Jesus-Romero Fuster formulas}\label{ss:MendesRomero}
The following Theorem~\ref{thm:MendesRomero} due to Mendes de Jesus and Romero Fuster \cite{MendesRomero} generalizes Theorem~\ref{thm:Arnold_additivity} in a different way from Proposition~\ref{prop:strange_J}.
This is the second key step for the proof of Theorems~\ref{thm:main_J}, \ref{thm:main_St}.

\begin{thm}[{\cite[Theorems 1, 2]{MendesRomero}}]\label{thm:MendesRomero}
Let $C_0,C_1$ be unoriented generic closed curves in general position.
Suppose that we are given a bridge $\Gamma$ between them satisfying $\Int\Gamma\cap(C_0\cup C_1)=\emptyset$ and $D(\Gamma)=\emptyset$ (see Figure~\ref{fig:sum_MendesRomero}, the left).
Then we have
\begin{enumerate}[(1)]
\item
	$J^{\pm}(C_0+_{\Gamma}C_1)=J^{\pm}(C_0)+J^{\pm}(C_1)\mp 2T^{\pm}_{\Gamma}(C_0,C_1)+2\,\ind_{v_0}(C_0)+2\,\ind_{v_1}(C_1)$,
\item
	$\St(C_0+_{\Gamma}C_1)=\St(C_0)+\St(C_1)-T^{\St}_{\Gamma}(C_0,C_1)+2\,\ind_{v_0}(C_0)\ind_{v_1}(C_1)$.
\end{enumerate}
\end{thm}
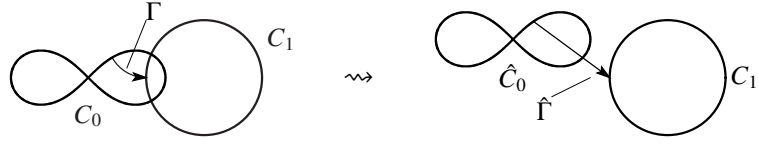
\begin{figure}
\begin{center}
\input{sum_MendesRomero}
\end{center}
\caption{$\Gamma$ in Theorem~\ref{thm:MendesRomero} and separation of $C_0$ from $C_1$ along $\Gamma$}
\label{fig:sum_MendesRomero}
\end{figure}

\begin{rem}
Under the assumption of Theorem~\ref{thm:MendesRomero} we have $n_{\Gamma}=\abs{\Int\Gamma\cap(C_0\cup C_1)}=0$ and $\ind^{R_i}_{v_i}(C_i)=\ind_{v_i}(C_i)$; the last equality holds because $b_{1-i}\in R_{v_i}$ in this case.
Thus Theorem~\ref{thm:main_J} generalizes Theorem~\ref{thm:MendesRomero} (1).
Similarly Theorem~\ref{thm:main_St} generalizes Theorem~\ref{thm:MendesRomero} (2) since under the assumption of Theorem~\ref{thm:MendesRomero}, the last term in Theorem~\ref{thm:main_St} is a sum over the empty-set.

Theorem~\ref{thm:MendesRomero} in fact slightly generalizes \cite[Theorems 1, 2]{MendesRomero}.
In \cite{MendesRomero} closed curves are in advance given orientations and bridges $\Gamma$ should be such that the orientations are compatible with respect to $\Gamma$.
Since we begin with unoriented curves, we can put $\Gamma$ anywhere (as long as it satisfies the assumptions).

In private communications the authors of \cite{MendesRomero} told us that their formula for $\St$ in \cite{MendesRomero} contains an error (in the ``index part'') and Theorem~\ref{thm:MendesRomero} corrects it.
The formula for $J^{\pm}$ in \cite{MendesRomero} is correct;
the difference in the sign of $T^{\pm}_{\Gamma}$ is due to our definition of $T^-_{\Gamma}$.
\end{rem}

\subsubsection{Proof of Theorem~\ref{thm:MendesRomero} (1)}
Performing a generic homotopy if necessary, we may assume that $\Gamma$ is a straight segment and the parallel translation of $C_0$ to the direction $-v_0$ is a generic separating homotopy.
Denote by $\hat{C}_0$ the resulting closed curve which is equivalent to $C_0$ and separated from $C_1$, and let $\hat{\Gamma}$ be the bridge obtained by extending $\Gamma$ so that its endpoint $\hat{b}_0\in\hat{C}_0$ corresponds to $b_0\in C_0$ and $\hat{b}_1=b_1\in C_1$.
See Figure~\ref{fig:sum_MendesRomero}, the right.

\begin{figure}
\begin{center}
\input{orientation_1}
\end{center}
\caption{The point $x_1$ (resp.\ $x_2,x_3,x_4$) contributes to $l^+_0$ (resp.\ $l^-_0,l^-_1,l^+_1$)}
\label{fig:orientation_1}
\end{figure}
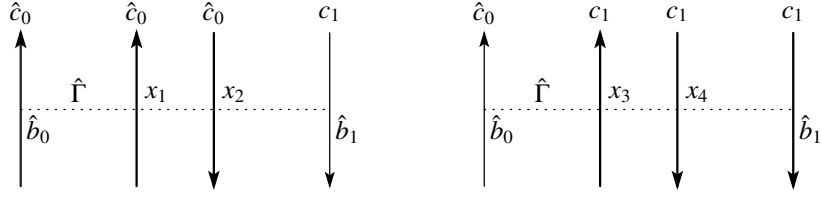

Denoting $\hat{C}_1:=C_1$ for coherence of notations, we define
\begin{equation}\label{eq:def_l^pm_i}
 l_i^{\pm}:=\abs{\{x\in\Int\hat{\Gamma}\cap\hat{C}_i\mid s(x)=\pm 1\}},
\end{equation}
see Figure~\ref{fig:orientation_1}.
These values are independent of the choice of compatible parameters.

\begin{lem}\label{lem:separation_parallel_translation_J}
Let $\hat{C}_0,C_1$ and $\hat{\Gamma}$ be as above.
Then
\[
 J^{\pm}(\hat{C}_0+_{\hat{\Gamma}}C_1)=J^{\pm}(C_0+_{\Gamma}C_1)\pm 2T^{\pm}_{\Gamma}(C_0,C_1)\pm 2(l^{\mp}_0+l^{\mp}_1).
\]
\end{lem}
\begin{proof}
The parallel translation of $C_0$ to the direction $-v_0$ is naturally extended to a generic homotopy $h$ from $C_0+_{\Gamma}C_1$ to $\hat{C}_0+_{\hat{\Gamma}}C_1$.
All the self-tangency points in the homotopy $h$ involve both $C_0$ and $C_1$ because $h$ is just a parallel translation of $C_0$.
The homotopy $h$ produces roughly $T^{\pm}_{\Gamma}(C_0,C_1)$ direct / inverse self-tangencies of $C_0+_{\Gamma}C_1$ (counted with signs).
But we have punched ``holes'' on $C_0$ and $C_1$ near $b_0$ and $b_1$ to make the generalized connected sum, and these holes may change the numbers of tangencies.

For example, let $x\in\Int\hat{\Gamma}\cap\hat{C}_i$ ($i=0,1$) be such that $s(x)=-1$.
Then $\hat{C}_i$ and $\hat{C}_{1-i}$ point the same side of $\hat{\Gamma}$ at $x$ and $\hat{b}_{1-i}$.
As we see in Figure~\ref{fig:change_T^pm} (the case $i=1$), an arc in $C_i$ near $x$ relates to one more positive direct tangency with $C_{1-i}$ in the homotopy $h$ than in the separating homotopy.
\begin{figure}
\begin{center}
\input{change_T_pm}
\end{center}
\caption{Tangencies before/after taking connected sum}
\label{fig:change_T^pm}
\end{figure}
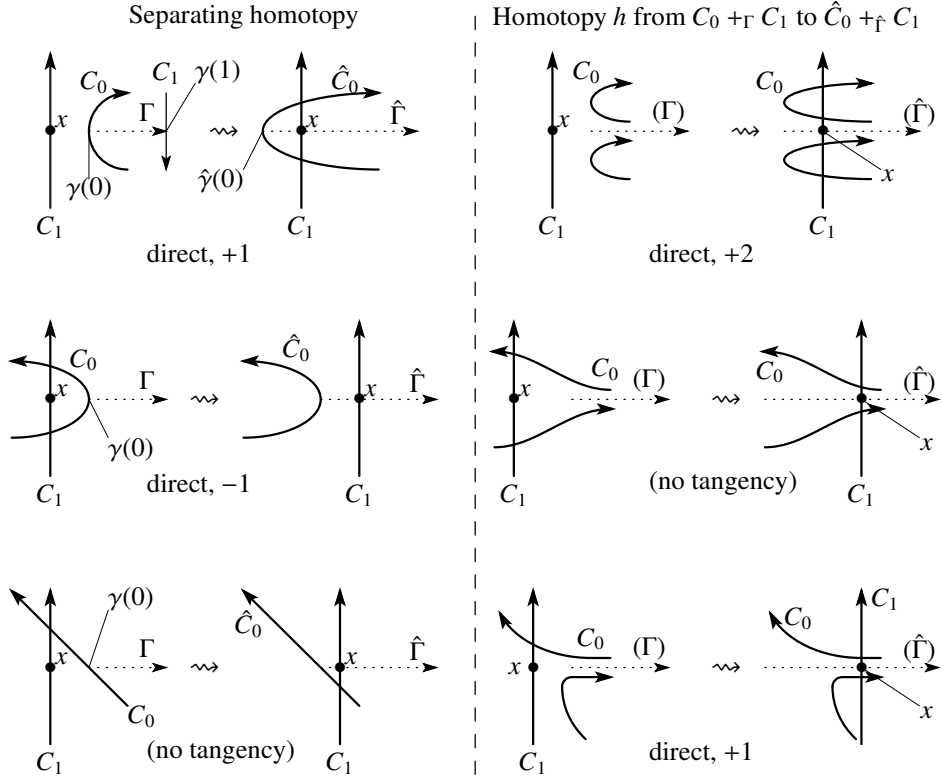
The number of such $x$'s is $l^-_i$.
Thus $T^+_{\Gamma}(C_0,C_1)+(l^-_0+l^-_1)$ direct self-tangencies occur in $h$.

Similarly we see that, for $x\in\Int\hat{\Gamma}\cap\hat{C}_i$ ($i=0,1$) such that $s(x)=+1$, an arc in $C_i$ near $x$ relates to one more positive inverse tangency with $C_{1-i}$ in the homotopy $h$ than in the separating homotopy.
Thus $T^-_{\Gamma}(C_0,C_1)+(l^+_0+l^+_1)$ inverse self-tangencies occur in $h$.
\end{proof}

\begin{lem}\label{lem:separation_parallel_translation_ind}
Let $\hat{C}_0$, $\hat{C}_1:=C_1$ and $\hat{\Gamma}$ be as above and $\hat{v}_i$ be a normal vector on $\hat{C}_i$ determined by $\hat{\Gamma}$.
Then
\[
 \ind_{\hat{v}_0}(\hat{C}_0)=l^+_0-l^-_0,\quad
 \ind_{v_1}(C_1)=l^+_1-l^-_1.
\]
\end{lem}
\begin{proof}
They follow from Lemma~\ref{lem:adjacent_indices} and the fact that $\hat{v}_i\in R^{\infty}_{1-i}$, whose index is zero.
\end{proof}

\begin{proof}[Proof of Theorem~\ref{thm:MendesRomero} (1)]
Let $\hat{C}_0$ and $\hat{\Gamma}$ be as above.
They satisfy all the assumptions in Proposition~\ref{prop:strange_J}, and since $n_{\hat{\Gamma}}=0$ we have
\[
 J^{\pm}(\hat{C}_0+_{\hat{\Gamma}}C_1)
 =J^{\pm}(\hat{C}_0)+J^{\pm}(C_1)+\ind_{\hat{v}_0}(\hat{C}_0)+\ind_{v_1}(C_1)
 \pm\abs{\Int\hat{\Gamma}\cap(\hat{C}_0\cup C_1)}.
\]
Because $\hat{C}_0$ is a parallel translation of $C_0$, we have $J^{\pm}(\hat{C}_0)=J^{\pm}(C_0)$ and $\ind_{\hat{v}_0}(\hat{C}_0)=\ind_{v_0}(C_0)$.
Using Lemma~\ref{lem:separation_parallel_translation_J} and $\abs{\Int\hat{\Gamma}\cap(\hat{C}_0\cup C_1)}=l^+_0+l^-_0+l^+_1+l^-_1$, we obtain
\begin{multline*}
 J^{\pm}(C_0+_{\Gamma}C_1)\pm 2T^{\pm}_{\Gamma}(C_0,C_1)\pm 2(l^{\mp}_0+l^{\mp}_1)\\
 =J^{\pm}(C_0)+J^{\pm}(C_1)+\ind_{v_0}(C_0)+\ind_{v_1}(C_1)\pm(l^+_0+l^-_0+l^+_1+l^-_1).
\end{multline*}
Therefore
\begin{multline*}
 J^{\pm}(C_0+_{\Gamma}C_1)=J^{\pm}(C_0)+J^{\pm}(C_1)\mp 2T^{\pm}_{\Gamma}(C_0,C_1)\\
 +\ind_{v_0}(C_0)+\ind_{v_1}(C_1)+(l^+_0-l^-_0)+(l^+_1-l^-_1),
\end{multline*}
and Lemma~\ref{lem:separation_parallel_translation_ind} completes the proof.
\end{proof}

\subsubsection{Proof of Theorem~\ref{thm:MendesRomero} (2)}
Let $\hat{C}_0$ and $\hat{\Gamma}$ be as above.
The idea is the same as above;
describe $\St(\hat{C}_0+_{\hat{\Gamma}}C_1)$ using $\St(C_0+_{\Gamma}C_1)$ and $T^{\St}_{\Gamma}(C_0,C_1)$ and substitute it to
\begin{equation}\label{eq:Arnold_additivity_St}
 \St(\hat{C}_0+_{\hat{\Gamma}}C_1)=\St(C_0)+\St(C_1),
\end{equation}
that follows from Theorem~\ref{thm:Arnold_additivity}.
To do this we need to count the number of triple points occurring in a homotopy $h$ from $C_0+_{\Gamma}C_1$ to $\hat{C}_0+_{\hat{\Gamma}}C_1$.

In contrast to the case of $J^{\pm}$, the ``holes'' used to make the generalized connected sum do not affect the number of triple points.
This is because, performing a small generic homotopy if necessary, we may assume that any double points of $C_i$ ($i=0,1$) do not get across small arcs on $C_{1-i}$ near $b_{1-i}$ in the separating homotopy of $C_0$ from $C_1$.

\begin{figure}
\begin{center}
\input{new_triple_point}
\end{center}
\caption{New triple points occurring in $h$}
\label{fig:new_triple_point}
\end{figure}
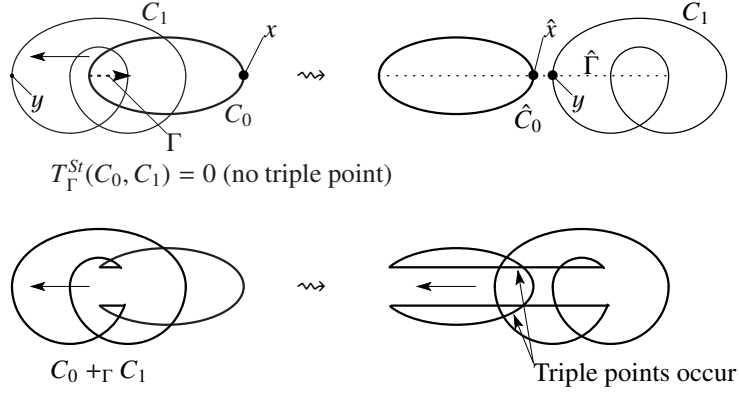

Instead, in the homotopy $h$ from $C_0+_{\Gamma}C_1$ to $\hat{C}_0+_{\Gamma}C_1$, an arc of $C_0$ near $x\in C_0$ that corresponds to some $\hat{x}\in\Int\hat{\Gamma}\cap\hat{C}_0$ makes two triple points with the parallels of $\Gamma$ (see Definition~\ref{def:generalized_connected_sum}) and each arc of $C_1$ near some $y\in\Int\hat{\Gamma}\cap C_1$ (see Figure~\ref{fig:new_triple_point}).
These triple points do not occur in the separating homotopy (i.e., parallel translation of $C_0$).

\begin{figure}
\begin{center}
\input{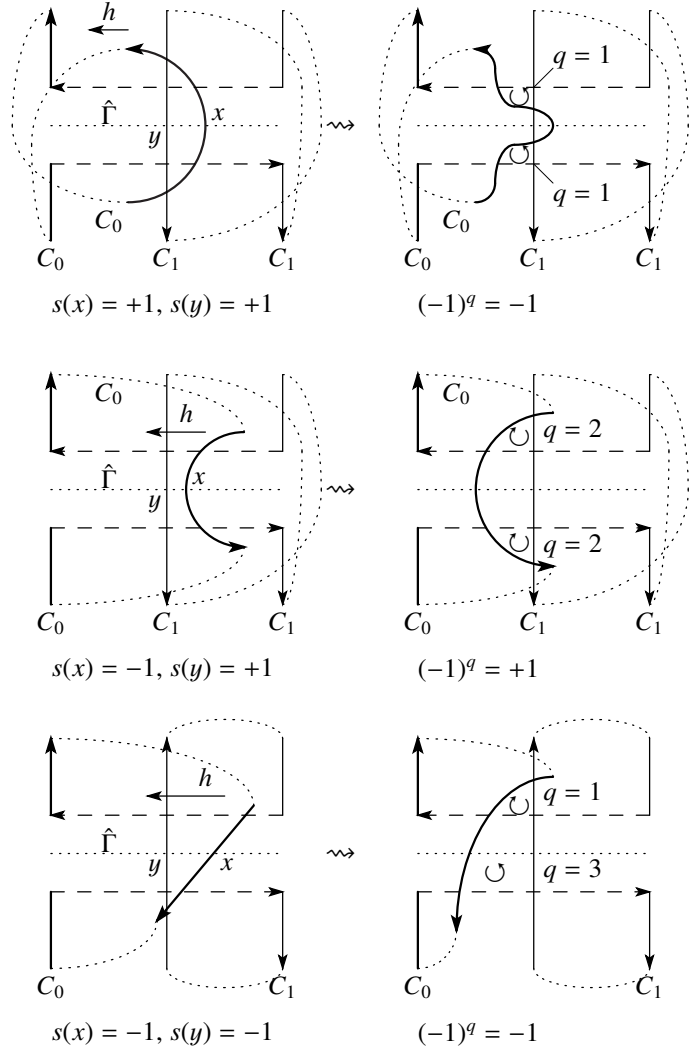}
\end{center}
\caption{New triple points after taking the generalized connected sum}
\label{fig:change_T^St}
\end{figure}

\begin{lem}\label{lem:sign_of_new_triple_points}
The signs of two newborn triangles that occur when $x\in C_0$ (corresponding to some $\hat{x}\in\Int\hat{\Gamma}\cap\hat{C}_0$) gets through an arc of $C_1$ near $y\in\Int\hat{\Gamma}\cap C_1$ in the homotopy $h$ are both $-s(x)s(y)$ (for $s(x)$ see Definition~\ref{def:sign_Gamma_cap_C}).
\end{lem}
\begin{proof}
Check all the cases.
Figure~\ref{fig:change_T^St} shows three typical examples:
In the first one $s(x)=s(y)=+1$ and in the separating homotopy an inverse tangency occurs.
In the second $s(x)=-1$, $s(y)=+1$ and in the separating homotopy a direct tangency occurs.
In the last one $s(x)=s(y)=-1$ and no tangency occurs.
All the other cases can be similarly checked.
\end{proof}

\begin{cor}\label{cor:number_of_triple_points_in_h}
Let $\hat{C}_0$ and $\hat{\Gamma}$ be as above.
Then the number of triple points (counted with signs) occurring in the homotopy $h$ from $C_0+_{\Gamma}C_1$ to $\hat{C}_0+_{\hat{\Gamma}}C_1$ is
\[
 T^{\St}_{\Gamma}(C_0,C_1)-2\,\ind_{v_0}(C_0)\ind_{v_1}(C_1).
\]
\end{cor}

\begin{proof}
Let $l^{\pm}_i$ be as in \eqref{eq:def_l^pm_i}.
By Lemma~\ref{lem:sign_of_new_triple_points}, each small arc near $x\in\Int\hat{\Gamma}\cap\hat{C}_0$ produces $2(l^-_1-l^+_1)s(x)$ triple points in the homotopy $h$.
Summing them up for all the points $x\in\Int\hat{\Gamma}\cap\hat{C}_0$, we see that the number of all the triple points occurring in the homotopy $h$ is
\[
 T^{\St}_{\Gamma}(C_0,C_1)+2(l^-_1-l^+_1)(l^+_0-l^-_0).
\]
Lemma~\ref{lem:adjacent_indices} together with the fact that $\hat{b}_1\in R^{\infty}_0$ and hence $\ind_{v_0}^{R^{\infty}_0}(\hat{C}_0)=0$ implies $l^+_0-l^-_0=\ind_{\hat{v}_0}(\hat{C}_0)=\ind_{v_0}(C_0)$.
Similarly we have $l^-_1-l^+_1=-\ind_{v_1}(C_1)$.
\end{proof}

\begin{proof}[Proof of Theorem~\ref{thm:MendesRomero} (2)]
By Corollary~\ref{cor:number_of_triple_points_in_h} we have
\[
 \St(\hat{C}_0+_{\hat{\Gamma}}C_1)=\St(C_0+_{\Gamma}C_1)+T^{\St}_{\Gamma}(C_0,C_1)-2\,\ind_{v_0}(C_0)\ind_{v_1}(C_1).
\]
Substitute it to \eqref{eq:Arnold_additivity_St}.
\end{proof}

\subsection{Proof of Theorem~\ref{thm:main_J}}
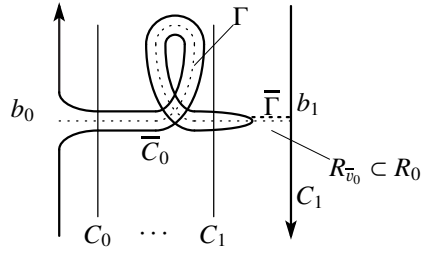
\begin{figure}
\begin{center}
\input{push_appendix_0}
\caption{Pushing-appendix from $C_0$}
\label{fig:push_appendix_0}
\end{center}
\end{figure}

Let $C_0,C_1$ be generic closed curves in general position and $\Gamma$ any bridge between them.
Pushing an appendix from $C_0$ along $\Gamma$ as in Figure~\ref{fig:push_appendix_0}, we obtain a generic curve $\overline{C}_0$ and a bridge $\overline{\Gamma}$ between $\overline{C}_0$ and $C_1$ such that $\overline{C}_0$, $\overline{\Gamma}$ and $C_1$ satisfy the assumptions in Theorem~\ref{thm:MendesRomero} and $C_0+_{\Gamma}C_1=\overline{C}_0+_{\overline{\Gamma}}C_1$.
Thus
\begin{equation}\label{eq:J_overline_plus}
 J^{\pm}(C_0+_{\Gamma}C_1)=J^{\pm}(\overline{C}_0)+J^{\pm}(C_1)\mp 2T^{\pm}_{\overline{\Gamma}}(\overline{C}_0,C_1)+2(\ind_{\overline{v}_0}(\overline{C}_0)+\ind_{v_1}(C_1)),
\end{equation}
where $\overline{\gamma}$ is a parameter of $\overline{\Gamma}$.
We need to compute $J^{\pm}(\overline{C}_0)$, $T^{\pm}_{\Gamma}(\overline{C}_0,C_1)$ and $\ind_{\overline{v}_0}(\overline{C}_0)$.

By Lemma~\ref{lem:compatible_ori} we can choose compatible orientations $c_i$ ($i=0,1$) of $C_i$ with respect to both $\Gamma$ and $v_i$.
Define $l^{\pm}_i$ ($i=0,1$) by
\[
 l_i^{\pm}:=\bigl\lvert\{x\in\Int\Gamma\cap C_i\mid s(x)=\pm 1\}\bigr\rvert
\]
(they are similar to \eqref{eq:def_l^pm_i}).
Then as we have seen in Figure~\ref{label:pushing-appendix_+}, at each point $x\in\Int\Gamma\cap C_0$ with $s(x)=+1$ (resp.\ $-1$), the pushing-appendix from $C_0$ produces a positive direct (resp.\ inverse) self-tangency.
Moreover at each double point of $\Gamma$ the pushing-appendix produces one direct self-tangency and one inverse self-tangency (see Figure~\ref{fig:tangency_from_double_point}).
Thus the pushing-appendix from $C_0$ along $\Gamma$ produces $l_0^\pm+n_{\Gamma}$ direct / inverse self-tangencies, and hence
\begin{equation}\label{eq:J_overline}
 J^{\pm}(\overline{C}_0)=J^{\pm}(C_0)\pm 2(l^{\pm}_0+n_{\Gamma}).
\end{equation}

Let $\overline{c}_0$ be a parameter of $\overline{C}_0$ that is naturally induced by $c_0$.
Then $\overline{c}_0$ is compatible with $\overline{v}_0$.
We notice that $R_0$ is the unique component of $\R^2\setminus C_0$ such that $R_{\overline{v}_0}\subset R_0$.
We have $\ind^{R_{\overline{v}_0}}(\overline{c}_0)=\ind^{R_0}(c_0)$ because the ``appendix part'' does not affect the index.
Moreover we have $\ind^{R_0}_{v_0}(C_0)=\ind_{v_0}(C_0)+l_0^--l_0^+$ by Lemma~\ref{lem:adjacent_indices}.
Since $c_0$ and $\overline{c}_0$ are compatible we obtain
\begin{align}\label{eq:ind_overline}
 \ind_{\overline{v}_0}(\overline{C}_0)&=\ind^{R_{\overline{v}_0}}(\overline{c}_0)=\ind^{R_0}(c_0)=\ind^{R_0}_{v_0}(C_0) \\
 &=\ind_{v_0}(C_0)+l_0^--l_0^+.\notag
\end{align}
Similarly by Lemma~\ref{lem:adjacent_indices}
\begin{equation}\label{eq:ind_R_1}
\ind^{R_1}_{v_1}(C_1)=\ind_{v_1}(C_1)+l^-_1-l^+_1
\end{equation}
since $c_1$ is compatible with $v_1$.
Using \eqref{eq:ind_overline} and \eqref{eq:ind_R_1} we have
\begin{multline}\label{eq:ind}
 2\bigl(\ind_{\overline{v}_0}(\overline{C}_0)+\ind_{v_1}(C_1)\bigr)\\
 =\bigl(\ind_{v_0}^{R_0}(C_0)+\ind_{v_0}(C_0)+l^-_0-l^+_0\bigr)
 +\bigl(\ind_{v_1}(C_1)+\ind^{R_1}_{v_1}(C_1)-l^-_1+l^+_1\bigr).
\end{multline}

To compute $T^{\pm}_{\overline{\Gamma}}(\overline{C}_0,C_1)$,
we may separate $\overline{C}_0$ from $C_1$ as follows;
\begin{enumerate}[(i)]
\item
	first pull the appendix back to $C_0$ along $\Gamma$,
\item
	separate $C_0$ from $C_1$, and then
\item
	push an appendix again from $C_0$ along $\Gamma$ to produce $\overline{C}_0$.
\end{enumerate}
Although this is not a generic separating homotopy in the sense of the beginning of this section, but as we can easily see, such a generic homotopy also passes through $T^{\pm}_{\overline{\Gamma}}(\overline{C}_0,C_1)$ times of tangencies.

In (i) all the tangencies between $C_0$ and $C_1$ are negative.
If $y\in\Int\Gamma\cap C_1$ is such that $s(y)=+1$ (resp.\ $s(y)=-1$), then $C_1$ and $C_0$ point the opposite sides (resp.\ the same side) of $\Gamma$ at $y$ and $b_0$ and an inverse (resp.\ a direct) tangency occurs between $C_0$ and $C_1$ at $y$.
For example, at the point $y=x_4$ in Figure~\ref{fig:orientation_1}, a negative inverse tangency occurs in (i).
Thus the numbers of direct / inverse tangencies in (i) are respectively $-l_1^-$ and $-l_1^+$.
The number of tangencies in (ii) is by definition $T^{\pm}_{\Gamma}(C_0,C_1)$.
In (iii) no tangency occurs between $C_0$ and $C_1$.
Thus
\begin{equation}\label{eq:T^pm_overline}
 T^{\pm}_{\overline{\Gamma}}(\overline{C}_0,C_1)=T^{\pm}_{\Gamma}(C_0,C_1)-l_1^{\mp}.
\end{equation}

Substituting \eqref{eq:J_overline}, \eqref{eq:ind} and \eqref{eq:T^pm_overline} to \eqref{eq:J_overline_plus} we obtain
\begin{multline*}
 J^{\pm}(C_0+_{\Gamma}C_1)
 =J^{\pm}(C_0)+J^{\pm}(C_1)\mp 2T^{\pm}_{\Gamma}(C_0,C_1)+\ind_{v_0}(C_0)+\ind_{v_1}(C_1)\\
 +\ind^{R_0}_{v_0}(C_0)+\ind^{R_1}_{v_1}(C_1)\pm 2(l^{\pm}_0+n_{\Gamma})+(l^-_0-l^+_0-l^-_1+l^+_1)\pm 2l^{\mp}_1.
\end{multline*}
Theorem~\ref{thm:main_J} thus follows because
\[
 \pm 2l^{\pm}_0+(l^-_0-l^+_0-l^-_1+l^+_1)\pm 2l^{\mp}_1=\pm(l^+_0+l^-_0+l^+_1+l^-_1)=\pm\abs{\Int\Gamma\cap(C_0\cup C_1)}.
\]

\subsection{Proof of Theorem~\ref{thm:main_St}}
By Theorem~\ref{thm:MendesRomero} (2) we have
\[
 \St(C_0+_{\Gamma}C_1)=\St(\overline{C}_0)+\St(C_1)-T^{\St}_{\overline{\Gamma}}(\overline{C}_0,C_1)+2\,\ind_{\overline{v}_0}(\overline{C}_0)\ind_{v_1}(C_1),
\]
where $\overline{C}_0$ and $\overline{\Gamma}$ are as above.
We notice that $\St(\overline{C}_0)=\St(C_0)$ because $\overline{C}_0$ can be obtained by a pushing-appendix from $C_0$ in which no triple point occurs.
Moreover $\ind_{\overline{v}_0}(\overline{C}_0)=\ind^{R_0}_{v_0}(C_0)$ by \eqref{eq:ind_overline}.
Thus
\begin{equation}\label{eq:MendesRomero_applied}
 \St(C_0+_{\Gamma}C_1)=\St(C_0)+\St(C_1)-T^{\St}_{\overline{\Gamma}}(\overline{C}_0,C_1)+2\,\ind^{R_0}_{v_0}(C_0)\ind_{v_1}(C_1).
\end{equation}
Suppose $\Int\Gamma\cap C_0=\{x_1,\dots,x_p\}$ and let $R^{(j)}_1\subset\R^2\setminus C_1$ be the component that contains $x_j$.
As we will show in Corollary~\ref{cor:triple_points_change_T^St},
\begin{equation}\label{T^St_S5}
 T^{\St}_{\overline{\Gamma}}(\overline{C}_0,C_1)=T^{\St}_{\Gamma}(C_0,C_1)-2\sum_{1\le j\le p}s(x_j)\ind^{R^{(j)}}_{v_1}(C_1)
\end{equation}
(here $\ind^{R^{(j)}}_{v_1}(C_1)=\ind^{R^{(j)}}(C_1)$ because $c_1$ is compatible with $v_1$).
By Lemma~\ref{lem:adjacent_indices} we have
\begin{equation}\label{eq:index_s}
 \ind^{R_0}_{v_0}(C_0)=\ind_{v_0}(C_0)-\sum_{1\le j\le p}s(x_j).
\end{equation}
Using \eqref{T^St_S5} and \eqref{eq:index_s} we obtain
\begin{multline}\label{eq:T^St+ind}
 -T^{\St}_{\overline{\Gamma}}(\overline{C}_0,C_1)+2\,\ind^{R_0}_{v_0}(C_0)\ind_{v_1}(C_1)\\
 =-T^{\St}_{\Gamma}(C_0,C_1)
  +2\,\ind_{v_0}(C_0)\ind_{v_1}(C_1)\\
  +2\sum_{1\le j\le p}s(x_j)\Bigl(\ind^{R^{(j)}_1}_{v_1}(C_1)-\ind_{v_1}(C_1)\Bigr).
\end{multline}
Again by Lemma~\ref{lem:adjacent_indices},
\begin{equation}\label{eq:sum_s}
 \ind^{R^{(j)}_1}_{v_1}(C_1)-\ind_{v_1}(C_1)=-\sum s(y),
\end{equation}
where $y$ runs over the points in $\Int\Gamma\cap C_1$ sitting between $x_j$ and $\gamma(1)$.
Substituting \eqref{eq:T^St+ind} and \eqref{eq:sum_s} into \eqref{eq:MendesRomero_applied} we complete the proof for $\St$.

%% file: sum_MendesRomero.tex
{\unitlength 0.1in%
\begin{picture}(37.0000,7.3000)(8.0000,-15.0000)%
%
\special{pn 13}%
\special{ar 1800 1200 300 300 0.0000000 6.2831853}%
\put(12.0000,-14.0000){\makebox(0,0){$C_0$}}%
\put(22.0000,-10.0000){\makebox(0,0){$C_1$}}%
%
\special{pn 8}%
\special{ar 1500 1000 200 200 1.5707963 2.6387494}%
%
\special{pn 8}%
\special{pa 1485 1199}%
\special{pa 1500 1200}%
\special{fp}%
\special{sh 1}%
\special{pa 1500 1200}%
\special{pa 1435 1176}%
\special{pa 1447 1196}%
\special{pa 1432 1216}%
\special{pa 1500 1200}%
\special{fp}%
\put(15.0000,-9.0000){\makebox(0,0)[lb]{$\Gamma$}}%
%
\special{pn 4}%
\special{pa 1500 900}%
\special{pa 1400 1150}%
\special{fp}%
\put(26.0000,-12.0000){\makebox(0,0){$\rightsquigarrow$}}%
%
\special{pn 13}%
\special{ar 4200 1200 300 300 0.0000000 6.2831853}%
\put(34.0000,-12.0000){\makebox(0,0){$\hat{C}_0$}}%
\put(46.0000,-12.0000){\makebox(0,0){$C_1$}}%
\put(36.0000,-13.0000){\makebox(0,0)[rt]{$\hat{\Gamma}$}}%
%
\special{pn 4}%
\special{pa 3600 1300}%
\special{pa 3800 1200}%
\special{fp}%
%
\special{pn 8}%
\special{pa 3510 910}%
\special{pa 3900 1200}%
\special{fp}%
\special{sh 1}%
\special{pa 3900 1200}%
\special{pa 3858 1144}%
\special{pa 3857 1168}%
\special{pa 3835 1176}%
\special{pa 3900 1200}%
\special{fp}%
\special{pn 13}%
\special{pa 1600 1200}%
\special{pa 1600 1189}%
\special{pa 1599 1188}%
\special{pa 1599 1181}%
\special{pa 1598 1180}%
\special{pa 1598 1175}%
\special{pa 1597 1174}%
\special{pa 1597 1170}%
\special{pa 1596 1169}%
\special{pa 1596 1166}%
\special{pa 1595 1165}%
\special{pa 1595 1163}%
\special{pa 1594 1162}%
\special{pa 1594 1159}%
\special{pa 1593 1158}%
\special{pa 1593 1156}%
\special{pa 1592 1155}%
\special{pa 1592 1154}%
\special{pa 1591 1153}%
\special{pa 1591 1151}%
\special{pa 1590 1150}%
\special{pa 1590 1149}%
\special{pa 1589 1148}%
\special{pa 1589 1146}%
\special{pa 1588 1145}%
\special{pa 1585 1138}%
\special{pa 1584 1137}%
\special{pa 1583 1134}%
\special{pa 1581 1132}%
\special{pa 1580 1129}%
\special{pa 1578 1127}%
\special{pa 1578 1126}%
\special{pa 1576 1124}%
\special{pa 1576 1123}%
\special{pa 1574 1121}%
\special{pa 1574 1120}%
\special{pa 1571 1117}%
\special{pa 1571 1116}%
\special{pa 1566 1111}%
\special{pa 1566 1110}%
\special{pa 1553 1097}%
\special{pa 1552 1097}%
\special{pa 1547 1092}%
\special{pa 1546 1092}%
\special{pa 1543 1089}%
\special{pa 1542 1089}%
\special{pa 1540 1087}%
\special{pa 1539 1087}%
\special{pa 1537 1085}%
\special{pa 1536 1085}%
\special{pa 1534 1083}%
\special{pa 1533 1083}%
\special{pa 1531 1081}%
\special{pa 1524 1078}%
\special{pa 1522 1076}%
\special{pa 1519 1075}%
\special{pa 1518 1074}%
\special{pa 1516 1074}%
\special{pa 1515 1073}%
\special{pa 1512 1072}%
\special{pa 1511 1071}%
\special{pa 1509 1071}%
\special{pa 1508 1070}%
\special{pa 1507 1070}%
\special{pa 1506 1069}%
\special{pa 1495 1066}%
\special{pa 1494 1065}%
\special{pa 1487 1064}%
\special{pa 1486 1063}%
\special{pa 1483 1063}%
\special{pa 1482 1062}%
\special{pa 1478 1062}%
\special{pa 1477 1061}%
\special{pa 1472 1061}%
\special{pa 1471 1060}%
\special{pa 1464 1060}%
\special{pa 1463 1059}%
\special{pa 1427 1059}%
\special{pa 1426 1060}%
\special{pa 1418 1060}%
\special{pa 1417 1061}%
\special{pa 1406 1062}%
\special{pa 1405 1063}%
\special{pa 1402 1063}%
\special{pa 1401 1064}%
\special{pa 1397 1064}%
\special{pa 1396 1065}%
\special{pa 1393 1065}%
\special{pa 1392 1066}%
\special{pa 1390 1066}%
\special{pa 1389 1067}%
\special{pa 1386 1067}%
\special{pa 1385 1068}%
\special{pa 1374 1071}%
\special{pa 1373 1072}%
\special{pa 1371 1072}%
\special{pa 1370 1073}%
\special{pa 1369 1073}%
\special{pa 1368 1074}%
\special{pa 1363 1075}%
\special{pa 1362 1076}%
\special{pa 1359 1077}%
\special{pa 1358 1078}%
\special{pa 1356 1078}%
\special{pa 1355 1079}%
\special{pa 1348 1082}%
\special{pa 1347 1083}%
\special{pa 1340 1086}%
\special{pa 1339 1087}%
\special{pa 1332 1090}%
\special{pa 1330 1092}%
\special{pa 1327 1093}%
\special{pa 1326 1094}%
\special{pa 1325 1094}%
\special{pa 1323 1096}%
\special{pa 1320 1097}%
\special{pa 1318 1099}%
\special{pa 1317 1099}%
\special{pa 1315 1101}%
\special{pa 1314 1101}%
\special{pa 1312 1103}%
\special{pa 1309 1104}%
\special{pa 1306 1107}%
\special{pa 1305 1107}%
\special{pa 1303 1109}%
\special{pa 1302 1109}%
\special{pa 1300 1111}%
\special{pa 1299 1111}%
\special{pa 1296 1114}%
\special{pa 1295 1114}%
\special{pa 1293 1116}%
\special{pa 1292 1116}%
\special{pa 1288 1120}%
\special{pa 1287 1120}%
\special{pa 1284 1123}%
\special{pa 1283 1123}%
\special{pa 1279 1127}%
\special{pa 1278 1127}%
\special{pa 1274 1131}%
\special{pa 1273 1131}%
\special{pa 1268 1136}%
\special{pa 1267 1136}%
\special{pa 1260 1143}%
\special{pa 1259 1143}%
\special{pa 1251 1151}%
\special{pa 1250 1151}%
\special{pa 1235 1166}%
\special{pa 1234 1166}%
\special{pa 1200 1200}%
\special{pa 1166 1166}%
\special{pa 1165 1166}%
\special{pa 1150 1151}%
\special{pa 1149 1151}%
\special{pa 1141 1143}%
\special{pa 1140 1143}%
\special{pa 1133 1136}%
\special{pa 1132 1136}%
\special{pa 1127 1131}%
\special{pa 1126 1131}%
\special{pa 1122 1127}%
\special{pa 1121 1127}%
\special{pa 1117 1123}%
\special{pa 1116 1123}%
\special{pa 1113 1120}%
\special{pa 1112 1120}%
\special{pa 1108 1116}%
\special{pa 1107 1116}%
\special{pa 1105 1114}%
\special{pa 1104 1114}%
\special{pa 1101 1111}%
\special{pa 1100 1111}%
\special{pa 1098 1109}%
\special{pa 1097 1109}%
\special{pa 1095 1107}%
\special{pa 1094 1107}%
\special{pa 1091 1104}%
\special{pa 1088 1103}%
\special{pa 1086 1101}%
\special{pa 1085 1101}%
\special{pa 1083 1099}%
\special{pa 1082 1099}%
\special{pa 1080 1097}%
\special{pa 1077 1096}%
\special{pa 1075 1094}%
\special{pa 1072 1093}%
\special{pa 1071 1092}%
\special{pa 1070 1092}%
\special{pa 1068 1090}%
\special{pa 1061 1087}%
\special{pa 1060 1086}%
\special{pa 1053 1083}%
\special{pa 1052 1082}%
\special{pa 1045 1079}%
\special{pa 1044 1078}%
\special{pa 1042 1078}%
\special{pa 1041 1077}%
\special{pa 1038 1076}%
\special{pa 1037 1075}%
\special{pa 1032 1074}%
\special{pa 1031 1073}%
\special{pa 1030 1073}%
\special{pa 1029 1072}%
\special{pa 1018 1069}%
\special{pa 1017 1068}%
\special{pa 1015 1068}%
\special{pa 1014 1067}%
\special{pa 1011 1067}%
\special{pa 1010 1066}%
\special{pa 1008 1066}%
\special{pa 1007 1065}%
\special{pa 1004 1065}%
\special{pa 1003 1064}%
\special{pa 999 1064}%
\special{pa 998 1063}%
\special{pa 995 1063}%
\special{pa 994 1062}%
\special{pa 983 1061}%
\special{pa 982 1060}%
\special{pa 974 1060}%
\special{pa 973 1059}%
\special{pa 937 1059}%
\special{pa 936 1060}%
\special{pa 929 1060}%
\special{pa 928 1061}%
\special{pa 923 1061}%
\special{pa 922 1062}%
\special{pa 918 1062}%
\special{pa 917 1063}%
\special{pa 910 1064}%
\special{pa 909 1065}%
\special{pa 906 1065}%
\special{pa 905 1066}%
\special{pa 894 1069}%
\special{pa 893 1070}%
\special{pa 892 1070}%
\special{pa 891 1071}%
\special{pa 889 1071}%
\special{pa 888 1072}%
\special{pa 885 1073}%
\special{pa 884 1074}%
\special{pa 882 1074}%
\special{pa 881 1075}%
\special{pa 878 1076}%
\special{pa 876 1078}%
\special{pa 869 1081}%
\special{pa 867 1083}%
\special{pa 866 1083}%
\special{pa 864 1085}%
\special{pa 863 1085}%
\special{pa 861 1087}%
\special{pa 860 1087}%
\special{pa 858 1089}%
\special{pa 857 1089}%
\special{pa 854 1092}%
\special{pa 853 1092}%
\special{pa 848 1097}%
\special{pa 847 1097}%
\special{pa 834 1110}%
\special{pa 834 1111}%
\special{pa 829 1116}%
\special{pa 829 1117}%
\special{pa 826 1120}%
\special{pa 826 1121}%
\special{pa 824 1123}%
\special{pa 824 1124}%
\special{pa 822 1126}%
\special{pa 822 1127}%
\special{pa 820 1129}%
\special{pa 819 1132}%
\special{pa 817 1134}%
\special{pa 814 1141}%
\special{pa 813 1142}%
\special{pa 812 1145}%
\special{pa 811 1146}%
\special{pa 811 1148}%
\special{pa 810 1149}%
\special{pa 810 1150}%
\special{pa 809 1151}%
\special{pa 809 1153}%
\special{pa 808 1154}%
\special{pa 808 1155}%
\special{pa 807 1156}%
\special{pa 807 1158}%
\special{pa 806 1159}%
\special{pa 806 1162}%
\special{pa 805 1163}%
\special{pa 805 1165}%
\special{pa 804 1166}%
\special{pa 804 1169}%
\special{pa 803 1170}%
\special{pa 803 1174}%
\special{pa 802 1175}%
\special{pa 802 1180}%
\special{pa 801 1181}%
\special{pa 801 1188}%
\special{pa 800 1189}%
\special{pa 800 1211}%
\special{pa 801 1212}%
\special{pa 801 1219}%
\special{pa 802 1220}%
\special{pa 802 1225}%
\special{pa 803 1226}%
\special{pa 803 1230}%
\special{pa 804 1231}%
\special{pa 804 1234}%
\special{pa 805 1235}%
\special{pa 805 1237}%
\special{pa 806 1238}%
\special{pa 806 1241}%
\special{pa 807 1242}%
\special{pa 807 1244}%
\special{pa 808 1245}%
\special{pa 808 1246}%
\special{pa 809 1247}%
\special{pa 809 1249}%
\special{pa 810 1250}%
\special{pa 810 1251}%
\special{pa 811 1252}%
\special{pa 811 1254}%
\special{pa 812 1255}%
\special{pa 815 1262}%
\special{pa 816 1263}%
\special{pa 817 1266}%
\special{pa 819 1268}%
\special{pa 820 1271}%
\special{pa 822 1273}%
\special{pa 822 1274}%
\special{pa 824 1276}%
\special{pa 824 1277}%
\special{pa 826 1279}%
\special{pa 826 1280}%
\special{pa 829 1283}%
\special{pa 829 1284}%
\special{pa 834 1289}%
\special{pa 834 1290}%
\special{pa 847 1303}%
\special{pa 848 1303}%
\special{pa 853 1308}%
\special{pa 854 1308}%
\special{pa 857 1311}%
\special{pa 858 1311}%
\special{pa 860 1313}%
\special{pa 861 1313}%
\special{pa 863 1315}%
\special{pa 864 1315}%
\special{pa 866 1317}%
\special{pa 867 1317}%
\special{pa 869 1319}%
\special{pa 876 1322}%
\special{pa 878 1324}%
\special{pa 881 1325}%
\special{pa 882 1326}%
\special{pa 884 1326}%
\special{pa 885 1327}%
\special{pa 888 1328}%
\special{pa 889 1329}%
\special{pa 891 1329}%
\special{pa 892 1330}%
\special{pa 893 1330}%
\special{pa 894 1331}%
\special{pa 905 1334}%
\special{pa 906 1335}%
\special{pa 913 1336}%
\special{pa 914 1337}%
\special{pa 917 1337}%
\special{pa 918 1338}%
\special{pa 922 1338}%
\special{pa 923 1339}%
\special{pa 928 1339}%
\special{pa 929 1340}%
\special{pa 936 1340}%
\special{pa 937 1341}%
\special{pa 973 1341}%
\special{pa 974 1340}%
\special{pa 982 1340}%
\special{pa 983 1339}%
\special{pa 994 1338}%
\special{pa 995 1337}%
\special{pa 998 1337}%
\special{pa 999 1336}%
\special{pa 1003 1336}%
\special{pa 1004 1335}%
\special{pa 1007 1335}%
\special{pa 1008 1334}%
\special{pa 1010 1334}%
\special{pa 1011 1333}%
\special{pa 1014 1333}%
\special{pa 1015 1332}%
\special{pa 1026 1329}%
\special{pa 1027 1328}%
\special{pa 1029 1328}%
\special{pa 1030 1327}%
\special{pa 1031 1327}%
\special{pa 1032 1326}%
\special{pa 1037 1325}%
\special{pa 1038 1324}%
\special{pa 1041 1323}%
\special{pa 1042 1322}%
\special{pa 1044 1322}%
\special{pa 1045 1321}%
\special{pa 1052 1318}%
\special{pa 1053 1317}%
\special{pa 1060 1314}%
\special{pa 1061 1313}%
\special{pa 1068 1310}%
\special{pa 1070 1308}%
\special{pa 1073 1307}%
\special{pa 1074 1306}%
\special{pa 1075 1306}%
\special{pa 1077 1304}%
\special{pa 1080 1303}%
\special{pa 1082 1301}%
\special{pa 1083 1301}%
\special{pa 1085 1299}%
\special{pa 1086 1299}%
\special{pa 1088 1297}%
\special{pa 1091 1296}%
\special{pa 1094 1293}%
\special{pa 1095 1293}%
\special{pa 1097 1291}%
\special{pa 1098 1291}%
\special{pa 1100 1289}%
\special{pa 1101 1289}%
\special{pa 1104 1286}%
\special{pa 1105 1286}%
\special{pa 1107 1284}%
\special{pa 1108 1284}%
\special{pa 1112 1280}%
\special{pa 1113 1280}%
\special{pa 1116 1277}%
\special{pa 1117 1277}%
\special{pa 1121 1273}%
\special{pa 1122 1273}%
\special{pa 1126 1269}%
\special{pa 1127 1269}%
\special{pa 1132 1264}%
\special{pa 1133 1264}%
\special{pa 1140 1257}%
\special{pa 1141 1257}%
\special{pa 1149 1249}%
\special{pa 1150 1249}%
\special{pa 1165 1234}%
\special{pa 1166 1234}%
\special{pa 1199 1201}%
\special{fp}%
\special{pa 1200 1200}%
\special{pa 1234 1234}%
\special{pa 1235 1234}%
\special{pa 1250 1249}%
\special{pa 1251 1249}%
\special{pa 1259 1257}%
\special{pa 1260 1257}%
\special{pa 1267 1264}%
\special{pa 1268 1264}%
\special{pa 1273 1269}%
\special{pa 1274 1269}%
\special{pa 1278 1273}%
\special{pa 1279 1273}%
\special{pa 1283 1277}%
\special{pa 1284 1277}%
\special{pa 1287 1280}%
\special{pa 1288 1280}%
\special{pa 1292 1284}%
\special{pa 1293 1284}%
\special{pa 1295 1286}%
\special{pa 1296 1286}%
\special{pa 1299 1289}%
\special{pa 1300 1289}%
\special{pa 1302 1291}%
\special{pa 1303 1291}%
\special{pa 1305 1293}%
\special{pa 1306 1293}%
\special{pa 1309 1296}%
\special{pa 1312 1297}%
\special{pa 1314 1299}%
\special{pa 1315 1299}%
\special{pa 1317 1301}%
\special{pa 1318 1301}%
\special{pa 1320 1303}%
\special{pa 1323 1304}%
\special{pa 1325 1306}%
\special{pa 1328 1307}%
\special{pa 1329 1308}%
\special{pa 1330 1308}%
\special{pa 1332 1310}%
\special{pa 1339 1313}%
\special{pa 1340 1314}%
\special{pa 1347 1317}%
\special{pa 1348 1318}%
\special{pa 1355 1321}%
\special{pa 1356 1322}%
\special{pa 1358 1322}%
\special{pa 1359 1323}%
\special{pa 1362 1324}%
\special{pa 1363 1325}%
\special{pa 1368 1326}%
\special{pa 1369 1327}%
\special{pa 1370 1327}%
\special{pa 1371 1328}%
\special{pa 1382 1331}%
\special{pa 1383 1332}%
\special{pa 1385 1332}%
\special{pa 1386 1333}%
\special{pa 1389 1333}%
\special{pa 1390 1334}%
\special{pa 1392 1334}%
\special{pa 1393 1335}%
\special{pa 1396 1335}%
\special{pa 1397 1336}%
\special{pa 1401 1336}%
\special{pa 1402 1337}%
\special{pa 1405 1337}%
\special{pa 1406 1338}%
\special{pa 1417 1339}%
\special{pa 1418 1340}%
\special{pa 1426 1340}%
\special{pa 1427 1341}%
\special{pa 1463 1341}%
\special{pa 1464 1340}%
\special{pa 1471 1340}%
\special{pa 1472 1339}%
\special{pa 1477 1339}%
\special{pa 1478 1338}%
\special{pa 1482 1338}%
\special{pa 1483 1337}%
\special{pa 1490 1336}%
\special{pa 1491 1335}%
\special{pa 1494 1335}%
\special{pa 1495 1334}%
\special{pa 1506 1331}%
\special{pa 1507 1330}%
\special{pa 1508 1330}%
\special{pa 1509 1329}%
\special{pa 1511 1329}%
\special{pa 1512 1328}%
\special{pa 1515 1327}%
\special{pa 1516 1326}%
\special{pa 1518 1326}%
\special{pa 1519 1325}%
\special{pa 1522 1324}%
\special{pa 1524 1322}%
\special{pa 1531 1319}%
\special{pa 1533 1317}%
\special{pa 1534 1317}%
\special{pa 1536 1315}%
\special{pa 1537 1315}%
\special{pa 1539 1313}%
\special{pa 1540 1313}%
\special{pa 1542 1311}%
\special{pa 1543 1311}%
\special{pa 1546 1308}%
\special{pa 1547 1308}%
\special{pa 1552 1303}%
\special{pa 1553 1303}%
\special{pa 1566 1290}%
\special{pa 1566 1289}%
\special{pa 1571 1284}%
\special{pa 1571 1283}%
\special{pa 1574 1280}%
\special{pa 1574 1279}%
\special{pa 1576 1277}%
\special{pa 1576 1276}%
\special{pa 1578 1274}%
\special{pa 1578 1273}%
\special{pa 1580 1271}%
\special{pa 1581 1268}%
\special{pa 1583 1266}%
\special{pa 1586 1259}%
\special{pa 1587 1258}%
\special{pa 1588 1255}%
\special{pa 1589 1254}%
\special{pa 1589 1252}%
\special{pa 1590 1251}%
\special{pa 1590 1250}%
\special{pa 1591 1249}%
\special{pa 1591 1247}%
\special{pa 1592 1246}%
\special{pa 1592 1245}%
\special{pa 1593 1244}%
\special{pa 1593 1242}%
\special{pa 1594 1241}%
\special{pa 1594 1238}%
\special{pa 1595 1237}%
\special{pa 1595 1235}%
\special{pa 1596 1234}%
\special{pa 1596 1231}%
\special{pa 1597 1230}%
\special{pa 1597 1226}%
\special{pa 1598 1225}%
\special{pa 1598 1220}%
\special{pa 1599 1219}%
\special{pa 1599 1212}%
\special{pa 1600 1211}%
\special{pa 1600 1201}%
\special{fp}%
\special{pn 13}%
\special{pa 3800 1000}%
\special{pa 3800 989}%
\special{pa 3799 988}%
\special{pa 3799 981}%
\special{pa 3798 980}%
\special{pa 3798 975}%
\special{pa 3797 974}%
\special{pa 3797 970}%
\special{pa 3796 969}%
\special{pa 3796 966}%
\special{pa 3795 965}%
\special{pa 3795 963}%
\special{pa 3794 962}%
\special{pa 3794 959}%
\special{pa 3793 958}%
\special{pa 3793 956}%
\special{pa 3792 955}%
\special{pa 3792 954}%
\special{pa 3791 953}%
\special{pa 3791 951}%
\special{pa 3790 950}%
\special{pa 3790 949}%
\special{pa 3789 948}%
\special{pa 3789 946}%
\special{pa 3788 945}%
\special{pa 3785 938}%
\special{pa 3784 937}%
\special{pa 3783 934}%
\special{pa 3781 932}%
\special{pa 3780 929}%
\special{pa 3778 927}%
\special{pa 3778 926}%
\special{pa 3776 924}%
\special{pa 3776 923}%
\special{pa 3774 921}%
\special{pa 3774 920}%
\special{pa 3771 917}%
\special{pa 3771 916}%
\special{pa 3766 911}%
\special{pa 3766 910}%
\special{pa 3753 897}%
\special{pa 3752 897}%
\special{pa 3747 892}%
\special{pa 3746 892}%
\special{pa 3743 889}%
\special{pa 3742 889}%
\special{pa 3740 887}%
\special{pa 3739 887}%
\special{pa 3737 885}%
\special{pa 3736 885}%
\special{pa 3734 883}%
\special{pa 3733 883}%
\special{pa 3731 881}%
\special{pa 3724 878}%
\special{pa 3722 876}%
\special{pa 3719 875}%
\special{pa 3718 874}%
\special{pa 3716 874}%
\special{pa 3715 873}%
\special{pa 3712 872}%
\special{pa 3711 871}%
\special{pa 3709 871}%
\special{pa 3708 870}%
\special{pa 3707 870}%
\special{pa 3706 869}%
\special{pa 3695 866}%
\special{pa 3694 865}%
\special{pa 3687 864}%
\special{pa 3686 863}%
\special{pa 3683 863}%
\special{pa 3682 862}%
\special{pa 3678 862}%
\special{pa 3677 861}%
\special{pa 3672 861}%
\special{pa 3671 860}%
\special{pa 3664 860}%
\special{pa 3663 859}%
\special{pa 3627 859}%
\special{pa 3626 860}%
\special{pa 3618 860}%
\special{pa 3617 861}%
\special{pa 3606 862}%
\special{pa 3605 863}%
\special{pa 3602 863}%
\special{pa 3601 864}%
\special{pa 3597 864}%
\special{pa 3596 865}%
\special{pa 3593 865}%
\special{pa 3592 866}%
\special{pa 3590 866}%
\special{pa 3589 867}%
\special{pa 3586 867}%
\special{pa 3585 868}%
\special{pa 3574 871}%
\special{pa 3573 872}%
\special{pa 3571 872}%
\special{pa 3570 873}%
\special{pa 3569 873}%
\special{pa 3568 874}%
\special{pa 3563 875}%
\special{pa 3562 876}%
\special{pa 3559 877}%
\special{pa 3558 878}%
\special{pa 3556 878}%
\special{pa 3555 879}%
\special{pa 3548 882}%
\special{pa 3547 883}%
\special{pa 3540 886}%
\special{pa 3539 887}%
\special{pa 3532 890}%
\special{pa 3530 892}%
\special{pa 3527 893}%
\special{pa 3526 894}%
\special{pa 3525 894}%
\special{pa 3523 896}%
\special{pa 3520 897}%
\special{pa 3518 899}%
\special{pa 3517 899}%
\special{pa 3515 901}%
\special{pa 3514 901}%
\special{pa 3512 903}%
\special{pa 3509 904}%
\special{pa 3506 907}%
\special{pa 3505 907}%
\special{pa 3503 909}%
\special{pa 3502 909}%
\special{pa 3500 911}%
\special{pa 3499 911}%
\special{pa 3496 914}%
\special{pa 3495 914}%
\special{pa 3493 916}%
\special{pa 3492 916}%
\special{pa 3488 920}%
\special{pa 3487 920}%
\special{pa 3484 923}%
\special{pa 3483 923}%
\special{pa 3479 927}%
\special{pa 3478 927}%
\special{pa 3474 931}%
\special{pa 3473 931}%
\special{pa 3468 936}%
\special{pa 3467 936}%
\special{pa 3460 943}%
\special{pa 3459 943}%
\special{pa 3451 951}%
\special{pa 3450 951}%
\special{pa 3435 966}%
\special{pa 3434 966}%
\special{pa 3400 1000}%
\special{pa 3366 966}%
\special{pa 3365 966}%
\special{pa 3350 951}%
\special{pa 3349 951}%
\special{pa 3341 943}%
\special{pa 3340 943}%
\special{pa 3333 936}%
\special{pa 3332 936}%
\special{pa 3327 931}%
\special{pa 3326 931}%
\special{pa 3322 927}%
\special{pa 3321 927}%
\special{pa 3317 923}%
\special{pa 3316 923}%
\special{pa 3313 920}%
\special{pa 3312 920}%
\special{pa 3308 916}%
\special{pa 3307 916}%
\special{pa 3305 914}%
\special{pa 3304 914}%
\special{pa 3301 911}%
\special{pa 3300 911}%
\special{pa 3298 909}%
\special{pa 3297 909}%
\special{pa 3295 907}%
\special{pa 3294 907}%
\special{pa 3291 904}%
\special{pa 3288 903}%
\special{pa 3286 901}%
\special{pa 3285 901}%
\special{pa 3283 899}%
\special{pa 3282 899}%
\special{pa 3280 897}%
\special{pa 3277 896}%
\special{pa 3275 894}%
\special{pa 3272 893}%
\special{pa 3271 892}%
\special{pa 3270 892}%
\special{pa 3268 890}%
\special{pa 3261 887}%
\special{pa 3260 886}%
\special{pa 3253 883}%
\special{pa 3252 882}%
\special{pa 3245 879}%
\special{pa 3244 878}%
\special{pa 3242 878}%
\special{pa 3241 877}%
\special{pa 3238 876}%
\special{pa 3237 875}%
\special{pa 3232 874}%
\special{pa 3231 873}%
\special{pa 3230 873}%
\special{pa 3229 872}%
\special{pa 3218 869}%
\special{pa 3217 868}%
\special{pa 3215 868}%
\special{pa 3214 867}%
\special{pa 3211 867}%
\special{pa 3210 866}%
\special{pa 3208 866}%
\special{pa 3207 865}%
\special{pa 3204 865}%
\special{pa 3203 864}%
\special{pa 3199 864}%
\special{pa 3198 863}%
\special{pa 3195 863}%
\special{pa 3194 862}%
\special{pa 3183 861}%
\special{pa 3182 860}%
\special{pa 3174 860}%
\special{pa 3173 859}%
\special{pa 3137 859}%
\special{pa 3136 860}%
\special{pa 3129 860}%
\special{pa 3128 861}%
\special{pa 3123 861}%
\special{pa 3122 862}%
\special{pa 3118 862}%
\special{pa 3117 863}%
\special{pa 3110 864}%
\special{pa 3109 865}%
\special{pa 3106 865}%
\special{pa 3105 866}%
\special{pa 3094 869}%
\special{pa 3093 870}%
\special{pa 3092 870}%
\special{pa 3091 871}%
\special{pa 3089 871}%
\special{pa 3088 872}%
\special{pa 3085 873}%
\special{pa 3084 874}%
\special{pa 3082 874}%
\special{pa 3081 875}%
\special{pa 3078 876}%
\special{pa 3076 878}%
\special{pa 3069 881}%
\special{pa 3067 883}%
\special{pa 3066 883}%
\special{pa 3064 885}%
\special{pa 3063 885}%
\special{pa 3061 887}%
\special{pa 3060 887}%
\special{pa 3058 889}%
\special{pa 3057 889}%
\special{pa 3054 892}%
\special{pa 3053 892}%
\special{pa 3048 897}%
\special{pa 3047 897}%
\special{pa 3034 910}%
\special{pa 3034 911}%
\special{pa 3029 916}%
\special{pa 3029 917}%
\special{pa 3026 920}%
\special{pa 3026 921}%
\special{pa 3024 923}%
\special{pa 3024 924}%
\special{pa 3022 926}%
\special{pa 3022 927}%
\special{pa 3020 929}%
\special{pa 3019 932}%
\special{pa 3017 934}%
\special{pa 3014 941}%
\special{pa 3013 942}%
\special{pa 3012 945}%
\special{pa 3011 946}%
\special{pa 3011 948}%
\special{pa 3010 949}%
\special{pa 3010 950}%
\special{pa 3009 951}%
\special{pa 3009 953}%
\special{pa 3008 954}%
\special{pa 3008 955}%
\special{pa 3007 956}%
\special{pa 3007 958}%
\special{pa 3006 959}%
\special{pa 3006 962}%
\special{pa 3005 963}%
\special{pa 3005 965}%
\special{pa 3004 966}%
\special{pa 3004 969}%
\special{pa 3003 970}%
\special{pa 3003 974}%
\special{pa 3002 975}%
\special{pa 3002 980}%
\special{pa 3001 981}%
\special{pa 3001 988}%
\special{pa 3000 989}%
\special{pa 3000 1011}%
\special{pa 3001 1012}%
\special{pa 3001 1019}%
\special{pa 3002 1020}%
\special{pa 3002 1025}%
\special{pa 3003 1026}%
\special{pa 3003 1030}%
\special{pa 3004 1031}%
\special{pa 3004 1034}%
\special{pa 3005 1035}%
\special{pa 3005 1037}%
\special{pa 3006 1038}%
\special{pa 3006 1041}%
\special{pa 3007 1042}%
\special{pa 3007 1044}%
\special{pa 3008 1045}%
\special{pa 3008 1046}%
\special{pa 3009 1047}%
\special{pa 3009 1049}%
\special{pa 3010 1050}%
\special{pa 3010 1051}%
\special{pa 3011 1052}%
\special{pa 3011 1054}%
\special{pa 3012 1055}%
\special{pa 3015 1062}%
\special{pa 3016 1063}%
\special{pa 3017 1066}%
\special{pa 3019 1068}%
\special{pa 3020 1071}%
\special{pa 3022 1073}%
\special{pa 3022 1074}%
\special{pa 3024 1076}%
\special{pa 3024 1077}%
\special{pa 3026 1079}%
\special{pa 3026 1080}%
\special{pa 3029 1083}%
\special{pa 3029 1084}%
\special{pa 3034 1089}%
\special{pa 3034 1090}%
\special{pa 3047 1103}%
\special{pa 3048 1103}%
\special{pa 3053 1108}%
\special{pa 3054 1108}%
\special{pa 3057 1111}%
\special{pa 3058 1111}%
\special{pa 3060 1113}%
\special{pa 3061 1113}%
\special{pa 3063 1115}%
\special{pa 3064 1115}%
\special{pa 3066 1117}%
\special{pa 3067 1117}%
\special{pa 3069 1119}%
\special{pa 3076 1122}%
\special{pa 3078 1124}%
\special{pa 3081 1125}%
\special{pa 3082 1126}%
\special{pa 3084 1126}%
\special{pa 3085 1127}%
\special{pa 3088 1128}%
\special{pa 3089 1129}%
\special{pa 3091 1129}%
\special{pa 3092 1130}%
\special{pa 3093 1130}%
\special{pa 3094 1131}%
\special{pa 3105 1134}%
\special{pa 3106 1135}%
\special{pa 3113 1136}%
\special{pa 3114 1137}%
\special{pa 3117 1137}%
\special{pa 3118 1138}%
\special{pa 3122 1138}%
\special{pa 3123 1139}%
\special{pa 3128 1139}%
\special{pa 3129 1140}%
\special{pa 3136 1140}%
\special{pa 3137 1141}%
\special{pa 3173 1141}%
\special{pa 3174 1140}%
\special{pa 3182 1140}%
\special{pa 3183 1139}%
\special{pa 3194 1138}%
\special{pa 3195 1137}%
\special{pa 3198 1137}%
\special{pa 3199 1136}%
\special{pa 3203 1136}%
\special{pa 3204 1135}%
\special{pa 3207 1135}%
\special{pa 3208 1134}%
\special{pa 3210 1134}%
\special{pa 3211 1133}%
\special{pa 3214 1133}%
\special{pa 3215 1132}%
\special{pa 3226 1129}%
\special{pa 3227 1128}%
\special{pa 3229 1128}%
\special{pa 3230 1127}%
\special{pa 3231 1127}%
\special{pa 3232 1126}%
\special{pa 3237 1125}%
\special{pa 3238 1124}%
\special{pa 3241 1123}%
\special{pa 3242 1122}%
\special{pa 3244 1122}%
\special{pa 3245 1121}%
\special{pa 3252 1118}%
\special{pa 3253 1117}%
\special{pa 3260 1114}%
\special{pa 3261 1113}%
\special{pa 3268 1110}%
\special{pa 3270 1108}%
\special{pa 3273 1107}%
\special{pa 3274 1106}%
\special{pa 3275 1106}%
\special{pa 3277 1104}%
\special{pa 3280 1103}%
\special{pa 3282 1101}%
\special{pa 3283 1101}%
\special{pa 3285 1099}%
\special{pa 3286 1099}%
\special{pa 3288 1097}%
\special{pa 3291 1096}%
\special{pa 3294 1093}%
\special{pa 3295 1093}%
\special{pa 3297 1091}%
\special{pa 3298 1091}%
\special{pa 3300 1089}%
\special{pa 3301 1089}%
\special{pa 3304 1086}%
\special{pa 3305 1086}%
\special{pa 3307 1084}%
\special{pa 3308 1084}%
\special{pa 3312 1080}%
\special{pa 3313 1080}%
\special{pa 3316 1077}%
\special{pa 3317 1077}%
\special{pa 3321 1073}%
\special{pa 3322 1073}%
\special{pa 3326 1069}%
\special{pa 3327 1069}%
\special{pa 3332 1064}%
\special{pa 3333 1064}%
\special{pa 3340 1057}%
\special{pa 3341 1057}%
\special{pa 3349 1049}%
\special{pa 3350 1049}%
\special{pa 3365 1034}%
\special{pa 3366 1034}%
\special{pa 3399 1001}%
\special{fp}%
\special{pa 3400 1000}%
\special{pa 3434 1034}%
\special{pa 3435 1034}%
\special{pa 3450 1049}%
\special{pa 3451 1049}%
\special{pa 3459 1057}%
\special{pa 3460 1057}%
\special{pa 3467 1064}%
\special{pa 3468 1064}%
\special{pa 3473 1069}%
\special{pa 3474 1069}%
\special{pa 3478 1073}%
\special{pa 3479 1073}%
\special{pa 3483 1077}%
\special{pa 3484 1077}%
\special{pa 3487 1080}%
\special{pa 3488 1080}%
\special{pa 3492 1084}%
\special{pa 3493 1084}%
\special{pa 3495 1086}%
\special{pa 3496 1086}%
\special{pa 3499 1089}%
\special{pa 3500 1089}%
\special{pa 3502 1091}%
\special{pa 3503 1091}%
\special{pa 3505 1093}%
\special{pa 3506 1093}%
\special{pa 3509 1096}%
\special{pa 3512 1097}%
\special{pa 3514 1099}%
\special{pa 3515 1099}%
\special{pa 3517 1101}%
\special{pa 3518 1101}%
\special{pa 3520 1103}%
\special{pa 3523 1104}%
\special{pa 3525 1106}%
\special{pa 3528 1107}%
\special{pa 3529 1108}%
\special{pa 3530 1108}%
\special{pa 3532 1110}%
\special{pa 3539 1113}%
\special{pa 3540 1114}%
\special{pa 3547 1117}%
\special{pa 3548 1118}%
\special{pa 3555 1121}%
\special{pa 3556 1122}%
\special{pa 3558 1122}%
\special{pa 3559 1123}%
\special{pa 3562 1124}%
\special{pa 3563 1125}%
\special{pa 3568 1126}%
\special{pa 3569 1127}%
\special{pa 3570 1127}%
\special{pa 3571 1128}%
\special{pa 3582 1131}%
\special{pa 3583 1132}%
\special{pa 3585 1132}%
\special{pa 3586 1133}%
\special{pa 3589 1133}%
\special{pa 3590 1134}%
\special{pa 3592 1134}%
\special{pa 3593 1135}%
\special{pa 3596 1135}%
\special{pa 3597 1136}%
\special{pa 3601 1136}%
\special{pa 3602 1137}%
\special{pa 3605 1137}%
\special{pa 3606 1138}%
\special{pa 3617 1139}%
\special{pa 3618 1140}%
\special{pa 3626 1140}%
\special{pa 3627 1141}%
\special{pa 3663 1141}%
\special{pa 3664 1140}%
\special{pa 3671 1140}%
\special{pa 3672 1139}%
\special{pa 3677 1139}%
\special{pa 3678 1138}%
\special{pa 3682 1138}%
\special{pa 3683 1137}%
\special{pa 3690 1136}%
\special{pa 3691 1135}%
\special{pa 3694 1135}%
\special{pa 3695 1134}%
\special{pa 3706 1131}%
\special{pa 3707 1130}%
\special{pa 3708 1130}%
\special{pa 3709 1129}%
\special{pa 3711 1129}%
\special{pa 3712 1128}%
\special{pa 3715 1127}%
\special{pa 3716 1126}%
\special{pa 3718 1126}%
\special{pa 3719 1125}%
\special{pa 3722 1124}%
\special{pa 3724 1122}%
\special{pa 3731 1119}%
\special{pa 3733 1117}%
\special{pa 3734 1117}%
\special{pa 3736 1115}%
\special{pa 3737 1115}%
\special{pa 3739 1113}%
\special{pa 3740 1113}%
\special{pa 3742 1111}%
\special{pa 3743 1111}%
\special{pa 3746 1108}%
\special{pa 3747 1108}%
\special{pa 3752 1103}%
\special{pa 3753 1103}%
\special{pa 3766 1090}%
\special{pa 3766 1089}%
\special{pa 3771 1084}%
\special{pa 3771 1083}%
\special{pa 3774 1080}%
\special{pa 3774 1079}%
\special{pa 3776 1077}%
\special{pa 3776 1076}%
\special{pa 3778 1074}%
\special{pa 3778 1073}%
\special{pa 3780 1071}%
\special{pa 3781 1068}%
\special{pa 3783 1066}%
\special{pa 3786 1059}%
\special{pa 3787 1058}%
\special{pa 3788 1055}%
\special{pa 3789 1054}%
\special{pa 3789 1052}%
\special{pa 3790 1051}%
\special{pa 3790 1050}%
\special{pa 3791 1049}%
\special{pa 3791 1047}%
\special{pa 3792 1046}%
\special{pa 3792 1045}%
\special{pa 3793 1044}%
\special{pa 3793 1042}%
\special{pa 3794 1041}%
\special{pa 3794 1038}%
\special{pa 3795 1037}%
\special{pa 3795 1035}%
\special{pa 3796 1034}%
\special{pa 3796 1031}%
\special{pa 3797 1030}%
\special{pa 3797 1026}%
\special{pa 3798 1025}%
\special{pa 3798 1020}%
\special{pa 3799 1019}%
\special{pa 3799 1012}%
\special{pa 3800 1011}%
\special{pa 3800 1001}%
\special{fp}%
\end{picture}}%

%% file: orientation_1.tex
{\unitlength 0.1in%
\begin{picture}(43.3500,9.6500)(6.9500,-16.0000)%
%
\special{pn 13}%
\special{pa 1000 1600}%
\special{pa 1000 800}%
\special{fp}%
\special{sh 1}%
\special{pa 1000 800}%
\special{pa 980 867}%
\special{pa 1000 853}%
\special{pa 1020 867}%
\special{pa 1000 800}%
\special{fp}%
%
\special{pn 8}%
\special{pa 2600 800}%
\special{pa 2600 1600}%
\special{fp}%
\special{sh 1}%
\special{pa 2600 1600}%
\special{pa 2620 1533}%
\special{pa 2600 1547}%
\special{pa 2580 1533}%
\special{pa 2600 1600}%
\special{fp}%
%
\special{pn 8}%
\special{pa 1000 1200}%
\special{pa 2600 1200}%
\special{dt 0.045}%
\put(10.0000,-7.0000){\makebox(0,0){$\hat{c}_0$}}%
\put(26.0000,-7.0000){\makebox(0,0){$c_1$}}%
%
\special{pn 13}%
\special{pa 1600 1600}%
\special{pa 1600 800}%
\special{fp}%
\special{sh 1}%
\special{pa 1600 800}%
\special{pa 1580 867}%
\special{pa 1600 853}%
\special{pa 1620 867}%
\special{pa 1600 800}%
\special{fp}%
%
\special{pn 13}%
\special{pa 2000 800}%
\special{pa 2000 1600}%
\special{fp}%
\special{sh 1}%
\special{pa 2000 1600}%
\special{pa 2020 1533}%
\special{pa 2000 1547}%
\special{pa 1980 1533}%
\special{pa 2000 1600}%
\special{fp}%
\put(16.4000,-11.6000){\makebox(0,0)[lb]{$x_1$}}%
\put(20.4000,-11.6000){\makebox(0,0)[lb]{$x_2$}}%
\put(13.0000,-11.0000){\makebox(0,0){$\hat{\Gamma}$}}%
%
\special{pn 8}%
\special{pa 3400 1600}%
\special{pa 3400 800}%
\special{fp}%
\special{sh 1}%
\special{pa 3400 800}%
\special{pa 3380 867}%
\special{pa 3400 853}%
\special{pa 3420 867}%
\special{pa 3400 800}%
\special{fp}%
%
\special{pn 13}%
\special{pa 5000 800}%
\special{pa 5000 1600}%
\special{fp}%
\special{sh 1}%
\special{pa 5000 1600}%
\special{pa 5020 1533}%
\special{pa 5000 1547}%
\special{pa 4980 1533}%
\special{pa 5000 1600}%
\special{fp}%
%
\special{pn 8}%
\special{pa 3400 1200}%
\special{pa 5000 1200}%
\special{dt 0.045}%
\put(37.0000,-11.0000){\makebox(0,0){$\hat{\Gamma}$}}%
%
\special{pn 13}%
\special{pa 4000 1600}%
\special{pa 4000 800}%
\special{fp}%
\special{sh 1}%
\special{pa 4000 800}%
\special{pa 3980 867}%
\special{pa 4000 853}%
\special{pa 4020 867}%
\special{pa 4000 800}%
\special{fp}%
%
\special{pn 13}%
\special{pa 4400 800}%
\special{pa 4400 1600}%
\special{fp}%
\special{sh 1}%
\special{pa 4400 1600}%
\special{pa 4420 1533}%
\special{pa 4400 1547}%
\special{pa 4380 1533}%
\special{pa 4400 1600}%
\special{fp}%
\put(40.4000,-11.6000){\makebox(0,0)[lb]{$x_3$}}%
\put(44.4000,-11.6000){\makebox(0,0)[lb]{$x_4$}}%
\put(34.0000,-7.0000){\makebox(0,0){$\hat{c}_0$}}%
\put(40.0000,-7.0000){\makebox(0,0){$c_1$}}%
\put(44.0000,-7.0000){\makebox(0,0){$c_1$}}%
\put(50.0000,-7.0000){\makebox(0,0){$c_1$}}%
\put(16.0000,-7.0000){\makebox(0,0){$\hat{c}_0$}}%
\put(20.0000,-7.0000){\makebox(0,0){$\hat{c}_0$}}%
\put(10.3000,-12.3000){\makebox(0,0)[lt]{$\hat{b}_0$}}%
\put(26.3000,-12.3000){\makebox(0,0)[lt]{$\hat{b}_1$}}%
\put(50.3000,-12.3000){\makebox(0,0)[lt]{$\hat{b}_1$}}%
\put(34.3000,-12.3000){\makebox(0,0)[lt]{$\hat{b}_0$}}%
\end{picture}}%

%% file: change_T_pm.tex
{\unitlength 0.1in%
\begin{picture}(48.2000,40.6500)(10.8000,-44.0000)%
%
\special{pn 13}%
\special{pa 1300 1400}%
\special{pa 1300 600}%
\special{fp}%
\special{sh 1}%
\special{pa 1300 600}%
\special{pa 1280 667}%
\special{pa 1300 653}%
\special{pa 1320 667}%
\special{pa 1300 600}%
\special{fp}%
%
\special{pn 13}%
\special{ar 1700 1000 200 200 1.5707963 4.7123890}%
%
\special{pn 13}%
\special{pa 1685 801}%
\special{pa 1700 800}%
\special{fp}%
\special{sh 1}%
\special{pa 1700 800}%
\special{pa 1632 784}%
\special{pa 1647 804}%
\special{pa 1635 824}%
\special{pa 1700 800}%
\special{fp}%
%
\special{pn 8}%
\special{pa 1500 1000}%
\special{pa 1900 1000}%
\special{dt 0.045}%
\special{sh 1}%
\special{pa 1900 1000}%
\special{pa 1833 980}%
\special{pa 1847 1000}%
\special{pa 1833 1020}%
\special{pa 1900 1000}%
\special{fp}%
\put(18.0000,-9.0000){\makebox(0,0){$\Gamma$}}%
\put(16.0000,-8.0000){\makebox(0,0)[rb]{$C_0$}}%
\put(13.0000,-15.0000){\makebox(0,0){$C_1$}}%
\put(13.0000,-10.0000){\makebox(0,0){$\bullet$}}%
\put(13.2000,-9.8000){\makebox(0,0)[lb]{$x$}}%
\put(22.0000,-10.0000){\makebox(0,0){$\rightsquigarrow$}}%
%
\special{pn 13}%
\special{pa 2600 1400}%
\special{pa 2600 600}%
\special{fp}%
\special{sh 1}%
\special{pa 2600 600}%
\special{pa 2580 667}%
\special{pa 2600 653}%
\special{pa 2620 667}%
\special{pa 2600 600}%
\special{fp}%
\put(26.0000,-15.0000){\makebox(0,0){$C_1$}}%
\put(26.0000,-10.0000){\makebox(0,0){$\bullet$}}%
\put(26.2000,-9.8000){\makebox(0,0)[lb]{$x$}}%
%
\special{pn 13}%
\special{ar 3000 1000 600 200 1.5707963 4.7123890}%
%
\special{pn 13}%
\special{pa 2978 800}%
\special{pa 3000 800}%
\special{fp}%
\special{sh 1}%
\special{pa 3000 800}%
\special{pa 2933 780}%
\special{pa 2947 800}%
\special{pa 2933 820}%
\special{pa 3000 800}%
\special{fp}%
\put(29.0000,-8.0000){\makebox(0,0)[rb]{$\hat{C}_0$}}%
%
\special{pn 8}%
\special{pa 2400 1000}%
\special{pa 3200 1000}%
\special{dt 0.045}%
\special{sh 1}%
\special{pa 3200 1000}%
\special{pa 3133 980}%
\special{pa 3147 1000}%
\special{pa 3133 1020}%
\special{pa 3200 1000}%
\special{fp}%
\put(31.0000,-9.0000){\makebox(0,0){$\hat{\Gamma}$}}%
\put(23.0000,-4.0000){\makebox(0,0){Separating homotopy}}%
\put(39.0000,-15.0000){\makebox(0,0){$C_1$}}%
%
\special{pn 13}%
\special{pa 3900 1400}%
\special{pa 3900 600}%
\special{fp}%
\special{sh 1}%
\special{pa 3900 600}%
\special{pa 3880 667}%
\special{pa 3900 653}%
\special{pa 3920 667}%
\special{pa 3900 600}%
\special{fp}%
\put(39.0000,-10.0000){\makebox(0,0){$\bullet$}}%
\put(39.2000,-9.8000){\makebox(0,0)[lb]{$x$}}%
%
\special{pn 13}%
\special{ar 4300 1150 200 100 1.5707963 4.7123890}%
%
\special{pn 13}%
\special{pa 4280 1050}%
\special{pa 4300 1050}%
\special{fp}%
\special{sh 1}%
\special{pa 4300 1050}%
\special{pa 4233 1030}%
\special{pa 4247 1050}%
\special{pa 4233 1070}%
\special{pa 4300 1050}%
\special{fp}%
%
\special{pn 13}%
\special{ar 4300 850 200 100 1.5707963 4.7123890}%
%
\special{pn 13}%
\special{pa 4280 750}%
\special{pa 4300 750}%
\special{fp}%
\special{sh 1}%
\special{pa 4300 750}%
\special{pa 4233 730}%
\special{pa 4247 750}%
\special{pa 4233 770}%
\special{pa 4300 750}%
\special{fp}%
\put(41.5000,-7.5000){\makebox(0,0)[rb]{$C_0$}}%
%
\special{pn 8}%
\special{pa 4100 1000}%
\special{pa 4600 1000}%
\special{dt 0.045}%
\special{sh 1}%
\special{pa 4600 1000}%
\special{pa 4533 980}%
\special{pa 4547 1000}%
\special{pa 4533 1020}%
\special{pa 4600 1000}%
\special{fp}%
\put(45.0000,-9.0000){\makebox(0,0){($\Gamma$)}}%
\put(49.0000,-10.0000){\makebox(0,0){$\rightsquigarrow$}}%
\put(53.0000,-15.0000){\makebox(0,0){$C_1$}}%
%
\special{pn 13}%
\special{pa 5300 1400}%
\special{pa 5300 600}%
\special{fp}%
\special{sh 1}%
\special{pa 5300 600}%
\special{pa 5280 667}%
\special{pa 5300 653}%
\special{pa 5320 667}%
\special{pa 5300 600}%
\special{fp}%
\put(53.0000,-10.0000){\makebox(0,0){$\bullet$}}%
\put(56.0000,-12.0000){\makebox(0,0)[lt]{$x$}}%
%
\special{pn 8}%
\special{pa 5100 1000}%
\special{pa 5800 1000}%
\special{dt 0.045}%
\special{sh 1}%
\special{pa 5800 1000}%
\special{pa 5733 980}%
\special{pa 5747 1000}%
\special{pa 5733 1020}%
\special{pa 5800 1000}%
\special{fp}%
\put(58.0000,-9.0000){\makebox(0,0){($\hat{\Gamma}$)}}%
\put(47.0000,-4.0000){\makebox(0,0){Homotopy $h$ from $C_0+_{\Gamma}C_1$ to $\hat{C}_0+_{\hat{\Gamma}}C_1$}}%
%
\special{pn 13}%
\special{pa 1300 2800}%
\special{pa 1300 2000}%
\special{fp}%
\special{sh 1}%
\special{pa 1300 2000}%
\special{pa 1280 2067}%
\special{pa 1300 2053}%
\special{pa 1320 2067}%
\special{pa 1300 2000}%
\special{fp}%
\put(13.0000,-29.0000){\makebox(0,0){$C_1$}}%
\put(13.0000,-24.0000){\makebox(0,0){$\bullet$}}%
\put(13.2000,-23.8000){\makebox(0,0)[lb]{$x$}}%
%
\special{pn 13}%
\special{ar 1100 2400 400 200 4.7123890 1.5707963}%
%
\special{pn 13}%
\special{pa 1120 2200}%
\special{pa 1100 2200}%
\special{fp}%
\special{sh 1}%
\special{pa 1100 2200}%
\special{pa 1167 2220}%
\special{pa 1153 2200}%
\special{pa 1167 2180}%
\special{pa 1100 2200}%
\special{fp}%
%
\special{pn 8}%
\special{pa 1500 2400}%
\special{pa 1900 2400}%
\special{dt 0.045}%
\special{sh 1}%
\special{pa 1900 2400}%
\special{pa 1833 2380}%
\special{pa 1847 2400}%
\special{pa 1833 2420}%
\special{pa 1900 2400}%
\special{fp}%
\put(18.0000,-23.0000){\makebox(0,0){$\Gamma$}}%
\put(14.0000,-22.5000){\makebox(0,0)[lb]{$C_0$}}%
\put(21.0000,-24.0000){\makebox(0,0){$\rightsquigarrow$}}%
\put(29.0000,-29.0000){\makebox(0,0){$C_1$}}%
%
\special{pn 13}%
\special{pa 2900 2800}%
\special{pa 2900 2000}%
\special{fp}%
\special{sh 1}%
\special{pa 2900 2000}%
\special{pa 2880 2067}%
\special{pa 2900 2053}%
\special{pa 2920 2067}%
\special{pa 2900 2000}%
\special{fp}%
\put(29.0000,-24.0000){\makebox(0,0){$\bullet$}}%
\put(29.2000,-23.8000){\makebox(0,0)[lb]{$x$}}%
%
\special{pn 13}%
\special{ar 2300 2400 400 200 4.7123890 1.5707963}%
%
\special{pn 13}%
\special{pa 2320 2200}%
\special{pa 2300 2200}%
\special{fp}%
\special{sh 1}%
\special{pa 2300 2200}%
\special{pa 2367 2220}%
\special{pa 2353 2200}%
\special{pa 2367 2180}%
\special{pa 2300 2200}%
\special{fp}%
\put(25.0000,-22.0000){\makebox(0,0)[lb]{$\hat{C}_0$}}%
%
\special{pn 8}%
\special{pa 2700 2400}%
\special{pa 3300 2400}%
\special{dt 0.045}%
\special{sh 1}%
\special{pa 3300 2400}%
\special{pa 3233 2380}%
\special{pa 3247 2400}%
\special{pa 3233 2420}%
\special{pa 3300 2400}%
\special{fp}%
\put(32.0000,-23.0000){\makebox(0,0){$\hat{\Gamma}$}}%
\put(37.0000,-29.0000){\makebox(0,0){$C_1$}}%
%
\special{pn 13}%
\special{pa 3700 2800}%
\special{pa 3700 2000}%
\special{fp}%
\special{sh 1}%
\special{pa 3700 2000}%
\special{pa 3680 2067}%
\special{pa 3700 2053}%
\special{pa 3720 2067}%
\special{pa 3700 2000}%
\special{fp}%
%
\special{pn 8}%
\special{pa 4000 2400}%
\special{pa 4500 2400}%
\special{dt 0.045}%
\special{sh 1}%
\special{pa 4500 2400}%
\special{pa 4433 2380}%
\special{pa 4447 2400}%
\special{pa 4433 2420}%
\special{pa 4500 2400}%
\special{fp}%
\put(44.0000,-23.0000){\makebox(0,0){($\Gamma$)}}%
\put(41.0000,-23.0000){\makebox(0,0)[lb]{$C_0$}}%
\put(37.0000,-24.0000){\makebox(0,0){$\bullet$}}%
\put(37.2000,-23.8000){\makebox(0,0)[lb]{$x$}}%
\put(48.0000,-24.0000){\makebox(0,0){$\rightsquigarrow$}}%
\put(55.0000,-24.0000){\makebox(0,0){$\bullet$}}%
%
\special{pn 13}%
\special{pa 5500 2800}%
\special{pa 5500 2000}%
\special{fp}%
\special{sh 1}%
\special{pa 5500 2000}%
\special{pa 5480 2067}%
\special{pa 5500 2053}%
\special{pa 5520 2067}%
\special{pa 5500 2000}%
\special{fp}%
\put(55.0000,-29.0000){\makebox(0,0){$C_1$}}%
\put(51.0000,-22.0000){\makebox(0,0)[rt]{$C_0$}}%
\put(58.0000,-23.0000){\makebox(0,0){($\hat{\Gamma}$)}}%
%
\special{pn 8}%
\special{pa 5000 2400}%
\special{pa 5900 2400}%
\special{dt 0.045}%
\special{sh 1}%
\special{pa 5900 2400}%
\special{pa 5833 2380}%
\special{pa 5847 2400}%
\special{pa 5833 2420}%
\special{pa 5900 2400}%
\special{fp}%
\put(13.0000,-43.0000){\makebox(0,0){$C_1$}}%
%
\special{pn 13}%
\special{pa 1300 4200}%
\special{pa 1300 3400}%
\special{fp}%
\special{sh 1}%
\special{pa 1300 3400}%
\special{pa 1280 3467}%
\special{pa 1300 3453}%
\special{pa 1320 3467}%
\special{pa 1300 3400}%
\special{fp}%
%
\special{pn 13}%
\special{pa 1700 4000}%
\special{pa 1100 3400}%
\special{fp}%
\special{sh 1}%
\special{pa 1100 3400}%
\special{pa 1133 3461}%
\special{pa 1138 3438}%
\special{pa 1161 3433}%
\special{pa 1100 3400}%
\special{fp}%
\put(17.0000,-40.0000){\makebox(0,0)[lt]{$C_0$}}%
%
\special{pn 8}%
\special{pa 1500 3800}%
\special{pa 1900 3800}%
\special{dt 0.045}%
\special{sh 1}%
\special{pa 1900 3800}%
\special{pa 1833 3780}%
\special{pa 1847 3800}%
\special{pa 1833 3820}%
\special{pa 1900 3800}%
\special{fp}%
\put(18.0000,-37.0000){\makebox(0,0){$\Gamma$}}%
\put(13.0000,-38.0000){\makebox(0,0){$\bullet$}}%
\put(13.2000,-37.8000){\makebox(0,0)[lb]{$x$}}%
\put(21.0000,-38.0000){\makebox(0,0){$\rightsquigarrow$}}%
%
\special{pn 13}%
\special{pa 2900 4000}%
\special{pa 2300 3400}%
\special{fp}%
\special{sh 1}%
\special{pa 2300 3400}%
\special{pa 2333 3461}%
\special{pa 2338 3438}%
\special{pa 2361 3433}%
\special{pa 2300 3400}%
\special{fp}%
%
\special{pn 13}%
\special{pa 2800 4200}%
\special{pa 2800 3400}%
\special{fp}%
\special{sh 1}%
\special{pa 2800 3400}%
\special{pa 2780 3467}%
\special{pa 2800 3453}%
\special{pa 2820 3467}%
\special{pa 2800 3400}%
\special{fp}%
\put(28.0000,-43.0000){\makebox(0,0){$C_1$}}%
\put(28.0000,-38.0000){\makebox(0,0){$\bullet$}}%
\put(28.2000,-37.8000){\makebox(0,0)[lb]{$x$}}%
%
\special{pn 8}%
\special{pa 2700 3800}%
\special{pa 3300 3800}%
\special{dt 0.045}%
\special{sh 1}%
\special{pa 3300 3800}%
\special{pa 3233 3780}%
\special{pa 3247 3800}%
\special{pa 3233 3820}%
\special{pa 3300 3800}%
\special{fp}%
\put(32.0000,-37.0000){\makebox(0,0){$\hat{\Gamma}$}}%
%
\special{pn 8}%
\special{pa 3500 400}%
\special{pa 3500 4400}%
\special{da 0.070}%
%
\special{pn 13}%
\special{pa 3600 2150}%
\special{pa 3590 2150}%
\special{fp}%
\special{sh 1}%
\special{pa 3590 2150}%
\special{pa 3657 2170}%
\special{pa 3643 2150}%
\special{pa 3657 2130}%
\special{pa 3590 2150}%
\special{fp}%
%
\special{pn 13}%
\special{pa 5000 2150}%
\special{pa 4990 2150}%
\special{fp}%
\special{sh 1}%
\special{pa 4990 2150}%
\special{pa 5057 2170}%
\special{pa 5043 2150}%
\special{pa 5057 2130}%
\special{pa 4990 2150}%
\special{fp}%
%
\special{pn 13}%
\special{pa 4200 2450}%
\special{pa 4210 2450}%
\special{fp}%
\special{sh 1}%
\special{pa 4210 2450}%
\special{pa 4143 2430}%
\special{pa 4157 2450}%
\special{pa 4143 2470}%
\special{pa 4210 2450}%
\special{fp}%
%
\special{pn 13}%
\special{pa 5600 2450}%
\special{pa 5610 2450}%
\special{fp}%
\special{sh 1}%
\special{pa 5610 2450}%
\special{pa 5543 2430}%
\special{pa 5557 2450}%
\special{pa 5543 2470}%
\special{pa 5610 2450}%
\special{fp}%
\put(38.0000,-43.0000){\makebox(0,0){$C_1$}}%
%
\special{pn 13}%
\special{pa 3800 4200}%
\special{pa 3800 3400}%
\special{fp}%
\special{sh 1}%
\special{pa 3800 3400}%
\special{pa 3780 3467}%
\special{pa 3800 3453}%
\special{pa 3820 3467}%
\special{pa 3800 3400}%
\special{fp}%
\put(38.0000,-38.0000){\makebox(0,0){$\bullet$}}%
\put(37.0000,-38.0000){\makebox(0,0){$x$}}%
%
\special{pn 13}%
\special{ar 4100 3400 500 350 1.5707963 2.7989911}%
%
\special{pn 13}%
\special{pa 3635 3529}%
\special{pa 3629 3518}%
\special{fp}%
\special{sh 1}%
\special{pa 3629 3518}%
\special{pa 3643 3586}%
\special{pa 3655 3565}%
\special{pa 3678 3567}%
\special{pa 3629 3518}%
\special{fp}%
%
\special{pn 13}%
\special{pa 4200 3750}%
\special{pa 4100 3750}%
\special{fp}%
%
\special{pn 13}%
\special{ar 4000 3900 50 50 3.1415927 4.7123890}%
%
\special{pn 13}%
\special{ar 4400 3900 450 400 2.3888851 3.1415927}%
%
\special{pn 13}%
\special{pa 5400 3850}%
\special{pa 5600 3850}%
\special{fp}%
\special{sh 1}%
\special{pa 5600 3850}%
\special{pa 5533 3830}%
\special{pa 5547 3850}%
\special{pa 5533 3870}%
\special{pa 5600 3850}%
\special{fp}%
%
\special{pn 8}%
\special{pa 4000 3800}%
\special{pa 4500 3800}%
\special{dt 0.045}%
\special{sh 1}%
\special{pa 4500 3800}%
\special{pa 4433 3780}%
\special{pa 4447 3800}%
\special{pa 4433 3820}%
\special{pa 4500 3800}%
\special{fp}%
\put(44.0000,-37.0000){\makebox(0,0){($\Gamma$)}}%
\put(48.0000,-38.0000){\makebox(0,0){$\rightsquigarrow$}}%
%
\special{pn 13}%
\special{ar 5400 3900 50 50 3.1415927 4.7123890}%
%
\special{pn 13}%
\special{pa 4000 3850}%
\special{pa 4200 3850}%
\special{fp}%
\special{sh 1}%
\special{pa 4200 3850}%
\special{pa 4133 3830}%
\special{pa 4147 3850}%
\special{pa 4133 3870}%
\special{pa 4200 3850}%
\special{fp}%
%
\special{pn 13}%
\special{ar 5800 3900 450 400 2.3888851 3.1415927}%
%
\special{pn 13}%
\special{pa 5500 4200}%
\special{pa 5500 3400}%
\special{fp}%
\special{sh 1}%
\special{pa 5500 3400}%
\special{pa 5480 3467}%
\special{pa 5500 3453}%
\special{pa 5520 3467}%
\special{pa 5500 3400}%
\special{fp}%
\put(55.0000,-38.0000){\makebox(0,0){$\bullet$}}%
\put(58.0000,-40.0000){\makebox(0,0)[lt]{$x$}}%
%
\special{pn 8}%
\special{pa 5000 3800}%
\special{pa 5900 3800}%
\special{dt 0.045}%
\special{sh 1}%
\special{pa 5900 3800}%
\special{pa 5833 3780}%
\special{pa 5847 3800}%
\special{pa 5833 3820}%
\special{pa 5900 3800}%
\special{fp}%
\put(58.0000,-37.0000){\makebox(0,0){($\hat{\Gamma}$)}}%
\put(24.0000,-35.0000){\makebox(0,0)[rt]{$\hat{C}_0$}}%
\put(41.0000,-36.5000){\makebox(0,0){$C_0$}}%
\put(51.0000,-36.0000){\makebox(0,0)[lb]{$C_0$}}%
%
\special{pn 13}%
\special{ar 5550 850 450 100 1.5707963 4.7123890}%
%
\special{pn 13}%
\special{pa 5525 750}%
\special{pa 5550 750}%
\special{fp}%
\special{sh 1}%
\special{pa 5550 750}%
\special{pa 5483 730}%
\special{pa 5497 750}%
\special{pa 5483 770}%
\special{pa 5550 750}%
\special{fp}%
%
\special{pn 13}%
\special{ar 5550 1150 450 100 1.5707963 4.7123890}%
%
\special{pn 13}%
\special{pa 5525 1050}%
\special{pa 5550 1050}%
\special{fp}%
\special{sh 1}%
\special{pa 5550 1050}%
\special{pa 5483 1030}%
\special{pa 5497 1050}%
\special{pa 5483 1070}%
\special{pa 5550 1050}%
\special{fp}%
\put(51.0000,-8.0000){\makebox(0,0)[rb]{$C_0$}}%
\special{pn 13}%
\special{pa 3600 2150}%
\special{pa 3615 2150}%
\special{pa 3620 2151}%
\special{pa 3630 2151}%
\special{pa 3635 2152}%
\special{pa 3640 2152}%
\special{pa 3645 2153}%
\special{pa 3650 2153}%
\special{pa 3700 2163}%
\special{pa 3705 2165}%
\special{pa 3710 2166}%
\special{pa 3715 2168}%
\special{pa 3720 2169}%
\special{pa 3725 2171}%
\special{pa 3730 2172}%
\special{pa 3740 2176}%
\special{pa 3745 2177}%
\special{pa 3785 2193}%
\special{pa 3790 2196}%
\special{pa 3805 2202}%
\special{pa 3810 2205}%
\special{pa 3820 2209}%
\special{pa 3825 2212}%
\special{pa 3830 2214}%
\special{pa 3835 2217}%
\special{pa 3840 2219}%
\special{pa 3845 2222}%
\special{pa 3850 2224}%
\special{pa 3855 2227}%
\special{pa 3860 2229}%
\special{pa 3865 2232}%
\special{pa 3870 2234}%
\special{pa 3880 2240}%
\special{pa 3885 2242}%
\special{pa 3890 2245}%
\special{pa 3895 2247}%
\special{pa 3905 2253}%
\special{pa 3910 2255}%
\special{pa 3915 2258}%
\special{pa 3920 2260}%
\special{pa 3930 2266}%
\special{pa 3935 2268}%
\special{pa 3940 2271}%
\special{pa 3945 2273}%
\special{pa 3950 2276}%
\special{pa 3955 2278}%
\special{pa 3960 2281}%
\special{pa 3965 2283}%
\special{pa 3970 2286}%
\special{pa 3975 2288}%
\special{pa 3980 2291}%
\special{pa 3990 2295}%
\special{pa 3995 2298}%
\special{pa 4010 2304}%
\special{pa 4015 2307}%
\special{pa 4055 2323}%
\special{pa 4060 2324}%
\special{pa 4070 2328}%
\special{pa 4075 2329}%
\special{pa 4080 2331}%
\special{pa 4085 2332}%
\special{pa 4090 2334}%
\special{pa 4095 2335}%
\special{pa 4100 2337}%
\special{pa 4150 2347}%
\special{pa 4155 2347}%
\special{pa 4160 2348}%
\special{pa 4165 2348}%
\special{pa 4170 2349}%
\special{pa 4180 2349}%
\special{pa 4185 2350}%
\special{pa 4200 2350}%
\special{fp}%
\special{pn 13}%
\special{pa 3600 2650}%
\special{pa 3615 2650}%
\special{pa 3620 2649}%
\special{pa 3630 2649}%
\special{pa 3635 2648}%
\special{pa 3640 2648}%
\special{pa 3645 2647}%
\special{pa 3650 2647}%
\special{pa 3700 2637}%
\special{pa 3705 2635}%
\special{pa 3710 2634}%
\special{pa 3715 2632}%
\special{pa 3720 2631}%
\special{pa 3725 2629}%
\special{pa 3730 2628}%
\special{pa 3740 2624}%
\special{pa 3745 2623}%
\special{pa 3785 2607}%
\special{pa 3790 2604}%
\special{pa 3805 2598}%
\special{pa 3810 2595}%
\special{pa 3820 2591}%
\special{pa 3825 2588}%
\special{pa 3830 2586}%
\special{pa 3835 2583}%
\special{pa 3840 2581}%
\special{pa 3845 2578}%
\special{pa 3850 2576}%
\special{pa 3855 2573}%
\special{pa 3860 2571}%
\special{pa 3865 2568}%
\special{pa 3870 2566}%
\special{pa 3880 2560}%
\special{pa 3885 2558}%
\special{pa 3890 2555}%
\special{pa 3895 2553}%
\special{pa 3905 2547}%
\special{pa 3910 2545}%
\special{pa 3915 2542}%
\special{pa 3920 2540}%
\special{pa 3930 2534}%
\special{pa 3935 2532}%
\special{pa 3940 2529}%
\special{pa 3945 2527}%
\special{pa 3950 2524}%
\special{pa 3955 2522}%
\special{pa 3960 2519}%
\special{pa 3965 2517}%
\special{pa 3970 2514}%
\special{pa 3975 2512}%
\special{pa 3980 2509}%
\special{pa 3990 2505}%
\special{pa 3995 2502}%
\special{pa 4010 2496}%
\special{pa 4015 2493}%
\special{pa 4055 2477}%
\special{pa 4060 2476}%
\special{pa 4070 2472}%
\special{pa 4075 2471}%
\special{pa 4080 2469}%
\special{pa 4085 2468}%
\special{pa 4090 2466}%
\special{pa 4095 2465}%
\special{pa 4100 2463}%
\special{pa 4150 2453}%
\special{pa 4155 2453}%
\special{pa 4160 2452}%
\special{pa 4165 2452}%
\special{pa 4170 2451}%
\special{pa 4180 2451}%
\special{pa 4185 2450}%
\special{pa 4200 2450}%
\special{fp}%
\special{pn 13}%
\special{pa 5000 2150}%
\special{pa 5015 2150}%
\special{pa 5020 2151}%
\special{pa 5030 2151}%
\special{pa 5035 2152}%
\special{pa 5040 2152}%
\special{pa 5045 2153}%
\special{pa 5050 2153}%
\special{pa 5100 2163}%
\special{pa 5105 2165}%
\special{pa 5110 2166}%
\special{pa 5115 2168}%
\special{pa 5120 2169}%
\special{pa 5125 2171}%
\special{pa 5130 2172}%
\special{pa 5140 2176}%
\special{pa 5145 2177}%
\special{pa 5185 2193}%
\special{pa 5190 2196}%
\special{pa 5205 2202}%
\special{pa 5210 2205}%
\special{pa 5220 2209}%
\special{pa 5225 2212}%
\special{pa 5230 2214}%
\special{pa 5235 2217}%
\special{pa 5240 2219}%
\special{pa 5245 2222}%
\special{pa 5250 2224}%
\special{pa 5255 2227}%
\special{pa 5260 2229}%
\special{pa 5265 2232}%
\special{pa 5270 2234}%
\special{pa 5280 2240}%
\special{pa 5285 2242}%
\special{pa 5290 2245}%
\special{pa 5295 2247}%
\special{pa 5305 2253}%
\special{pa 5310 2255}%
\special{pa 5315 2258}%
\special{pa 5320 2260}%
\special{pa 5330 2266}%
\special{pa 5335 2268}%
\special{pa 5340 2271}%
\special{pa 5345 2273}%
\special{pa 5350 2276}%
\special{pa 5355 2278}%
\special{pa 5360 2281}%
\special{pa 5365 2283}%
\special{pa 5370 2286}%
\special{pa 5375 2288}%
\special{pa 5380 2291}%
\special{pa 5390 2295}%
\special{pa 5395 2298}%
\special{pa 5410 2304}%
\special{pa 5415 2307}%
\special{pa 5455 2323}%
\special{pa 5460 2324}%
\special{pa 5470 2328}%
\special{pa 5475 2329}%
\special{pa 5480 2331}%
\special{pa 5485 2332}%
\special{pa 5490 2334}%
\special{pa 5495 2335}%
\special{pa 5500 2337}%
\special{pa 5550 2347}%
\special{pa 5555 2347}%
\special{pa 5560 2348}%
\special{pa 5565 2348}%
\special{pa 5570 2349}%
\special{pa 5580 2349}%
\special{pa 5585 2350}%
\special{pa 5600 2350}%
\special{fp}%
\special{pn 13}%
\special{pa 5000 2650}%
\special{pa 5015 2650}%
\special{pa 5020 2649}%
\special{pa 5030 2649}%
\special{pa 5035 2648}%
\special{pa 5040 2648}%
\special{pa 5045 2647}%
\special{pa 5050 2647}%
\special{pa 5100 2637}%
\special{pa 5105 2635}%
\special{pa 5110 2634}%
\special{pa 5115 2632}%
\special{pa 5120 2631}%
\special{pa 5125 2629}%
\special{pa 5130 2628}%
\special{pa 5140 2624}%
\special{pa 5145 2623}%
\special{pa 5185 2607}%
\special{pa 5190 2604}%
\special{pa 5205 2598}%
\special{pa 5210 2595}%
\special{pa 5220 2591}%
\special{pa 5225 2588}%
\special{pa 5230 2586}%
\special{pa 5235 2583}%
\special{pa 5240 2581}%
\special{pa 5245 2578}%
\special{pa 5250 2576}%
\special{pa 5255 2573}%
\special{pa 5260 2571}%
\special{pa 5265 2568}%
\special{pa 5270 2566}%
\special{pa 5280 2560}%
\special{pa 5285 2558}%
\special{pa 5290 2555}%
\special{pa 5295 2553}%
\special{pa 5305 2547}%
\special{pa 5310 2545}%
\special{pa 5315 2542}%
\special{pa 5320 2540}%
\special{pa 5330 2534}%
\special{pa 5335 2532}%
\special{pa 5340 2529}%
\special{pa 5345 2527}%
\special{pa 5350 2524}%
\special{pa 5355 2522}%
\special{pa 5360 2519}%
\special{pa 5365 2517}%
\special{pa 5370 2514}%
\special{pa 5375 2512}%
\special{pa 5380 2509}%
\special{pa 5390 2505}%
\special{pa 5395 2502}%
\special{pa 5410 2496}%
\special{pa 5415 2493}%
\special{pa 5455 2477}%
\special{pa 5460 2476}%
\special{pa 5470 2472}%
\special{pa 5475 2471}%
\special{pa 5480 2469}%
\special{pa 5485 2468}%
\special{pa 5490 2466}%
\special{pa 5495 2465}%
\special{pa 5500 2463}%
\special{pa 5550 2453}%
\special{pa 5555 2453}%
\special{pa 5560 2452}%
\special{pa 5565 2452}%
\special{pa 5570 2451}%
\special{pa 5580 2451}%
\special{pa 5585 2450}%
\special{pa 5600 2450}%
\special{fp}%
%
\special{pn 13}%
\special{ar 5500 3400 500 350 1.5707963 2.7989911}%
%
\special{pn 13}%
\special{pa 5035 3529}%
\special{pa 5029 3518}%
\special{fp}%
\special{sh 1}%
\special{pa 5029 3518}%
\special{pa 5043 3586}%
\special{pa 5055 3565}%
\special{pa 5078 3567}%
\special{pa 5029 3518}%
\special{fp}%
%
\special{pn 13}%
\special{pa 5600 3750}%
\special{pa 5500 3750}%
\special{fp}%
\put(55.5000,-35.0000){\makebox(0,0)[lb]{$C_1$}}%
\put(18.0000,-17.0000){\makebox(0,0)[lb]{direct, $+1$}}%
\put(44.0000,-17.0000){\makebox(0,0)[lb]{direct, $+2$}}%
\put(18.0000,-29.0000){\makebox(0,0)[lb]{direct, $-1$}}%
\put(44.0000,-29.0000){\makebox(0,0)[lb]{(no tangency)}}%
\put(18.0000,-43.0000){\makebox(0,0)[lb]{(no tangency)}}%
\put(44.0000,-43.0000){\makebox(0,0)[lb]{direct, $+1$}}%
%
\special{pn 4}%
\special{pa 5600 1200}%
\special{pa 5300 1000}%
\special{fp}%
\put(58.0000,-26.0000){\makebox(0,0)[lt]{$x$}}%
%
\special{pn 4}%
\special{pa 5800 2600}%
\special{pa 5500 2400}%
\special{fp}%
%
\special{pn 4}%
\special{pa 5800 4000}%
\special{pa 5500 3800}%
\special{fp}%
%
\special{pn 8}%
\special{pa 1900 800}%
\special{pa 1900 1200}%
\special{fp}%
\special{sh 1}%
\special{pa 1900 1200}%
\special{pa 1920 1133}%
\special{pa 1900 1147}%
\special{pa 1880 1133}%
\special{pa 1900 1200}%
\special{fp}%
\put(19.0000,-7.0000){\makebox(0,0){$C_1$}}%
\put(15.0000,-13.0000){\makebox(0,0){$\gamma(0)$}}%
%
\special{pn 4}%
\special{pa 1500 1250}%
\special{pa 1500 1000}%
\special{fp}%
\put(20.5000,-7.5000){\makebox(0,0)[lb]{$\gamma(1)$}}%
%
\special{pn 4}%
\special{pa 2050 750}%
\special{pa 1900 1000}%
\special{fp}%
\put(16.0000,-26.0000){\makebox(0,0)[lt]{$\gamma(0)$}}%
\put(16.0000,-35.0000){\makebox(0,0)[lb]{$\gamma(0)$}}%
%
\special{pn 4}%
\special{pa 1600 3500}%
\special{pa 1500 3800}%
\special{fp}%
%
\special{pn 4}%
\special{pa 1600 2600}%
\special{pa 1500 2400}%
\special{fp}%
\put(23.0000,-12.0000){\makebox(0,0)[rt]{$\hat{\gamma}(0)$}}%
%
\special{pn 4}%
\special{pa 2300 1200}%
\special{pa 2400 1000}%
\special{fp}%
\end{picture}}%

%% file: new_triple_point.tex
{\unitlength 0.1in%
\begin{picture}(37.0000,19.1500)(7.0000,-20.0000)%
\put(18.0000,-6.5000){\makebox(0,0)[lt]{$C_0$}}%
\put(13.7000,-2.3000){\makebox(0,0)[lb]{$C_1$}}%
\put(9.0000,-11.0000){\makebox(0,0)[lb]{$T^{\St}_{\Gamma}(C_0,C_1)=0$ (no triple point)}}%
%
\special{pn 13}%
\special{ar 1500 500 400 200 0.0000000 6.2831853}%
\put(15.0000,-8.0000){\makebox(0,0)[lt]{$\Gamma$}}%
%
\special{pn 4}%
\special{pa 1500 800}%
\special{pa 1200 500}%
\special{fp}%
%
\special{pn 13}%
\special{ar 1500 1600 400 200 3.6864566 2.5967287}%
%
\special{pn 8}%
\special{ar 1000 500 300 300 6.2831853 3.1415927}%
%
\special{pn 8}%
\special{ar 1300 500 300 300 6.2831853 3.1415927}%
%
\special{pn 8}%
\special{ar 1150 500 150 150 3.1415927 6.2831853}%
%
\special{pn 8}%
\special{ar 1150 500 450 300 3.1415927 6.2831853}%
%
\special{pn 13}%
\special{pa 1100 500}%
\special{pa 1300 500}%
\special{dt 0.045}%
\special{sh 1}%
\special{pa 1300 500}%
\special{pa 1233 480}%
\special{pa 1247 500}%
\special{pa 1233 520}%
\special{pa 1300 500}%
\special{fp}%
%
\special{pn 4}%
\special{pa 1100 400}%
\special{pa 800 400}%
\special{fp}%
\special{sh 1}%
\special{pa 800 400}%
\special{pa 867 420}%
\special{pa 853 400}%
\special{pa 867 380}%
\special{pa 800 400}%
\special{fp}%
%
\special{pn 13}%
\special{ar 1150 1600 150 150 3.1415927 5.5884470}%
%
\special{pn 13}%
\special{pa 1155 1500}%
\special{pa 1260 1500}%
\special{fp}%
%
\special{pn 13}%
\special{ar 1000 1600 300 300 0.3217506 3.1415927}%
%
\special{pn 13}%
\special{pa 2665 1700}%
\special{pa 3790 1700}%
\special{fp}%
%
\special{pn 13}%
\special{ar 1150 1600 450 300 3.1415927 6.2831853}%
%
\special{pn 13}%
\special{ar 1300 1600 300 300 6.2831853 3.1415927}%
\put(22.5000,-16.0000){\makebox(0,0){$\rightsquigarrow$}}%
%
\special{pn 13}%
\special{ar 3000 1600 400 200 3.6864566 2.5967287}%
%
\special{pn 13}%
\special{ar 3650 1600 450 300 3.1415927 6.2831853}%
%
\special{pn 13}%
\special{ar 3500 1600 300 300 0.3217506 3.1415927}%
%
\special{pn 13}%
\special{ar 3800 1600 300 300 6.2831853 3.1415927}%
%
\special{pn 13}%
\special{ar 3650 1600 150 150 3.1415927 5.5884470}%
%
\special{pn 13}%
\special{pa 2655 1500}%
\special{pa 3760 1500}%
\special{fp}%
%
\special{pn 13}%
\special{pa 1155 1700}%
\special{pa 1280 1700}%
\special{fp}%
\put(9.0000,-21.0000){\makebox(0,0)[lb]{$C_0+_{\Gamma}C_1$}}%
%
\special{pn 4}%
\special{pa 1100 1600}%
\special{pa 800 1600}%
\special{fp}%
\special{sh 1}%
\special{pa 800 1600}%
\special{pa 867 1620}%
\special{pa 853 1600}%
\special{pa 867 1580}%
\special{pa 800 1600}%
\special{fp}%
%
\special{pn 4}%
\special{pa 3100 1600}%
\special{pa 2800 1600}%
\special{fp}%
\special{sh 1}%
\special{pa 2800 1600}%
\special{pa 2867 1620}%
\special{pa 2853 1600}%
\special{pa 2867 1580}%
\special{pa 2800 1600}%
\special{fp}%
\put(34.0000,-20.0000){\makebox(0,0)[lt]{Triple points occur}}%
%
\special{pn 4}%
\special{pa 3400 2000}%
\special{pa 3300 1750}%
\special{fp}%
\special{sh 1}%
\special{pa 3300 1750}%
\special{pa 3306 1819}%
\special{pa 3320 1800}%
\special{pa 3343 1804}%
\special{pa 3300 1750}%
\special{fp}%
%
\special{pn 4}%
\special{pa 3400 2000}%
\special{pa 3320 1520}%
\special{fp}%
\special{sh 1}%
\special{pa 3320 1520}%
\special{pa 3311 1589}%
\special{pa 3329 1573}%
\special{pa 3351 1582}%
\special{pa 3320 1520}%
\special{fp}%
\put(22.5000,-5.0000){\makebox(0,0){$\rightsquigarrow$}}%
%
\special{pn 13}%
\special{ar 3000 500 400 200 0.0000000 6.2831853}%
%
\special{pn 8}%
\special{ar 3950 500 450 300 3.1415927 6.2831853}%
%
\special{pn 8}%
\special{ar 3800 500 300 300 6.2831853 3.1415927}%
%
\special{pn 8}%
\special{ar 3950 500 150 150 3.1415927 6.2831853}%
%
\special{pn 8}%
\special{ar 4100 500 300 300 6.2831853 3.1415927}%
\put(33.0000,-6.5000){\makebox(0,0)[lt]{$\hat{C}_0$}}%
\put(41.7000,-2.3000){\makebox(0,0)[lb]{$C_1$}}%
%
\special{pn 8}%
\special{pa 2600 500}%
\special{pa 4100 500}%
\special{dt 0.045}%
\put(37.0000,-4.3000){\makebox(0,0){$\hat{\Gamma}$}}%
\put(34.0000,-5.0000){\makebox(0,0){$\bullet$}}%
\put(34.5000,-3.0000){\makebox(0,0)[lb]{$\hat{x}$}}%
\put(19.0000,-5.0000){\makebox(0,0){$\bullet$}}%
\put(20.0000,-3.0000){\makebox(0,0)[lb]{$x$}}%
%
\special{pn 4}%
\special{pa 2000 300}%
\special{pa 1900 500}%
\special{fp}%
%
\special{pn 4}%
\special{pa 3450 300}%
\special{pa 3400 500}%
\special{fp}%
\put(36.0000,-6.0000){\makebox(0,0)[lt]{$y$}}%
%
\special{pn 4}%
\special{pa 3600 600}%
\special{pa 3500 500}%
\special{fp}%
\put(35.0000,-5.0000){\makebox(0,0){$\bullet$}}%
%
\special{pn 4}%
\special{sh 1}%
\special{ar 700 500 8 8 0 6.2831853}%
\special{sh 1}%
\special{ar 700 500 8 8 0 6.2831853}%
\put(8.0000,-6.0000){\makebox(0,0)[lt]{$y$}}%
%
\special{pn 4}%
\special{pa 800 600}%
\special{pa 700 500}%
\special{fp}%
\end{picture}}%

%% file: push_appendix_0.tex
{\unitlength 0.1in%
\begin{picture}(16.5000,12.3000)(9.5000,-20.0000)%
%
\special{pn 13}%
\special{pa 1200 1250}%
\special{pa 1200 800}%
\special{fp}%
\special{sh 1}%
\special{pa 1200 800}%
\special{pa 1180 867}%
\special{pa 1200 853}%
\special{pa 1220 867}%
\special{pa 1200 800}%
\special{fp}%
%
\special{pn 13}%
\special{pa 1200 1550}%
\special{pa 1200 2000}%
\special{fp}%
%
\special{pn 8}%
\special{pa 1200 1400}%
\special{pa 1700 1400}%
\special{dt 0.045}%
%
\special{pn 8}%
\special{pa 1900 1400}%
\special{pa 2400 1400}%
\special{dt 0.045}%
%
\special{pn 13}%
\special{pa 2400 800}%
\special{pa 2400 2000}%
\special{fp}%
\special{sh 1}%
\special{pa 2400 2000}%
\special{pa 2420 1933}%
\special{pa 2400 1947}%
\special{pa 2380 1933}%
\special{pa 2400 2000}%
\special{fp}%
%
\special{pn 13}%
\special{ar 1400 1250 200 100 1.5707963 3.1415927}%
%
\special{pn 13}%
\special{ar 1400 1550 200 100 3.1415927 4.7123890}%
%
\special{pn 8}%
\special{pn 8}%
\special{pa 1900 1000}%
\special{pa 1900 1008}%
\special{fp}%
\special{pa 1899 1044}%
\special{pa 1898 1052}%
\special{fp}%
\special{pa 1895 1088}%
\special{pa 1894 1095}%
\special{fp}%
\special{pa 1889 1131}%
\special{pa 1888 1138}%
\special{fp}%
\special{pa 1880 1173}%
\special{pa 1878 1180}%
\special{fp}%
\special{pa 1869 1215}%
\special{pa 1866 1222}%
\special{fp}%
\special{pa 1854 1255}%
\special{pa 1851 1262}%
\special{fp}%
\special{pa 1835 1294}%
\special{pa 1832 1301}%
\special{fp}%
\special{pa 1812 1331}%
\special{pa 1808 1337}%
\special{fp}%
\special{pa 1783 1364}%
\special{pa 1777 1369}%
\special{fp}%
\special{pa 1748 1388}%
\special{pa 1742 1391}%
\special{fp}%
\special{pa 1708 1400}%
\special{pa 1700 1400}%
\special{fp}%
%
\special{pn 8}%
\special{pn 8}%
\special{pa 1700 1000}%
\special{pa 1700 991}%
\special{fp}%
\special{pa 1710 956}%
\special{pa 1714 950}%
\special{fp}%
\special{pa 1734 925}%
\special{pa 1740 920}%
\special{fp}%
\special{pa 1772 904}%
\special{pa 1779 902}%
\special{fp}%
\special{pa 1817 902}%
\special{pa 1825 903}%
\special{fp}%
\special{pa 1858 919}%
\special{pa 1864 923}%
\special{fp}%
\special{pa 1886 948}%
\special{pa 1890 956}%
\special{fp}%
\special{pa 1900 991}%
\special{pa 1900 1000}%
\special{fp}%
%
\special{pn 8}%
\special{pn 8}%
\special{pa 1900 1400}%
\special{pa 1892 1400}%
\special{fp}%
\special{pa 1858 1391}%
\special{pa 1851 1388}%
\special{fp}%
\special{pa 1822 1368}%
\special{pa 1816 1363}%
\special{fp}%
\special{pa 1792 1337}%
\special{pa 1788 1331}%
\special{fp}%
\special{pa 1769 1301}%
\special{pa 1765 1295}%
\special{fp}%
\special{pa 1749 1263}%
\special{pa 1746 1256}%
\special{fp}%
\special{pa 1734 1223}%
\special{pa 1732 1215}%
\special{fp}%
\special{pa 1722 1181}%
\special{pa 1720 1174}%
\special{fp}%
\special{pa 1712 1139}%
\special{pa 1711 1131}%
\special{fp}%
\special{pa 1706 1096}%
\special{pa 1705 1088}%
\special{fp}%
\special{pa 1702 1052}%
\special{pa 1701 1044}%
\special{fp}%
\special{pa 1700 1008}%
\special{pa 1700 1000}%
\special{fp}%
%
\special{pn 13}%
\special{pa 1400 1350}%
\special{pa 1700 1350}%
\special{fp}%
%
\special{pn 13}%
\special{ar 1700 1000 150 350 6.2831853 1.5707963}%
%
\special{pn 13}%
\special{ar 1800 1000 50 50 3.1415927 6.2831853}%
%
\special{pn 13}%
\special{ar 1900 1010 150 340 1.5707963 3.1549252}%
%
\special{pn 13}%
\special{pa 1400 1450}%
\special{pa 1700 1450}%
\special{fp}%
%
\special{pn 13}%
\special{ar 1700 1000 250 450 6.2831853 1.5707963}%
%
\special{pn 13}%
\special{ar 1800 1000 150 150 3.1415927 6.2831853}%
%
\special{pn 13}%
\special{ar 1900 1000 250 450 1.5707963 3.1415927}%
%
\special{pn 13}%
\special{ar 1900 1400 300 50 4.7123890 1.5707963}%
\put(17.0000,-15.5000){\makebox(0,0){$\overline{C}_0$}}%
\put(21.0000,-9.0000){\makebox(0,0)[lb]{$\Gamma$}}%
%
\special{pn 4}%
\special{pa 2100 900}%
\special{pa 1900 1110}%
\special{fp}%
\put(25.0000,-18.0000){\makebox(0,0){$C_1$}}%
%
\special{pn 13}%
\special{pa 2200 1380}%
\special{pa 2400 1380}%
\special{dt 0.045}%
\put(23.0000,-13.0000){\makebox(0,0){$\overline{\Gamma}$}}%
%
\special{pn 8}%
\special{pa 1400 1900}%
\special{pa 1400 900}%
\special{fp}%
%
\special{pn 8}%
\special{pa 2000 1900}%
\special{pa 2000 900}%
\special{fp}%
\put(14.0000,-20.0000){\makebox(0,0){$C_0$}}%
\put(20.0000,-20.0000){\makebox(0,0){$C_1$}}%
\put(17.0000,-20.0000){\makebox(0,0){$\cdots$}}%
\put(9.5000,-14.0000){\makebox(0,0)[lb]{$b_0$}}%
\put(24.3000,-13.7000){\makebox(0,0)[lb]{$b_1$}}%
%
\special{pn 4}%
\special{pa 2300 1450}%
\special{pa 2600 1600}%
\special{fp}%
\put(26.0000,-16.0000){\makebox(0,0)[lt]{$R_{\overline{v}_0}\subset R_0$}}%
\end{picture}}%

%% file: 5_T.tex
\section{Computations of $T^{\pm}$ and $T^{\St}_{\Gamma}$}\label{s:T}
Here we give algorithms to compute $T^{\pm}_{\Gamma}(C_0,C_1)$ and $T^{\St}_{\Gamma}(C_0,C_1)$ for any generic closed curves $C_0,C_1$ in general position and any bridge $\Gamma$ between them.
This completes our generalized sum formulas.

\subsection{$T^{\pm}$}
It is enough to compute $T^{\pm}(c_0,c_1)$ for any $c_0,c_1\in\GI$ in general position.
We call a curve equivalent to $K_1$ a \emph{simple closed curve}.
\begin{defn}
For $c\in\GI$, let $S(c)$ be the disjoint union of simple closed curves obtained by splicing all the double points of $c$ as shown in Figure~\ref{fig:splicing}.
\end{defn}
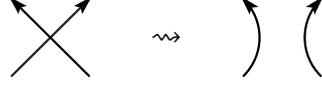
\begin{figure}
\begin{center}
\input{splicing}
\end{center}
\caption{Splice of a double point}
\label{fig:splicing}
\end{figure}


\begin{lem}[{\cite[Lemma~3]{MendesRomero}}]\label{lem:splice_not_change_T}
Let $c_0,c_1\in\GI$ be in general position.
Suppose $S(c_0)=c_0^1\sqcup\dots\sqcup c_0^k$, where each simple closed curve $c_0^i$ inherits the orientation of $c_0$, and the splicing is done in sufficiently small neighborhood of double points of $c_0$ so that $\bigcup_jc^j_0\cap c_1=c_0\cap c_1$.
Then
\[
 T^{\pm}(c_0,c_1)=\sum_{1\le j\le k}T^{\pm}(c_0^j,c_1).
\]
\end{lem}
\begin{proof}
A generic separating homotopy of $c_0$ from $c_1$ induces generic separating homotopies of $c_0^1,\dots,c_0^k$ from $c_1$.
Some triple points occurring in the original homotopy may be replaced by pairs of positive and negative tangencies between $c_0^j$ and $c_1$ as in Figure~\ref{fig:new_tangencies}, but they cancel with each other in $\sum_iT^{\pm}(c_0^i,c_1)$.
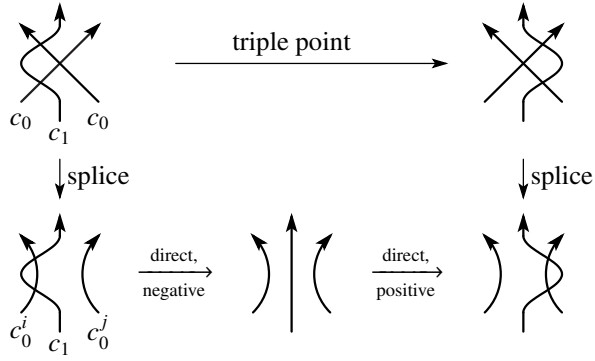
\begin{figure}
\begin{center}
\input{new_tangencies}
\end{center}
\caption{New tangencies}
\label{fig:new_tangencies}
\end{figure}
\end{proof}

Thus the case both $c_0,c_1$ are simple closed curves is essential in the computation of $T^{\pm}$, because in general we have by Lemma~\ref{lem:splice_not_change_T}
\[
 T^{\pm}(c_0,c_1)=\sum_{j=1}^k\sum_{m=1}^lT^{\pm}(c_0^j,c_1^m),
\]
where $S(c_0)=c_0^1\sqcup\dots\sqcup c_0^k$ and $S(c_1)=c_1^1\sqcup\dots\sqcup c_1^l$ and each simple closed curve $c_i^j$ ($i=0,1$) inherits the orientation of $c_i$.

Recall the rotation number $r$ from Definition~\ref{def:rotation_number}.
If $c$ is a simple closed curve then $r(c)=\pm 1$.

\begin{prop}\label{prop:T^pm}
Let $c_0,c_1$ be non-separated oriented simple closed curves in general position and let $D_i$ be the disk bounded by $C_i$.
If $r(c_0)=r(c_1)$ then
\[
 T^+(c_0,c_1)=-\frac{1}{2}\abs{C_0\cap C_1}+\abs{\pi_0(D_0\cap D_1)},\quad
 T^-(c_0,c_1)=-\abs{\pi_0(D_0\cap D_1)}.
\]
If $r(c_0)=-r(c_1)$ then
\[
 T^+(c_0,c_1)=-\abs{\pi_0(D_0\cap D_1},\quad
 T^-(c_0,c_1)=-\frac{1}{2}\abs{C_0\cap C_1}+\abs{\pi_0(D_0\cap D_1)}.
\]
\end{prop}

\begin{rem}
By Proposition~\ref{prop:T^pm} we have $T^+(c_0,c_1)+T^-(c_0,c_1)=-\abs{C_0\cap C_1}/2$ for simple closed curves $c_0,c_1$.
In fact this is true for any generic closed curves $c_0,c_1$ in general position.
This can be proved by choosing a bridge $\Gamma$ so that $c_0,c_1$ are compatible orientations, and by using Theorem~\ref{thm:main_J} and the fact that in general $J^+(C)-J^-(C)=\abs{D(C)}$ \cite[\S4]{Arnold}.
\end{rem}

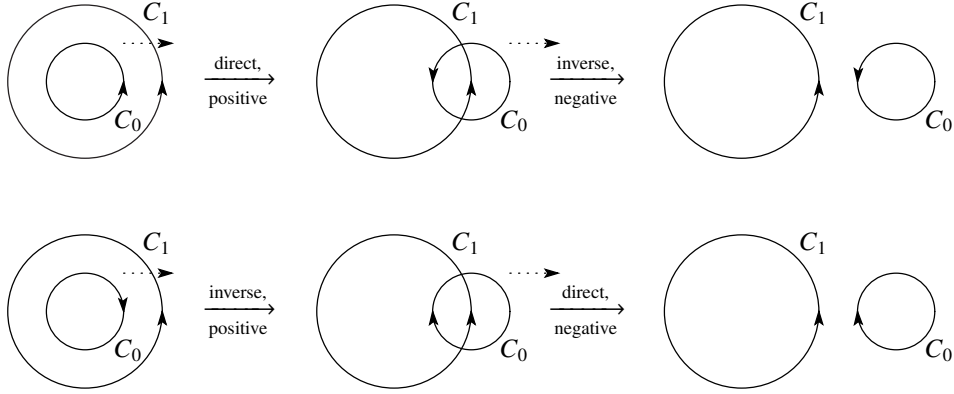
\begin{figure}
\begin{center}
\input{contained}
\end{center}
\caption{The case $C_0\subset D_1$}
\label{fig:contained}
\end{figure}

First consider the case that $C_0\cap C_1=\emptyset$, in which case $C_0\subset D_1$ or $C_1\subset D_0$.
In this case Proposition~\ref{prop:T^pm} claims that $T^{\pm}(c_0,c_1)=\pm r(c_0)r(c_1)$.
Figure~\ref{fig:contained} would be enough for the proof of this case.
Only the cases $(r(c_0),r(c_1))=(+1,+1),(-1,+1)$ are shown, but for the remaining cases we only need to reverse the orientations, that does not change the types (direct or inverse) of tangencies.

For the case $C_0\cap C_1\ne\emptyset$, we note that in this case $C_1$ divides $D_0$ into finitely many polygons, each of which is a $2k$-gon for some $k\ge 1$.
Notice also that the components of $D_0\cap D_1$ appear alternately in $D_0$, namely if two components of $D_0\setminus(C_1\cap D_0)$ are adjacent to each other, then one component is in $D_0\cap D_1$ and the other is not.

\begin{lem}\label{lem:adjacent_to_2-gon}
There is at least one component $R$ of $D_0\cap D_1$ satisfying one of the following;
\begin{enumerate}[(i)]
\item
	all but possibly one components of $D_0\setminus(C_1\cap D_0)$ adjacent to $R$ in $D_0$ are $2$-gons, or
\item
	$R$ is a $2$-gon (see Figure~\ref{fig:polygons}).
\end{enumerate}
\end{lem}
\begin{figure}
\begin{center}
\input{polygons}
\end{center}
\caption{Polygons}
\label{fig:polygons}
\end{figure}
\begin{proof}
If $D_0\cap D_1$ is connected, then $C_0,C_1$ are as in Figure~\ref{fig:star}, and hence $D_0\cap D_1$ is a $2$-gon or otherwise all the adjacent components to $D_0\cap D_1$ in $D_0$ are $2$-gons.
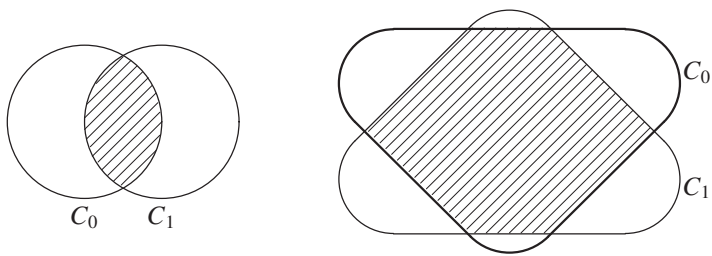
\begin{figure}
\begin{center}
\input{star}
\end{center}
\caption{The cases $D_0\cap D_1$ are connected}
\label{fig:star}
\end{figure}

Suppose $\abs{\pi_0(D_0\cap D_1)}\ge 2$, and choose any component $R_0$ of $D_0\cap D_1$.
If $R_0$ is not a component as claimed,
then choose an adjacent component $R_1$ to $R_0$ in $D_0$ that is not a $2$-gon,
and by a surgery along $\overline{R}_0\cap\overline{R}_1$ we `cut' $C_0$ into two simple closed curves (see Figure~\ref{fig:no_2-gon}).
\begin{figure}
\begin{center}
\input{no_2-gon}
\end{center}
\caption{A surgery}
\label{fig:no_2-gon}
\end{figure}
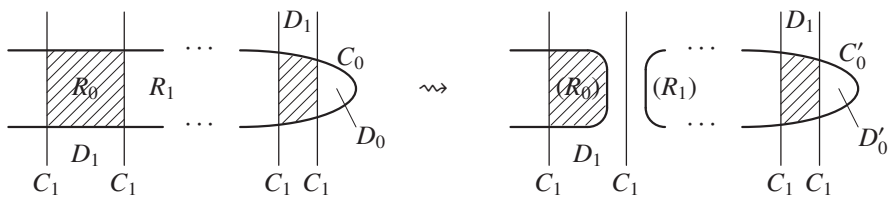
Let $C'_0$ be one of these simple closed curves that bounds the disk $D'_0$ with $D'_0\cap R_1\ne\emptyset$.
Then $D'_0\cap D_1$ has at least one component (otherwise $R_1$ is a $2$-gon) but fewer than $D_0\cap D_1$.
Repeating similar surgeries to $C'_0$ we can find a component $R_2\subset D_0\cap D_1$ and $R_3\subset D_0$ adjacent to $R_2$ such that the surgery along $\overline{R}_2\cap\overline{R}_3$ produces a simple closed curve $C''_0$ such that $R':=D''_0\cap D_1$ is connected, where $D''_0$ is the disk bounded by $C''_0$ (see Figure~\ref{fig:claimed_region}).
All the adjacent components to $R'$ in $D''_0$ are thus $2$-gon, and in the original $D_0$ all the adjacent components to $R'$ except for $R_3$ are $2$-gons.
\begin{figure}
\begin{center}
\input{claimed_region}
\end{center}
\caption{$R_2,R_3$}
\label{fig:claimed_region}
\end{figure}
\end{proof}

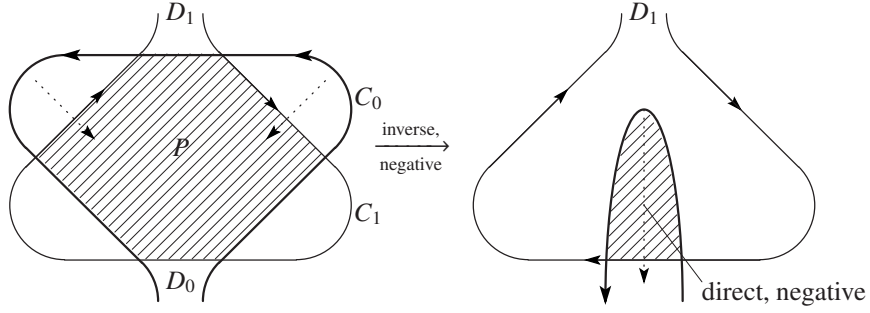
\begin{figure}
\begin{center}
\input{star_resolution}
\end{center}
\caption{Transformation of $C_0$ along $P$ (the case $r(c_0)=-r(c_1)$)}
\label{fig:star_resolution}
\end{figure}

Choose a component $P$ of $D_0\cap D_1$ as claimed in Lemma~\ref{lem:adjacent_to_2-gon}.
If $P$ is a $2k$-gon, then we can transform $C_0$ as in Figure~\ref{fig:star_resolution} through $k$ negative tangencies so that the disk $D'_0$ bounded by the resulting curve $C'_0$ satisfies $D'_0\cap D_1=(D_0\cap D_1)\setminus P$.
If $r(c_0)=r(c_1)$ then $k-1$ of the tangencies are direct and the remaining one is inverse,
while if $r(c_0)=-r(c_1)$ then $k-1$ of the tangencies are inverse and the remaining one is direct.
Repeating the similar transformations we see the following.

\begin{prop}\label{prop:T_k}
Let $c_0,c_1$ be oriented simple closed curves in general position and $D_i$ ($i=0,1$) the disks bounded by $C_i:=c_i(S^1)$.
Suppose $D_0\cap D_1=P_1\sqcup\dots\sqcup P_s$ where each $P_i$ is a $2k_i$-gon.
Then if $r(c_0)=r(c_1)$ we have
\[
 T^+(c_0,c_1)=-\sum_{i=1}^s(k_i-1),\quad
 T^-(c_0,c_1)=-s,
\]
and if $r(c_0)=-r(c_1)$ we have
\[
 T^+(c_0,c_1)=-s,\quad
 T^-(c_0,c_1)=-\sum_{i=1}^s(k_i-1).
\]
\end{prop}

\begin{proof}[Proof of Proposition~\ref{prop:T^pm}, the case $C_0\cap C_1\ne\emptyset$]
Let $P_i$ ($1\le i\le s$) be as in Proposition~\ref{prop:T_k}.
Then $s=\abs{\pi_0(D_0\cap D_1)}$.
Each point in $C_0\cap C_1$ is a vertex of exactly one $P_i$ for some $i$,
and hence $\abs{C_0\cap C_1}=\sum_{1\le i\le s}2k_i$.
Thus the claim follows from Proposition~\ref{prop:T_k}.
\end{proof}

\subsection{$T^{\St}_{\Gamma}$}
Let $C_0,C_1$ be generic closed curves in general position and $\Gamma$ a bridge between them.
A triple point occurs in a generic separating homotopy of $C_0$ from $C_1$ when a double point of $C_i$ gets across $C_{1-i}$ ($i=0,1$).
To compute $T^{\St}_{\Gamma}(C_0,C_1)$ we need to determine the sign $(-1)^q$ of the newborn triangle at each triple point (see Definition~\ref{def:signs} and Example~\ref{ex:signs_triplepoint}).

\begin{defn}[{\cite[\S6]{Arnold}}]
Let $c\in\GI$ and let $x=c(t_0)\in C\setminus D(C)$, where $C:=c(S^1)$.
For a double point $p=c(t_1)=c(t_2)\in D(C)$, $t_0<t_1<t_2<1+t_0$, define $e^c_x(p)$ as $+1$ (resp.\ $-1$) if $\ang{c'(t_1),c'(t_2)}$ is a positive (resp.\ negative) basis of $\R^2$.
\end{defn}

We can check that $e^{-c}_x(p)=-e^c_x(p)$.

Below we choose compatible parameters $c_0,c_1$ of $C_0,C_1$ with respect to $\Gamma$.

\begin{lem}\label{lem:sign_through_triplepoint}
When a double point $p$ of $c_0$ gets across $c_1$ from the left in a separating homotopy, the sign $(-1)^q$ of the corresponding newborn triangle is equal to $e^{c_0}_{b_0}(p)$, where $b_0\in C_0:=c_0(S^1)$ is one of the endpoints of $\Gamma$.
\end{lem}
\begin{proof}
Check all the cases.
See Figure~\ref{fig:all_the_newborn_triangles}.
\end{proof}
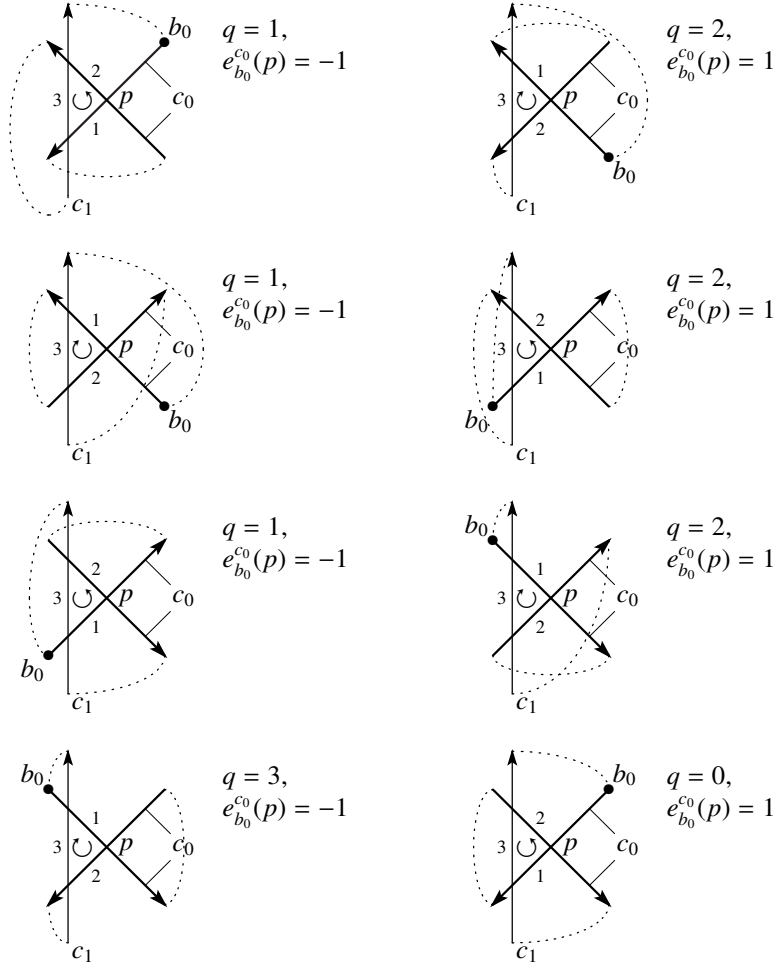
\begin{figure}
\begin{center}
\input{all_the_newborn_triangles}
\end{center}
\caption{All the possible newborn triangles after a double point $p$ gets across $c_1$ from the left}
\label{fig:all_the_newborn_triangles}
\end{figure}

Notice that if a double point $p$ of $c_0$ gets across $c_1$ from the right, then the sign of the newborn triangle is $-e^{c_0}_{b_0}(p)$.

\begin{prop}\label{prop:T^St}
Let $C_0,C_1$ be generic closed curves in general position and $\Gamma$ a bridge between them.
Let $p_1,\dots,p_k$ (resp.\ $q_1,\dots,q_l$) be the double points of $C_0$ (resp.\ $C_1$) and $R_1^i$ (resp.\ $R_0^j$) the component of $\R^2\setminus C_1$ (resp.\ $\R^2\setminus C_0$) that contains $p_i$ (resp.\ $q_j$).
Then
\[
 T^{\St}_{\Gamma}(C_0,C_1)=\sum_{1\le j\le k}e^{c_0}_{b_0}(p_j)\ind^{R_1^j}(c_1)+\sum_{1\le m\le l}e^{c_1}_{b_1}(q_m)\ind^{R_0^m}(c_0),
\]
where $b_i\in C_i$ ($i=0,1$) are the endpoints of $\Gamma$.
\end{prop}

The right hand side of the formula in Proposition~\ref{prop:T^St} does not depend on the compatible parameters $c_0,c_1$ since reversing the orientations change the signs of both $e$ and $\ind$.

\begin{proof}[Proof of Proposition~\ref{prop:T^St}]
Suppose that $p_j$ gets across $c_1$ from the right (resp.\ left) $l^+$ times (resp.\ $l^-$ times) in a generic separating homotopy.
Then $p_j$ contributes to $T^{\St}_{\Gamma}(C_0,C_1)$ by $(l^+-l^-)e^{c_0}_{b_0}(p_j)$ by Lemma~\ref{lem:sign_through_triplepoint}.
On the other hand, we have seen in the proof of Lemma~\ref{lem:adjacent_indices} that $l^--l^+=\ind^{R_1^{\infty}}(c_1)-\ind^{R_1^j}(c_1)$.
Since $\ind^{R_1^{\infty}}(c_1)=0$ we have $l^+-l^-=\ind^{R_1^j}(c_1)$.
The same observation on $q_m$ completes the proof.
\end{proof}

\begin{rem}
Proposition~\ref{prop:T^St} relates to the ``pushing away formula'' in \cite[\S6]{Arnold}, which can be seen as the case that $R_1^i$ are common for any $i\ge 1$ and all $R_0^j$ are the unbounded component of $\R^2\setminus C_0$.
\end{rem}

\begin{cor}\label{cor:triple_points_change_T^St}
Let $C_i$ ($i=0,1$), $\Gamma$, $\overline{C}_0$ and $\overline{\Gamma}$ be as in \S\ref{ss:MendesRomero}.
Suppose $\Int\Gamma\cap C_0=\{x_1,\dots,x_p\}$ and let $R^{(j)}_1\subset\R^2\setminus C_1$ be the component that contains $x_j$.
Then
\[
 T^{\St}_{\overline{\Gamma}}(\overline{C}_0,C_1)=T^{\St}_{\Gamma}(C_0,C_1)-2\sum_{1\le j\le p}s(x_j)\ind^{R^{(j)}}(C_1).
\]
\end{cor}

\begin{proof}
To make use of Proposition~\ref{prop:T^St} we need to find double points of $\overline{C}_0$ that do not come from those of $C_0$.
There exist two kinds of such double points:
\begin{enumerate}[(i)]
\item
	Corresponding to each $x_j$, the curve $\overline{C}_0$ has two more double points $y_j,z_j$ than $C_0$.
\item
	Corresponding to each double points of $\Gamma$, the curve $\overline{C}_0$ has four more double points $y_j,z_j$ than $C_0$.
\end{enumerate}
Firstly we calculate the contribution of double points of type (i) to $T^{\St}_{\overline{\Gamma}}(\overline{C}_0,C_1)$.
As shown in Figure~\ref{fig:e(y)=e(z)=-s(x)},
we have $e^{\overline{c}_0}_{\overline{b}_0}(y_j)=e^{\overline{c}_0}_{\overline{b}_0}(z_j)=-s(x_j)$.
\begin{figure}
\begin{center}
\input{e_y=e_z=-s_x}
\end{center}
\caption{$e^{\overline{c}_0}_{\overline{b}_0}(y_j)=e^{\overline{c}_0}_{\overline{b}_0}(z_j)=-s(x_j)$; the case $s(x_j)=+1$}
\label{fig:e(y)=e(z)=-s(x)}
\end{figure}
Since $y_j,z_j\in R^{(j)}_1$, the contribution of double points of type (i) to $T^{\St}_{\overline{\Gamma}}(\overline{C}_0,C_1)$ is $-2\sum_{1\le j\le p}s(x_j)\ind^{R^{(j)}}(C_1)$.

All the contributions of type (ii) double points cancel out, because
\begin{itemize}
\item
	all the four double points are included in a common region of $\R^2\setminus C_1$, and
\item
	$e=+1$ for two of these double points and $e=-1$ for remaining two, as we can see in Figure~\ref{fig:for_e_cancel}.\qedhere
\end{itemize}
\begin{figure}
\begin{center}
\input{four_e_cancel}
\end{center}
\caption{Any double point of $\Gamma$ does not contribute to $T^{\St}_{\overline{\Gamma}}(\overline{C}_0,C_1)$}
\label{fig:for_e_cancel}
\end{figure}
\end{proof}

%% file: splicing.tex
{\unitlength 0.1in%
\begin{picture}(16.0000,4.0000)(6.0000,-12.0000)%
%
\special{pn 13}%
\special{pa 600 1200}%
\special{pa 1000 800}%
\special{fp}%
\special{sh 1}%
\special{pa 1000 800}%
\special{pa 939 833}%
\special{pa 962 838}%
\special{pa 967 861}%
\special{pa 1000 800}%
\special{fp}%
%
\special{pn 13}%
\special{pa 1000 1200}%
\special{pa 600 800}%
\special{fp}%
\special{sh 1}%
\special{pa 600 800}%
\special{pa 633 861}%
\special{pa 638 838}%
\special{pa 661 833}%
\special{pa 600 800}%
\special{fp}%
\put(14.0000,-10.0000){\makebox(0,0){$\rightsquigarrow$}}%
%
\special{pn 13}%
\special{ar 1600 1000 283 283 5.4977871 0.7853982}%
%
\special{pn 13}%
\special{pa 1810 811}%
\special{pa 1800 800}%
\special{fp}%
\special{sh 1}%
\special{pa 1800 800}%
\special{pa 1830 863}%
\special{pa 1836 839}%
\special{pa 1860 836}%
\special{pa 1800 800}%
\special{fp}%
%
\special{pn 13}%
\special{ar 2400 1000 283 283 2.3561945 3.9269908}%
%
\special{pn 13}%
\special{pa 2190 811}%
\special{pa 2200 800}%
\special{fp}%
\special{sh 1}%
\special{pa 2200 800}%
\special{pa 2140 836}%
\special{pa 2164 839}%
\special{pa 2170 863}%
\special{pa 2200 800}%
\special{fp}%
\end{picture}}%

%% file: new_tangencies.tex
{\unitlength 0.1in%
\begin{picture}(31.8500,17.1500)(2.1500,-24.1500)%
%
\special{pn 13}%
\special{pa 600 1200}%
\special{pa 1000 800}%
\special{fp}%
\special{sh 1}%
\special{pa 1000 800}%
\special{pa 939 833}%
\special{pa 962 838}%
\special{pa 967 861}%
\special{pa 1000 800}%
\special{fp}%
%
\special{pn 13}%
\special{pa 1000 1200}%
\special{pa 600 800}%
\special{fp}%
\special{sh 1}%
\special{pa 600 800}%
\special{pa 633 861}%
\special{pa 638 838}%
\special{pa 661 833}%
\special{pa 600 800}%
\special{fp}%
\put(20.0000,-9.0000){\makebox(0,0){triple point}}%
%
\special{pn 13}%
\special{ar 400 2100 283 283 5.4977871 0.7853982}%
%
\special{pn 13}%
\special{pa 610 1911}%
\special{pa 600 1900}%
\special{fp}%
\special{sh 1}%
\special{pa 600 1900}%
\special{pa 630 1963}%
\special{pa 636 1939}%
\special{pa 660 1936}%
\special{pa 600 1900}%
\special{fp}%
%
\special{pn 13}%
\special{ar 1200 2100 283 283 2.3561945 3.9269908}%
%
\special{pn 13}%
\special{pa 990 1911}%
\special{pa 1000 1900}%
\special{fp}%
\special{sh 1}%
\special{pa 1000 1900}%
\special{pa 940 1936}%
\special{pa 964 1939}%
\special{pa 970 1963}%
\special{pa 1000 1900}%
\special{fp}%
\special{pn 13}%
\special{pa 800 800}%
\special{pa 800 805}%
\special{pa 799 810}%
\special{pa 795 820}%
\special{pa 789 830}%
\special{pa 781 840}%
\special{pa 771 850}%
\special{pa 759 860}%
\special{pa 731 880}%
\special{pa 723 885}%
\special{pa 716 890}%
\special{pa 684 910}%
\special{pa 677 915}%
\special{pa 669 920}%
\special{pa 641 940}%
\special{pa 629 950}%
\special{pa 619 960}%
\special{pa 611 970}%
\special{pa 605 980}%
\special{pa 601 990}%
\special{pa 600 995}%
\special{pa 600 1005}%
\special{pa 601 1010}%
\special{pa 605 1020}%
\special{pa 611 1030}%
\special{pa 619 1040}%
\special{pa 629 1050}%
\special{pa 641 1060}%
\special{pa 669 1080}%
\special{pa 677 1085}%
\special{pa 684 1090}%
\special{pa 716 1110}%
\special{pa 723 1115}%
\special{pa 731 1120}%
\special{pa 759 1140}%
\special{pa 771 1150}%
\special{pa 781 1160}%
\special{pa 789 1170}%
\special{pa 795 1180}%
\special{pa 799 1190}%
\special{pa 800 1195}%
\special{pa 800 1200}%
\special{fp}%
%
\special{pn 13}%
\special{pa 800 800}%
\special{pa 800 700}%
\special{fp}%
\special{sh 1}%
\special{pa 800 700}%
\special{pa 780 767}%
\special{pa 800 753}%
\special{pa 820 767}%
\special{pa 800 700}%
\special{fp}%
%
\special{pn 13}%
\special{pa 800 1300}%
\special{pa 800 1200}%
\special{fp}%
%
\special{pn 13}%
\special{pa 3000 1200}%
\special{pa 3400 800}%
\special{fp}%
\special{sh 1}%
\special{pa 3400 800}%
\special{pa 3339 833}%
\special{pa 3362 838}%
\special{pa 3367 861}%
\special{pa 3400 800}%
\special{fp}%
%
\special{pn 13}%
\special{pa 3400 1200}%
\special{pa 3000 800}%
\special{fp}%
\special{sh 1}%
\special{pa 3000 800}%
\special{pa 3033 861}%
\special{pa 3038 838}%
\special{pa 3061 833}%
\special{pa 3000 800}%
\special{fp}%
\special{pn 13}%
\special{pa 3200 800}%
\special{pa 3200 805}%
\special{pa 3201 810}%
\special{pa 3205 820}%
\special{pa 3211 830}%
\special{pa 3219 840}%
\special{pa 3229 850}%
\special{pa 3241 860}%
\special{pa 3269 880}%
\special{pa 3277 885}%
\special{pa 3284 890}%
\special{pa 3316 910}%
\special{pa 3323 915}%
\special{pa 3331 920}%
\special{pa 3359 940}%
\special{pa 3371 950}%
\special{pa 3381 960}%
\special{pa 3389 970}%
\special{pa 3395 980}%
\special{pa 3399 990}%
\special{pa 3400 995}%
\special{pa 3400 1005}%
\special{pa 3399 1010}%
\special{pa 3395 1020}%
\special{pa 3389 1030}%
\special{pa 3381 1040}%
\special{pa 3371 1050}%
\special{pa 3359 1060}%
\special{pa 3331 1080}%
\special{pa 3323 1085}%
\special{pa 3316 1090}%
\special{pa 3284 1110}%
\special{pa 3277 1115}%
\special{pa 3269 1120}%
\special{pa 3241 1140}%
\special{pa 3229 1150}%
\special{pa 3219 1160}%
\special{pa 3211 1170}%
\special{pa 3205 1180}%
\special{pa 3201 1190}%
\special{pa 3200 1195}%
\special{pa 3200 1200}%
\special{fp}%
%
\special{pn 13}%
\special{pa 3200 800}%
\special{pa 3200 700}%
\special{fp}%
\special{sh 1}%
\special{pa 3200 700}%
\special{pa 3180 767}%
\special{pa 3200 753}%
\special{pa 3220 767}%
\special{pa 3200 700}%
\special{fp}%
%
\special{pn 13}%
\special{pa 3200 1300}%
\special{pa 3200 1200}%
\special{fp}%
\special{pn 13}%
\special{pa 800 1900}%
\special{pa 800 1905}%
\special{pa 799 1910}%
\special{pa 795 1920}%
\special{pa 789 1930}%
\special{pa 781 1940}%
\special{pa 771 1950}%
\special{pa 759 1960}%
\special{pa 731 1980}%
\special{pa 723 1985}%
\special{pa 716 1990}%
\special{pa 684 2010}%
\special{pa 677 2015}%
\special{pa 669 2020}%
\special{pa 641 2040}%
\special{pa 629 2050}%
\special{pa 619 2060}%
\special{pa 611 2070}%
\special{pa 605 2080}%
\special{pa 601 2090}%
\special{pa 600 2095}%
\special{pa 600 2105}%
\special{pa 601 2110}%
\special{pa 605 2120}%
\special{pa 611 2130}%
\special{pa 619 2140}%
\special{pa 629 2150}%
\special{pa 641 2160}%
\special{pa 669 2180}%
\special{pa 677 2185}%
\special{pa 684 2190}%
\special{pa 716 2210}%
\special{pa 723 2215}%
\special{pa 731 2220}%
\special{pa 759 2240}%
\special{pa 771 2250}%
\special{pa 781 2260}%
\special{pa 789 2270}%
\special{pa 795 2280}%
\special{pa 799 2290}%
\special{pa 800 2295}%
\special{pa 800 2300}%
\special{fp}%
%
\special{pn 13}%
\special{pa 800 1900}%
\special{pa 800 1800}%
\special{fp}%
\special{sh 1}%
\special{pa 800 1800}%
\special{pa 780 1867}%
\special{pa 800 1853}%
\special{pa 820 1867}%
\special{pa 800 1800}%
\special{fp}%
%
\special{pn 13}%
\special{pa 800 2400}%
\special{pa 800 2300}%
\special{fp}%
\put(14.0000,-21.0000){\makebox(0,0){$\xrightarrow[\text{negative}]{\text{direct},}$}}%
%
\special{pn 13}%
\special{ar 1600 2100 283 283 5.4977871 0.7853982}%
%
\special{pn 13}%
\special{pa 1810 1911}%
\special{pa 1800 1900}%
\special{fp}%
\special{sh 1}%
\special{pa 1800 1900}%
\special{pa 1830 1963}%
\special{pa 1836 1939}%
\special{pa 1860 1936}%
\special{pa 1800 1900}%
\special{fp}%
%
\special{pn 13}%
\special{ar 2400 2100 283 283 2.3561945 3.9269908}%
%
\special{pn 13}%
\special{pa 2190 1911}%
\special{pa 2200 1900}%
\special{fp}%
\special{sh 1}%
\special{pa 2200 1900}%
\special{pa 2140 1936}%
\special{pa 2164 1939}%
\special{pa 2170 1963}%
\special{pa 2200 1900}%
\special{fp}%
%
\special{pn 13}%
\special{pa 2000 2400}%
\special{pa 2000 1800}%
\special{fp}%
\special{sh 1}%
\special{pa 2000 1800}%
\special{pa 1980 1867}%
\special{pa 2000 1853}%
\special{pa 2020 1867}%
\special{pa 2000 1800}%
\special{fp}%
\put(26.0000,-21.0000){\makebox(0,0){$\xrightarrow[\text{positive}]{\text{direct},}$}}%
%
\special{pn 13}%
\special{ar 2800 2100 283 283 5.4977871 0.7853982}%
%
\special{pn 13}%
\special{pa 3010 1911}%
\special{pa 3000 1900}%
\special{fp}%
\special{sh 1}%
\special{pa 3000 1900}%
\special{pa 3030 1963}%
\special{pa 3036 1939}%
\special{pa 3060 1936}%
\special{pa 3000 1900}%
\special{fp}%
%
\special{pn 13}%
\special{ar 3600 2100 283 283 2.3561945 3.9269908}%
%
\special{pn 13}%
\special{pa 3390 1911}%
\special{pa 3400 1900}%
\special{fp}%
\special{sh 1}%
\special{pa 3400 1900}%
\special{pa 3340 1936}%
\special{pa 3364 1939}%
\special{pa 3370 1963}%
\special{pa 3400 1900}%
\special{fp}%
\special{pn 13}%
\special{pa 3200 1900}%
\special{pa 3200 1905}%
\special{pa 3201 1910}%
\special{pa 3205 1920}%
\special{pa 3211 1930}%
\special{pa 3219 1940}%
\special{pa 3229 1950}%
\special{pa 3241 1960}%
\special{pa 3269 1980}%
\special{pa 3277 1985}%
\special{pa 3284 1990}%
\special{pa 3316 2010}%
\special{pa 3323 2015}%
\special{pa 3331 2020}%
\special{pa 3359 2040}%
\special{pa 3371 2050}%
\special{pa 3381 2060}%
\special{pa 3389 2070}%
\special{pa 3395 2080}%
\special{pa 3399 2090}%
\special{pa 3400 2095}%
\special{pa 3400 2105}%
\special{pa 3399 2110}%
\special{pa 3395 2120}%
\special{pa 3389 2130}%
\special{pa 3381 2140}%
\special{pa 3371 2150}%
\special{pa 3359 2160}%
\special{pa 3331 2180}%
\special{pa 3323 2185}%
\special{pa 3316 2190}%
\special{pa 3284 2210}%
\special{pa 3277 2215}%
\special{pa 3269 2220}%
\special{pa 3241 2240}%
\special{pa 3229 2250}%
\special{pa 3219 2260}%
\special{pa 3211 2270}%
\special{pa 3205 2280}%
\special{pa 3201 2290}%
\special{pa 3200 2295}%
\special{pa 3200 2300}%
\special{fp}%
%
\special{pn 13}%
\special{pa 3200 1900}%
\special{pa 3200 1800}%
\special{fp}%
\special{sh 1}%
\special{pa 3200 1800}%
\special{pa 3180 1867}%
\special{pa 3200 1853}%
\special{pa 3220 1867}%
\special{pa 3200 1800}%
\special{fp}%
%
\special{pn 13}%
\special{pa 3200 2400}%
\special{pa 3200 2300}%
\special{fp}%
\put(10.0000,-16.0000){\makebox(0,0){splice}}%
\put(34.0000,-16.0000){\makebox(0,0){splice}}%
\put(6.0000,-13.0000){\makebox(0,0){$c_0$}}%
\put(10.0000,-13.0000){\makebox(0,0){$c_0$}}%
\put(8.0000,-13.8000){\makebox(0,0){$c_1$}}%
\put(6.0000,-24.0000){\makebox(0,0){$c^i_0$}}%
\put(10.0000,-24.0000){\makebox(0,0){$c^j_0$}}%
\put(8.0000,-24.8000){\makebox(0,0){$c_1$}}%
%
\special{pn 8}%
\special{pa 1400 1000}%
\special{pa 2800 1000}%
\special{fp}%
\special{sh 1}%
\special{pa 2800 1000}%
\special{pa 2733 980}%
\special{pa 2747 1000}%
\special{pa 2733 1020}%
\special{pa 2800 1000}%
\special{fp}%
%
\special{pn 8}%
\special{pa 800 1500}%
\special{pa 800 1700}%
\special{fp}%
\special{sh 1}%
\special{pa 800 1700}%
\special{pa 820 1633}%
\special{pa 800 1647}%
\special{pa 780 1633}%
\special{pa 800 1700}%
\special{fp}%
%
\special{pn 8}%
\special{pa 3200 1500}%
\special{pa 3200 1700}%
\special{fp}%
\special{sh 1}%
\special{pa 3200 1700}%
\special{pa 3220 1633}%
\special{pa 3200 1647}%
\special{pa 3180 1633}%
\special{pa 3200 1700}%
\special{fp}%
\end{picture}}%

%% file: contained.tex
{\unitlength 0.1in%
\begin{picture}(48.0000,20.3000)(6.0000,-24.0000)%
%
\special{pn 8}%
\special{ar 1000 800 400 400 0.0000000 6.2831853}%
%
\special{pn 8}%
\special{pa 1400 815}%
\special{pa 1400 800}%
\special{fp}%
\special{sh 1}%
\special{pa 1400 800}%
\special{pa 1380 867}%
\special{pa 1400 853}%
\special{pa 1420 867}%
\special{pa 1400 800}%
\special{fp}%
%
\special{pn 8}%
\special{ar 1000 800 200 200 0.0000000 6.2831853}%
%
\special{pn 8}%
\special{pa 1199 815}%
\special{pa 1200 800}%
\special{fp}%
\special{sh 1}%
\special{pa 1200 800}%
\special{pa 1176 865}%
\special{pa 1196 853}%
\special{pa 1216 868}%
\special{pa 1200 800}%
\special{fp}%
%
\special{pn 8}%
\special{pa 1200 600}%
\special{pa 1450 600}%
\special{dt 0.045}%
\special{sh 1}%
\special{pa 1450 600}%
\special{pa 1383 580}%
\special{pa 1397 600}%
\special{pa 1383 620}%
\special{pa 1450 600}%
\special{fp}%
\put(18.0000,-8.0000){\makebox(0,0){$\xrightarrow[\text{positive}]{\text{direct},}$}}%
%
\special{pn 8}%
\special{ar 2600 800 400 400 0.0000000 6.2831853}%
%
\special{pn 8}%
\special{pa 3000 815}%
\special{pa 3000 800}%
\special{fp}%
\special{sh 1}%
\special{pa 3000 800}%
\special{pa 2980 867}%
\special{pa 3000 853}%
\special{pa 3020 867}%
\special{pa 3000 800}%
\special{fp}%
%
\special{pn 8}%
\special{ar 3000 800 200 200 0.0000000 6.2831853}%
%
\special{pn 8}%
\special{pa 2801 785}%
\special{pa 2800 800}%
\special{fp}%
\special{sh 1}%
\special{pa 2800 800}%
\special{pa 2824 735}%
\special{pa 2804 747}%
\special{pa 2784 732}%
\special{pa 2800 800}%
\special{fp}%
\put(36.0000,-8.0000){\makebox(0,0){$\xrightarrow[\text{negative}]{\text{inverse},}$}}%
%
\special{pn 8}%
\special{pa 3200 600}%
\special{pa 3450 600}%
\special{dt 0.045}%
\special{sh 1}%
\special{pa 3450 600}%
\special{pa 3383 580}%
\special{pa 3397 600}%
\special{pa 3383 620}%
\special{pa 3450 600}%
\special{fp}%
%
\special{pn 8}%
\special{ar 4400 800 400 400 0.0000000 6.2831853}%
%
\special{pn 8}%
\special{pa 4800 815}%
\special{pa 4800 800}%
\special{fp}%
\special{sh 1}%
\special{pa 4800 800}%
\special{pa 4780 867}%
\special{pa 4800 853}%
\special{pa 4820 867}%
\special{pa 4800 800}%
\special{fp}%
%
\special{pn 8}%
\special{ar 5200 800 200 200 0.0000000 6.2831853}%
%
\special{pn 8}%
\special{pa 5001 785}%
\special{pa 5000 800}%
\special{fp}%
\special{sh 1}%
\special{pa 5000 800}%
\special{pa 5024 735}%
\special{pa 5004 747}%
\special{pa 4984 732}%
\special{pa 5000 800}%
\special{fp}%
\put(11.5000,-9.5000){\makebox(0,0)[lt]{$C_0$}}%
\put(13.0000,-5.0000){\makebox(0,0)[lb]{$C_1$}}%
\put(29.0000,-5.0000){\makebox(0,0)[lb]{$C_1$}}%
\put(31.5000,-9.5000){\makebox(0,0)[lt]{$C_0$}}%
\put(53.5000,-9.5000){\makebox(0,0)[lt]{$C_0$}}%
\put(47.0000,-5.0000){\makebox(0,0)[lb]{$C_1$}}%
%
\special{pn 8}%
\special{ar 1000 2000 400 400 0.0000000 6.2831853}%
%
\special{pn 8}%
\special{pa 1400 2015}%
\special{pa 1400 2000}%
\special{fp}%
\special{sh 1}%
\special{pa 1400 2000}%
\special{pa 1380 2067}%
\special{pa 1400 2053}%
\special{pa 1420 2067}%
\special{pa 1400 2000}%
\special{fp}%
%
\special{pn 8}%
\special{ar 1000 2000 200 200 0.0000000 6.2831853}%
%
\special{pn 8}%
\special{pa 1199 1985}%
\special{pa 1200 2000}%
\special{fp}%
\special{sh 1}%
\special{pa 1200 2000}%
\special{pa 1216 1932}%
\special{pa 1196 1947}%
\special{pa 1176 1935}%
\special{pa 1200 2000}%
\special{fp}%
\put(13.0000,-17.0000){\makebox(0,0)[lb]{$C_1$}}%
\put(11.5000,-21.5000){\makebox(0,0)[lt]{$C_0$}}%
\put(18.0000,-20.0000){\makebox(0,0){$\xrightarrow[\text{positive}]{\text{inverse},}$}}%
%
\special{pn 8}%
\special{ar 2600 2000 400 400 0.0000000 6.2831853}%
%
\special{pn 8}%
\special{pa 3000 2015}%
\special{pa 3000 2000}%
\special{fp}%
\special{sh 1}%
\special{pa 3000 2000}%
\special{pa 2980 2067}%
\special{pa 3000 2053}%
\special{pa 3020 2067}%
\special{pa 3000 2000}%
\special{fp}%
%
\special{pn 8}%
\special{ar 3000 2000 200 200 0.0000000 6.2831853}%
%
\special{pn 8}%
\special{pa 2801 2015}%
\special{pa 2800 2000}%
\special{fp}%
\special{sh 1}%
\special{pa 2800 2000}%
\special{pa 2784 2068}%
\special{pa 2804 2053}%
\special{pa 2824 2065}%
\special{pa 2800 2000}%
\special{fp}%
%
\special{pn 8}%
\special{ar 4400 2000 400 400 0.0000000 6.2831853}%
%
\special{pn 8}%
\special{pa 4800 2015}%
\special{pa 4800 2000}%
\special{fp}%
\special{sh 1}%
\special{pa 4800 2000}%
\special{pa 4780 2067}%
\special{pa 4800 2053}%
\special{pa 4820 2067}%
\special{pa 4800 2000}%
\special{fp}%
\put(36.0000,-20.0000){\makebox(0,0){$\xrightarrow[\text{negative}]{\text{direct},}$}}%
%
\special{pn 8}%
\special{ar 5200 2000 200 200 0.0000000 6.2831853}%
%
\special{pn 8}%
\special{pa 5001 2015}%
\special{pa 5000 2000}%
\special{fp}%
\special{sh 1}%
\special{pa 5000 2000}%
\special{pa 4984 2068}%
\special{pa 5004 2053}%
\special{pa 5024 2065}%
\special{pa 5000 2000}%
\special{fp}%
\put(53.5000,-21.5000){\makebox(0,0)[lt]{$C_0$}}%
\put(47.0000,-17.0000){\makebox(0,0)[lb]{$C_1$}}%
\put(29.0000,-17.0000){\makebox(0,0)[lb]{$C_1$}}%
\put(31.5000,-21.5000){\makebox(0,0)[lt]{$C_0$}}%
%
\special{pn 8}%
\special{pa 1200 1800}%
\special{pa 1450 1800}%
\special{dt 0.045}%
\special{sh 1}%
\special{pa 1450 1800}%
\special{pa 1383 1780}%
\special{pa 1397 1800}%
\special{pa 1383 1820}%
\special{pa 1450 1800}%
\special{fp}%
%
\special{pn 8}%
\special{pa 3200 1800}%
\special{pa 3450 1800}%
\special{dt 0.045}%
\special{sh 1}%
\special{pa 3450 1800}%
\special{pa 3383 1780}%
\special{pa 3397 1800}%
\special{pa 3383 1820}%
\special{pa 3450 1800}%
\special{fp}%
\end{picture}}%

%% file: polygons.tex
{\unitlength 0.1in%
\begin{picture}(34.0000,12.7000)(6.0000,-18.7000)%
%
\special{pn 8}%
\special{ar 1000 800 447 447 0.4636476 2.6779450}%
%
\special{pn 8}%
\special{ar 1000 1400 447 447 3.6052403 5.8195377}%
%
\special{pn 4}%
\special{pa 1160 980}%
\special{pa 900 1240}%
\special{fp}%
\special{pa 1200 1000}%
\special{pa 960 1240}%
\special{fp}%
\special{pa 1240 1020}%
\special{pa 1010 1250}%
\special{fp}%
\special{pa 1270 1050}%
\special{pa 1080 1240}%
\special{fp}%
\special{pa 1300 1080}%
\special{pa 1170 1210}%
\special{fp}%
\special{pa 1110 970}%
\special{pa 860 1220}%
\special{fp}%
\special{pa 1060 960}%
\special{pa 820 1200}%
\special{fp}%
\special{pa 1010 950}%
\special{pa 780 1180}%
\special{fp}%
\special{pa 940 960}%
\special{pa 740 1160}%
\special{fp}%
\special{pa 870 970}%
\special{pa 710 1130}%
\special{fp}%
\put(10.0000,-15.0000){\makebox(0,0){a $2$-gon}}%
\put(18.0000,-18.0000){\makebox(0,0)[lb]{a hexagon in $D_0\cap D_1$ all but possibly one}}%
\special{pn 8}%
\special{pa 2400 600}%
\special{pa 2415 600}%
\special{pa 2420 601}%
\special{pa 2430 601}%
\special{pa 2435 602}%
\special{pa 2440 602}%
\special{pa 2460 606}%
\special{pa 2465 606}%
\special{pa 2470 608}%
\special{pa 2490 612}%
\special{pa 2495 614}%
\special{pa 2500 615}%
\special{pa 2505 617}%
\special{pa 2510 618}%
\special{pa 2530 626}%
\special{pa 2535 627}%
\special{pa 2540 629}%
\special{pa 2545 632}%
\special{pa 2560 638}%
\special{pa 2565 641}%
\special{pa 2575 645}%
\special{pa 2580 648}%
\special{pa 2585 650}%
\special{pa 2600 659}%
\special{pa 2605 661}%
\special{pa 2640 682}%
\special{pa 2645 686}%
\special{pa 2655 692}%
\special{pa 2660 696}%
\special{pa 2670 702}%
\special{pa 2675 706}%
\special{pa 2680 709}%
\special{pa 2685 713}%
\special{pa 2690 716}%
\special{pa 2695 720}%
\special{pa 2700 723}%
\special{pa 2710 731}%
\special{pa 2715 734}%
\special{pa 2735 750}%
\special{pa 2740 753}%
\special{pa 2765 773}%
\special{pa 2770 776}%
\special{pa 2830 824}%
\special{pa 2835 827}%
\special{pa 2860 847}%
\special{pa 2865 850}%
\special{pa 2885 866}%
\special{pa 2890 869}%
\special{pa 2900 877}%
\special{pa 2905 880}%
\special{pa 2910 884}%
\special{pa 2915 887}%
\special{pa 2920 891}%
\special{pa 2925 894}%
\special{pa 2930 898}%
\special{pa 2940 904}%
\special{pa 2945 908}%
\special{pa 2955 914}%
\special{pa 2960 918}%
\special{pa 2995 939}%
\special{pa 3000 941}%
\special{pa 3015 950}%
\special{pa 3020 952}%
\special{pa 3025 955}%
\special{pa 3035 959}%
\special{pa 3040 962}%
\special{pa 3055 968}%
\special{pa 3060 971}%
\special{pa 3065 973}%
\special{pa 3070 974}%
\special{pa 3090 982}%
\special{pa 3095 983}%
\special{pa 3100 985}%
\special{pa 3105 986}%
\special{pa 3110 988}%
\special{pa 3130 992}%
\special{pa 3135 994}%
\special{pa 3140 994}%
\special{pa 3160 998}%
\special{pa 3165 998}%
\special{pa 3170 999}%
\special{pa 3180 999}%
\special{pa 3185 1000}%
\special{pa 3200 1000}%
\special{fp}%
\special{pn 8}%
\special{pa 2400 1600}%
\special{pa 2415 1600}%
\special{pa 2420 1599}%
\special{pa 2430 1599}%
\special{pa 2435 1598}%
\special{pa 2440 1598}%
\special{pa 2460 1594}%
\special{pa 2465 1594}%
\special{pa 2470 1592}%
\special{pa 2490 1588}%
\special{pa 2495 1586}%
\special{pa 2500 1585}%
\special{pa 2505 1583}%
\special{pa 2510 1582}%
\special{pa 2530 1574}%
\special{pa 2535 1573}%
\special{pa 2540 1571}%
\special{pa 2545 1568}%
\special{pa 2560 1562}%
\special{pa 2565 1559}%
\special{pa 2575 1555}%
\special{pa 2580 1552}%
\special{pa 2585 1550}%
\special{pa 2600 1541}%
\special{pa 2605 1539}%
\special{pa 2640 1518}%
\special{pa 2645 1514}%
\special{pa 2655 1508}%
\special{pa 2660 1504}%
\special{pa 2670 1498}%
\special{pa 2675 1494}%
\special{pa 2680 1491}%
\special{pa 2685 1487}%
\special{pa 2690 1484}%
\special{pa 2695 1480}%
\special{pa 2700 1477}%
\special{pa 2710 1469}%
\special{pa 2715 1466}%
\special{pa 2735 1450}%
\special{pa 2740 1447}%
\special{pa 2765 1427}%
\special{pa 2770 1424}%
\special{pa 2830 1376}%
\special{pa 2835 1373}%
\special{pa 2860 1353}%
\special{pa 2865 1350}%
\special{pa 2885 1334}%
\special{pa 2890 1331}%
\special{pa 2900 1323}%
\special{pa 2905 1320}%
\special{pa 2910 1316}%
\special{pa 2915 1313}%
\special{pa 2920 1309}%
\special{pa 2925 1306}%
\special{pa 2930 1302}%
\special{pa 2940 1296}%
\special{pa 2945 1292}%
\special{pa 2955 1286}%
\special{pa 2960 1282}%
\special{pa 2995 1261}%
\special{pa 3000 1259}%
\special{pa 3015 1250}%
\special{pa 3020 1248}%
\special{pa 3025 1245}%
\special{pa 3035 1241}%
\special{pa 3040 1238}%
\special{pa 3055 1232}%
\special{pa 3060 1229}%
\special{pa 3065 1227}%
\special{pa 3070 1226}%
\special{pa 3090 1218}%
\special{pa 3095 1217}%
\special{pa 3100 1215}%
\special{pa 3105 1214}%
\special{pa 3110 1212}%
\special{pa 3130 1208}%
\special{pa 3135 1206}%
\special{pa 3140 1206}%
\special{pa 3160 1202}%
\special{pa 3165 1202}%
\special{pa 3170 1201}%
\special{pa 3180 1201}%
\special{pa 3185 1200}%
\special{pa 3200 1200}%
\special{fp}%
%
\special{pn 8}%
\special{ar 2400 1100 500 500 1.5707963 4.7123890}%
\special{pn 13}%
\special{pa 2400 1200}%
\special{pa 2415 1200}%
\special{pa 2420 1201}%
\special{pa 2430 1201}%
\special{pa 2435 1202}%
\special{pa 2440 1202}%
\special{pa 2460 1206}%
\special{pa 2465 1206}%
\special{pa 2470 1208}%
\special{pa 2490 1212}%
\special{pa 2495 1214}%
\special{pa 2500 1215}%
\special{pa 2505 1217}%
\special{pa 2510 1218}%
\special{pa 2530 1226}%
\special{pa 2535 1227}%
\special{pa 2540 1229}%
\special{pa 2545 1232}%
\special{pa 2560 1238}%
\special{pa 2565 1241}%
\special{pa 2575 1245}%
\special{pa 2580 1248}%
\special{pa 2585 1250}%
\special{pa 2600 1259}%
\special{pa 2605 1261}%
\special{pa 2640 1282}%
\special{pa 2645 1286}%
\special{pa 2655 1292}%
\special{pa 2660 1296}%
\special{pa 2670 1302}%
\special{pa 2675 1306}%
\special{pa 2680 1309}%
\special{pa 2685 1313}%
\special{pa 2690 1316}%
\special{pa 2695 1320}%
\special{pa 2700 1323}%
\special{pa 2710 1331}%
\special{pa 2715 1334}%
\special{pa 2735 1350}%
\special{pa 2740 1353}%
\special{pa 2765 1373}%
\special{pa 2770 1376}%
\special{pa 2830 1424}%
\special{pa 2835 1427}%
\special{pa 2860 1447}%
\special{pa 2865 1450}%
\special{pa 2885 1466}%
\special{pa 2890 1469}%
\special{pa 2900 1477}%
\special{pa 2905 1480}%
\special{pa 2910 1484}%
\special{pa 2915 1487}%
\special{pa 2920 1491}%
\special{pa 2925 1494}%
\special{pa 2930 1498}%
\special{pa 2940 1504}%
\special{pa 2945 1508}%
\special{pa 2955 1514}%
\special{pa 2960 1518}%
\special{pa 2995 1539}%
\special{pa 3000 1541}%
\special{pa 3015 1550}%
\special{pa 3020 1552}%
\special{pa 3025 1555}%
\special{pa 3035 1559}%
\special{pa 3040 1562}%
\special{pa 3055 1568}%
\special{pa 3060 1571}%
\special{pa 3065 1573}%
\special{pa 3070 1574}%
\special{pa 3090 1582}%
\special{pa 3095 1583}%
\special{pa 3100 1585}%
\special{pa 3105 1586}%
\special{pa 3110 1588}%
\special{pa 3130 1592}%
\special{pa 3135 1594}%
\special{pa 3140 1594}%
\special{pa 3160 1598}%
\special{pa 3165 1598}%
\special{pa 3170 1599}%
\special{pa 3180 1599}%
\special{pa 3185 1600}%
\special{pa 3200 1600}%
\special{fp}%
\special{pn 13}%
\special{pa 2400 1000}%
\special{pa 2415 1000}%
\special{pa 2420 999}%
\special{pa 2430 999}%
\special{pa 2435 998}%
\special{pa 2440 998}%
\special{pa 2460 994}%
\special{pa 2465 994}%
\special{pa 2470 992}%
\special{pa 2490 988}%
\special{pa 2495 986}%
\special{pa 2500 985}%
\special{pa 2505 983}%
\special{pa 2510 982}%
\special{pa 2530 974}%
\special{pa 2535 973}%
\special{pa 2540 971}%
\special{pa 2545 968}%
\special{pa 2560 962}%
\special{pa 2565 959}%
\special{pa 2575 955}%
\special{pa 2580 952}%
\special{pa 2585 950}%
\special{pa 2600 941}%
\special{pa 2605 939}%
\special{pa 2640 918}%
\special{pa 2645 914}%
\special{pa 2655 908}%
\special{pa 2660 904}%
\special{pa 2670 898}%
\special{pa 2675 894}%
\special{pa 2680 891}%
\special{pa 2685 887}%
\special{pa 2690 884}%
\special{pa 2695 880}%
\special{pa 2700 877}%
\special{pa 2710 869}%
\special{pa 2715 866}%
\special{pa 2735 850}%
\special{pa 2740 847}%
\special{pa 2765 827}%
\special{pa 2770 824}%
\special{pa 2830 776}%
\special{pa 2835 773}%
\special{pa 2860 753}%
\special{pa 2865 750}%
\special{pa 2885 734}%
\special{pa 2890 731}%
\special{pa 2900 723}%
\special{pa 2905 720}%
\special{pa 2910 716}%
\special{pa 2915 713}%
\special{pa 2920 709}%
\special{pa 2925 706}%
\special{pa 2930 702}%
\special{pa 2940 696}%
\special{pa 2945 692}%
\special{pa 2955 686}%
\special{pa 2960 682}%
\special{pa 2995 661}%
\special{pa 3000 659}%
\special{pa 3015 650}%
\special{pa 3020 648}%
\special{pa 3025 645}%
\special{pa 3035 641}%
\special{pa 3040 638}%
\special{pa 3055 632}%
\special{pa 3060 629}%
\special{pa 3065 627}%
\special{pa 3070 626}%
\special{pa 3090 618}%
\special{pa 3095 617}%
\special{pa 3100 615}%
\special{pa 3105 614}%
\special{pa 3110 612}%
\special{pa 3130 608}%
\special{pa 3135 606}%
\special{pa 3140 606}%
\special{pa 3160 602}%
\special{pa 3165 602}%
\special{pa 3170 601}%
\special{pa 3180 601}%
\special{pa 3185 600}%
\special{pa 3200 600}%
\special{fp}%
%
\special{pn 13}%
\special{ar 3200 1100 500 500 4.7123890 1.5707963}%
%
\special{pn 8}%
\special{pa 3200 1000}%
\special{pa 3800 1000}%
\special{fp}%
%
\special{pn 8}%
\special{pa 3200 1200}%
\special{pa 3800 1200}%
\special{fp}%
%
\special{pn 13}%
\special{pa 1800 1200}%
\special{pa 2400 1200}%
\special{fp}%
%
\special{pn 13}%
\special{pa 1800 1000}%
\special{pa 2400 1000}%
\special{fp}%
%
\special{pn 4}%
\special{pa 2940 900}%
\special{pa 2590 1250}%
\special{fp}%
\special{pa 2970 930}%
\special{pa 2630 1270}%
\special{fp}%
\special{pa 3010 950}%
\special{pa 2670 1290}%
\special{fp}%
\special{pa 3050 970}%
\special{pa 2700 1320}%
\special{fp}%
\special{pa 3100 980}%
\special{pa 2740 1340}%
\special{fp}%
\special{pa 3150 990}%
\special{pa 2770 1370}%
\special{fp}%
\special{pa 3200 1000}%
\special{pa 2840 1360}%
\special{fp}%
\special{pa 3260 1000}%
\special{pa 3020 1240}%
\special{fp}%
\special{pa 3320 1000}%
\special{pa 3110 1210}%
\special{fp}%
\special{pa 3380 1000}%
\special{pa 3180 1200}%
\special{fp}%
\special{pa 3440 1000}%
\special{pa 3240 1200}%
\special{fp}%
\special{pa 3500 1000}%
\special{pa 3300 1200}%
\special{fp}%
\special{pa 3560 1000}%
\special{pa 3360 1200}%
\special{fp}%
\special{pa 3620 1000}%
\special{pa 3420 1200}%
\special{fp}%
\special{pa 3680 1000}%
\special{pa 3480 1200}%
\special{fp}%
\special{pa 3690 1050}%
\special{pa 3540 1200}%
\special{fp}%
\special{pa 3700 1100}%
\special{pa 3600 1200}%
\special{fp}%
\special{pa 3690 1170}%
\special{pa 3660 1200}%
\special{fp}%
\special{pa 2900 880}%
\special{pa 2550 1230}%
\special{fp}%
\special{pa 2870 850}%
\special{pa 2510 1210}%
\special{fp}%
\special{pa 2830 830}%
\special{pa 2460 1200}%
\special{fp}%
\special{pa 2800 800}%
\special{pa 2400 1200}%
\special{fp}%
\special{pa 2600 940}%
\special{pa 2340 1200}%
\special{fp}%
\special{pa 2500 980}%
\special{pa 2280 1200}%
\special{fp}%
\special{pa 2430 990}%
\special{pa 2220 1200}%
\special{fp}%
\special{pa 2360 1000}%
\special{pa 2160 1200}%
\special{fp}%
\special{pa 2300 1000}%
\special{pa 2100 1200}%
\special{fp}%
\special{pa 2240 1000}%
\special{pa 2040 1200}%
\special{fp}%
\special{pa 2180 1000}%
\special{pa 1980 1200}%
\special{fp}%
\special{pa 2120 1000}%
\special{pa 1920 1200}%
\special{fp}%
\special{pa 2060 1000}%
\special{pa 1910 1150}%
\special{fp}%
\special{pa 2000 1000}%
\special{pa 1900 1100}%
\special{fp}%
\special{pa 1940 1000}%
\special{pa 1910 1030}%
\special{fp}%
\put(18.0000,-20.0000){\makebox(0,0)[lb]{adjacent regions in $D_0$ are $2$-gons}}%
\put(36.0000,-8.0000){\makebox(0,0)[lb]{$C_0$}}%
\put(20.0000,-8.0000){\makebox(0,0)[rb]{$C_1$}}%
\put(32.0000,-8.0000){\makebox(0,0){($2$-gon)}}%
\put(32.0000,-14.0000){\makebox(0,0){($2$-gon)}}%
%
\special{pn 13}%
\special{pn 13}%
\special{pa 1800 1200}%
\special{pa 1786 1200}%
\special{fp}%
\special{pa 1754 1197}%
\special{pa 1741 1196}%
\special{fp}%
\special{pa 1710 1189}%
\special{pa 1697 1186}%
\special{fp}%
\special{pa 1668 1175}%
\special{pa 1657 1170}%
\special{fp}%
\special{pa 1631 1154}%
\special{pa 1623 1146}%
\special{fp}%
\special{pa 1606 1125}%
\special{pa 1602 1114}%
\special{fp}%
\special{pa 1603 1084}%
\special{pa 1607 1073}%
\special{fp}%
\special{pa 1624 1053}%
\special{pa 1632 1046}%
\special{fp}%
\special{pa 1657 1030}%
\special{pa 1667 1025}%
\special{fp}%
\special{pa 1697 1014}%
\special{pa 1709 1011}%
\special{fp}%
\special{pa 1741 1004}%
\special{pa 1754 1003}%
\special{fp}%
\special{pa 1786 1000}%
\special{pa 1800 1000}%
\special{fp}%
%
\special{pn 8}%
\special{pn 8}%
\special{pa 3800 1000}%
\special{pa 3808 1000}%
\special{fp}%
\special{pa 3845 1003}%
\special{pa 3853 1004}%
\special{fp}%
\special{pa 3888 1010}%
\special{pa 3895 1012}%
\special{fp}%
\special{pa 3924 1021}%
\special{pa 3930 1024}%
\special{fp}%
\special{pa 3958 1039}%
\special{pa 3963 1042}%
\special{fp}%
\special{pa 3985 1062}%
\special{pa 3989 1067}%
\special{fp}%
\special{pa 4000 1095}%
\special{pa 4000 1104}%
\special{fp}%
\special{pa 3989 1133}%
\special{pa 3985 1138}%
\special{fp}%
\special{pa 3964 1158}%
\special{pa 3958 1162}%
\special{fp}%
\special{pa 3930 1176}%
\special{pa 3924 1179}%
\special{fp}%
\special{pa 3893 1188}%
\special{pa 3887 1190}%
\special{fp}%
\special{pa 3853 1196}%
\special{pa 3845 1197}%
\special{fp}%
\special{pa 3808 1200}%
\special{pa 3800 1200}%
\special{fp}%
\end{picture}}%

%% file: star.tex
{\unitlength 0.1in%
\begin{picture}(35.0000,12.6600)(4.0000,-18.8300)%
%
\special{pn 13}%
\special{pa 2400 715}%
\special{pa 3600 715}%
\special{fp}%
%
\special{pn 13}%
\special{pa 2200 1200}%
\special{pa 2800 1800}%
\special{fp}%
%
\special{pn 13}%
\special{pa 3200 1800}%
\special{pa 3800 1200}%
\special{fp}%
%
\special{pn 13}%
\special{ar 2400 1000 283 283 2.3561945 4.7123890}%
%
\special{pn 13}%
\special{ar 3000 1600 283 283 0.7853982 2.3561945}%
%
\special{pn 13}%
\special{ar 3600 1000 283 283 4.7123890 0.7853982}%
%
\special{pn 8}%
\special{ar 3000 900 283 283 3.9269908 5.4977871}%
%
\special{pn 8}%
\special{pa 2200 1300}%
\special{pa 2800 700}%
\special{fp}%
%
\special{pn 8}%
\special{ar 2400 1500 283 283 1.5707963 3.9269908}%
%
\special{pn 8}%
\special{pa 2400 1785}%
\special{pa 3600 1785}%
\special{fp}%
%
\special{pn 8}%
\special{ar 3600 1500 283 283 5.4977871 1.5707963}%
%
\special{pn 8}%
\special{pa 3200 700}%
\special{pa 3800 1300}%
\special{fp}%
%
\special{pn 4}%
\special{pa 3370 870}%
\special{pa 2630 1610}%
\special{fp}%
\special{pa 3400 900}%
\special{pa 2660 1640}%
\special{fp}%
\special{pa 3430 930}%
\special{pa 2690 1670}%
\special{fp}%
\special{pa 3460 960}%
\special{pa 2720 1700}%
\special{fp}%
\special{pa 3490 990}%
\special{pa 2750 1730}%
\special{fp}%
\special{pa 3520 1020}%
\special{pa 2780 1760}%
\special{fp}%
\special{pa 3550 1050}%
\special{pa 2820 1780}%
\special{fp}%
\special{pa 3580 1080}%
\special{pa 2880 1780}%
\special{fp}%
\special{pa 3610 1110}%
\special{pa 2940 1780}%
\special{fp}%
\special{pa 3640 1140}%
\special{pa 3000 1780}%
\special{fp}%
\special{pa 3670 1170}%
\special{pa 3060 1780}%
\special{fp}%
\special{pa 3700 1200}%
\special{pa 3120 1780}%
\special{fp}%
\special{pa 3730 1230}%
\special{pa 3180 1780}%
\special{fp}%
\special{pa 3340 840}%
\special{pa 2600 1580}%
\special{fp}%
\special{pa 3310 810}%
\special{pa 2570 1550}%
\special{fp}%
\special{pa 3280 780}%
\special{pa 2540 1520}%
\special{fp}%
\special{pa 3250 750}%
\special{pa 2510 1490}%
\special{fp}%
\special{pa 3220 720}%
\special{pa 2480 1460}%
\special{fp}%
\special{pa 3170 710}%
\special{pa 2450 1430}%
\special{fp}%
\special{pa 3110 710}%
\special{pa 2420 1400}%
\special{fp}%
\special{pa 3050 710}%
\special{pa 2390 1370}%
\special{fp}%
\special{pa 2990 710}%
\special{pa 2360 1340}%
\special{fp}%
\special{pa 2930 710}%
\special{pa 2330 1310}%
\special{fp}%
\special{pa 2870 710}%
\special{pa 2300 1280}%
\special{fp}%
\special{pa 2810 710}%
\special{pa 2270 1250}%
\special{fp}%
\put(39.0000,-10.0000){\makebox(0,0)[lb]{$C_0$}}%
\put(39.0000,-16.0000){\makebox(0,0)[lb]{$C_1$}}%
%
\special{pn 8}%
\special{ar 800 1200 400 400 0.0000000 6.2831853}%
%
\special{pn 8}%
\special{ar 1200 1200 400 400 0.0000000 6.2831853}%
%
\special{pn 4}%
\special{pa 1170 1050}%
\special{pa 840 1380}%
\special{fp}%
\special{pa 1180 1100}%
\special{pa 870 1410}%
\special{fp}%
\special{pa 1190 1150}%
\special{pa 890 1450}%
\special{fp}%
\special{pa 1200 1200}%
\special{pa 920 1480}%
\special{fp}%
\special{pa 1190 1270}%
\special{pa 950 1510}%
\special{fp}%
\special{pa 1170 1350}%
\special{pa 990 1530}%
\special{fp}%
\special{pa 1150 1010}%
\special{pa 830 1330}%
\special{fp}%
\special{pa 1130 970}%
\special{pa 810 1290}%
\special{fp}%
\special{pa 1100 940}%
\special{pa 800 1240}%
\special{fp}%
\special{pa 1070 910}%
\special{pa 800 1180}%
\special{fp}%
\special{pa 1040 880}%
\special{pa 810 1110}%
\special{fp}%
\special{pa 1000 860}%
\special{pa 860 1000}%
\special{fp}%
\put(8.0000,-17.0000){\makebox(0,0){$C_0$}}%
\put(12.0000,-17.0000){\makebox(0,0){$C_1$}}%
\end{picture}}%

%% file: no_2-gon.tex
{\unitlength 0.1in%
\begin{picture}(44.0000,8.5000)(14.0000,-32.3500)%
\put(36.0000,-28.0000){\makebox(0,0){$\rightsquigarrow$}}%
%
\special{pn 13}%
\special{pa 1400 2600}%
\special{pa 2200 2600}%
\special{fp}%
\put(18.0000,-28.0000){\makebox(0,0){$R_0$}}%
%
\special{pn 13}%
\special{pa 1400 3000}%
\special{pa 2200 3000}%
\special{fp}%
%
\special{pn 8}%
\special{pa 1600 3200}%
\special{pa 1600 2400}%
\special{fp}%
%
\special{pn 8}%
\special{pa 2000 3200}%
\special{pa 2000 2400}%
\special{fp}%
%
\special{pn 8}%
\special{pa 3000 3200}%
\special{pa 3000 2400}%
\special{fp}%
\put(24.0000,-26.0000){\makebox(0,0){$\cdots$}}%
\put(24.0000,-30.0000){\makebox(0,0){$\cdots$}}%
%
\special{pn 13}%
\special{ar 2600 2800 600 200 4.7123890 1.5707963}%
%
\special{pn 4}%
\special{pa 1990 2610}%
\special{pa 1610 2990}%
\special{fp}%
\special{pa 2000 2660}%
\special{pa 1660 3000}%
\special{fp}%
\special{pa 2000 2720}%
\special{pa 1720 3000}%
\special{fp}%
\special{pa 2000 2780}%
\special{pa 1780 3000}%
\special{fp}%
\special{pa 2000 2840}%
\special{pa 1840 3000}%
\special{fp}%
\special{pa 2000 2900}%
\special{pa 1900 3000}%
\special{fp}%
\special{pa 2000 2960}%
\special{pa 1960 3000}%
\special{fp}%
\special{pa 1940 2600}%
\special{pa 1600 2940}%
\special{fp}%
\special{pa 1880 2600}%
\special{pa 1600 2880}%
\special{fp}%
\special{pa 1820 2600}%
\special{pa 1600 2820}%
\special{fp}%
\special{pa 1760 2600}%
\special{pa 1600 2760}%
\special{fp}%
\special{pa 1700 2600}%
\special{pa 1600 2700}%
\special{fp}%
\special{pa 1640 2600}%
\special{pa 1600 2640}%
\special{fp}%
\put(31.0000,-27.0000){\makebox(0,0)[lb]{$C_0$}}%
\put(16.0000,-33.0000){\makebox(0,0){$C_1$}}%
\put(20.0000,-33.0000){\makebox(0,0){$C_1$}}%
\put(28.0000,-33.0000){\makebox(0,0){$C_1$}}%
\put(18.0000,-31.5000){\makebox(0,0){$D_1$}}%
%
\special{pn 8}%
\special{pa 2800 3200}%
\special{pa 2800 2400}%
\special{fp}%
\put(22.0000,-28.0000){\makebox(0,0){$R_1$}}%
\put(32.0000,-30.0000){\makebox(0,0)[lt]{$D_0$}}%
%
\special{pn 4}%
\special{pa 3200 3000}%
\special{pa 3100 2800}%
\special{fp}%
\put(30.0000,-33.0000){\makebox(0,0){$C_1$}}%
%
\special{pn 4}%
\special{pa 3000 2740}%
\special{pa 2800 2940}%
\special{fp}%
\special{pa 3000 2800}%
\special{pa 2820 2980}%
\special{fp}%
\special{pa 3000 2860}%
\special{pa 2890 2970}%
\special{fp}%
\special{pa 3000 2680}%
\special{pa 2800 2880}%
\special{fp}%
\special{pa 2970 2650}%
\special{pa 2800 2820}%
\special{fp}%
\special{pa 2920 2640}%
\special{pa 2800 2760}%
\special{fp}%
\special{pa 2870 2630}%
\special{pa 2800 2700}%
\special{fp}%
%
\special{pn 13}%
\special{pa 4000 2600}%
\special{pa 4400 2600}%
\special{fp}%
%
\special{pn 13}%
\special{pa 4000 3000}%
\special{pa 4400 3000}%
\special{fp}%
%
\special{pn 8}%
\special{pa 4200 3200}%
\special{pa 4200 2400}%
\special{fp}%
%
\special{pn 13}%
\special{ar 4400 2900 100 100 6.2831853 1.5707963}%
%
\special{pn 13}%
\special{ar 4400 2700 100 100 4.7123890 6.2831853}%
%
\special{pn 13}%
\special{pa 4500 2900}%
\special{pa 4500 2700}%
\special{fp}%
%
\special{pn 4}%
\special{pa 4490 2710}%
\special{pa 4210 2990}%
\special{fp}%
\special{pa 4480 2660}%
\special{pa 4200 2940}%
\special{fp}%
\special{pa 4450 2630}%
\special{pa 4200 2880}%
\special{fp}%
\special{pa 4410 2610}%
\special{pa 4200 2820}%
\special{fp}%
\special{pa 4360 2600}%
\special{pa 4200 2760}%
\special{fp}%
\special{pa 4300 2600}%
\special{pa 4200 2700}%
\special{fp}%
\special{pa 4240 2600}%
\special{pa 4200 2640}%
\special{fp}%
\special{pa 4490 2770}%
\special{pa 4270 2990}%
\special{fp}%
\special{pa 4490 2830}%
\special{pa 4330 2990}%
\special{fp}%
\special{pa 4490 2890}%
\special{pa 4390 2990}%
\special{fp}%
\put(43.5000,-28.0000){\makebox(0,0){($R_0$)}}%
%
\special{pn 13}%
\special{pa 4700 2900}%
\special{pa 4700 2700}%
\special{fp}%
%
\special{pn 13}%
\special{ar 4800 2700 100 100 3.1415927 4.7123890}%
%
\special{pn 13}%
\special{ar 4800 2900 100 100 1.5707963 3.1415927}%
\put(50.0000,-26.0000){\makebox(0,0){$\cdots$}}%
\put(50.0000,-30.0000){\makebox(0,0){$\cdots$}}%
%
\special{pn 13}%
\special{ar 5200 2800 600 200 4.7123890 1.5707963}%
%
\special{pn 8}%
\special{pa 5400 3200}%
\special{pa 5400 2400}%
\special{fp}%
%
\special{pn 8}%
\special{pa 5600 3200}%
\special{pa 5600 2400}%
\special{fp}%
%
\special{pn 4}%
\special{pa 5600 2740}%
\special{pa 5400 2940}%
\special{fp}%
\special{pa 5600 2800}%
\special{pa 5420 2980}%
\special{fp}%
\special{pa 5600 2860}%
\special{pa 5490 2970}%
\special{fp}%
\special{pa 5600 2680}%
\special{pa 5400 2880}%
\special{fp}%
\special{pa 5570 2650}%
\special{pa 5400 2820}%
\special{fp}%
\special{pa 5520 2640}%
\special{pa 5400 2760}%
\special{fp}%
\special{pa 5470 2630}%
\special{pa 5400 2700}%
\special{fp}%
\put(58.0000,-30.0000){\makebox(0,0)[lt]{$D'_0$}}%
%
\special{pn 4}%
\special{pa 5800 3000}%
\special{pa 5700 2800}%
\special{fp}%
\put(57.0000,-27.0000){\makebox(0,0)[lb]{$C'_0$}}%
\put(48.5000,-28.0000){\makebox(0,0){($R_1$)}}%
%
\special{pn 8}%
\special{pa 4600 3200}%
\special{pa 4600 2400}%
\special{fp}%
\put(42.0000,-33.0000){\makebox(0,0){$C_1$}}%
\put(46.0000,-33.0000){\makebox(0,0){$C_1$}}%
\put(54.0000,-33.0000){\makebox(0,0){$C_1$}}%
\put(56.0000,-33.0000){\makebox(0,0){$C_1$}}%
\put(44.0000,-31.5000){\makebox(0,0){$D_1$}}%
\put(55.0000,-24.5000){\makebox(0,0){$D_1$}}%
\put(29.0000,-24.5000){\makebox(0,0){$D_1$}}%
\end{picture}}%

%% file: claimed_region.tex
{\unitlength 0.1in%
\begin{picture}(29.5500,8.3500)(7.4500,-12.3500)%
%
\special{pn 13}%
\special{ar 800 800 1200 200 4.7123890 1.5707963}%
%
\special{pn 8}%
\special{pa 1000 1200}%
\special{pa 1000 400}%
\special{fp}%
%
\special{pn 8}%
\special{pa 1400 1200}%
\special{pa 1400 400}%
\special{fp}%
%
\special{pn 8}%
\special{pa 1600 1200}%
\special{pa 1600 400}%
\special{fp}%
\put(9.0000,-8.0000){\makebox(0,0){$R_2$}}%
\put(12.0000,-8.0000){\makebox(0,0){$R_3$}}%
%
\special{pn 8}%
\special{pa 800 1000}%
\special{pa 800 600}%
\special{ip}%
%
\special{pn 4}%
\special{pa 990 630}%
\special{pa 800 820}%
\special{fp}%
\special{pa 1000 680}%
\special{pa 800 880}%
\special{fp}%
\special{pa 1000 740}%
\special{pa 800 940}%
\special{fp}%
\special{pa 1000 800}%
\special{pa 820 980}%
\special{fp}%
\special{pa 1000 860}%
\special{pa 880 980}%
\special{fp}%
\special{pa 1000 920}%
\special{pa 940 980}%
\special{fp}%
\special{pa 940 620}%
\special{pa 800 760}%
\special{fp}%
\special{pa 880 620}%
\special{pa 800 700}%
\special{fp}%
%
\special{pn 4}%
\special{pa 1600 740}%
\special{pa 1400 940}%
\special{fp}%
\special{pa 1600 800}%
\special{pa 1450 950}%
\special{fp}%
\special{pa 1600 860}%
\special{pa 1510 950}%
\special{fp}%
\special{pa 1600 680}%
\special{pa 1400 880}%
\special{fp}%
\special{pa 1560 660}%
\special{pa 1400 820}%
\special{fp}%
\special{pa 1510 650}%
\special{pa 1400 760}%
\special{fp}%
\special{pa 1450 650}%
\special{pa 1400 700}%
\special{fp}%
\put(10.0000,-13.0000){\makebox(0,0){$C_1$}}%
\put(14.0000,-13.0000){\makebox(0,0){$C_1$}}%
\put(16.0000,-13.0000){\makebox(0,0){$C_1$}}%
\put(15.0000,-5.0000){\makebox(0,0){$D_1$}}%
\put(22.0000,-8.0000){\makebox(0,0){$\rightsquigarrow$}}%
%
\special{pn 13}%
\special{ar 2500 700 100 100 4.7123890 6.2831853}%
%
\special{pn 13}%
\special{ar 2500 900 100 100 6.2831853 1.5707963}%
%
\special{pn 13}%
\special{pa 2600 900}%
\special{pa 2600 700}%
\special{fp}%
%
\special{pn 4}%
\special{pa 2500 1000}%
\special{pa 2500 600}%
\special{ip}%
%
\special{pn 4}%
\special{pa 2590 710}%
\special{pa 2500 800}%
\special{fp}%
\special{pa 2590 770}%
\special{pa 2500 860}%
\special{fp}%
\special{pa 2590 830}%
\special{pa 2500 920}%
\special{fp}%
\special{pa 2590 890}%
\special{pa 2500 980}%
\special{fp}%
\special{pa 2580 660}%
\special{pa 2500 740}%
\special{fp}%
\special{pa 2550 630}%
\special{pa 2500 680}%
\special{fp}%
%
\special{pn 13}%
\special{ar 2900 700 100 100 3.1415927 4.7123890}%
%
\special{pn 13}%
\special{ar 2900 900 100 100 1.5707963 3.1415927}%
%
\special{pn 13}%
\special{pa 2800 900}%
\special{pa 2800 700}%
\special{fp}%
%
\special{pn 8}%
\special{pa 2700 1200}%
\special{pa 2700 400}%
\special{fp}%
%
\special{pn 8}%
\special{pa 3100 1200}%
\special{pa 3100 400}%
\special{fp}%
%
\special{pn 8}%
\special{pa 3300 1200}%
\special{pa 3300 400}%
\special{fp}%
%
\special{pn 13}%
\special{ar 2900 800 800 200 4.7123890 1.5707963}%
%
\special{pn 4}%
\special{pa 3300 720}%
\special{pa 3100 920}%
\special{fp}%
\special{pa 3300 780}%
\special{pa 3100 980}%
\special{fp}%
\special{pa 3300 840}%
\special{pa 3160 980}%
\special{fp}%
\special{pa 3300 900}%
\special{pa 3220 980}%
\special{fp}%
\special{pa 3300 660}%
\special{pa 3100 860}%
\special{fp}%
\special{pa 3270 630}%
\special{pa 3100 800}%
\special{fp}%
\special{pa 3220 620}%
\special{pa 3100 740}%
\special{fp}%
\special{pa 3160 620}%
\special{pa 3100 680}%
\special{fp}%
\put(19.0000,-7.0000){\makebox(0,0)[lb]{$C_0$}}%
\put(37.0000,-7.0000){\makebox(0,0)[lb]{$C''_0$}}%
\put(27.0000,-13.0000){\makebox(0,0){$C_1$}}%
\put(31.0000,-13.0000){\makebox(0,0){$C_1$}}%
\put(33.0000,-13.0000){\makebox(0,0){$C_1$}}%
\put(18.0000,-10.0000){\makebox(0,0)[lt]{$D_0$}}%
%
\special{pn 4}%
\special{pa 1800 1000}%
\special{pa 1700 800}%
\special{fp}%
\put(36.0000,-10.0000){\makebox(0,0)[lt]{$D''_0$}}%
%
\special{pn 4}%
\special{pa 3600 1000}%
\special{pa 3500 800}%
\special{fp}%
\put(32.0000,-5.0000){\makebox(0,0){$D_1$}}%
\put(32.0000,-8.0000){\makebox(0,0){$R'$}}%
\put(15.0000,-8.0000){\makebox(0,0){$R'$}}%
\end{picture}}%

%% file: star_resolution.tex
\unitlength 0.1in
\begin{picture}( 41.6600, 15.6500)(  5.1700,-20.0000)
%
{\color[named]{Black}{%
\special{pn 13}%
\special{pa 600 1200}%
\special{pa 1200 1800}%
\special{fp}%
}}%
%
{\color[named]{Black}{%
\special{pn 13}%
\special{pa 1600 1800}%
\special{pa 2200 1200}%
\special{fp}%
}}%
%
{\color[named]{Black}{%
\special{pn 13}%
\special{ar 800 1000 284 284  2.3561945  4.7123890}%
}}%
%
{\color[named]{Black}{%
\special{pn 13}%
\special{ar 2000 1000 284 284  4.7123890  0.7853982}%
}}%
%
{\color[named]{Black}{%
\special{pn 13}%
\special{pa 2016 718}%
\special{pa 2000 718}%
\special{fp}%
\special{sh 1}%
\special{pa 2000 718}%
\special{pa 2068 738}%
\special{pa 2054 718}%
\special{pa 2068 698}%
\special{pa 2000 718}%
\special{fp}%
}}%
%
{\color[named]{Black}{%
\special{pn 8}%
\special{pa 600 1300}%
\special{pa 1200 700}%
\special{fp}%
}}%
%
{\color[named]{Black}{%
\special{pn 8}%
\special{ar 800 1500 284 284  1.5707963  3.9269908}%
}}%
%
{\color[named]{Black}{%
\special{pn 8}%
\special{pa 800 1786}%
\special{pa 2000 1786}%
\special{fp}%
}}%
%
{\color[named]{Black}{%
\special{pn 8}%
\special{ar 2000 1500 284 284  5.4977871  1.5707963}%
}}%
%
{\color[named]{Black}{%
\special{pn 8}%
\special{pa 1600 700}%
\special{pa 2200 1300}%
\special{fp}%
}}%
%
{\color[named]{Black}{%
\special{pn 4}%
\special{pa 1770 870}%
\special{pa 1030 1610}%
\special{fp}%
\special{pa 1800 900}%
\special{pa 1060 1640}%
\special{fp}%
\special{pa 1830 930}%
\special{pa 1090 1670}%
\special{fp}%
\special{pa 1860 960}%
\special{pa 1120 1700}%
\special{fp}%
\special{pa 1890 990}%
\special{pa 1150 1730}%
\special{fp}%
\special{pa 1920 1020}%
\special{pa 1180 1760}%
\special{fp}%
\special{pa 1950 1050}%
\special{pa 1220 1780}%
\special{fp}%
\special{pa 1980 1080}%
\special{pa 1280 1780}%
\special{fp}%
\special{pa 2010 1110}%
\special{pa 1340 1780}%
\special{fp}%
\special{pa 2040 1140}%
\special{pa 1400 1780}%
\special{fp}%
\special{pa 2070 1170}%
\special{pa 1460 1780}%
\special{fp}%
\special{pa 2100 1200}%
\special{pa 1520 1780}%
\special{fp}%
\special{pa 2130 1230}%
\special{pa 1580 1780}%
\special{fp}%
\special{pa 1740 840}%
\special{pa 1000 1580}%
\special{fp}%
\special{pa 1710 810}%
\special{pa 970 1550}%
\special{fp}%
\special{pa 1680 780}%
\special{pa 940 1520}%
\special{fp}%
\special{pa 1650 750}%
\special{pa 910 1490}%
\special{fp}%
\special{pa 1620 720}%
\special{pa 880 1460}%
\special{fp}%
\special{pa 1570 710}%
\special{pa 850 1430}%
\special{fp}%
\special{pa 1510 710}%
\special{pa 820 1400}%
\special{fp}%
\special{pa 1450 710}%
\special{pa 790 1370}%
\special{fp}%
\special{pa 1390 710}%
\special{pa 760 1340}%
\special{fp}%
\special{pa 1330 710}%
\special{pa 730 1310}%
\special{fp}%
\special{pa 1270 710}%
\special{pa 700 1280}%
\special{fp}%
\special{pa 1210 710}%
\special{pa 670 1250}%
\special{fp}%
}}%
\put(23.0000,-10.0000){\makebox(0,0)[lb]{$C_0$}}%
\put(23.0000,-16.0000){\makebox(0,0)[lb]{$C_1$}}%
%
{\color[named]{Black}{%
\special{pn 13}%
\special{ar 1000 2000 284 284  5.4977871  6.2831853}%
}}%
%
{\color[named]{Black}{%
\special{pn 13}%
\special{ar 1800 2000 284 284  3.1415927  3.9269908}%
}}%
%
{\color[named]{Black}{%
\special{pn 8}%
\special{ar 1000 500 284 284  6.2831853  0.7853982}%
}}%
%
{\color[named]{Black}{%
\special{pn 8}%
\special{ar 1800 500 284 284  2.3561945  3.1415927}%
}}%
%
{\color[named]{Black}{%
\special{pn 13}%
\special{pa 2000 716}%
\special{pa 800 716}%
\special{fp}%
\special{sh 1}%
\special{pa 800 716}%
\special{pa 868 736}%
\special{pa 854 716}%
\special{pa 868 696}%
\special{pa 800 716}%
\special{fp}%
}}%
%
{\color[named]{Black}{%
\special{pn 8}%
\special{pa 600 1300}%
\special{pa 1000 900}%
\special{fp}%
\special{sh 1}%
\special{pa 1000 900}%
\special{pa 940 934}%
\special{pa 962 938}%
\special{pa 968 962}%
\special{pa 1000 900}%
\special{fp}%
}}%
%
{\color[named]{Black}{%
\special{pn 8}%
\special{pa 1600 700}%
\special{pa 1900 1000}%
\special{fp}%
\special{sh 1}%
\special{pa 1900 1000}%
\special{pa 1868 940}%
\special{pa 1862 962}%
\special{pa 1840 968}%
\special{pa 1900 1000}%
\special{fp}%
}}%
%
{\color[named]{Black}{%
\special{pn 8}%
\special{pa 650 850}%
\special{pa 950 1150}%
\special{dt 0.045}%
\special{sh 1}%
\special{pa 950 1150}%
\special{pa 918 1090}%
\special{pa 912 1112}%
\special{pa 890 1118}%
\special{pa 950 1150}%
\special{fp}%
}}%
%
{\color[named]{Black}{%
\special{pn 8}%
\special{pa 2150 850}%
\special{pa 1850 1150}%
\special{dt 0.045}%
\special{sh 1}%
\special{pa 1850 1150}%
\special{pa 1912 1118}%
\special{pa 1888 1112}%
\special{pa 1884 1090}%
\special{pa 1850 1150}%
\special{fp}%
}}%
\put(26.0000,-12.0000){\makebox(0,0){$\xrightarrow[\text{negative}]{\text{inverse},}$}}%
%
{\color[named]{Black}{%
\special{pn 8}%
\special{ar 3200 1500 284 284  1.5707963  3.9269908}%
}}%
%
{\color[named]{Black}{%
\special{pn 8}%
\special{pa 3000 1300}%
\special{pa 3600 700}%
\special{fp}%
}}%
%
{\color[named]{Black}{%
\special{pn 8}%
\special{pa 3000 1300}%
\special{pa 3400 900}%
\special{fp}%
\special{sh 1}%
\special{pa 3400 900}%
\special{pa 3340 934}%
\special{pa 3362 938}%
\special{pa 3368 962}%
\special{pa 3400 900}%
\special{fp}%
}}%
%
{\color[named]{Black}{%
\special{pn 8}%
\special{ar 3400 500 284 284  6.2831853  0.7853982}%
}}%
%
{\color[named]{Black}{%
\special{pn 8}%
\special{ar 4200 500 284 284  2.3561945  3.1415927}%
}}%
%
{\color[named]{Black}{%
\special{pn 8}%
\special{pa 4000 700}%
\special{pa 4300 1000}%
\special{fp}%
\special{sh 1}%
\special{pa 4300 1000}%
\special{pa 4268 940}%
\special{pa 4262 962}%
\special{pa 4240 968}%
\special{pa 4300 1000}%
\special{fp}%
}}%
%
{\color[named]{Black}{%
\special{pn 8}%
\special{pa 4000 700}%
\special{pa 4600 1300}%
\special{fp}%
}}%
%
{\color[named]{Black}{%
\special{pn 8}%
\special{ar 4400 1500 284 284  5.4977871  1.5707963}%
}}%
%
{\color[named]{Black}{%
\special{pn 8}%
\special{pa 3200 1786}%
\special{pa 4400 1786}%
\special{fp}%
}}%
%
{\color[named]{Black}{%
\special{pn 13}%
\special{ar 3800 2000 200 1000  3.1415927  6.2831853}%
}}%
%
{\color[named]{Black}{%
\special{pn 13}%
\special{pa 3600 1976}%
\special{pa 3600 2000}%
\special{fp}%
\special{sh 1}%
\special{pa 3600 2000}%
\special{pa 3620 1934}%
\special{pa 3600 1948}%
\special{pa 3580 1934}%
\special{pa 3600 2000}%
\special{fp}%
}}%
%
{\color[named]{Black}{%
\special{pn 4}%
\special{pa 3960 1440}%
\special{pa 3620 1780}%
\special{fp}%
\special{pa 3950 1390}%
\special{pa 3610 1730}%
\special{fp}%
\special{pa 3940 1340}%
\special{pa 3620 1660}%
\special{fp}%
\special{pa 3930 1290}%
\special{pa 3620 1600}%
\special{fp}%
\special{pa 3920 1240}%
\special{pa 3630 1530}%
\special{fp}%
\special{pa 3910 1190}%
\special{pa 3640 1460}%
\special{fp}%
\special{pa 3890 1150}%
\special{pa 3650 1390}%
\special{fp}%
\special{pa 3880 1100}%
\special{pa 3660 1320}%
\special{fp}%
\special{pa 3850 1070}%
\special{pa 3680 1240}%
\special{fp}%
\special{pa 3830 1030}%
\special{pa 3700 1160}%
\special{fp}%
\special{pa 3970 1490}%
\special{pa 3680 1780}%
\special{fp}%
\special{pa 3970 1550}%
\special{pa 3740 1780}%
\special{fp}%
\special{pa 3980 1600}%
\special{pa 3800 1780}%
\special{fp}%
\special{pa 3980 1660}%
\special{pa 3860 1780}%
\special{fp}%
\special{pa 3990 1710}%
\special{pa 3920 1780}%
\special{fp}%
}}%
%
{\color[named]{Black}{%
\special{pn 8}%
\special{pa 4400 1786}%
\special{pa 3500 1786}%
\special{fp}%
\special{sh 1}%
\special{pa 3500 1786}%
\special{pa 3568 1806}%
\special{pa 3554 1786}%
\special{pa 3568 1766}%
\special{pa 3500 1786}%
\special{fp}%
}}%
%
{\color[named]{Black}{%
\special{pn 8}%
\special{pa 3800 1000}%
\special{pa 3800 1900}%
\special{dt 0.045}%
\special{sh 1}%
\special{pa 3800 1900}%
\special{pa 3820 1834}%
\special{pa 3800 1848}%
\special{pa 3780 1834}%
\special{pa 3800 1900}%
\special{fp}%
}}%
\put(41.0000,-19.0000){\makebox(0,0)[lt]{direct, negative}}%
%
{\color[named]{Black}{%
\special{pn 4}%
\special{pa 4100 1900}%
\special{pa 3800 1500}%
\special{fp}%
}}%
\put(14.0000,-12.0000){\makebox(0,0){$P$}}%
\put(14.0000,-19.0000){\makebox(0,0){$D_0$}}%
\put(14.0000,-5.0000){\makebox(0,0){$D_1$}}%
\put(38.0000,-5.0000){\makebox(0,0){$D_1$}}%
\end{picture}%

%% file: all_the_newborn_triangles.tex
{\unitlength 0.1in%
\begin{picture}(35.4500,49.2000)(2.5500,-53.2000)%
%
\special{pn 13}%
\special{pa 1200 600}%
\special{pa 600 1200}%
\special{fp}%
\special{sh 1}%
\special{pa 600 1200}%
\special{pa 661 1167}%
\special{pa 638 1162}%
\special{pa 633 1139}%
\special{pa 600 1200}%
\special{fp}%
%
\special{pn 13}%
\special{pa 1200 1200}%
\special{pa 600 600}%
\special{fp}%
\special{sh 1}%
\special{pa 600 600}%
\special{pa 633 661}%
\special{pa 638 638}%
\special{pa 661 633}%
\special{pa 600 600}%
\special{fp}%
%
\special{pn 8}%
\special{pa 700 1400}%
\special{pa 700 400}%
\special{fp}%
\special{sh 1}%
\special{pa 700 400}%
\special{pa 680 467}%
\special{pa 700 453}%
\special{pa 720 467}%
\special{pa 700 400}%
\special{fp}%
\put(7.2000,-14.2000){\makebox(0,0)[lt]{$c_1$}}%
%
\special{pn 8}%
\special{pn 8}%
\special{pa 700 400}%
\special{pa 708 400}%
\special{fp}%
\special{pa 747 401}%
\special{pa 756 401}%
\special{fp}%
\special{pa 795 403}%
\special{pa 803 405}%
\special{fp}%
\special{pa 842 408}%
\special{pa 850 409}%
\special{fp}%
\special{pa 889 415}%
\special{pa 897 416}%
\special{fp}%
\special{pa 935 423}%
\special{pa 943 425}%
\special{fp}%
\special{pa 981 435}%
\special{pa 989 437}%
\special{fp}%
\special{pa 1026 448}%
\special{pa 1034 451}%
\special{fp}%
\special{pa 1070 466}%
\special{pa 1078 469}%
\special{fp}%
\special{pa 1113 487}%
\special{pa 1120 491}%
\special{fp}%
\special{pa 1151 514}%
\special{pa 1157 519}%
\special{fp}%
\special{pa 1183 548}%
\special{pa 1187 555}%
\special{fp}%
\special{pa 1200 592}%
\special{pa 1200 600}%
\special{fp}%
%
\special{pn 8}%
\special{pn 8}%
\special{pa 1200 1200}%
\special{pa 1199 1208}%
\special{fp}%
\special{pa 1181 1235}%
\special{pa 1175 1240}%
\special{fp}%
\special{pa 1144 1258}%
\special{pa 1137 1261}%
\special{fp}%
\special{pa 1103 1273}%
\special{pa 1096 1276}%
\special{fp}%
\special{pa 1060 1285}%
\special{pa 1052 1286}%
\special{fp}%
\special{pa 1016 1292}%
\special{pa 1008 1293}%
\special{fp}%
\special{pa 972 1297}%
\special{pa 964 1298}%
\special{fp}%
\special{pa 927 1300}%
\special{pa 919 1300}%
\special{fp}%
\special{pa 882 1300}%
\special{pa 873 1300}%
\special{fp}%
\special{pa 837 1298}%
\special{pa 829 1297}%
\special{fp}%
\special{pa 792 1293}%
\special{pa 784 1292}%
\special{fp}%
\special{pa 748 1286}%
\special{pa 740 1285}%
\special{fp}%
\special{pa 705 1276}%
\special{pa 697 1274}%
\special{fp}%
\special{pa 663 1261}%
\special{pa 656 1259}%
\special{fp}%
\special{pa 625 1240}%
\special{pa 619 1235}%
\special{fp}%
\special{pa 601 1208}%
\special{pa 600 1200}%
\special{fp}%
\put(12.0000,-6.0000){\makebox(0,0){$\bullet$}}%
\put(12.2000,-5.8000){\makebox(0,0)[lb]{$b_0$}}%
\put(13.0000,-9.0000){\makebox(0,0){$c_0$}}%
%
\special{pn 4}%
\special{pa 1230 830}%
\special{pa 1100 700}%
\special{fp}%
%
\special{pn 8}%
\special{pn 8}%
\special{pa 600 1500}%
\special{pa 592 1500}%
\special{fp}%
\special{pa 557 1489}%
\special{pa 550 1485}%
\special{fp}%
\special{pa 520 1463}%
\special{pa 514 1457}%
\special{fp}%
\special{pa 491 1428}%
\special{pa 487 1421}%
\special{fp}%
\special{pa 468 1388}%
\special{pa 465 1381}%
\special{fp}%
\special{pa 450 1347}%
\special{pa 447 1339}%
\special{fp}%
\special{pa 435 1303}%
\special{pa 433 1296}%
\special{fp}%
\special{pa 423 1259}%
\special{pa 421 1251}%
\special{fp}%
\special{pa 414 1214}%
\special{pa 412 1206}%
\special{fp}%
\special{pa 407 1169}%
\special{pa 406 1161}%
\special{fp}%
\special{pa 403 1123}%
\special{pa 402 1115}%
\special{fp}%
\special{pa 400 1077}%
\special{pa 400 1069}%
\special{fp}%
\special{pa 400 1031}%
\special{pa 400 1023}%
\special{fp}%
\special{pa 402 985}%
\special{pa 403 977}%
\special{fp}%
\special{pa 406 940}%
\special{pa 407 932}%
\special{fp}%
\special{pa 412 894}%
\special{pa 414 886}%
\special{fp}%
\special{pa 421 849}%
\special{pa 422 842}%
\special{fp}%
\special{pa 432 805}%
\special{pa 435 797}%
\special{fp}%
\special{pa 447 761}%
\special{pa 449 754}%
\special{fp}%
\special{pa 465 719}%
\special{pa 468 712}%
\special{fp}%
\special{pa 487 679}%
\special{pa 491 673}%
\special{fp}%
\special{pa 514 643}%
\special{pa 520 638}%
\special{fp}%
\special{pa 549 615}%
\special{pa 557 611}%
\special{fp}%
\special{pa 592 600}%
\special{pa 600 600}%
\special{fp}%
%
\special{pn 8}%
\special{pn 8}%
\special{pa 700 1400}%
\special{pa 700 1410}%
\special{fp}%
\special{pa 686 1450}%
\special{pa 682 1457}%
\special{fp}%
\special{pa 658 1482}%
\special{pa 651 1486}%
\special{fp}%
\special{pa 610 1500}%
\special{pa 600 1500}%
\special{fp}%
%
\special{pn 4}%
\special{pa 1230 970}%
\special{pa 1100 1100}%
\special{fp}%
\put(8.5000,-10.5000){\makebox(0,0){$\sb 1$}}%
\put(8.5000,-7.5000){\makebox(0,0){$\sb 2$}}%
\put(6.5000,-9.0000){\makebox(0,0){$\sb 3$}}%
\put(15.0000,-6.0000){\makebox(0,0)[lb]{$q=1$,}}%
%
\special{pn 8}%
\special{pa 3000 1400}%
\special{pa 3000 400}%
\special{fp}%
\special{sh 1}%
\special{pa 3000 400}%
\special{pa 2980 467}%
\special{pa 3000 453}%
\special{pa 3020 467}%
\special{pa 3000 400}%
\special{fp}%
%
\special{pn 13}%
\special{pa 3500 1200}%
\special{pa 2900 600}%
\special{fp}%
\special{sh 1}%
\special{pa 2900 600}%
\special{pa 2933 661}%
\special{pa 2938 638}%
\special{pa 2961 633}%
\special{pa 2900 600}%
\special{fp}%
%
\special{pn 13}%
\special{pa 3500 600}%
\special{pa 2900 1200}%
\special{fp}%
\special{sh 1}%
\special{pa 2900 1200}%
\special{pa 2961 1167}%
\special{pa 2938 1162}%
\special{pa 2933 1139}%
\special{pa 2900 1200}%
\special{fp}%
%
\special{pn 8}%
\special{pn 8}%
\special{pa 3000 400}%
\special{pa 3008 400}%
\special{fp}%
\special{pa 3045 401}%
\special{pa 3054 401}%
\special{fp}%
\special{pa 3090 404}%
\special{pa 3098 405}%
\special{fp}%
\special{pa 3135 409}%
\special{pa 3143 411}%
\special{fp}%
\special{pa 3179 417}%
\special{pa 3187 418}%
\special{fp}%
\special{pa 3223 426}%
\special{pa 3230 428}%
\special{fp}%
\special{pa 3266 438}%
\special{pa 3273 439}%
\special{fp}%
\special{pa 3308 451}%
\special{pa 3316 453}%
\special{fp}%
\special{pa 3350 467}%
\special{pa 3358 470}%
\special{fp}%
\special{pa 3391 485}%
\special{pa 3399 489}%
\special{fp}%
\special{pa 3432 506}%
\special{pa 3439 510}%
\special{fp}%
\special{pa 3470 530}%
\special{pa 3477 534}%
\special{fp}%
\special{pa 3507 555}%
\special{pa 3513 560}%
\special{fp}%
\special{pa 3541 583}%
\special{pa 3547 588}%
\special{fp}%
\special{pa 3574 614}%
\special{pa 3580 620}%
\special{fp}%
\special{pa 3604 647}%
\special{pa 3609 654}%
\special{fp}%
\special{pa 3631 684}%
\special{pa 3635 690}%
\special{fp}%
\special{pa 3654 722}%
\special{pa 3658 729}%
\special{fp}%
\special{pa 3673 762}%
\special{pa 3676 770}%
\special{fp}%
\special{pa 3687 804}%
\special{pa 3689 811}%
\special{fp}%
\special{pa 3696 847}%
\special{pa 3697 855}%
\special{fp}%
\special{pa 3700 892}%
\special{pa 3700 900}%
\special{fp}%
%
\special{pn 8}%
\special{pn 8}%
\special{pa 3700 900}%
\special{pa 3700 909}%
\special{fp}%
\special{pa 3698 947}%
\special{pa 3696 956}%
\special{fp}%
\special{pa 3690 993}%
\special{pa 3688 1001}%
\special{fp}%
\special{pa 3678 1036}%
\special{pa 3676 1044}%
\special{fp}%
\special{pa 3661 1077}%
\special{pa 3658 1083}%
\special{fp}%
\special{pa 3640 1114}%
\special{pa 3636 1120}%
\special{fp}%
\special{pa 3615 1145}%
\special{pa 3611 1149}%
\special{fp}%
\special{pa 3587 1170}%
\special{pa 3582 1174}%
\special{fp}%
\special{pa 3551 1190}%
\special{pa 3544 1193}%
\special{fp}%
\special{pa 3509 1200}%
\special{pa 3500 1200}%
\special{fp}%
%
\special{pn 8}%
\special{pn 8}%
\special{pa 2900 600}%
\special{pa 2901 593}%
\special{fp}%
\special{pa 2916 568}%
\special{pa 2921 564}%
\special{fp}%
\special{pa 2945 547}%
\special{pa 2952 544}%
\special{fp}%
\special{pa 2981 532}%
\special{pa 2988 529}%
\special{fp}%
\special{pa 3020 520}%
\special{pa 3027 518}%
\special{fp}%
\special{pa 3062 511}%
\special{pa 3069 510}%
\special{fp}%
\special{pa 3105 505}%
\special{pa 3113 504}%
\special{fp}%
\special{pa 3149 501}%
\special{pa 3157 501}%
\special{fp}%
\special{pa 3193 500}%
\special{pa 3201 500}%
\special{fp}%
\special{pa 3238 501}%
\special{pa 3246 501}%
\special{fp}%
\special{pa 3282 504}%
\special{pa 3290 505}%
\special{fp}%
\special{pa 3325 509}%
\special{pa 3333 510}%
\special{fp}%
\special{pa 3368 517}%
\special{pa 3375 519}%
\special{fp}%
\special{pa 3409 528}%
\special{pa 3416 531}%
\special{fp}%
\special{pa 3446 543}%
\special{pa 3453 546}%
\special{fp}%
\special{pa 3479 564}%
\special{pa 3485 568}%
\special{fp}%
\special{pa 3499 593}%
\special{pa 3500 600}%
\special{fp}%
%
\special{pn 8}%
\special{pn 8}%
\special{pa 3000 1400}%
\special{pa 2993 1399}%
\special{fp}%
\special{pa 2965 1388}%
\special{pa 2960 1384}%
\special{fp}%
\special{pa 2942 1362}%
\special{pa 2938 1357}%
\special{fp}%
\special{pa 2924 1330}%
\special{pa 2921 1324}%
\special{fp}%
\special{pa 2911 1293}%
\special{pa 2910 1285}%
\special{fp}%
\special{pa 2903 1252}%
\special{pa 2903 1244}%
\special{fp}%
\special{pa 2900 1208}%
\special{pa 2900 1200}%
\special{fp}%
\put(31.5000,-10.5000){\makebox(0,0){$\sb 2$}}%
\put(31.5000,-7.5000){\makebox(0,0){$\sb 1$}}%
\put(29.5000,-9.0000){\makebox(0,0){$\sb 3$}}%
\put(15.0000,-8.0000){\makebox(0,0)[lb]{$e^{c_0}_{b_0}(p)=-1$}}%
\put(38.0000,-6.0000){\makebox(0,0)[lb]{$q=2$,}}%
\put(38.0000,-8.0000){\makebox(0,0)[lb]{$e^{c_0}_{b_0}(p)=1$}}%
\put(35.0000,-12.0000){\makebox(0,0){$\bullet$}}%
\put(35.2000,-12.2000){\makebox(0,0)[lt]{$b_0$}}%
\put(30.2000,-14.2000){\makebox(0,0)[lt]{$c_1$}}%
\put(36.0000,-9.0000){\makebox(0,0){$c_0$}}%
%
\special{pn 4}%
\special{pa 3530 830}%
\special{pa 3400 700}%
\special{fp}%
%
\special{pn 4}%
\special{pa 3530 970}%
\special{pa 3400 1100}%
\special{fp}%
\put(7.7500,-9.0000){\makebox(0,0){$\circlearrowleft$}}%
\put(30.7500,-9.0000){\makebox(0,0){$\circlearrowright$}}%
\put(10.0000,-9.0000){\makebox(0,0){$p$}}%
\put(33.0000,-9.0000){\makebox(0,0){$p$}}%
%
\special{pn 8}%
\special{pa 700 2700}%
\special{pa 700 1700}%
\special{fp}%
\special{sh 1}%
\special{pa 700 1700}%
\special{pa 680 1767}%
\special{pa 700 1753}%
\special{pa 720 1767}%
\special{pa 700 1700}%
\special{fp}%
%
\special{pn 13}%
\special{pa 1200 2500}%
\special{pa 600 1900}%
\special{fp}%
\special{sh 1}%
\special{pa 600 1900}%
\special{pa 633 1961}%
\special{pa 638 1938}%
\special{pa 661 1933}%
\special{pa 600 1900}%
\special{fp}%
%
\special{pn 13}%
\special{pa 600 2500}%
\special{pa 1200 1900}%
\special{fp}%
\special{sh 1}%
\special{pa 1200 1900}%
\special{pa 1139 1933}%
\special{pa 1162 1938}%
\special{pa 1167 1961}%
\special{pa 1200 1900}%
\special{fp}%
\put(12.0000,-25.0000){\makebox(0,0){$\bullet$}}%
\put(12.2000,-25.2000){\makebox(0,0)[lt]{$b_0$}}%
\put(8.5000,-20.5000){\makebox(0,0){$\sb 1$}}%
\put(8.5000,-23.5000){\makebox(0,0){$\sb 2$}}%
%
\special{pn 8}%
\special{pn 8}%
\special{pa 600 2500}%
\special{pa 592 2499}%
\special{fp}%
\special{pa 566 2482}%
\special{pa 561 2477}%
\special{fp}%
\special{pa 543 2447}%
\special{pa 540 2440}%
\special{fp}%
\special{pa 527 2406}%
\special{pa 525 2398}%
\special{fp}%
\special{pa 516 2362}%
\special{pa 514 2355}%
\special{fp}%
\special{pa 508 2319}%
\special{pa 507 2311}%
\special{fp}%
\special{pa 503 2274}%
\special{pa 502 2266}%
\special{fp}%
\special{pa 500 2228}%
\special{pa 500 2219}%
\special{fp}%
\special{pa 500 2181}%
\special{pa 500 2173}%
\special{fp}%
\special{pa 502 2135}%
\special{pa 503 2127}%
\special{fp}%
\special{pa 507 2089}%
\special{pa 508 2082}%
\special{fp}%
\special{pa 514 2045}%
\special{pa 516 2037}%
\special{fp}%
\special{pa 525 2002}%
\special{pa 527 1994}%
\special{fp}%
\special{pa 540 1960}%
\special{pa 543 1953}%
\special{fp}%
\special{pa 561 1923}%
\special{pa 566 1918}%
\special{fp}%
\special{pa 593 1901}%
\special{pa 600 1900}%
\special{fp}%
%
\special{pn 8}%
\special{pn 8}%
\special{pa 1200 1900}%
\special{pa 1200 1908}%
\special{fp}%
\special{pa 1199 1945}%
\special{pa 1199 1953}%
\special{fp}%
\special{pa 1197 1990}%
\special{pa 1196 1997}%
\special{fp}%
\special{pa 1193 2034}%
\special{pa 1192 2042}%
\special{fp}%
\special{pa 1187 2078}%
\special{pa 1186 2086}%
\special{fp}%
\special{pa 1180 2122}%
\special{pa 1178 2130}%
\special{fp}%
\special{pa 1171 2166}%
\special{pa 1169 2174}%
\special{fp}%
\special{pa 1161 2210}%
\special{pa 1158 2217}%
\special{fp}%
\special{pa 1149 2253}%
\special{pa 1146 2260}%
\special{fp}%
\special{pa 1135 2296}%
\special{pa 1132 2303}%
\special{fp}%
\special{pa 1119 2337}%
\special{pa 1116 2345}%
\special{fp}%
\special{pa 1101 2378}%
\special{pa 1097 2385}%
\special{fp}%
\special{pa 1081 2418}%
\special{pa 1077 2425}%
\special{fp}%
\special{pa 1059 2457}%
\special{pa 1055 2464}%
\special{fp}%
\special{pa 1035 2494}%
\special{pa 1030 2501}%
\special{fp}%
\special{pa 1008 2530}%
\special{pa 1003 2536}%
\special{fp}%
\special{pa 979 2564}%
\special{pa 973 2570}%
\special{fp}%
\special{pa 947 2596}%
\special{pa 941 2601}%
\special{fp}%
\special{pa 913 2624}%
\special{pa 906 2629}%
\special{fp}%
\special{pa 876 2649}%
\special{pa 869 2653}%
\special{fp}%
\special{pa 837 2670}%
\special{pa 830 2673}%
\special{fp}%
\special{pa 795 2685}%
\special{pa 788 2688}%
\special{fp}%
\special{pa 752 2696}%
\special{pa 744 2697}%
\special{fp}%
\special{pa 708 2700}%
\special{pa 700 2700}%
\special{fp}%
%
\special{pn 8}%
\special{pn 8}%
\special{pa 700 1700}%
\special{pa 708 1700}%
\special{fp}%
\special{pa 745 1701}%
\special{pa 753 1701}%
\special{fp}%
\special{pa 790 1704}%
\special{pa 798 1705}%
\special{fp}%
\special{pa 834 1709}%
\special{pa 842 1710}%
\special{fp}%
\special{pa 879 1717}%
\special{pa 886 1718}%
\special{fp}%
\special{pa 922 1726}%
\special{pa 930 1728}%
\special{fp}%
\special{pa 966 1738}%
\special{pa 974 1739}%
\special{fp}%
\special{pa 1009 1752}%
\special{pa 1016 1754}%
\special{fp}%
\special{pa 1051 1768}%
\special{pa 1058 1770}%
\special{fp}%
\special{pa 1092 1786}%
\special{pa 1099 1789}%
\special{fp}%
\special{pa 1132 1806}%
\special{pa 1139 1810}%
\special{fp}%
\special{pa 1170 1830}%
\special{pa 1176 1834}%
\special{fp}%
\special{pa 1206 1855}%
\special{pa 1213 1860}%
\special{fp}%
\special{pa 1241 1883}%
\special{pa 1247 1888}%
\special{fp}%
\special{pa 1274 1914}%
\special{pa 1280 1920}%
\special{fp}%
\special{pa 1304 1947}%
\special{pa 1309 1953}%
\special{fp}%
\special{pa 1331 1983}%
\special{pa 1335 1990}%
\special{fp}%
\special{pa 1354 2021}%
\special{pa 1358 2028}%
\special{fp}%
\special{pa 1373 2062}%
\special{pa 1376 2070}%
\special{fp}%
\special{pa 1387 2104}%
\special{pa 1389 2112}%
\special{fp}%
\special{pa 1396 2148}%
\special{pa 1397 2155}%
\special{fp}%
\special{pa 1400 2192}%
\special{pa 1400 2200}%
\special{fp}%
%
\special{pn 8}%
\special{pn 8}%
\special{pa 1400 2200}%
\special{pa 1400 2208}%
\special{fp}%
\special{pa 1398 2245}%
\special{pa 1397 2253}%
\special{fp}%
\special{pa 1391 2289}%
\special{pa 1389 2296}%
\special{fp}%
\special{pa 1380 2331}%
\special{pa 1378 2338}%
\special{fp}%
\special{pa 1365 2371}%
\special{pa 1361 2378}%
\special{fp}%
\special{pa 1343 2409}%
\special{pa 1339 2416}%
\special{fp}%
\special{pa 1317 2443}%
\special{pa 1313 2448}%
\special{fp}%
\special{pa 1287 2470}%
\special{pa 1280 2475}%
\special{fp}%
\special{pa 1249 2491}%
\special{pa 1242 2493}%
\special{fp}%
\special{pa 1208 2500}%
\special{pa 1200 2500}%
\special{fp}%
\put(13.0000,-22.0000){\makebox(0,0){$c_0$}}%
%
\special{pn 4}%
\special{pa 1230 2130}%
\special{pa 1100 2000}%
\special{fp}%
%
\special{pn 4}%
\special{pa 1230 2270}%
\special{pa 1100 2400}%
\special{fp}%
\put(6.5000,-22.0000){\makebox(0,0){$\sb 3$}}%
\put(7.7500,-22.0000){\makebox(0,0){$\circlearrowright$}}%
\put(7.2000,-27.2000){\makebox(0,0)[lt]{$c_1$}}%
\put(10.0000,-22.0000){\makebox(0,0){$p$}}%
\put(15.0000,-19.0000){\makebox(0,0)[lb]{$q=1$,}}%
\put(15.0000,-21.0000){\makebox(0,0)[lb]{$e^{c_0}_{b_0}(p)=-1$}}%
%
\special{pn 8}%
\special{pa 3000 2700}%
\special{pa 3000 1700}%
\special{fp}%
\special{sh 1}%
\special{pa 3000 1700}%
\special{pa 2980 1767}%
\special{pa 3000 1753}%
\special{pa 3020 1767}%
\special{pa 3000 1700}%
\special{fp}%
%
\special{pn 13}%
\special{pa 3500 2500}%
\special{pa 2900 1900}%
\special{fp}%
\special{sh 1}%
\special{pa 2900 1900}%
\special{pa 2933 1961}%
\special{pa 2938 1938}%
\special{pa 2961 1933}%
\special{pa 2900 1900}%
\special{fp}%
%
\special{pn 13}%
\special{pa 2900 2500}%
\special{pa 3500 1900}%
\special{fp}%
\special{sh 1}%
\special{pa 3500 1900}%
\special{pa 3439 1933}%
\special{pa 3462 1938}%
\special{pa 3467 1961}%
\special{pa 3500 1900}%
\special{fp}%
\put(29.0000,-25.0000){\makebox(0,0){$\bullet$}}%
\put(28.8000,-25.2000){\makebox(0,0)[rt]{$b_0$}}%
%
\special{pn 8}%
\special{pn 8}%
\special{pa 3500 1900}%
\special{pa 3507 1901}%
\special{fp}%
\special{pa 3533 1916}%
\special{pa 3537 1921}%
\special{fp}%
\special{pa 3554 1948}%
\special{pa 3557 1954}%
\special{fp}%
\special{pa 3569 1984}%
\special{pa 3571 1990}%
\special{fp}%
\special{pa 3580 2021}%
\special{pa 3582 2028}%
\special{fp}%
\special{pa 3589 2062}%
\special{pa 3590 2070}%
\special{fp}%
\special{pa 3595 2105}%
\special{pa 3596 2112}%
\special{fp}%
\special{pa 3598 2148}%
\special{pa 3599 2156}%
\special{fp}%
\special{pa 3600 2193}%
\special{pa 3600 2201}%
\special{fp}%
\special{pa 3599 2238}%
\special{pa 3599 2246}%
\special{fp}%
\special{pa 3596 2282}%
\special{pa 3595 2290}%
\special{fp}%
\special{pa 3591 2326}%
\special{pa 3589 2334}%
\special{fp}%
\special{pa 3583 2368}%
\special{pa 3581 2375}%
\special{fp}%
\special{pa 3572 2408}%
\special{pa 3569 2416}%
\special{fp}%
\special{pa 3557 2447}%
\special{pa 3553 2453}%
\special{fp}%
\special{pa 3537 2479}%
\special{pa 3532 2484}%
\special{fp}%
\special{pa 3507 2499}%
\special{pa 3500 2500}%
\special{fp}%
%
\special{pn 8}%
\special{pn 8}%
\special{pa 2800 2200}%
\special{pa 2800 2192}%
\special{fp}%
\special{pa 2801 2156}%
\special{pa 2802 2148}%
\special{fp}%
\special{pa 2804 2113}%
\special{pa 2805 2105}%
\special{fp}%
\special{pa 2810 2072}%
\special{pa 2811 2064}%
\special{fp}%
\special{pa 2818 2030}%
\special{pa 2819 2023}%
\special{fp}%
\special{pa 2829 1990}%
\special{pa 2831 1983}%
\special{fp}%
\special{pa 2843 1953}%
\special{pa 2846 1947}%
\special{fp}%
\special{pa 2863 1921}%
\special{pa 2868 1916}%
\special{fp}%
\special{pa 2893 1901}%
\special{pa 2900 1900}%
\special{fp}%
%
\special{pn 8}%
\special{pn 8}%
\special{pa 3000 2700}%
\special{pa 2992 2700}%
\special{fp}%
\special{pa 2959 2689}%
\special{pa 2952 2686}%
\special{fp}%
\special{pa 2928 2666}%
\special{pa 2923 2662}%
\special{fp}%
\special{pa 2902 2636}%
\special{pa 2898 2631}%
\special{fp}%
\special{pa 2880 2600}%
\special{pa 2876 2593}%
\special{fp}%
\special{pa 2861 2560}%
\special{pa 2858 2552}%
\special{fp}%
\special{pa 2845 2518}%
\special{pa 2843 2510}%
\special{fp}%
\special{pa 2833 2476}%
\special{pa 2831 2468}%
\special{fp}%
\special{pa 2823 2432}%
\special{pa 2821 2424}%
\special{fp}%
\special{pa 2815 2388}%
\special{pa 2813 2380}%
\special{fp}%
\special{pa 2809 2344}%
\special{pa 2808 2336}%
\special{fp}%
\special{pa 2804 2299}%
\special{pa 2803 2291}%
\special{fp}%
\special{pa 2801 2254}%
\special{pa 2801 2246}%
\special{fp}%
\special{pa 2800 2208}%
\special{pa 2800 2200}%
\special{fp}%
%
\special{pn 8}%
\special{pn 8}%
\special{pa 2900 2500}%
\special{pa 2900 2492}%
\special{fp}%
\special{pa 2900 2454}%
\special{pa 2900 2446}%
\special{fp}%
\special{pa 2901 2408}%
\special{pa 2901 2400}%
\special{fp}%
\special{pa 2901 2362}%
\special{pa 2902 2353}%
\special{fp}%
\special{pa 2903 2316}%
\special{pa 2903 2307}%
\special{fp}%
\special{pa 2904 2270}%
\special{pa 2905 2262}%
\special{fp}%
\special{pa 2906 2224}%
\special{pa 2907 2216}%
\special{fp}%
\special{pa 2908 2178}%
\special{pa 2909 2170}%
\special{fp}%
\special{pa 2911 2132}%
\special{pa 2912 2124}%
\special{fp}%
\special{pa 2914 2087}%
\special{pa 2915 2079}%
\special{fp}%
\special{pa 2918 2041}%
\special{pa 2919 2033}%
\special{fp}%
\special{pa 2922 1996}%
\special{pa 2923 1988}%
\special{fp}%
\special{pa 2927 1951}%
\special{pa 2928 1943}%
\special{fp}%
\special{pa 2933 1906}%
\special{pa 2934 1898}%
\special{fp}%
\special{pa 2940 1862}%
\special{pa 2941 1854}%
\special{fp}%
\special{pa 2948 1818}%
\special{pa 2949 1810}%
\special{fp}%
\special{pa 2958 1775}%
\special{pa 2960 1768}%
\special{fp}%
\special{pa 2971 1735}%
\special{pa 2974 1728}%
\special{fp}%
\special{pa 2993 1702}%
\special{pa 3000 1700}%
\special{fp}%
\put(31.5000,-23.5000){\makebox(0,0){$\sb 1$}}%
\put(31.5000,-20.5000){\makebox(0,0){$\sb 2$}}%
\put(29.5000,-22.0000){\makebox(0,0){$\sb 3$}}%
\put(30.7500,-22.0000){\makebox(0,0){$\circlearrowleft$}}%
\put(36.0000,-22.0000){\makebox(0,0){$c_0$}}%
%
\special{pn 4}%
\special{pa 3530 2130}%
\special{pa 3400 2000}%
\special{fp}%
%
\special{pn 4}%
\special{pa 3530 2270}%
\special{pa 3400 2400}%
\special{fp}%
\put(30.2000,-27.2000){\makebox(0,0)[lt]{$c_1$}}%
\put(38.0000,-19.0000){\makebox(0,0)[lb]{$q=2$,}}%
\put(38.0000,-21.0000){\makebox(0,0)[lb]{$e^{c_0}_{b_0}(p)=1$}}%
\put(33.0000,-22.0000){\makebox(0,0){$p$}}%
%
\special{pn 8}%
\special{pa 700 4000}%
\special{pa 700 3000}%
\special{fp}%
\special{sh 1}%
\special{pa 700 3000}%
\special{pa 680 3067}%
\special{pa 700 3053}%
\special{pa 720 3067}%
\special{pa 700 3000}%
\special{fp}%
%
\special{pn 13}%
\special{pa 600 3800}%
\special{pa 1200 3200}%
\special{fp}%
\special{sh 1}%
\special{pa 1200 3200}%
\special{pa 1139 3233}%
\special{pa 1162 3238}%
\special{pa 1167 3261}%
\special{pa 1200 3200}%
\special{fp}%
%
\special{pn 13}%
\special{pa 600 3200}%
\special{pa 1200 3800}%
\special{fp}%
\special{sh 1}%
\special{pa 1200 3800}%
\special{pa 1167 3739}%
\special{pa 1162 3762}%
\special{pa 1139 3767}%
\special{pa 1200 3800}%
\special{fp}%
\put(6.0000,-38.0000){\makebox(0,0){$\bullet$}}%
\put(5.8000,-38.2000){\makebox(0,0)[rt]{$b_0$}}%
%
\special{pn 8}%
\special{pn 8}%
\special{pa 600 3200}%
\special{pa 601 3193}%
\special{fp}%
\special{pa 616 3168}%
\special{pa 621 3164}%
\special{fp}%
\special{pa 645 3147}%
\special{pa 652 3144}%
\special{fp}%
\special{pa 681 3132}%
\special{pa 688 3129}%
\special{fp}%
\special{pa 720 3120}%
\special{pa 727 3118}%
\special{fp}%
\special{pa 762 3111}%
\special{pa 769 3110}%
\special{fp}%
\special{pa 805 3105}%
\special{pa 813 3104}%
\special{fp}%
\special{pa 849 3101}%
\special{pa 857 3101}%
\special{fp}%
\special{pa 893 3100}%
\special{pa 901 3100}%
\special{fp}%
\special{pa 938 3101}%
\special{pa 946 3101}%
\special{fp}%
\special{pa 982 3104}%
\special{pa 990 3105}%
\special{fp}%
\special{pa 1025 3109}%
\special{pa 1033 3110}%
\special{fp}%
\special{pa 1068 3117}%
\special{pa 1075 3119}%
\special{fp}%
\special{pa 1109 3128}%
\special{pa 1116 3131}%
\special{fp}%
\special{pa 1146 3143}%
\special{pa 1153 3146}%
\special{fp}%
\special{pa 1179 3164}%
\special{pa 1185 3168}%
\special{fp}%
\special{pa 1199 3193}%
\special{pa 1200 3200}%
\special{fp}%
%
\special{pn 8}%
\special{pn 8}%
\special{pa 1200 3800}%
\special{pa 1200 3808}%
\special{fp}%
\special{pa 1190 3841}%
\special{pa 1186 3846}%
\special{fp}%
\special{pa 1167 3871}%
\special{pa 1163 3876}%
\special{fp}%
\special{pa 1138 3897}%
\special{pa 1132 3900}%
\special{fp}%
\special{pa 1103 3919}%
\special{pa 1096 3922}%
\special{fp}%
\special{pa 1063 3937}%
\special{pa 1056 3941}%
\special{fp}%
\special{pa 1022 3953}%
\special{pa 1014 3956}%
\special{fp}%
\special{pa 979 3966}%
\special{pa 972 3968}%
\special{fp}%
\special{pa 936 3976}%
\special{pa 928 3978}%
\special{fp}%
\special{pa 892 3985}%
\special{pa 883 3986}%
\special{fp}%
\special{pa 846 3991}%
\special{pa 838 3992}%
\special{fp}%
\special{pa 801 3996}%
\special{pa 793 3997}%
\special{fp}%
\special{pa 755 3999}%
\special{pa 746 3999}%
\special{fp}%
\special{pa 708 4000}%
\special{pa 700 4000}%
\special{fp}%
%
\special{pn 8}%
\special{pn 8}%
\special{pa 500 3500}%
\special{pa 500 3492}%
\special{fp}%
\special{pa 501 3454}%
\special{pa 501 3445}%
\special{fp}%
\special{pa 503 3408}%
\special{pa 504 3400}%
\special{fp}%
\special{pa 508 3363}%
\special{pa 509 3355}%
\special{fp}%
\special{pa 514 3318}%
\special{pa 515 3310}%
\special{fp}%
\special{pa 522 3274}%
\special{pa 523 3266}%
\special{fp}%
\special{pa 532 3230}%
\special{pa 534 3222}%
\special{fp}%
\special{pa 544 3187}%
\special{pa 547 3179}%
\special{fp}%
\special{pa 559 3145}%
\special{pa 562 3138}%
\special{fp}%
\special{pa 577 3106}%
\special{pa 581 3099}%
\special{fp}%
\special{pa 598 3070}%
\special{pa 602 3064}%
\special{fp}%
\special{pa 624 3038}%
\special{pa 628 3033}%
\special{fp}%
\special{pa 653 3014}%
\special{pa 659 3011}%
\special{fp}%
\special{pa 692 3000}%
\special{pa 700 3000}%
\special{fp}%
%
\special{pn 8}%
\special{pn 8}%
\special{pa 600 3800}%
\special{pa 593 3799}%
\special{fp}%
\special{pa 567 3784}%
\special{pa 563 3779}%
\special{fp}%
\special{pa 546 3752}%
\special{pa 543 3746}%
\special{fp}%
\special{pa 531 3717}%
\special{pa 529 3711}%
\special{fp}%
\special{pa 520 3679}%
\special{pa 518 3673}%
\special{fp}%
\special{pa 511 3638}%
\special{pa 510 3631}%
\special{fp}%
\special{pa 505 3596}%
\special{pa 504 3588}%
\special{fp}%
\special{pa 502 3553}%
\special{pa 501 3545}%
\special{fp}%
\special{pa 500 3508}%
\special{pa 500 3500}%
\special{fp}%
\put(8.5000,-36.5000){\makebox(0,0){$\sb 1$}}%
\put(8.5000,-33.5000){\makebox(0,0){$\sb 2$}}%
\put(6.5000,-35.0000){\makebox(0,0){$\sb 3$}}%
\put(7.7500,-35.0000){\makebox(0,0){$\circlearrowleft$}}%
\put(10.0000,-35.0000){\makebox(0,0){$p$}}%
\put(15.0000,-32.0000){\makebox(0,0)[lb]{$q=1$,}}%
\put(15.0000,-34.0000){\makebox(0,0)[lb]{$e^{c_0}_{b_0}(p)=-1$}}%
\put(7.2000,-40.2000){\makebox(0,0)[lt]{$c_1$}}%
\put(13.0000,-35.0000){\makebox(0,0){$c_0$}}%
%
\special{pn 4}%
\special{pa 1230 3430}%
\special{pa 1100 3300}%
\special{fp}%
%
\special{pn 4}%
\special{pa 1230 3570}%
\special{pa 1100 3700}%
\special{fp}%
%
\special{pn 8}%
\special{pa 3000 4000}%
\special{pa 3000 3000}%
\special{fp}%
\special{sh 1}%
\special{pa 3000 3000}%
\special{pa 2980 3067}%
\special{pa 3000 3053}%
\special{pa 3020 3067}%
\special{pa 3000 3000}%
\special{fp}%
%
\special{pn 13}%
\special{pa 2900 3800}%
\special{pa 3500 3200}%
\special{fp}%
\special{sh 1}%
\special{pa 3500 3200}%
\special{pa 3439 3233}%
\special{pa 3462 3238}%
\special{pa 3467 3261}%
\special{pa 3500 3200}%
\special{fp}%
%
\special{pn 13}%
\special{pa 2900 3200}%
\special{pa 3500 3800}%
\special{fp}%
\special{sh 1}%
\special{pa 3500 3800}%
\special{pa 3467 3739}%
\special{pa 3462 3762}%
\special{pa 3439 3767}%
\special{pa 3500 3800}%
\special{fp}%
\put(29.0000,-32.0000){\makebox(0,0){$\bullet$}}%
\put(28.8000,-31.8000){\makebox(0,0)[rb]{$b_0$}}%
%
\special{pn 8}%
\special{pn 8}%
\special{pa 3500 3800}%
\special{pa 3499 3807}%
\special{fp}%
\special{pa 3484 3832}%
\special{pa 3480 3836}%
\special{fp}%
\special{pa 3455 3853}%
\special{pa 3448 3856}%
\special{fp}%
\special{pa 3416 3869}%
\special{pa 3410 3871}%
\special{fp}%
\special{pa 3378 3881}%
\special{pa 3370 3882}%
\special{fp}%
\special{pa 3336 3889}%
\special{pa 3328 3890}%
\special{fp}%
\special{pa 3294 3895}%
\special{pa 3287 3896}%
\special{fp}%
\special{pa 3250 3899}%
\special{pa 3242 3899}%
\special{fp}%
\special{pa 3205 3900}%
\special{pa 3197 3900}%
\special{fp}%
\special{pa 3159 3899}%
\special{pa 3151 3899}%
\special{fp}%
\special{pa 3115 3896}%
\special{pa 3108 3895}%
\special{fp}%
\special{pa 3073 3891}%
\special{pa 3065 3889}%
\special{fp}%
\special{pa 3030 3882}%
\special{pa 3023 3881}%
\special{fp}%
\special{pa 2990 3872}%
\special{pa 2984 3869}%
\special{fp}%
\special{pa 2953 3857}%
\special{pa 2945 3853}%
\special{fp}%
\special{pa 2921 3836}%
\special{pa 2916 3832}%
\special{fp}%
\special{pa 2901 3807}%
\special{pa 2900 3800}%
\special{fp}%
%
\special{pn 8}%
\special{pn 8}%
\special{pa 3500 3200}%
\special{pa 3500 3208}%
\special{fp}%
\special{pa 3499 3245}%
\special{pa 3499 3253}%
\special{fp}%
\special{pa 3497 3290}%
\special{pa 3496 3298}%
\special{fp}%
\special{pa 3493 3335}%
\special{pa 3492 3343}%
\special{fp}%
\special{pa 3487 3379}%
\special{pa 3486 3387}%
\special{fp}%
\special{pa 3480 3423}%
\special{pa 3479 3431}%
\special{fp}%
\special{pa 3471 3466}%
\special{pa 3470 3474}%
\special{fp}%
\special{pa 3461 3510}%
\special{pa 3459 3517}%
\special{fp}%
\special{pa 3449 3552}%
\special{pa 3447 3560}%
\special{fp}%
\special{pa 3435 3595}%
\special{pa 3432 3602}%
\special{fp}%
\special{pa 3419 3637}%
\special{pa 3416 3644}%
\special{fp}%
\special{pa 3401 3678}%
\special{pa 3397 3685}%
\special{fp}%
\special{pa 3380 3718}%
\special{pa 3377 3726}%
\special{fp}%
\special{pa 3358 3757}%
\special{pa 3354 3764}%
\special{fp}%
\special{pa 3334 3795}%
\special{pa 3330 3802}%
\special{fp}%
\special{pa 3307 3831}%
\special{pa 3302 3837}%
\special{fp}%
\special{pa 3278 3865}%
\special{pa 3272 3871}%
\special{fp}%
\special{pa 3246 3896}%
\special{pa 3240 3901}%
\special{fp}%
\special{pa 3212 3924}%
\special{pa 3206 3929}%
\special{fp}%
\special{pa 3175 3949}%
\special{pa 3168 3953}%
\special{fp}%
\special{pa 3136 3970}%
\special{pa 3129 3973}%
\special{fp}%
\special{pa 3095 3986}%
\special{pa 3087 3988}%
\special{fp}%
\special{pa 3052 3996}%
\special{pa 3045 3997}%
\special{fp}%
\special{pa 3008 4000}%
\special{pa 3000 4000}%
\special{fp}%
%
\special{pn 8}%
\special{pn 8}%
\special{pa 2900 3200}%
\special{pa 2900 3192}%
\special{fp}%
\special{pa 2902 3156}%
\special{pa 2903 3148}%
\special{fp}%
\special{pa 2910 3115}%
\special{pa 2911 3108}%
\special{fp}%
\special{pa 2921 3077}%
\special{pa 2924 3071}%
\special{fp}%
\special{pa 2938 3044}%
\special{pa 2942 3038}%
\special{fp}%
\special{pa 2960 3016}%
\special{pa 2965 3013}%
\special{fp}%
\special{pa 2993 3001}%
\special{pa 3000 3000}%
\special{fp}%
\put(29.5000,-35.0000){\makebox(0,0){$\sb 3$}}%
\put(31.5000,-36.5000){\makebox(0,0){$\sb 2$}}%
\put(31.5000,-33.5000){\makebox(0,0){$\sb 1$}}%
\put(30.7500,-35.0000){\makebox(0,0){$\circlearrowright$}}%
\put(30.2000,-40.2000){\makebox(0,0)[lt]{$c_1$}}%
\put(36.0000,-35.0000){\makebox(0,0){$c_0$}}%
%
\special{pn 4}%
\special{pa 3530 3430}%
\special{pa 3400 3300}%
\special{fp}%
%
\special{pn 4}%
\special{pa 3530 3570}%
\special{pa 3400 3700}%
\special{fp}%
\put(38.0000,-32.0000){\makebox(0,0)[lb]{$q=2$,}}%
\put(38.0000,-34.0000){\makebox(0,0)[lb]{$e^{c_0}_{b_0}(p)=1$}}%
%
\special{pn 8}%
\special{pa 700 5300}%
\special{pa 700 4300}%
\special{fp}%
\special{sh 1}%
\special{pa 700 4300}%
\special{pa 680 4367}%
\special{pa 700 4353}%
\special{pa 720 4367}%
\special{pa 700 4300}%
\special{fp}%
%
\special{pn 13}%
\special{pa 600 4500}%
\special{pa 1200 5100}%
\special{fp}%
\special{sh 1}%
\special{pa 1200 5100}%
\special{pa 1167 5039}%
\special{pa 1162 5062}%
\special{pa 1139 5067}%
\special{pa 1200 5100}%
\special{fp}%
%
\special{pn 13}%
\special{pa 1200 4500}%
\special{pa 600 5100}%
\special{fp}%
\special{sh 1}%
\special{pa 600 5100}%
\special{pa 661 5067}%
\special{pa 638 5062}%
\special{pa 633 5039}%
\special{pa 600 5100}%
\special{fp}%
\put(6.0000,-45.0000){\makebox(0,0){$\bullet$}}%
\put(5.8000,-44.8000){\makebox(0,0)[rb]{$b_0$}}%
%
\special{pn 8}%
\special{pn 8}%
\special{pa 1200 4500}%
\special{pa 1207 4501}%
\special{fp}%
\special{pa 1232 4516}%
\special{pa 1236 4521}%
\special{fp}%
\special{pa 1253 4545}%
\special{pa 1256 4552}%
\special{fp}%
\special{pa 1269 4582}%
\special{pa 1271 4588}%
\special{fp}%
\special{pa 1280 4620}%
\special{pa 1282 4627}%
\special{fp}%
\special{pa 1289 4662}%
\special{pa 1290 4670}%
\special{fp}%
\special{pa 1295 4705}%
\special{pa 1296 4713}%
\special{fp}%
\special{pa 1299 4749}%
\special{pa 1299 4758}%
\special{fp}%
\special{pa 1300 4796}%
\special{pa 1300 4804}%
\special{fp}%
\special{pa 1299 4842}%
\special{pa 1299 4850}%
\special{fp}%
\special{pa 1296 4887}%
\special{pa 1295 4894}%
\special{fp}%
\special{pa 1290 4929}%
\special{pa 1289 4938}%
\special{fp}%
\special{pa 1282 4973}%
\special{pa 1280 4979}%
\special{fp}%
\special{pa 1271 5012}%
\special{pa 1269 5018}%
\special{fp}%
\special{pa 1257 5047}%
\special{pa 1253 5054}%
\special{fp}%
\special{pa 1237 5079}%
\special{pa 1232 5084}%
\special{fp}%
\special{pa 1207 5099}%
\special{pa 1200 5100}%
\special{fp}%
%
\special{pn 8}%
\special{pn 8}%
\special{pa 700 5300}%
\special{pa 693 5299}%
\special{fp}%
\special{pa 665 5288}%
\special{pa 660 5284}%
\special{fp}%
\special{pa 642 5263}%
\special{pa 638 5258}%
\special{fp}%
\special{pa 625 5231}%
\special{pa 622 5225}%
\special{fp}%
\special{pa 612 5194}%
\special{pa 610 5188}%
\special{fp}%
\special{pa 604 5153}%
\special{pa 603 5145}%
\special{fp}%
\special{pa 600 5108}%
\special{pa 600 5100}%
\special{fp}%
%
\special{pn 8}%
\special{pn 8}%
\special{pa 600 4500}%
\special{pa 600 4492}%
\special{fp}%
\special{pa 603 4455}%
\special{pa 603 4448}%
\special{fp}%
\special{pa 610 4414}%
\special{pa 611 4407}%
\special{fp}%
\special{pa 621 4376}%
\special{pa 624 4371}%
\special{fp}%
\special{pa 638 4343}%
\special{pa 642 4338}%
\special{fp}%
\special{pa 660 4316}%
\special{pa 665 4313}%
\special{fp}%
\special{pa 693 4301}%
\special{pa 700 4300}%
\special{fp}%
\put(8.5000,-49.5000){\makebox(0,0){$\sb 2$}}%
\put(8.5000,-46.5000){\makebox(0,0){$\sb 1$}}%
\put(6.5000,-48.0000){\makebox(0,0){$\sb 3$}}%
\put(7.2000,-53.2000){\makebox(0,0)[lt]{$c_1$}}%
\put(13.0000,-48.0000){\makebox(0,0){$c_0$}}%
%
\special{pn 4}%
\special{pa 1230 4730}%
\special{pa 1100 4600}%
\special{fp}%
%
\special{pn 4}%
\special{pa 1230 4870}%
\special{pa 1100 5000}%
\special{fp}%
\put(10.0000,-48.0000){\makebox(0,0){$p$}}%
\put(7.7500,-48.0000){\makebox(0,0){$\circlearrowright$}}%
\put(15.0000,-45.0000){\makebox(0,0)[lb]{$q=3$,}}%
\put(15.0000,-47.0000){\makebox(0,0)[lb]{$e^{c_0}_{b_0}(p)=-1$}}%
%
\special{pn 8}%
\special{pa 3000 5300}%
\special{pa 3000 4300}%
\special{fp}%
\special{sh 1}%
\special{pa 3000 4300}%
\special{pa 2980 4367}%
\special{pa 3000 4353}%
\special{pa 3020 4367}%
\special{pa 3000 4300}%
\special{fp}%
%
\special{pn 13}%
\special{pa 3500 4500}%
\special{pa 2900 5100}%
\special{fp}%
\special{sh 1}%
\special{pa 2900 5100}%
\special{pa 2961 5067}%
\special{pa 2938 5062}%
\special{pa 2933 5039}%
\special{pa 2900 5100}%
\special{fp}%
%
\special{pn 13}%
\special{pa 2900 4500}%
\special{pa 3500 5100}%
\special{fp}%
\special{sh 1}%
\special{pa 3500 5100}%
\special{pa 3467 5039}%
\special{pa 3462 5062}%
\special{pa 3439 5067}%
\special{pa 3500 5100}%
\special{fp}%
\put(35.0000,-45.0000){\makebox(0,0){$\bullet$}}%
\put(35.2000,-44.8000){\makebox(0,0)[lb]{$b_0$}}%
%
\special{pn 8}%
\special{pn 8}%
\special{pa 2900 5100}%
\special{pa 2893 5099}%
\special{fp}%
\special{pa 2868 5084}%
\special{pa 2864 5079}%
\special{fp}%
\special{pa 2847 5055}%
\special{pa 2844 5048}%
\special{fp}%
\special{pa 2831 5018}%
\special{pa 2829 5012}%
\special{fp}%
\special{pa 2820 4980}%
\special{pa 2818 4973}%
\special{fp}%
\special{pa 2811 4938}%
\special{pa 2810 4930}%
\special{fp}%
\special{pa 2805 4895}%
\special{pa 2804 4887}%
\special{fp}%
\special{pa 2801 4851}%
\special{pa 2801 4842}%
\special{fp}%
\special{pa 2800 4804}%
\special{pa 2800 4796}%
\special{fp}%
\special{pa 2801 4758}%
\special{pa 2801 4750}%
\special{fp}%
\special{pa 2804 4713}%
\special{pa 2805 4706}%
\special{fp}%
\special{pa 2810 4671}%
\special{pa 2811 4662}%
\special{fp}%
\special{pa 2818 4627}%
\special{pa 2820 4621}%
\special{fp}%
\special{pa 2829 4588}%
\special{pa 2831 4582}%
\special{fp}%
\special{pa 2843 4553}%
\special{pa 2847 4546}%
\special{fp}%
\special{pa 2863 4521}%
\special{pa 2868 4516}%
\special{fp}%
\special{pa 2893 4501}%
\special{pa 2900 4500}%
\special{fp}%
%
\special{pn 8}%
\special{pn 8}%
\special{pa 3500 5100}%
\special{pa 3500 5108}%
\special{fp}%
\special{pa 3491 5139}%
\special{pa 3487 5145}%
\special{fp}%
\special{pa 3470 5168}%
\special{pa 3466 5172}%
\special{fp}%
\special{pa 3444 5192}%
\special{pa 3439 5196}%
\special{fp}%
\special{pa 3413 5213}%
\special{pa 3407 5216}%
\special{fp}%
\special{pa 3379 5231}%
\special{pa 3373 5233}%
\special{fp}%
\special{pa 3343 5245}%
\special{pa 3337 5248}%
\special{fp}%
\special{pa 3305 5259}%
\special{pa 3298 5260}%
\special{fp}%
\special{pa 3266 5269}%
\special{pa 3258 5271}%
\special{fp}%
\special{pa 3225 5279}%
\special{pa 3217 5280}%
\special{fp}%
\special{pa 3183 5286}%
\special{pa 3175 5287}%
\special{fp}%
\special{pa 3140 5292}%
\special{pa 3132 5293}%
\special{fp}%
\special{pa 3096 5296}%
\special{pa 3088 5297}%
\special{fp}%
\special{pa 3052 5299}%
\special{pa 3044 5299}%
\special{fp}%
\special{pa 3008 5300}%
\special{pa 3000 5300}%
\special{fp}%
%
\special{pn 8}%
\special{pn 8}%
\special{pa 3000 4300}%
\special{pa 3008 4300}%
\special{fp}%
\special{pa 3047 4301}%
\special{pa 3056 4301}%
\special{fp}%
\special{pa 3093 4303}%
\special{pa 3102 4304}%
\special{fp}%
\special{pa 3139 4308}%
\special{pa 3147 4309}%
\special{fp}%
\special{pa 3184 4314}%
\special{pa 3192 4315}%
\special{fp}%
\special{pa 3230 4322}%
\special{pa 3238 4324}%
\special{fp}%
\special{pa 3274 4333}%
\special{pa 3281 4335}%
\special{fp}%
\special{pa 3317 4345}%
\special{pa 3325 4348}%
\special{fp}%
\special{pa 3357 4360}%
\special{pa 3364 4363}%
\special{fp}%
\special{pa 3395 4377}%
\special{pa 3402 4381}%
\special{fp}%
\special{pa 3432 4399}%
\special{pa 3438 4404}%
\special{fp}%
\special{pa 3462 4424}%
\special{pa 3466 4428}%
\special{fp}%
\special{pa 3486 4454}%
\special{pa 3490 4459}%
\special{fp}%
\special{pa 3500 4492}%
\special{pa 3500 4500}%
\special{fp}%
\put(31.5000,-49.5000){\makebox(0,0){$\sb 1$}}%
\put(31.5000,-46.5000){\makebox(0,0){$\sb 2$}}%
\put(29.5000,-48.0000){\makebox(0,0){$\sb 3$}}%
\put(30.2000,-53.2000){\makebox(0,0)[lt]{$c_1$}}%
\put(36.0000,-48.0000){\makebox(0,0){$c_0$}}%
%
\special{pn 4}%
\special{pa 3530 4730}%
\special{pa 3400 4600}%
\special{fp}%
%
\special{pn 4}%
\special{pa 3530 4870}%
\special{pa 3400 5000}%
\special{fp}%
\put(30.7500,-48.0000){\makebox(0,0){$\circlearrowleft$}}%
\put(38.0000,-45.0000){\makebox(0,0)[lb]{$q=0$,}}%
\put(38.0000,-47.0000){\makebox(0,0)[lb]{$e^{c_0}_{b_0}(p)=1$}}%
\put(33.0000,-35.0000){\makebox(0,0){$p$}}%
\put(33.0000,-48.0000){\makebox(0,0){$p$}}%
\end{picture}}%

%% file: e_y=e_z=-s_x.tex
{\unitlength 0.1in%
\begin{picture}(29.5500,10.6500)(4.4500,-13.0000)%
%
\special{pn 13}%
\special{pa 600 1200}%
\special{pa 600 400}%
\special{fp}%
\special{sh 1}%
\special{pa 600 400}%
\special{pa 580 467}%
\special{pa 600 453}%
\special{pa 620 467}%
\special{pa 600 400}%
\special{fp}%
%
\special{pn 13}%
\special{pa 1000 1200}%
\special{pa 1000 400}%
\special{fp}%
\special{sh 1}%
\special{pa 1000 400}%
\special{pa 980 467}%
\special{pa 1000 453}%
\special{pa 1020 467}%
\special{pa 1000 400}%
\special{fp}%
%
\special{pn 8}%
\special{pa 600 800}%
\special{pa 1400 800}%
\special{dt 0.045}%
\put(13.0000,-9.0000){\makebox(0,0){$\Gamma$}}%
\put(10.3000,-7.7000){\makebox(0,0)[lb]{$x_j$}}%
\put(6.3000,-7.7000){\makebox(0,0)[lb]{$b_0$}}%
\put(6.0000,-13.0000){\makebox(0,0){$C_0$}}%
\put(10.0000,-13.0000){\makebox(0,0){$C_0$}}%
%
\special{pn 13}%
\special{ar 2500 600 100 100 1.5707963 3.1415927}%
%
\special{pn 13}%
\special{pa 2400 600}%
\special{pa 2400 400}%
\special{fp}%
\special{sh 1}%
\special{pa 2400 400}%
\special{pa 2380 467}%
\special{pa 2400 453}%
\special{pa 2420 467}%
\special{pa 2400 400}%
\special{fp}%
%
\special{pn 13}%
\special{pa 2400 1200}%
\special{pa 2400 1000}%
\special{fp}%
%
\special{pn 13}%
\special{ar 2500 1000 100 100 3.1415927 4.7123890}%
%
\special{pn 13}%
\special{ar 2500 800 700 100 4.7123890 1.5707963}%
%
\special{pn 13}%
\special{pa 2800 1200}%
\special{pa 2800 400}%
\special{fp}%
\special{sh 1}%
\special{pa 2800 400}%
\special{pa 2780 467}%
\special{pa 2800 453}%
\special{pa 2820 467}%
\special{pa 2800 400}%
\special{fp}%
\put(32.3000,-7.7000){\makebox(0,0)[lb]{$\overline{b}_0$}}%
%
\special{pn 8}%
\special{pn 8}%
\special{pa 2400 400}%
\special{pa 2400 392}%
\special{fp}%
\special{pa 2408 361}%
\special{pa 2411 355}%
\special{fp}%
\special{pa 2428 331}%
\special{pa 2432 327}%
\special{fp}%
\special{pa 2457 310}%
\special{pa 2464 307}%
\special{fp}%
\special{pa 2495 300}%
\special{pa 2503 300}%
\special{fp}%
\special{pa 2536 307}%
\special{pa 2542 309}%
\special{fp}%
\special{pa 2568 326}%
\special{pa 2572 330}%
\special{fp}%
\special{pa 2589 355}%
\special{pa 2592 361}%
\special{fp}%
\special{pa 2600 392}%
\special{pa 2600 400}%
\special{fp}%
%
\special{pn 8}%
\special{pn 8}%
\special{pa 2800 1200}%
\special{pa 2800 1208}%
\special{fp}%
\special{pa 2792 1239}%
\special{pa 2789 1246}%
\special{fp}%
\special{pa 2772 1270}%
\special{pa 2768 1274}%
\special{fp}%
\special{pa 2742 1291}%
\special{pa 2735 1294}%
\special{fp}%
\special{pa 2704 1300}%
\special{pa 2696 1300}%
\special{fp}%
\special{pa 2665 1294}%
\special{pa 2658 1291}%
\special{fp}%
\special{pa 2633 1274}%
\special{pa 2629 1270}%
\special{fp}%
\special{pa 2611 1245}%
\special{pa 2608 1239}%
\special{fp}%
\special{pa 2600 1208}%
\special{pa 2600 1200}%
\special{fp}%
%
\special{pn 8}%
\special{pa 2600 1200}%
\special{pa 2600 400}%
\special{dt 0.045}%
\put(32.0000,-8.0000){\makebox(0,0){$\bullet$}}%
%
\special{pn 8}%
\special{pa 3100 1200}%
\special{pa 3100 1000}%
\special{fp}%
\special{sh 1}%
\special{pa 3100 1000}%
\special{pa 3080 1067}%
\special{pa 3100 1053}%
\special{pa 3120 1067}%
\special{pa 3100 1000}%
\special{fp}%
%
\special{pn 8}%
\special{pa 3100 1200}%
\special{pa 3300 1200}%
\special{fp}%
\special{sh 1}%
\special{pa 3300 1200}%
\special{pa 3233 1180}%
\special{pa 3247 1200}%
\special{pa 3233 1220}%
\special{pa 3300 1200}%
\special{fp}%
\put(33.0000,-11.5000){\makebox(0,0){$\sb 2$}}%
\put(31.5000,-10.0000){\makebox(0,0){$\sb 1$}}%
%
\special{pn 4}%
\special{pa 2800 900}%
\special{pa 3000 1100}%
\special{fp}%
%
\special{pn 8}%
\special{pa 3200 500}%
\special{pa 3000 500}%
\special{fp}%
\special{sh 1}%
\special{pa 3000 500}%
\special{pa 3067 520}%
\special{pa 3053 500}%
\special{pa 3067 480}%
\special{pa 3000 500}%
\special{fp}%
%
\special{pn 8}%
\special{pa 3200 500}%
\special{pa 3200 300}%
\special{fp}%
\special{sh 1}%
\special{pa 3200 300}%
\special{pa 3180 367}%
\special{pa 3200 353}%
\special{pa 3220 367}%
\special{pa 3200 300}%
\special{fp}%
\put(30.0000,-4.5000){\makebox(0,0){$\sb 1$}}%
\put(31.5000,-3.0000){\makebox(0,0){$\sb 2$}}%
%
\special{pn 4}%
\special{pa 2800 700}%
\special{pa 3000 600}%
\special{fp}%
\put(27.7000,-9.3000){\makebox(0,0)[rt]{$y_j$}}%
\put(27.7000,-6.7000){\makebox(0,0)[rb]{$z_j$}}%
\put(19.0000,-8.0000){\makebox(0,0){$\rightsquigarrow$}}%
\put(34.0000,-12.0000){\makebox(0,0)[lb]{$e_{\overline{b}_0}(y_j)=-1$}}%
\put(34.0000,-5.0000){\makebox(0,0)[lb]{$e_{\overline{b}_0}(z_j)=-1$}}%
\end{picture}}%

%% file: four_e_cancel.tex
{\unitlength 0.1in%
\begin{picture}(38.5500,9.0000)(1.4500,-10.3500)%
%
\special{pn 13}%
\special{pa 300 900}%
\special{pa 300 300}%
\special{fp}%
\special{sh 1}%
\special{pa 300 300}%
\special{pa 280 367}%
\special{pa 300 353}%
\special{pa 320 367}%
\special{pa 300 300}%
\special{fp}%
\put(3.0000,-10.0000){\makebox(0,0){$C_0$}}%
%
\special{pn 13}%
\special{pa 300 600}%
\special{pa 500 600}%
\special{dt 0.045}%
%
\special{pn 13}%
\special{pn 13}%
\special{pa 700 400}%
\special{pa 712 400}%
\special{fp}%
\special{pa 743 405}%
\special{pa 755 407}%
\special{fp}%
\special{pa 784 418}%
\special{pa 795 424}%
\special{fp}%
\special{pa 820 440}%
\special{pa 830 448}%
\special{fp}%
\special{pa 852 470}%
\special{pa 860 480}%
\special{fp}%
\special{pa 876 505}%
\special{pa 882 516}%
\special{fp}%
\special{pa 892 545}%
\special{pa 895 557}%
\special{fp}%
\special{pa 900 588}%
\special{pa 900 600}%
\special{fp}%
\special{pn 13}%
\special{pn 13}%
\special{pa 700 800}%
\special{pa 713 799}%
\special{fp}%
\special{pa 744 794}%
\special{pa 756 791}%
\special{fp}%
\special{pa 785 779}%
\special{pa 796 772}%
\special{fp}%
\special{pa 824 757}%
\special{pa 834 749}%
\special{fp}%
\special{pa 860 731}%
\special{pa 870 723}%
\special{fp}%
\special{pa 895 704}%
\special{pa 905 696}%
\special{fp}%
\special{pa 930 677}%
\special{pa 940 669}%
\special{fp}%
\special{pa 966 651}%
\special{pa 976 643}%
\special{fp}%
\special{pa 1004 628}%
\special{pa 1015 621}%
\special{fp}%
\special{pa 1044 609}%
\special{pa 1056 606}%
\special{fp}%
\special{pa 1087 601}%
\special{pa 1100 600}%
\special{fp}%
%
\special{pn 13}%
\special{pa 1100 300}%
\special{pa 1100 900}%
\special{fp}%
\special{sh 1}%
\special{pa 1100 900}%
\special{pa 1120 833}%
\special{pa 1100 847}%
\special{pa 1080 833}%
\special{pa 1100 900}%
\special{fp}%
\put(11.0000,-10.0000){\makebox(0,0){$C_1$}}%
\put(6.0000,-5.0000){\makebox(0,0){$\Gamma$}}%
\put(15.0000,-6.0000){\makebox(0,0){$\rightsquigarrow$}}%
%
\special{pn 13}%
\special{pa 1900 550}%
\special{pa 1900 300}%
\special{fp}%
\special{sh 1}%
\special{pa 1900 300}%
\special{pa 1880 367}%
\special{pa 1900 353}%
\special{pa 1920 367}%
\special{pa 1900 300}%
\special{fp}%
%
\special{pn 13}%
\special{pa 500 600}%
\special{pa 900 600}%
\special{fp}%
%
\special{pn 13}%
\special{pa 700 400}%
\special{pa 700 800}%
\special{fp}%
%
\special{pn 13}%
\special{pn 13}%
\special{pa 2350 400}%
\special{pa 2364 401}%
\special{fp}%
\special{pa 2396 407}%
\special{pa 2408 412}%
\special{fp}%
\special{pa 2435 427}%
\special{pa 2445 434}%
\special{fp}%
\special{pa 2467 456}%
\special{pa 2474 465}%
\special{fp}%
\special{pa 2488 492}%
\special{pa 2493 504}%
\special{fp}%
\special{pa 2499 536}%
\special{pa 2500 550}%
\special{fp}%
%
\special{pn 13}%
\special{pn 13}%
\special{pa 2250 400}%
\special{pa 2257 390}%
\special{fp}%
\special{pa 2275 371}%
\special{pa 2283 363}%
\special{fp}%
\special{pa 2306 346}%
\special{pa 2317 340}%
\special{fp}%
\special{pa 2345 329}%
\special{pa 2356 325}%
\special{fp}%
\special{pa 2387 320}%
\special{pa 2401 320}%
\special{fp}%
\special{pa 2432 323}%
\special{pa 2444 325}%
\special{fp}%
\special{pa 2473 335}%
\special{pa 2485 341}%
\special{fp}%
\special{pa 2509 357}%
\special{pa 2518 364}%
\special{fp}%
\special{pa 2537 383}%
\special{pa 2544 392}%
\special{fp}%
\special{pa 2560 417}%
\special{pa 2565 427}%
\special{fp}%
\special{pa 2575 456}%
\special{pa 2577 469}%
\special{fp}%
\special{pa 2580 500}%
\special{pa 2579 514}%
\special{fp}%
\special{pa 2575 544}%
\special{pa 2571 555}%
\special{fp}%
\special{pa 2559 584}%
\special{pa 2553 595}%
\special{fp}%
\special{pa 2536 618}%
\special{pa 2529 625}%
\special{fp}%
\special{pa 2509 643}%
\special{pa 2500 650}%
\special{fp}%
%
\special{pn 13}%
\special{pa 1900 900}%
\special{pa 1900 650}%
\special{fp}%
%
\special{pn 13}%
\special{pa 1900 550}%
\special{pa 2100 550}%
\special{dt 0.045}%
%
\special{pn 13}%
\special{pa 1900 650}%
\special{pa 2100 650}%
\special{dt 0.045}%
\special{pn 13}%
\special{pn 13}%
\special{pa 2370 800}%
\special{pa 2383 797}%
\special{fp}%
\special{pa 2412 783}%
\special{pa 2423 776}%
\special{fp}%
\special{pa 2449 756}%
\special{pa 2459 747}%
\special{fp}%
\special{pa 2482 724}%
\special{pa 2491 715}%
\special{fp}%
\special{pa 2512 691}%
\special{pa 2521 681}%
\special{fp}%
\special{pa 2543 657}%
\special{pa 2552 647}%
\special{fp}%
\special{pa 2574 624}%
\special{pa 2583 615}%
\special{fp}%
\special{pa 2607 592}%
\special{pa 2617 584}%
\special{fp}%
\special{pa 2644 566}%
\special{pa 2656 561}%
\special{fp}%
\special{pa 2687 551}%
\special{pa 2700 550}%
\special{fp}%
%
\special{pn 13}%
\special{pn 13}%
\special{pa 2350 900}%
\special{pa 2339 899}%
\special{fp}%
\special{pa 2312 893}%
\special{pa 2302 888}%
\special{fp}%
\special{pa 2282 874}%
\special{pa 2276 868}%
\special{fp}%
\special{pa 2262 848}%
\special{pa 2257 838}%
\special{fp}%
\special{pa 2251 811}%
\special{pa 2250 800}%
\special{fp}%
\special{pn 13}%
\special{pn 13}%
\special{pa 2370 900}%
\special{pa 2383 897}%
\special{fp}%
\special{pa 2412 883}%
\special{pa 2423 876}%
\special{fp}%
\special{pa 2449 856}%
\special{pa 2459 847}%
\special{fp}%
\special{pa 2482 824}%
\special{pa 2491 815}%
\special{fp}%
\special{pa 2512 791}%
\special{pa 2521 781}%
\special{fp}%
\special{pa 2543 757}%
\special{pa 2552 747}%
\special{fp}%
\special{pa 2574 724}%
\special{pa 2583 715}%
\special{fp}%
\special{pa 2607 692}%
\special{pa 2617 684}%
\special{fp}%
\special{pa 2644 666}%
\special{pa 2656 661}%
\special{fp}%
\special{pa 2687 651}%
\special{pa 2700 650}%
\special{fp}%
%
\special{pn 13}%
\special{ar 2700 600 50 50 4.7123890 1.5707963}%
%
\special{pn 13}%
\special{pa 2800 300}%
\special{pa 2800 900}%
\special{fp}%
\special{sh 1}%
\special{pa 2800 900}%
\special{pa 2820 833}%
\special{pa 2800 847}%
\special{pa 2780 833}%
\special{pa 2800 900}%
\special{fp}%
%
\special{pn 13}%
\special{pa 2100 650}%
\special{pa 2500 650}%
\special{fp}%
\special{sh 1}%
\special{pa 2500 650}%
\special{pa 2433 630}%
\special{pa 2447 650}%
\special{pa 2433 670}%
\special{pa 2500 650}%
\special{fp}%
%
\special{pn 13}%
\special{pa 2500 550}%
\special{pa 2100 550}%
\special{fp}%
\special{sh 1}%
\special{pa 2100 550}%
\special{pa 2167 570}%
\special{pa 2153 550}%
\special{pa 2167 530}%
\special{pa 2100 550}%
\special{fp}%
%
\special{pn 13}%
\special{pa 2250 400}%
\special{pa 2250 800}%
\special{fp}%
\special{sh 1}%
\special{pa 2250 800}%
\special{pa 2270 733}%
\special{pa 2250 747}%
\special{pa 2230 733}%
\special{pa 2250 800}%
\special{fp}%
%
\special{pn 13}%
\special{pa 2350 800}%
\special{pa 2350 400}%
\special{fp}%
\special{sh 1}%
\special{pa 2350 400}%
\special{pa 2330 467}%
\special{pa 2350 453}%
\special{pa 2370 467}%
\special{pa 2350 400}%
\special{fp}%
\put(28.0000,-10.0000){\makebox(0,0){$C_1$}}%
\put(19.0000,-10.0000){\makebox(0,0){$\overline{C}_0$}}%
%
\special{pn 13}%
\special{pa 3200 800}%
\special{pa 4000 800}%
\special{fp}%
\special{sh 1}%
\special{pa 4000 800}%
\special{pa 3933 780}%
\special{pa 3947 800}%
\special{pa 3933 820}%
\special{pa 4000 800}%
\special{fp}%
%
\special{pn 13}%
\special{pa 4000 500}%
\special{pa 3200 500}%
\special{fp}%
\special{sh 1}%
\special{pa 3200 500}%
\special{pa 3267 520}%
\special{pa 3253 500}%
\special{pa 3267 480}%
\special{pa 3200 500}%
\special{fp}%
%
\special{pn 13}%
\special{pa 3400 300}%
\special{pa 3400 1000}%
\special{fp}%
\special{sh 1}%
\special{pa 3400 1000}%
\special{pa 3420 933}%
\special{pa 3400 947}%
\special{pa 3380 933}%
\special{pa 3400 1000}%
\special{fp}%
%
\special{pn 13}%
\special{pa 3800 1000}%
\special{pa 3800 300}%
\special{fp}%
\special{sh 1}%
\special{pa 3800 300}%
\special{pa 3780 367}%
\special{pa 3800 353}%
\special{pa 3820 367}%
\special{pa 3800 300}%
\special{fp}%
\put(35.0000,-7.5000){\makebox(0,0){$\sb 1$}}%
\put(39.0000,-8.5000){\makebox(0,0){$\sb 1$}}%
\put(34.5000,-6.0000){\makebox(0,0){$\sb 1$}}%
\put(33.5000,-9.0000){\makebox(0,0){$\sb 2$}}%
\put(38.5000,-7.0000){\makebox(0,0){$\sb 2$}}%
\put(37.0000,-4.5000){\makebox(0,0){$\sb 2$}}%
\put(38.5000,-4.0000){\makebox(0,0){$\sb 1$}}%
\put(33.0000,-5.5000){\makebox(0,0){$\sb 2$}}%
\put(38.0000,-11.0000){\makebox(0,0){$e=+1$}}%
\put(38.0000,-2.0000){\makebox(0,0){$e=+1$}}%
\put(33.0000,-11.0000){\makebox(0,0){$e=-1$}}%
\put(33.0000,-2.0000){\makebox(0,0){$e=-1$}}%
\put(27.5000,-6.0000){\makebox(0,0){$\bullet$}}%
%
\special{pn 8}%
\special{pa 2750 600}%
\special{pa 2900 400}%
\special{fp}%
\put(29.0000,-4.0000){\makebox(0,0)[lb]{$\overline{b}_0$}}%
\end{picture}}%

%% file: Arnold_Inv_Sum.bbl
\begin{thebibliography}{99}
\bibitem{Arnold}
V.~I.~Arnold,
\emph{Plane curves, their invariants, perestroikas and classifications},
Adv.\ Soviet Math., \textbf{21}, Singularities and bifurcations, 33--91, Amer.\ Math.\ Soc., Providence, RI, 1994

\bibitem{MendesRomero}
C.~Mendes de Jesus and M.~C.~Romero Fuster,
\emph{Bridges, channels and Arnold's invariants for generic plane curves},
Topology Appl.\ \textbf{125} (2002), no.\ 3, 505--524

\bibitem{Shumakovich95}
A.~Shumakovich,
\emph{Explicit formulas for the strangeness of plane curves} (Russian. Russian summary),
Algebra i Analiz \textbf{7} (1995), no.\ 3, 165--199;
translation in St.\ Petersburg Math.\ J.\ \textbf{7} (1996), no.\ 3, 445--472

\bibitem{Whitney}
H.~Whitney,
\emph{On regular closed curves in the plane},
Compositio Math.\ \textbf{4} (1937), 276--284
\end{thebibliography}
